\documentclass{amsart}
\usepackage{amssymb,latexsym,amsmath,amscd,graphicx,graphics,epic,eepic}

\newcommand{\zed}{\mathbb{Z}}

\newcommand{\C}{\mathbb{C}}
\newcommand{\fil}{\mathcal{F}}
\newcommand{\Hom}{\text{Hom}}
\newcommand{\im}{\text{Im}~}

\theoremstyle{plain}
\newtheorem{theorem}{Theorem}[section]
\newtheorem{lemma}[theorem]{Lemma}
\newtheorem{proposition}[theorem]{Proposition}
\newtheorem{corollary}[theorem]{Corollary}

\theoremstyle{definition}
\newtheorem{definition}[theorem]{Definition}

\newtheorem{acknowledgments}{Acknowledgments\ignorespaces}

\theoremstyle{remark}
\newtheorem{remark}[theorem]{Remark}

\numberwithin{equation}{section}

\newcommand{\id}{\text{id}}

\begin{document}

\title{On the quantum filtration of the Khovanov-Rozansky cohomology}

\author{Hao Wu}

\address{Department of mathematics and Statistics, Lederle Graduate Research Tower, 710 North Pleasant Street, University of Massachusetts, Amherst, MA 01003-9305, USA}

\email{wu@math.umass.edu}



\begin{abstract}
We prove the quantum filtration on the Khovanov-Rozansky link cohomology $H_p$ with a general degree $(n+1)$ monic potential polynomial $p(x)$ is invariant under Reidemeister moves, and construct a spectral sequence converging to $H_p$ that is invariant under Reidemeister moves, whose $E_1$ term is isomorphic to the Khovanov-Rozansky $\mathfrak{sl}(n)$-cohomology $H_n$. Then we define a generalization of the Rasmussen invariant, and study some of its properties. We also discuss relations between upper bounds of the self-linking number of transversal links in standard contact $S^3$.
\end{abstract}

\maketitle

\section{Introduction and Statements of Results}

For $n\geq2$, Khovanov and Rozansky \cite{KR1} defined an invariant for links in $S^3$ known as the Khovanov-Rozansky $\mathfrak{sl}(n)$-cohomology $H_n$, which is a $\zed\oplus\zed$-graded vector space, where one of the $\zed$-grading is the cohomological grading, and the other comes from the degree of polynomials, and is called the quantum grading. Its construction is based on matrix factorizations with potential from $x^{n+1}$ associated to certain planar diagrams. When $n=2$, their theory is equivalent to the Khovanov homology defined in \cite{K1}. When using a non-homogeneous potential, the cohomological grading remains a $\zed$-grading, while the quantum grading degenerates into a filtration, which we call the quantum filtration, and denote by $\mathcal{F}$. The works of Lee \cite{Lee2} and Rasmussen \cite{Ras1} show that it is interesting to study Khovanov-Rozansky cohomologies defined using a non-homogeneous polynomial. In \cite{Gornik}, Gornik studied the Khovanov-Rozansky cohomology with the potential $x^{n+1}-(n+1)\beta^nx$. Recently, Rasmussen \cite{Ras2} proved that, for any polynomial $p(x)$, the Khovanov-Rozansky cohomology $H_p$ with potential $p(x)$ is invariant as a $\mathbb{Z}$-graded vector space under braid-like Reidemeister moves, where the $\mathbb{Z}$-grading is the cohomological grading. However, it was not known if the quantum filtration on $H_p$ is invariant under Reidemeister moves except in the special case of the $\mathfrak{sl}(n)$-cohomology $H_n$, i.e. when $p(x)=x^{n+1}$. The main goal of this paper is to prove the invariance of the quantum filtration under Reidemeister moves. This makes $H_p$ a much more interesting link invariant.

All the links in this paper are in $S^3$. For minor technical convenience, we work over the base field $\mathbb{C}$ and consider only monic potential polynomials (so that we can cite the results of Khovanov and Rozansky \cite{KR1} and Gornik \cite{Gornik} without reformulate and re-prove them.) Note that the arguments and calculations in \cite{KR1} remain true over $\C$.

\begin{theorem}\label{fil-inv}
For any
\[
p(x)=x^{n+1}+\sum_{i=0}^{n}c_i x^i \in \mathbb{C}[x],
\]
the Khovanov-Rozansky cohomology $H_p$ with potential $p(x)$ is invariant under Reidemeister moves as a finite dimensional $\zed$-graded and filtered $\mathbb{C}$-vector space, where the $\zed$-grading is the cohomological grading, and the filtration is the quantum filtration $\mathcal{F}$.
\end{theorem}

Our method also gives us the following relation between the $\mathfrak{sl}(n)$-cohomology and the $H_p$ cohomology, which is a generalization of the corresponding result in \cite{Gornik}. 

\begin{theorem}\label{spec}
For any
\[
p(x)=x^{n+1}+\sum_{i=0}^{n}c_i x^i \in \mathbb{C}[x],
\]
there is a spectral sequence $\{E_k\}$ converging to $H_p$ with $E_1=H_n$. For each $k\geq1$, the $E_k$ term of this spectral sequence with its two gradings is invariant under Reidemeister moves.
\end{theorem}

The spectral sequence $\{E_k\}$ behaves well under link cobodisms.

\begin{theorem}\label{cobo-homo}
For $k\geq1$, a link cobodism $S$ from $L_0$ to $L_1$ induces, up to multiplication by a non-zero scalar, a homomorphism 
\[
\Psi_{S,k}: E_k(L_0) \rightarrow E_k(L_1)
\]
that shifts the first $\zed$-grading of $E_k^{\ast,\ast}$ by $-(n-1)\chi(S)$.
\end{theorem}

It's clear from the construction that the quantum filtration $\mathcal{F}$ on $H_p(L)$ is bounded from above and below for any link $L$. We define 
\begin{eqnarray*}
g^{max}_p(L) & = & \min\{k~|~\mathcal{F}^kH_p(L)=H_p(L)\}, \\
g^{min}_p(L) & = & \min\{k~|~\fil^k H_p(L)\neq 0\},
\end{eqnarray*}
i.e. the highest and lowest filtration levels of $H_p(L)$. These are numerical link invariants. In view of Theorem \ref{spec}, we have the following easy corollary:

\begin{corollary}\label{dominate}
$g^{max}_{x^{n+1}}(L)\geq g^{max}_p(L)$ for any link $L$ and any monic polynomial $p(x)$ of degree $n+1$.
\end{corollary}
\begin{proof}
Let $\{E_k\}$ be the spectral sequence in Theorem \ref{spec}. Note that $g^{max}_{x^{n+1}}(L) = \max\{i~|~E_0^{i\ast}(L)\neq0\}$ and $g^{max}_p(L) = \max\{i~|~E_\infty^{i\ast}(L)\neq0\}$. This implies the corollary.
\end{proof}

We also generalize the Rasmussen invariant $s$ defined in \cite{Ras1}.

\begin{definition}
For any knot $K$, and $n\geq2$, define the $n$-th order Rasmussen invariant $s_n$ by $s_n(K)=\frac{1}{2}(g^{max}_p(K)+g^{min}_p(K))$, where $p(x)=x^{n+1}-(n+1)x$.
\end{definition}

By the definitions of $s_n$ and $s$, it is easy to see that $s_2(K)=-s(\overline{K})$ for any knot $K$, where $\overline{K}$ is the mirror image of $K$. By generalizing the arguments in \cite{Ras1}, we prove:

\begin{theorem}\label{ras-genus}
Let $p(x)=x^{n+1}-(n+1)x$. Then, for any link $L$ in $S^3$,
\[
g^{max}_p(L) \geq (n-1)\chi_s(L),
\]
where $\chi_s(L)$ is the slice Euler characteristic of $L$, i.e. the maximal Euler characteristic, $\chi(S)$, of an oriented smoothly embedded compact surface $F$ without closed components in $D^4$ bounded by $L$.

Specially, $|s_n(K)|\leq 2(n-1)g_s(K)$ for any knot $K$ in $S^3$, where $g_s(K)$ is the slice genus of $K$.
\end{theorem}

\begin{remark}
After an earlier version of this paper was posted on arxiv.org, the author learned that Andrew Lobb \cite{Lobb} gave a different proof of Theorem \ref{ras-genus} without using the invariance of $H_p$.
\end{remark}

It is proved in \cite{Wu} that, for any transversal link $L$ in the standard contact $S^3$, 
\[
sl(L)\leq -\limsup_{n\rightarrow\infty}\frac{g^{max}_{x^{n+1}}(\overline{L})}{n-1},
\]
where $sl$ is the self-linking number, and $\overline{L}$ is the mirror image of $L$. We call the right hand side the asymptotic Khovanov-Rozansky bound of the self-linking number. By Corollary \ref{dominate} and Theorem \ref{ras-genus}, we can see that it is sharper that the upper bound of $sl(L)$ given by $-\chi_s(L)$ established by Rudolph in \cite{Ru}. Also, as explained in \cite{Wu}, the asymptotic Khovanov-Rozansky bound is sharper than the upper bound $e_P(L)$ of the self-linking number established by Franks and Williams \cite{FW} and Morton \cite{Mo} using the HOMFLY polynomial $P(L)$. These discussions plus Ferrand's results from Section 7 of \cite{Fer} give us:

\begin{proposition}
The asymptotic Khovanov-Rozansky bound of the self-linking number is sharper than both of the other two upper bounds, and can be arbitrarily sharper than either of them. That is, for any $N>0$, there are links $L_1$ and $L_2$, such that
\begin{eqnarray*}
-\chi_s(L_1)+ \limsup_{n\rightarrow\infty}\frac{g^{max}_{x^{n+1}}(\overline{L}_1)}{n-1}& > & N \\
e_P(L_2)+ \limsup_{n\rightarrow\infty}\frac{g^{max}_{x^{n+1}}(\overline{L}_2)}{n-1}& > & N
\end{eqnarray*}
\end{proposition}
\begin{proof}
Straightforward.
\end{proof}

The proofs of Theorems \ref{fil-inv}, \ref{spec} and \ref{cobo-homo} are based on generalizations of methods developed by Khovanov and Rozansky in \cite{KR1,KR2}. We assume the readers are somewhat familiar with these papers. Some understanding of Rasmussen's work in \cite{Ras2} would be useful too. In the treatment of slice Euler characteristic and the generalized Rasmussen invariants, we will use techniques developed by Gornik \cite{Gornik} and Rasmussen \cite{Ras1}.

The rest of this paper is organized as following. In section \ref{kr-def}, we recall the construction of the Khovanov-Rozansky cohomology. Then we establish the MOY decompositions in section \ref{moy}. Theorems \ref{fil-inv}, \ref{spec} and \ref{cobo-homo} are proved in section \ref{reidemeister}. Finally, in section \ref{rasmussen}, we discuss the generalized Rasmussen invariants and prove Theorem \ref{ras-genus}. 

\begin{acknowledgments}
The author would like to thank Mikhail Khovanov for many helpful discussions and Bojan Gornik for answering some questions about \cite{Gornik}. He would like to thank Michael G. Sullivan for his help on spectral sequences. He would also like to thank Marco Mackaay for helpful comments on an earlier version of this paper.
\end{acknowledgments}

\section{The Khovanov-Rozansky Cohomology}\label{kr-def}

\subsection{Filtered modules over a graded commutative $\C$-algebra}
Let $V$ be a $\C$-linear space. A filtration $\fil$ on $V$ is a sequence of subspaces $\{\fil^i V\}_{i\in \zed}$ such that 
\[
\cdots\subset\fil^i V\subset\fil^{i+1} V \subset \fil^{i+2} V \subset\cdots,
\]
\[
\bigcup_{i\in\zed}\fil^i V=V, ~\text{and}~\bigcap_{i\in\zed}\fil^i V=0.
\]
We call $(V,\fil)$ a filtered $\C$-linear space, and, when there is no danger of confusion, we drop $\fil$ from the notation.

\begin{definition}
The filtration $\fil$ is bounded from below if $\exists i\in\zed$ such that $\fil^i V=0$, and bounded from above if $\exists j\in\zed$ such that $\fil^j V=V$. 
\end{definition}

\begin{definition}
Let $(V, \fil)$ be a finite dimensional filtered $\C$-linear space. Then the filtered dimension of $V$ is defined to be the Laurent polynomial
\[
\sum_{k=-\infty}^{\infty} q^k \dim {(\fil^kV/\fil^{k-1}V)}.
\]
\end{definition}

\begin{definition}
A commutative $\C$-algebra $R$ is said to be graded if there is a decomposition of $\C$-linear space
\[
R = \bigoplus_{i=0}^{\infty} R_i,
\]
such that $xy\in R_{i+j}$ whenever $x\in R_i$ and $y\in R_j$. An element of $R_j$ is called a homogeneous element of degree $j$. For a non-zero element $p\in R$, $\deg p =j$ if $p \in \bigoplus_{i=0}^{j} R_i$ and $p \notin \bigoplus_{i=0}^{j-1} R_i$. Define $\deg 0 =-\infty$.

Unless otherwise specified, all $\C$-algebras in the rest of this section will be commutative and graded.
\end{definition}

A filtered module over a graded commutative $\C$-algebra $R$ is a pair $(M,\fil)$, where $M$ is an $R$-module, and $\fil$ is a filtration of the underlying $\C$-space of $M$ such that $p\cdot\fil^i M\subset\fil^{i+\deg{p}}M$ for any $p\in R$. When there is no danger of confusion, we drop $\fil$ from the notation. For an element $x$ of $M$, we say $\deg{x}=i$ if $x\in\fil^iM$, but $x\notin\fil^{i-1}M$. ( Again, our convention is that $\deg0=-\infty$.) Clearly, the grading of $R$ induces a filtration on $R$, which makes it into a filtered $R$-module.

If $M$ and $N$ are filtered $R$-modules, then $M\oplus N$, $M\otimes N$ and $\Hom_R(M,N)$ are also filtered with filtrations given by
\begin{eqnarray*}
\fil^i(M\oplus N) & = & \fil^iM\oplus\fil^iN \subset M\oplus N\\
\fil^i(M\otimes N) & = & \sum_{j+k=i}\fil^jM\otimes\fil^kN \subset M\otimes N\\
\fil^i\Hom_R(M,N) & = & \{f\in\Hom_R(M,N)~|~f(\fil^kM)\subset\fil^{k+i}M\}
\end{eqnarray*}

\begin{definition}
A homomorphism $f:M\rightarrow N$ of $R$-modules is called a filtered homomorphism if $f\in\fil^0\Hom_R(M,N)$. More generally, we say $\deg{f}=i$ if $f\in\fil^i\Hom_R(M,N)$, but $f\notin\fil^{i-1}\Hom_R(M,N)$. (Again, $\deg0=-\infty$.) 

$f:M\rightarrow N$ is an isomorphism of filtered modules if $f$ is an isomorphism of the underlying $R$-modules, and both $f$ and $f^{-1}$ are filtered maps.
\end{definition}

We also need a notion of filtration shift.

\begin{definition}
Let $M$ be a filtered $R$-module. Then $M\{k\}$ is the filtered $R$-module obtained by shifting the filtration by $k$, i.e. $\fil^{i+k}(M\{k\})=\fil^iM$.
\end{definition}

Following are some technical results needed for later use.

\begin{definition}\label{special-basis}
Let $M$ be a finitely generated free filtered $R$-module. $(\alpha_1,\cdots,\alpha_m)$ is called a special basis for $M$ if

\begin{enumerate}
	\item $(\alpha_1,\cdots,\alpha_m)$ is an $R$-basis of the underlying $R$-module $M$;
	
	\item The canonical isomorphism
	\[
	M \cong \bigoplus_{i=1}^m R\cdot\alpha_i
	\]
	is an isomorphism of filtered $R$-modules, where the filtration on $R\cdot\alpha_i$ is given by 
	\[
	\fil^{j}(R\cdot\alpha_i) = (R\cdot\alpha_i)\cap \fil^j M.
	\]
\end{enumerate}
\end{definition}

\begin{lemma}\label{special-sum}
Assume that $M$ is a finitely generated free filtered $R$-module, and $(\alpha_1,\cdots,\alpha_m)$ is a special basis for $M$.
\begin{enumerate}
	\item For any $i\in\zed$ and $f_1,\cdots,f_m\in R$, $\sum_{k=1}^{m}f_k\alpha_k\in\fil^iM$ if and only if $f_k\alpha_k\in\fil^iM$, $k=1,\cdots,m$;
	\item For any $g\in M^{\ast}=\Hom_R(M,R)$, $\deg{(g)}=\max\{\deg g(\alpha_k)-\deg \alpha_k~|~k=1,\cdots,m\}$;
	\item If $N$ is also a finitely generated free filtered $R$-module with a special basis, then $M\oplus N$ and $M\otimes_R N$ have special bases too.
\end{enumerate}
\end{lemma}
\begin{proof}
Straightforward.
\end{proof}

\begin{definition}\label{faithful}
Let $(M,\fil)$ be a filtered $R$-module. $\fil$ is called faithful if 
\[
\deg{(r\cdot m)}=\deg{r}+\deg{m},
\] 
for any non-zero elements $r\in R,~m\in M$. In this case, $M$ is said to be a faithfully filtered $R$-module.
\end{definition}

\begin{lemma}\label{faithful-sum}
If $M$ and $N$ are both faithfully filtered $R$-modules, then so is $M\oplus N$.
\end{lemma}
\begin{proof}
Straightforward.
\end{proof}

\begin{proposition}\label{hom-tensor}
If $M$ and $N$ are finitely generated free faithfully filtered $R$-modules with spacial bases, then the natural isomorphism $\Hom_R(M,N)\xrightarrow{\cong} N\otimes M^\ast$ is an isomorphism of filtered modules.
\end{proposition}
\begin{proof}
Let $F:\Hom_R(M,N)\xrightarrow{\cong} N\otimes M^\ast$ be the natural isomorphism. Then $F^{-1}:N\otimes M^\ast\rightarrow\Hom_R(M,N)$ is given by 
\[
F^{-1}(n\otimes f)(m)=f(m)\cdot n, ~\forall~ n\in N,~ f\in M^\ast, ~m\in M.
\] 
It's easy to see that $F^{-1}$ is a filtered homomorphism. It remains to check that $F$ is also filtered. 

Let $(\alpha_1,\cdots,\alpha_m)$ and $(\beta_1,\cdots,\beta_n)$ be special basis for $M$ and $N$, and $g:M\rightarrow N$ a $R$-homomorphism of degree $d$. Under these basis, we have
\[
g(\alpha_i)=\sum_{j=1}^{n}g_{ij}\beta_j, ~i=1,\cdots,m.
\]
Define $f_{ij}:M\rightarrow R$ by 
\[
f_{ij}(\alpha_k)=
\left\{
  \begin{array}{l}
    g_{ij}, \text{ if } k=i, \\
    0, \text{ otherwise.}
  \end{array}
\right.
\]
Then 
\[
F(g)=\sum_{i=1}^m\sum_{j=1}^n\beta_j\otimes f_{ij}.
\]
Note that 
\[
\sum_{j=1}^{n}g_{ij}\beta_j=g(\alpha_i)\in\fil^{d+\deg{\alpha_i}}N.
\]
So 
\[
g_{ij}\beta_j\in\fil^{d+\deg{\alpha_i}}N, ~\forall~i,j,
\] 
and 
\[
\deg{g_{ij}}=\deg{g_{ij}\beta_j}-\deg{\beta_j}\leq d+\deg{\alpha_i}-\deg{\beta_j}.
\]
Then
\begin{eqnarray*}
\deg{f_{ij}} & = & \max\{\deg{f_{ij}(\alpha_k)}-\deg{\alpha_k}~|~k=1,\cdots,m\} \\
& = & \deg{g_{ij}}-\deg{\alpha_i}\leq d-\deg{\beta_j}.
\end{eqnarray*}
Thus, $\deg{(\beta_j\otimes f_{ij})}\leq d$, $\forall~i,j$, and, therefore, $\deg{F(g)}\leq d$. This shows $F$ is also filtered.
\end{proof}

\subsection{Filtered matrix factorizations}
Again, let $R$ be a graded commutative $\C$-algebra. Fix an integer $n\geq2$. Let $w$ be an element of $R$ with $\deg{w}\leq2n+2$. A filtered matrix factorization over $R$ with potential $w$ is a collection of two filtered free $R$-modules $M^0$, $M^1$ and two $R$-module homomorphisms $d^0:M^0\rightarrow M^1$, $d^1:M^1\rightarrow M^0$, called differential maps, s.t.,
\[
d^0\in\fil^{n+1}\Hom_R(M^0,M^1), \hspace{1cm} d^1\in\fil^{n+1}\Hom_R(M^1,M^0)
\]
\[
d^1d^0=w\cdot\id_{M^0}, \hspace{1cm} d^0d^1=w\cdot\id_{M^1}.
\]
We usually write such a matrix factorization $M$ as
\[
M^0 \xrightarrow{d_0} M^1 \xrightarrow{d_1} M^0.
\]
Following \cite{KR1}, we denote by $M\left\langle1\right\rangle$ the filtered matrix factorization
\[
M^1 \xrightarrow{-d_1} M^0 \xrightarrow{-d_0} M^1,
\]
by $M^{\ast}$ the filtered matrix factorization
\[
(M^0)^{\ast} \xrightarrow{(d_1)^{\ast}} (M^1)^{\ast} \xrightarrow{(d_0)^{\ast}} (M^0)^{\ast},
\]
by $M_-$ the filtered matrix factorization
\[
M^0 \xrightarrow{-d_0} M^1 \xrightarrow{d_1} M^0,
\]
and by $M_{\bullet}$ the filtered matrix factorization $(M^{\ast})_-$. Note that $M\left\langle1\right\rangle$, $M^{\ast}$ have potential $w$, and $M_-$, $M_{\bullet}$ have potential $-w$.

For filtered matrix factorizations $M$ with potential $w_1$ and $N$ with potential $w_2$, the tensor product $M\otimes N$ is the filtered matrix factorization with 
\begin{eqnarray*}
(M\otimes N)^0 & = & (M^0\otimes N^0)\oplus (M^1\otimes N^1), \\
(M\otimes N)^1 & = & (M^1\otimes N^0)\oplus (M^1\otimes N^0),
\end{eqnarray*}
and the differential given by signed Leibniz rule. It has potential $w_1+w_2$.

If filtered matrix factorizations $M$ and $N$ both have potential $w$, then the free $R$-module $\Hom_R(M,N)$ is a $\zed_2$-complex
\[
\Hom_R^0(M,N)\xrightarrow{d}\Hom_R^1(M,N)\xrightarrow{d}\Hom_R^0(M,N),
\]
where 
\begin{eqnarray*}
\Hom_R^0(M,N) & = & \Hom_R(M^0,N^0) \oplus \Hom_R(M^1,N^1), \\
\Hom_R^1(M,N) & = & \Hom_R(M^0,N^1) \oplus \Hom_R(M^1,N^0),
\end{eqnarray*}
and
\[
d(f)(m)=d_N(f(m))-(-1)^if(d_M(m)), \text{ for }f\in \Hom_R^i(M,N).
\]
Note that $\Hom_R(M,N)$ comes with a natural filtration, and $d\in\fil^{n+1}\Hom_R(\Hom_R(M,N),\Hom_R(M,N))$. For such $M$ and $N$, a homomorphism of matrix factorization $f:M\rightarrow N$ is a cycle $f\in\Hom_R^0(M,N)=\Hom_R(M^0,N^0) \oplus \Hom_R(M^1,N^1)$, i.e. $df=0$. We say that two homomorphisms $f$ and $g$ are homotopic, denoted by $f\sim g$, if $\exists h\in\Hom_R^1(M,N)$ such that $f-g=dh$. Denote by $\Hom_{MF}(M,N)$ the $R$-module of homomorphisms of matrix factorization from $M$ to $N$. Define $\Hom_{HMF}(M,N)=\Hom_{MF}(M,N)/\sim$, i.e. the $R$-module of homotopy classes of matrix factorization homomorphisms from $M$ to $N$. Clearly, 
\[
\Hom_{HMF}(M,N)=H^0(\Hom_R(M,N),d).
\]
Note that the filtration on $\Hom_R^0(M,N)$ descends onto $\Hom_{MF}(M,N)$ and $\Hom_{HMF}(M,N)$. We  define
\[
\Hom_{mf}(M,N)=\fil^0\Hom_{MF}(M,N),
\]
and 
\[
\Hom_{hmf}(M,N)=\fil^0\Hom_{HMF}(M,N).
\]
An element of $\Hom_{mf}(M,N)$ is called a filtered homomorphism of matrix factorizations. An element of $\Hom_{hmf}(M,N)$ is a homotopy class containing a filtered homomorphism. A map $f:M\rightarrow N$ is said to be an isomorphism of filtered matrix factorizations if $f$ has an inverse $f^{-1}$, and both $f$ and $f^{-1}$ are filtered homomorphisms of matrix factorizations.

Let $a$ and $b$ be elements of $R$ with $a\neq0$ and $\deg{(ab)}\leq2n+2$. Denote by $(a,b)_R$ the filtered matrix factorization 
\[
R \xrightarrow{a} R\{n+1-\deg{a}\} \xrightarrow{b} R,
\]
where $a,~b$ act on $R$ by multiplication. This matrix factorization has potential $ab$. When the graded $\C$-algebra $R$ is clear from context, we drop it from the notation. For $a_1\cdots,a_k,b_1,\cdots,b_k\in R$ with $a_i\neq0$ and $\deg{a_ib_i}\leq2n+2$, $i=1,\cdots,k$, denote by
\[
\left(%
\begin{array}{cc}
  a_1 & b_1 \\
  a_2 & b_2 \\
  \vdots & \vdots \\
  a_k & b_k \\
\end{array}%
\right)_R
\]
the tensor product of
$(a_1,b_1)_R,~(a_2,b_2)_R,\cdots,~(a_k,b_k)_R$. This is a filtered matrix factorization with potential
\[
w=a_1b_1+a_2b_2+\cdots+a_kb_k.
\]
Again, when $R$ is clear from context, we drop it from the notation. It's straightforward to check that
\[
\left(%
\begin{array}{cc}
  a_1 & b_1 \\
  a_2 & b_2 \\
  \vdots & \vdots \\
  a_k & b_k \\
\end{array}%
\right)_{\bullet}
=
\left(%
\begin{array}{cc}
  a_1 & -b_1 \\
  a_2 & -b_2 \\
  \vdots & \vdots \\
  a_k & -b_k \\
\end{array}%
\right)\left\langle k\right\rangle\{-k(n+1)+\sum_{i=1}^k\deg a_i\},
\]
which has potential $-w$, where $\left\langle k\right\rangle$ means $\underbrace{\left\langle 1\right\rangle\left\langle 1\right\rangle\cdots\left\langle 1\right\rangle}_{k \text{ times}}$.

\begin{proposition}\label{hmf-tensor}
Suppose $M=\left(%
\begin{array}{cc}
  a_1 & b_1 \\
  a_2 & b_2 \\
  \vdots & \vdots \\
  a_k & b_k \\
\end{array}%
\right)$ and
$N=\left(%
\begin{array}{cc}
  \alpha_1 & \beta_1 \\
  \alpha_2 & \beta_2 \\
  \vdots & \vdots \\
  \alpha_l & \beta_l \\
\end{array}%
\right)\left\langle \varepsilon \right\rangle$ are two filtered matrix factorizations with potential $w$ constructed as above, where $\varepsilon=0,1$. Then there is an isomorphism of filtered matrix factorization of potential $0$ 
\[
\Hom_R(M,N) \xrightarrow{\cong} N\otimes M_{\bullet}.
\]
Specially, 
\[
\Hom_{HMF}(M,N)=H^0(N\otimes M_{\bullet}).
\]
\end{proposition}
\begin{proof}
First, note that $N\otimes M_{\bullet}$ is a filtered matrix factorization with potential $w-w=0$. So it is indeed a $\zed_2$-complex. From the construction of $M$ and $N$, we know that each of them admits a direct sum decomposition in which every factor is a copy of $R$ with some filtration shift. Then by Lemmas \ref{special-sum} and \ref{faithful-sum}, the filtrations on $M$ and $N$ are faithful, and $M$ and $N$ admit special bases. So the natural isomorphism $F:\Hom_R(M,N) \xrightarrow{\cong} N\otimes M_{\bullet}$ is an isomorphism of filtered $R$-modules. It's easy to check that $F$ commutes with the differential maps of these two chain complexes. This proves the proposition.
\end{proof}

\begin{lemma}\label{negative-differential}
Let $M$ be the filtered matrix factorization $M^0 \xrightarrow{d_0} M^1 \xrightarrow{d_1} M^0$ over $R$ with potential $w$. For $c\in\C\setminus\{0\}$, define $M_c$ to be the matrix factorization $M^0 \xrightarrow{cd_0} M^1 \xrightarrow{c^{-1}d_1} M^0$. Then the potential of $M_c$ is $w$, and there is an isomorphism of filtered matrix factorizations $M\xrightarrow{\cong}M_c$.
\end{lemma}
\begin{proof}
Define $f:M\rightarrow M_c$ and $g:M_c\rightarrow M$ by the following diagram. It is clear that $f$ and $g$ are filtered homomorphisms and are inverses of each other.
\[
\begin{CD}
M^0 @>{d_0}>> M^1 @>{d_1}>> M^0\\
@VV{c^{-1}\id}V @VV{\id}V @VV{c^{-1}\id}V\\
M^0 @>{cd_0}>> M^1 @>{c^{-1}d_1}>> M^0\\ 
@VV{c\id}V @VV{\id}V @VV{c\id}V \\
M^0 @>{d_0}>> M^1 @>{d_1}>> M^0
\end{CD}
\]
\end{proof}

The next lemma is a slight generalization of Lemma 3.5 of \cite{Ras2}

\begin{lemma}\label{general-twist}
Let $M$ be the filtered matrix factorization $M^0 \xrightarrow{d_0} M^1 \xrightarrow{d_1} M^0$ over $R$ with potential $w$, and $H_i:M^i\rightarrow M^i$, $i\in\zed_2$, filtered homomorphisms with $H_i^2=0$. Define $\tilde{d}_i:M^i\rightarrow M^{i+1}$, $i\in\zed_2$, by 
\[
\tilde{d}_i= (\id_{M^{i+1}}-H_{i+1})\circ d_i \circ (\id_{M^{i}}+H_{i}),
\]
and $\widetilde{M}$ by 
\[
M^0 \xrightarrow{\tilde{d}_0} M^1 \xrightarrow{\tilde{d}_1} M^0.
\]
Then $\widetilde{M}$ is also a filtered matrix factorization over $R$ with potential $w$. And there is an isomorphism of filtered matrix factorizations
$M\xrightarrow{\cong}\widetilde{M}$.
\end{lemma}
\begin{proof}
Since $H_i$, $i\in\zed_2$, are filtered maps. It's clear that $\tilde{d}_i\in\fil^{n+1}\Hom_R(M^i,M^{i+1})$. Using $H_i^2=0$, it's easy to check that $\tilde{d}^1\tilde{d}^0=w\cdot\id_{M^0}$ and $\tilde{d}^0\tilde{d}^1=w\cdot\id_{M^1}$. So $\widetilde{M}$ is a filtered matrix factorization over $R$ with potential $w$. Consider the following diagram.

\[
\begin{CD}
M^0 @>{d_0}>> M^1 @>{d_1}>> M^0\\
@VV{\id - H_0}V @VV{\id-H_1}V @VV{\id-H_0}V\\
M^0 @>{\tilde{d}_0}>> M^1 @>{\tilde{d}_1}>> M^0\\ 
@VV{\id + H_0}V @VV{\id+H_1}V @VV{\id+H_0}V \\
M^0 @>{d_0}>> M^1 @>{d_1}>> M^0
\end{CD}
\]

Note that each square in the above diagram is commutative, and every vertical homomorphism is filtered. This gives us a pair of filtered homomorphisms \[
F:M\rightarrow\widetilde{M}, \hspace{1cm} G:\widetilde{M}\rightarrow M.
\]
Again, using $H_i^2=0$, it's easy to check that $G\circ F=\id_M$ and $F\circ G=\id_{\widetilde{M}}$. This proves the lemma.
\end{proof}

\begin{corollary}[Twisting]\label{twist}
Suppose $a_1,a_2,b_1,b_2,k\in R$ satisfy $a_i\neq0$, $\deg a_ib_i\leq 2n+2$, $i=1,2$, and $\deg k \leq \deg a_1 +\deg a_2 -2n-2$. Then there is an isomorphism of filtered matrix factorizations
\[
\left(%
\begin{array}{cc}
  a_1 & b_1 \\
  a_2 & b_2 
\end{array}%
\right)_R
\xrightarrow{\cong}
\left(%
\begin{array}{cc}
  a_1+kb_2 & b_1 \\
  a_2-kb_1 & b_2 
\end{array}%
\right)_R.
\]
\end{corollary}
\begin{proof}
The filtered matrix factorization 
\[
\left(%
\begin{array}{cc}
  a_1 & b_1 \\
  a_2 & b_2 
\end{array}%
\right)
\]
is given by
{\tiny
\[
\left.%
\begin{array}{c}
  R \\
  \oplus \\
  R\{2n+2-\deg a_1 -\deg a_2\}
\end{array}%
\right.
\xrightarrow{\left(%
\begin{array}{cc}
  a_1 & b_2 \\
  a_2 & -b_1 
\end{array}%
\right)}
\left.%
\begin{array}{c}
  R\{n+1-\deg a_1\} \\
  \oplus \\
  R\{n+1 -\deg a_2\}
\end{array}%
\right.
\xrightarrow{\left(%
\begin{array}{cc}
  b_1 & b_2 \\
  a_2 & -a_1 
\end{array}%
\right)}
\left.%
\begin{array}{c}
  R \\
  \oplus \\
  R\{2n+2-\deg a_1 -\deg a_2\}
\end{array}%
\right..
\]
}

The filtered matrix factorization 
\[
\left(%
\begin{array}{cc}
  a_1+kb_2 & b_1 \\
  a_2-kb_1 & b_2 
\end{array}%
\right)
\]
is given by
{\tiny
\[
\left.%
\begin{array}{c}
  R \\
  \oplus \\
  R\{2n+2-\deg a_1 -\deg a_2\}
\end{array}%
\right.
\xrightarrow{\left(%
\begin{array}{cc}
  a_1+kb_2 & b_2 \\
  a_2-kb_1 & -b_1 
\end{array}%
\right)}
\left.%
\begin{array}{c}
  R\{n+1-\deg a_1\} \\
  \oplus \\
  R\{n+1 -\deg a_2\}
\end{array}%
\right.
\xrightarrow{\left(%
\begin{array}{cc}
  b_1 & b_2 \\
  a_2-kb_1 & -a_1-kb_2
\end{array}%
\right)}
\left.%
\begin{array}{c}
  R \\
  \oplus \\
  R\{2n+2-\deg a_1 -\deg a_2\}
\end{array}%
\right..
\]
}
Define 
\[
H_0: \left.%
\begin{array}{c}
  R \\
  \oplus \\
  R\{2n+2-\deg a_1 -\deg a_2\}
\end{array}%
\right.
\rightarrow
\left.%
\begin{array}{c}
  R \\
  \oplus \\
  R\{2n+2-\deg a_1 -\deg a_2\}
\end{array}%
\right.
\]
and
\[
H_1:\left.%
\begin{array}{c}
  R\{n+1-\deg a_1\} \\
  \oplus \\
  R\{n+1 -\deg a_2\}
\end{array}%
\right.
\rightarrow
\left.%
\begin{array}{c}
  R\{n+1-\deg a_1\} \\
  \oplus \\
  R\{n+1 -\deg a_2\}
\end{array}%
\right.
\]
by
\[
H_0=\left(%
\begin{array}{cc}
  0 & 0 \\
  k & 0
\end{array}%
\right),
\hspace{1cm}
H_1=\left(%
\begin{array}{cc}
  0 & 0 \\
  0 & 0
\end{array}%
\right).
\]
Apply Lemma \ref{general-twist} to $M=\left(%
\begin{array}{cc}
  a_1 & b_1 \\
  a_2 & b_2 
\end{array}%
\right)$ and $H_1,~H_2$. Then the corollary follows.
\end{proof}

\begin{corollary}[Row operation]\label{row-op} 
Suppose $a_1,a_2,b_1,b_2,c\in R$ satisfy $a_i\neq0$, $\deg a_ib_i\leq 2n+2$, $i=1,2$, and $\deg c \leq \deg a_1 -\deg a_2$. Then there is an isomorphism of filtered matrix factorizations
\[
\left(%
\begin{array}{cc}
  a_1 & b_1 \\
  a_2 & b_2 
\end{array}%
\right)_R
\xrightarrow{\cong}
\left(%
\begin{array}{cc}
  a_1+ca_2 & b_1 \\
  a_2 & b_2-cb_1 
\end{array}%
\right)_R.
\]
\end{corollary}
\begin{proof}
The filtered matrix factorization 
\[
\left(%
\begin{array}{cc}
  a_1 & b_1 \\
  a_2 & b_2 
\end{array}%
\right)
\]
is given by
{\tiny
\[
\left.%
\begin{array}{c}
  R \\
  \oplus \\
  R\{2n+2-\deg a_1 -\deg a_2\}
\end{array}%
\right.
\xrightarrow{\left(%
\begin{array}{cc}
  a_1 & b_2 \\
  a_2 & -b_1 
\end{array}%
\right)}
\left.%
\begin{array}{c}
  R\{n+1-\deg a_1\} \\
  \oplus \\
  R\{n+1 -\deg a_2\}
\end{array}%
\right.
\xrightarrow{\left(%
\begin{array}{cc}
  b_1 & b_2 \\
  a_2 & -a_1 
\end{array}%
\right)}
\left.%
\begin{array}{c}
  R \\
  \oplus \\
  R\{2n+2-\deg a_1 -\deg a_2\}
\end{array}%
\right..
\]
}

The filtered matrix factorization 
\[
\left(%
\begin{array}{cc}
  a_1+ca_2 & b_1 \\
  a_2 & b_2-cb_1 
\end{array}%
\right)
\]
is given by
{\tiny
\[
\left.%
\begin{array}{c}
  R \\
  \oplus \\
  R\{2n+2-\deg a_1 -\deg a_2\}
\end{array}%
\right.
\xrightarrow{\left(%
\begin{array}{cc}
  a_1+ca_2 & b_2-cb_1 \\
  a_2 & -b_1 
\end{array}%
\right)}
\left.%
\begin{array}{c}
  R\{n+1-\deg a_1\} \\
  \oplus \\
  R\{n+1 -\deg a_2\}
\end{array}%
\right.
\xrightarrow{\left(%
\begin{array}{cc}
  b_1 & b_2-cb_1 \\
  a_2 & -a_1-ca_2
\end{array}%
\right)}
\left.%
\begin{array}{c}
  R \\
  \oplus \\
  R\{2n+2-\deg a_1 -\deg a_2\}
\end{array}%
\right..
\]
}
Define 
\[
H_0: \left.%
\begin{array}{c}
  R \\
  \oplus \\
  R\{2n+2-\deg a_1 -\deg a_2\}
\end{array}%
\right.
\rightarrow
\left.%
\begin{array}{c}
  R \\
  \oplus \\
  R\{2n+2-\deg a_1 -\deg a_2\}
\end{array}%
\right.
\]
and
\[
H_1:\left.%
\begin{array}{c}
  R\{n+1-\deg a_1\} \\
  \oplus \\
  R\{n+1 -\deg a_2\}
\end{array}%
\right.
\rightarrow
\left.%
\begin{array}{c}
  R\{n+1-\deg a_1\} \\
  \oplus \\
  R\{n+1 -\deg a_2\}
\end{array}%
\right.
\]
by
\[
H_0=\left(%
\begin{array}{cc}
  0 & 0 \\
  0 & 0
\end{array}%
\right),
\hspace{1cm}
H_1=\left(%
\begin{array}{cc}
  0 & -c \\
  0 & 0
\end{array}%
\right).
\]
Apply Lemma \ref{general-twist} to $M=\left(%
\begin{array}{cc}
  a_1 & b_1 \\
  a_2 & b_2 
\end{array}%
\right)$ and $H_1,~H_2$. Then the corollary follows.
\end{proof}

Let $R$ be a graded commutative $\C$-algebra with the grading given by 
\[
R = \bigoplus_{i=0}^{\infty} R_i.
\]
Assume that $b$ is a homogeneous element of $R$. Then $R'=R/bR$ is also a graded $\C$-algebra with grading given by
\[
R' = \bigoplus_{i=0}^{\infty} R_i/(R_i\cap bR).
\]
The following is the filtered version of Proposition 9 of \cite{KR1}.

\begin{proposition}\label{variable-exclusion}
Let $R$ be a graded $\C$-algebra, and $b$ a homogeneous element of $R$ which is not a zero divisor. Define $R'=R/bR$ and denote by $P:R\rightarrow R'$ the canonical projection. Let 
\[
M=\left(%
\begin{array}{cc}
  a_1 & b_1 \\
  a_2 & b_2 \\
  \vdots & \vdots \\
  a_k & b_k
\end{array}%
\right)_R
\]
be a matrix factorization of potential $0$ such that $\exists ~i$ with $1\leq i \leq k$, $b_i=b$ and $\deg a_i+\deg b =2n+2$. Then 
\[
M'=\left(%
\begin{array}{cc}
  P(a_1) & P(b_1) \\
  P(a_2) & P(b_2) \\
  \vdots & \vdots \\
  P(a_{i-1}) & P(b_{i-1}) \\
  P(a_{i+1}) & P(b_{i+1}) \\
  \vdots & \vdots \\
  P(a_k) & P(b_k)
\end{array}%
\right)_{R'}
\]
is also a matrix factorization of potential $0$. And there is an isomorphism of filtered $\C$-spaces
\[
H(M') \cong H(M),
\]
where $H(M')$ and $H(M)$ are the cohomologies of $M'$ and $M$ considered as $\zed_2$-complexes.
\end{proposition}

\begin{proof}
Let 
\[
\hat{M}=\left(%
\begin{array}{cc}
  a_1 & b_1 \\
  a_2 & b_2 \\
  \vdots & \vdots \\
  a_{i-1} & b_{i-1} \\
  a_{i+1} & b_{i+1} \\
  \vdots & \vdots \\
  a_k & b_k
\end{array}%
\right)_R
\]
Then
\begin{eqnarray*}
M & \cong & \hat{M}\otimes_R (a_i,b)_R \\
  & = & \hat{M}_0 \xrightarrow{a_i} \hat{M}_1\{n+1-\deg a_i\} \xrightarrow{b} \hat{M}_0,
\end{eqnarray*}
where $\hat{M}_0$ and $\hat{M}_1$ are two copies of $\hat{M}$. Clearly, $P:R\rightarrow R'$ is a graded homomorphism and induces a filtered projection $\hat{P}:\hat{M}\rightarrow M'$. Define $F:M\rightarrow M'$ by $F|_{\hat{M}_1\{n+1-\deg a_i\}}=0$ and $F|_{\hat{M}_0}=\hat{P}$. It is easy to check that $F$ is a filtered surjective chain map with kernel  
\[
\ker F =N= b\hat{M}_0 \xrightarrow{a_i} \hat{M}_1\{n+1-\deg a_i\} \xrightarrow{b} b\hat{M}_0,
\] 
which is a subcomplex of $M$.

Since $b$ is not a zero divisor, the homomorphism $R$-modules $\varphi_b:bR\rightarrow R$ given by $\varphi_b(br)=r$ is well defined. Denote by $1_j^{\epsilon}$ the element $1$ of the copy of $R$ of $\zed_2$-grading $\epsilon$ in the matrix factorization
\[
(a_j,b_j) ~=~ R \xrightarrow{a_j} R\{n+1-\deg a_j\} \xrightarrow{b_j} R.
\]
Define a homomorphism $h:N\rightarrow N$ of $R$-modules by 
\begin{eqnarray*}
h|_{\hat{M}_1\{n+1-\deg a_i\}} & = & 0, \\
h(br 1_1^{\epsilon_1} \otimes \cdots \otimes 1_{i-1}^{\epsilon_{i-1}} \otimes 1_i^{0} \otimes 1_{i+1}^{\epsilon_{i+1}} \otimes \cdots \otimes 1_k^{\epsilon_k}) & = & (-1)^{\sum_{j=1}^{i-1}\epsilon_j} r 1_1^{\epsilon_1} \otimes \cdots \otimes 1_{i-1}^{\epsilon_{i-1}} \otimes 1_i^{1} \otimes 1_{i+1}^{\epsilon_{i+1}} \otimes \cdots \otimes 1_k^{\epsilon_k}.
\end{eqnarray*}
Then $d|_N h+hd|_N=\id_N$, where $d$ is the differential map of $M$. In particular, this means that $H(N)=0$, which implies that $F$ induces an isomorphism $F_\ast:H(M)\rightarrow H(M')$. Since $F$ is filtered, so is $F_\ast$. To prove the proposition, we need to show that $F_\ast^{-1}$ is also filtered. For every cocycle $\alpha$ in $M'$, there is a cochain $\beta$ in $\hat{M}_0$ such that $F(\beta)=\hat{P}(\beta)=\alpha$ and $\deg\beta=\deg\alpha$. So $F(d \beta)=d' F(\beta)=d' \alpha=0$, where $d'$ is the differential map of $M'$. Therefore, $d \beta \in \ker F=N$, and $d \beta = d h(d \beta)+hd(d \beta)=d h(d \beta)$.  $h\in\fil^{-n-1}\Hom_R(N,N)$ since $\deg a_i+\deg b =2n+2$. But $d\in\fil^{n+1}\Hom_R(M,M)$. So $\deg h(d \beta) \leq \deg\beta$. Now we have $d (\beta-h(d \beta))=0$, $F(\beta-h(d \beta))=\alpha$ and $\deg (\beta-h(d \beta)) \leq \deg\beta=\deg\alpha$. Note that $F_\ast^{-1}([\alpha])=[\beta-h(d \beta)]$. This implies that $F_\ast^{-1}$ is also filtered.
\end{proof}

\begin{remark}\label{direct-sum-variable-exclusion}
Consider the special case of Proposition \ref{variable-exclusion} when 
\begin{enumerate}
	\item $R'$ is a graded subalgebra of $R$, and there exists a grading preserving direct sum decomposition $R=R'\oplus bR$;
	
	\item $a_1,\cdots,a_{i-1},a_{i+1},\cdots,a_k,b_1,\cdots,b_{i-1},b_{i+1},\cdots,b_k \in R'$.
\end{enumerate}
Then, in the above proof, $\hat{M}_0=M'\oplus b\hat{M}_0$ as filtered $R'$-modules. Let $J:M'\rightarrow \hat{M}_0$ be the inclusion map from this decomposition. We can choose that $\beta=J(\alpha)\in M_0'$. From (2), we have $d\beta =J(d'\alpha) + a_i\beta$, where  $a_i\beta\in\hat{M}_1\{n+1-\deg a_i\}$. But $d'\alpha=0$ since $\alpha$ is a cocycle in $M'$. So $h(d \beta)=h(a_i\beta)=0$. Therefore, $\beta-h(d \beta)=\beta\in \hat{M}_0$ is a cocycle in $M$, which is mapped to $\alpha$ by the projection $F$. This fact will be used in later computations.
\end{remark}

\subsection{Matrix factorizations associated to planer diagrams}\label{pd}
To construct the Khovanov-Rozansky cohomology, one
considers the planar diagrams $\Gamma$ with the following
properties:
\begin{enumerate}
    \item $\Gamma$ consists of two types of edges: regular edges
    and wide edges. These edges intersect only at their endpoints.
    \item Regular edges are disjoint from each other.
    Wide edges are disjoint from each other.
    \item Each regular edge is oriented, and contains at least one
    marked point. Open endpoints of regular edges are marked.
    \item Each wide edge has exactly two regular edges entering at one endpoint,
    and exactly two regular edges exiting from the other endpoint.
\end{enumerate}
Following \cite{Ras2}, we call such a planar diagram a graph.

In the rest of this section, we fix an integer $n\geq2$ and a polynomial
\[
p(x)=x^{n+1}+\sum_{i=0}^{n}c_i x^i \in \mathbb{C}[x].
\]
Let $\Gamma$ be a graph, and $\{x_1, \cdots, x_m\}$ the set of markings on $\Gamma$. Let $R=\C[x_1,\cdots,x_m]$ be the $\mathbb{Z}$-graded polynomial ring with $\deg x_i=2$. This grading induces a filtration $\fil$ on $R$, and makes it into a filtered $R$-module.

\begin{figure}[ht]

\setlength{\unitlength}{1pt}

\begin{picture}(420,20)(-210,0)

\linethickness{.5pt}

\put(-15,0){\vector(1,0){30}}

\put(-20,5){$x_j$}

\put(15,5){$x_i$}

\end{picture}

\caption{$L^i_j$}\label{lij}

\end{figure}

For an oriented regular arc $L^i_j$ in $\Gamma$ from the point
marked by $x_j$ to the point marked by $x_i$ with no marked
interior points, let $C_p(L^i_j)$ be the matrix factorization
$(\pi_{ij},x_i-x_j)_{R}$ given by
\[
R \xrightarrow{\pi_{ij}} R\{1-n\} \xrightarrow{x_i-x_j} R,
\]
where $\pi_{ij}=\frac{p(x_i)-p(x_j)}{x_i-x_j}$ is a polynomial of degree $2n$. The
purpose of the filtration shift here is to make $C_p(L^i_j)$ a filtered matrix factorization.

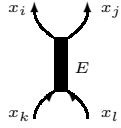
\begin{figure}[ht]

\setlength{\unitlength}{1pt}

\begin{picture}(420,45)(-210,0)

\linethickness{.5pt}

\put(-2.5,10){\vector(1,1){0}}

\qbezier(-10,0)(-10,5)(-2.5,10)

\qbezier(-10,40)(-10,35)(-2.5,30)

\put(-10,44){\vector(0,1){0}}

\put(2.5,10){\vector(-1,1){0}}

\qbezier(10,0)(10,5)(2.5,10)

\qbezier(10,40)(10,35)(2.5,30)

\put(10,44){\vector(0,1){0}}

\put(-20,40){\tiny{$x_i$}}

\put(15,40){\tiny{$x_j$}}

\put(-20,0){\tiny{$x_k$}}

\put(15,0){\tiny{$x_l$}}

\linethickness{5pt}

\put(0,10){\vector(0,1){20}}

\put(5,17){\tiny{$E$}}

\end{picture}

\caption{Wide edge $E$}\label{wideedge}

\end{figure}

For a wide edge $E$, let $x_i$, $x_j$, $x_k$, $x_l$ be the closest
markings to $E$ as depicted in Figure \ref{wideedge}. (It is
possible that $x_i=x_k$ or $x_j=x_l$.) Let $g$ be the unique
two-variable polynomial such that $g(x+y,xy)=p(x)+p(y)$, and
\[
u_{ijkl} = u(x_i,x_j,x_k,x_l) =
\frac{g(x_i+x_j,x_ix_j)-g(x_k+x_l,x_ix_j)}{x_i+x_j-x_k-x_l},
\]
\[
v_{ijkl} = v(x_i,x_j,x_k,x_l) =
\frac{g(x_k+x_l,x_ix_j)-g(x_k+x_l,x_kx_l)}{x_ix_j-x_kx_l}.
\]
Note that $u_{ijkl}$ and $v_{ijkl}$ are polynomials in $x_i$, $x_j$, $x_k$ and $x_l$ of degrees $2n$ and $2n-2$, respectively.

Define $C_p(E)$ to be the matrix factorization
\[
\left(%
\begin{array}{cc}
  u_{ijkl} & x_i+x_j-x_k-x_l \\
  v_{ijkl} & x_ix_j-x_kx_l \\
\end{array}%
\right)_{R}\{-1\},
\]
which is the tensor product of the matrix factorizations
\[
R\{-1\} \xrightarrow{u_{ijkl}} R\{-n\} \xrightarrow{x_i+x_j-x_k-x_l}
R\{-1\}
\]
and
\[
R \xrightarrow{v_{ijkl}} R\{3-n\} \xrightarrow{x_ix_j-x_kx_l} R.
\]
Again, the filtration shifts here make $C_p(E)$ a filtered matrix factorization.

We define that
\[
C_p(\Gamma) = (\bigotimes_{L^i_j}C_p(L^i_j)) ~\bigotimes~
(\bigotimes_{E}C_p(E)),
\]
where $L^i_j$ runs through all the regular arcs starting and ending at marked points with no marked interior points, and $E$ runs through all wide edges. Note that the potential of $C_p(\Gamma)$ is $w(\Gamma)=\sum\pm p(x)$, where $x$ runs through all the open endpoints, and the sign is "$+$" if the corresponding open endpoint is an exit, and is "$-$" if it's an entrance. 

When $\Gamma$ is closed, i.e. with no open endpoints, one can see that $w(\Gamma)=0$, and $C_p(\Gamma)$ is a $\zed_2$-chain complex. In this case, we denote by $H_p(\Gamma)$ the cohomology of this chain complex. Note that the filtration $\fil$ on $C_p(\Gamma)$ descends onto $H_p(\Gamma)$. We call it the quantum filtration, and, again, denote it by $\fil$. More specifically, an element of $H_p(\Gamma)$ is in $\fil^kH_p(\Gamma)$ if and only if it is represented by a cocycle in $\fil^kC_p(\Gamma)$. The $\zed_2$-grading also descend onto $H_p(\Gamma)$. For $i\in \zed_2$, denote by $H_p^i(\Gamma)$ the $R$-submodule of $H_p(\Gamma)$ consists of elements with $\zed_2$-degree $i$. Note that, by Proposition 22 of \cite{KR1}, $H_p(\Gamma)$ is independent of the marked points on $\Gamma$, i.e., given two different markings of the same closed graph, both of which satisfy the conditions of at the beginning of this subsection, then there is a canonical isomorphism of the two cohomologies computed by these two sets of markings, which preserves the $\zed_2$-grading. By Proposition \ref{variable-exclusion}, this isomorphism also preserves the quantum filtration.

When $p(x)=x^{n+1}$, let $C_n(\Gamma)=C_{x^{n+1}}(\Gamma)$, and, when $\Gamma$ is closed, $H_n(\Gamma)=H_{x^{n+1}}(\Gamma)$. In this case, the $\zed$-grading on $C_n(\Gamma)$ descends onto $H_n(\Gamma)$, which is called the quantum grading, and induces the quantum filtration $\fil$ on $H_n(\Gamma)$. We denote by $H_n^{i,k}(\Gamma)$ the subspace of $H_n(\Gamma)$ consists of homogeneous elements of $\zed_2$-degree $i$ and quantum degree $k$. (Note that $C_n(\Gamma)$ is the matrix factorization $C(\Gamma)$ defined in \cite{KR1}.)

For any $f\in R$, the map $\tilde{m}(f):C_p(\Gamma)\rightarrow C_p(\Gamma)$ given by $\tilde{m}(f)(X)=f\cdot X$ commutes with the differential of $C_p(\Gamma)$. So, when $\Gamma$ is closed, it induces a homomorphism $m(f):H_p(\Gamma)\rightarrow H_p(\Gamma)$. The following is a slight generalization of Proposition 2.3 of \cite{Gornik}.

\begin{lemma}\label{m(p')}
Let $\Gamma$ be a closed graph, and $x_i$ an marking on it. Then $m(p'(x_i)):H_p(\Gamma)\rightarrow H_p(\Gamma)$ is the zero homomorphism.
\end{lemma}
\begin{proof}
From Proposition 2 of \cite{KR1}, we know that, if $a_1,\cdots,a_k,b_1,\cdots,b_k\in R$ satisfy $a_1b_1+\cdots a_kb_k=0$, and $f$ is an element of the ideal $I$ of $R$ generated by $a_1,\cdots,a_k,b_1,\cdots,b_k$, then 
\[
m(f): \left(%
\begin{array}{cc}
  a_1 & b_1 \\
  \vdots & \vdots \\
  a_k & b_k \\
\end{array}%
\right)
\rightarrow
\left(%
\begin{array}{cc}
  a_1 & b_1 \\
  \vdots & \vdots \\
  a_k & b_k \\
\end{array}%
\right)
\]
is homotopic to the zero homomorphism. So we only need to check that $p'(x_i)$ is in the ideal $I$ corresponding to $C_p(\Gamma)$. If $x_i$ is the only marking on a closed circle, then $C_p(\Gamma)$ contains the factor $(\pi_{ii},0)$. In particular, $p'(x_i)=\pi_{ii}\in I$. If $x_i$ is connected to another marking $x_j$ by a regular arc, then $C_p(\Gamma)$ contains the factor $(\pi_{ij},\pm(x_i-x_j))$. And 
\[
p'(x_i)=\frac{\partial}{\partial x_i} (p(x_i)-p(x_j)) = \frac{\partial}{\partial x_i}((x_i-x_j)\pi_{ij}) = (x_i-x_j)\frac{\partial}{\partial x_i}\pi_{ij} + \pi_{ij} \in I. 
\]
If $x_i$ is next to a wide edge, let $x_j$, $x_k$, $x_l$ be the other three markings near this wide edge. Without loss of generality, assume $x_i$ and $x_j$ are the two exits of the wide edge. Then $C_p(\Gamma)$ contains the factor 
\[
\left(%
\begin{array}{cc}
  u_{ijkl} & x_i+x_j-x_k-x_l \\
  v_{ijkl} & x_ix_j-x_kx_l \\
\end{array}%
\right)\{-1\}.
\]
So
\begin{eqnarray*}
p'(x_i) & = & \frac{\partial}{\partial x_i} (p(x_i)+p(x_j)-p(x_k)-p(x_l)) \\
        & = & \frac{\partial}{\partial x_i}(u_{ijkl}(x_i+x_j-x_k-x_l)+v_{ijkl}(x_ix_j-x_kx_l)) \\
        & = & (x_i+x_j-x_k-x_l)\frac{\partial}{\partial x_i}u_{ijkl} + u_{ijkl} + (x_ix_j-x_kx_l)\frac{\partial}{\partial x_i}v_{ijkl} +x_jv_{ijkl} \in I
\end{eqnarray*}
Note that one of these three premises is always true. This proves the lemma.
\end{proof}

Given a closed graph $\Gamma$, modify $\Gamma$ by replacing each wide edge in it with a pair of parallel regular edges, i.e., change from the right side of Figure \ref{maps} to the left side of Figure \ref{maps}. We get a new graph consists of disjoint circles. Denote by $i(\Gamma)$ the $\zed_2$-number of circles in this new graph. The following proposition is a generalization of Proposition 3.2 of \cite{Gornik}.

\begin{proposition}\label{fil-grade}
Suppose that $\Gamma$ is a closed graph with $i(\Gamma)=i\in \zed_2$, and $p(x)$ is a monic polynomial of degree $n+1$ in $\C[x]$. Then $H_p^{i-1}(\Gamma)=0$, and, for each $k\in \zed$, there is a $\C$-linear isomorphism 
\[
\phi_k: H_n^{i,k}(\Gamma) \rightarrow \fil^k H_p^{i}(\Gamma)/\fil^{k-1} H_p^{i}(\Gamma).
\]
\end{proposition}

\begin{remark} 
Proposition \ref{fil-grade} plays a key role in the construction of the spectral sequence in Theorem \ref{spec}. It also implies that the filtered dimension of $H_p(\Gamma)$ is equal to the graded dimension of $H_n(\Gamma)$ as defined in \cite{KR1}, and, therefore, allows us to use many results from \cite{KR1} without repeating the actual computation.
\end{remark}

\begin{proof}[Proof of Proposition \ref{fil-grade}]
Denote the $\zed_2$-complex $C_p(\Gamma)$ by 
\[
M^0\xrightarrow{d_0}M^1\xrightarrow{d_1}M^0.
\]
Then $M^0$, $M^1$ are $\zed$-graded. And, the grading is bounded from below. Also, since all elements of $R$ have even degrees, from the construction of $C_p(\Gamma)$, we can see that the degrees of elements of $M^0$ have the same parity, and the degrees of elements of $M^0$ have the same parity, while these two parities differ by $n+1$. For $j\in\zed_2$, $d_j$ has a decomposition
\[
d_j=\sum_{l=0}^{n}d_j^{(l)},
\]
where $d_j^{(l)}$ is a homogeneous $R$-homomorphism of degree $n+1-2j$. Consider the homogeneous parts of $d_i\circ d_{i-1}=0$ and $d_{i-1}\circ d_i=0$ of degree $2(n+1-l)$. It's easy to see that
\[
\sum_{j+k=l,~ 0\leq j,k \leq n} d_i^{(j)}\circ d_{i-1}^{(k)}=0, \text{ and } \sum_{j+k=l,~ 0\leq j,k \leq n} d_{i-1}^{(j)}\circ d_i^{(k)}=0.
\]
It is also clear that the $\zed_2$-complex $C_n(\Gamma)$ is given by 
\[
M^0\xrightarrow{d_0^{(0)}}M^1\xrightarrow{d_1^{(0)}}M^0.
\]

First, we prove that $H_p^{i-1}(\Gamma)=0$. This is proved in \cite{KR1} for the special case when $p(x)=x^{n+1}$, i.e. we have $H_n^{i-1}(\Gamma)=0$. Now we prove it for a general $p(x)$. 

Let $\alpha\in M^{i-1}$ with $d_{i-1}\alpha=0$. Write 
\[
\alpha=\sum_{k=-\infty}^{\infty}\alpha_k,
\]
where $\alpha_k$ is a homogeneous element of degree $g-2k$, $\alpha_k=0$ for $k<0$, and, since the $\zed$-grading on $M^{i-1}$ is bounded from below, $\alpha_k=0$ for $k>>1$.

We construct by induction a sequence $\{\beta_k\}_{-\infty}^{\infty}\subset M^i$, such that $\beta_k$ is a homogeneous element degree $g-2k-n-1$, $\beta_k=0$ for $k<0$, and 
\[
\alpha_k=\sum_{l=0}^{n}d_i^{(l)}\beta_{k-l}.
\]
Of course, we have $\beta_k=0$ for $k>>1$ since the grading on $M^i$ is bounded from below. By definition of $\alpha_k$ and $\beta_k$, the above equation is true for $k<0$. Assume that we have found $\{\beta_k\}_{-\infty}^{m-1}\subset M^i$ so that this equation is true for $k<m$. Let's find $\beta_m$. Consider the homogeneous part of $d_{i-1}\alpha=0$ of degree $n+1+g-2m$. We have that
\begin{eqnarray*}
0 & = & \sum_{l=0}^{n}d_{i-1}^{(l)}\alpha_{m-l} ~=~ d_{i-1}^{(0)}\alpha_{m}+\sum_{l=1}^{n}d_{i-1}^{(l)}\alpha_{m-l} \\
  & = & d_{i-1}^{(0)}\alpha_{m}+\sum_{l=1}^{n}d_{i-1}^{(l)}\sum_{j=0}^n d_i^j\beta_{m-l-j} \\
  & = & d_{i-1}^{(0)}\alpha_{m} + \sum_{k=1}^{2n}(\sum_{l+j=k,~1\leq l\leq n,~ 0\leq j\leq n}d_{i-1}^{(l)}d_i^{(j)}\beta_{m-k}) \\
  & = & d_{i-1}^{(0)}\alpha_{m} - \sum_{k=1}^n d_{i-1}^{(0)}d_i^{(k)}\beta_{m-k} ~=~ d_{i-1}^{(0)}(\alpha_{m} - \sum_{k=1}^n d_i^{(k)}\beta_{m-k}).    
\end{eqnarray*}
Since $H_n^{i-1}(\Gamma)=0$, i.e. $\ker{d_{i-1}^{(0)}}=\im{d_i^{(0)}}$, there exists a homogeneous element $\beta_m$ of $M^i$ of degree $g-2m-n-1$ such that 
\[
\alpha_{m} - \sum_{k=1}^n d_i^{(k)}\beta_{m-k} = d_i^{(0)}\beta_m.
\] 
This completes the induction.

Now , we have 
\begin{eqnarray*}
\alpha & = & \sum_{k=-\infty}^{\infty}\alpha_k = \sum_{k=-\infty}^{\infty}\sum_{l=0}^{n}d_i^{(l)}\beta_{k-l} \\
       & = & \sum_{l=0}^{n}d_i^{(l)}(\sum_{k=-\infty}^{\infty}\beta_{k-l}) = \sum_{l=0}^{n}d_i^{(l)}(\sum_{j=-\infty}^{\infty}\beta_{j}) \\
       & = & d_i(\sum_{j=-\infty}^{\infty}\beta_{j}).
\end{eqnarray*}
This shows $\ker{d_{i-1}}=\im{d_i}$, i.e. $H_p^{i-1}(\Gamma)=0$.

Next we establish that $H_n^{i,k}(\Gamma) \cong \fil^k H_p^{i}(\Gamma)/\fil^{k-1} H_p^{i}(\Gamma)$. Let $(\ker{d_i^{(0)}})_k$ be the subspace of $\ker{d_i^{(0)}}$ consisting of elements of degree $k$. Define $\tilde{\phi}_k:(\ker{d_i^{(0)}})_k \rightarrow \fil^k H_p^{i}(\Gamma)/\fil^{k-1} H_p^{i}(\Gamma)$ as following.

For $\alpha\in (\ker{d_i^{(0)}})_k$, construct by induction a sequence $\{\alpha_l\}_{0}^{\infty}\subset M^i$, such that $\alpha_0=\alpha$, $\alpha_l$ is a homogeneous element of degree $k-2l$, and 
\[
\sum_{j=0}^{l}d_i^{(j)}\alpha_{l-j}=0, ~\forall ~l\in\zed,
\]
where we use the convention that $d_i^{(j)}=0$ for $j>n$. Note that $\alpha_l=0$ for $l>>1$ since the grading on $M^i$ is bounded from below. The above equation is true for $l=0$. Now assume that, for some $m\geq1$, we have found $\{\alpha_l\}_{0}^{m-1}$ such that the above equation is true for $l<m$. Let's find an $\alpha_m$. Consider
\begin{eqnarray*}
d_{i-1}^{(0)}(\sum_{j=1}^{m}d_i^{(j)}\alpha_{m-j}) & = & \sum_{j=1}^{m}d_{i-1}^{(0)}d_i^{(j)}\alpha_{m-j} \\
                                                   & = & -\sum_{j=1}^{m}\sum_{l=0}^{j-1}d_{i-1}^{(j-l)}d_i^{(l)}\alpha_{m-j} \\
                      (\text{set }p=m-j+l)         & = & -\sum_{p=0}^{m-1}\sum_{l=0}^p d_{i-1}^{(m-p)}d_i^{(l)}\alpha_{p-l} \\
                                                   & = & -\sum_{p=0}^{m-1}d_{i-1}^{(m-p)} (\sum_{l=0}^p d_i^{(l)}\alpha_{p-l}) \\
             (\text{by induction hypothesis})      & = & 0
\end{eqnarray*}
But $H_n^{i-1}(\Gamma)=0$, i.e. $\ker{d_{i-1}^{(0)}}=\im{d_i^{(0)}}$. So there exists a homogeneous element $\alpha_m$ of degree $k-2m$ such that 
\[
d_i^{(0)}\alpha_m = -\sum_{j=1}^{m}d_i^{(j)}\alpha_{m-j}.
\]
This completes the induction. 

It is clear from the above construction of $\{\alpha_l\}_{0}^{\infty}$ that 
\[
d_i(\sum_{l=0}^{\infty}\alpha_l)=0.
\]

Now define 
\[
\tilde{\phi}_k(\alpha)=[\sum_{l=0}^{\infty}\alpha_l]\in\fil^k H_p^{i}(\Gamma)/\fil^{k-1} H_p^{i}(\Gamma).
\] 
Note the sum here is in fact a finite sum. We need to check that $\tilde{\phi}_k(\alpha)$ is independent of the choice of the sequence $\{\alpha_l\}_{0}^{\infty}$. This is easy. Let $\{\alpha'_l\}_{0}^{\infty}$ another such sequence. Then
\[
(\sum_{l=0}^{\infty}\alpha_l)-(\sum_{l=0}^{\infty}\alpha'_l) = (\sum_{l=1}^{\infty}\alpha_l)-(\sum_{l=1}^{\infty}\alpha'_l) \in \fil^{k-1} \ker d_i.
\]
So
\[
[\sum_{l=0}^{\infty}\alpha_l] = [\sum_{l=0}^{\infty}\alpha'_l] \in \fil^k H_p^{i}(\Gamma)/\fil^{k-1} H_p^{i}(\Gamma).
\]
Thus, $\tilde{\phi}_k:(\ker{d_i^{(0)}})_k \rightarrow \fil^k H_p^{i}(\Gamma)/\fil^{k-1} H_p^{i}(\Gamma)$ is a well defined homomorphism. Next, we compute $\ker{\tilde{\phi}_k}$. Let $\alpha\in \ker{\tilde{\phi}_k}$, i.e., for the above constructed sequence $\{\alpha_l\}_{0}^{\infty}$, we have 
\[
\sum_{l=0}^{\infty}\alpha_l = d_{i-1}\beta+\gamma,
\]
where $\gamma$ is a cocycle in $\fil^{k-1}M^i$, and $\beta\in M^{i-1}$.  This equation implies that $d_{i-1}\beta\in\fil^kM^i$. We claim that we can choose $\beta$ so that $\deg\beta\leq k-n-1$. Assume that $\deg\beta=g> k-n-1$. Let $\beta_0$ be the top homogeneous part of $\beta$. Comparing the top homogeneous part in the above equation, we have $d_{i-1}^{(0)}\beta_0=0$. So there exists a homogeneous element $\theta\in M^i$ of degree $g-n-1$ such that $d_i^{(0)}\theta=\beta_0$. Let $\beta'=\beta-d_i\theta$. Then $\beta'$ also satisfies the above equation, and $\deg \beta' \leq \deg \beta -2$. Repeat this process. Within finite steps, we can find a $\hat{\beta}$ with $\deg\hat{\beta}\leq k-n-1$ and
\[
\sum_{l=0}^{\infty}\alpha_l = d_{i-1}\hat{\beta}+\gamma.
\]
Let $\hat{\beta}_0$ be the homogeneous part of $\hat{\beta}$ of degree $k-n-1$. From the above equation, one can see that $\alpha=\alpha_0=d_{i-1}^{(0)}\hat{\beta}_0$. This shows $\alpha\in(\im{d_{i-1}^{(0)}})_k$. So $\ker{\tilde{\phi}_k}\subset(\im{d_{i-1}^{(0)}})_k$. On the other hand, if $\alpha\in(\im{d_{i-1}^{(0)}})_k$, then there is a homogeneous element $\mu\in M^{i-1}$ of degree $k-n-1$ with $d_{i-1}^{(0)}\mu=\alpha$. Then 
\[
\sum_{l=0}^{\infty}\alpha_l = d_{i-1}\mu+\lambda,
\]
where $\lambda=-d_{i-1}\mu+\sum_{l=0}^{\infty}\alpha_l$ is a cocycle in $\fil^{k-1}M^i$. This shows that $\alpha\in \ker{\tilde{\phi}_k}$. Thus, $\ker{\tilde{\phi}_k}=(\im{d_{i-1}^{(0)}})_k$, and $\tilde{\phi}_k$ induces an isomorphism
\[
\phi_k: H_n^{i,k}(\Gamma)=(\ker{d_i^{(0)}})_k/(\im{d_{i-1}^{(0)}})_k \rightarrow \fil^k H_p^{i}(\Gamma)/\fil^{k-1} H_p^{i}(\Gamma).
\]
\end{proof}

\subsection{Graph cobodisms}\label{graph-cobodisms}

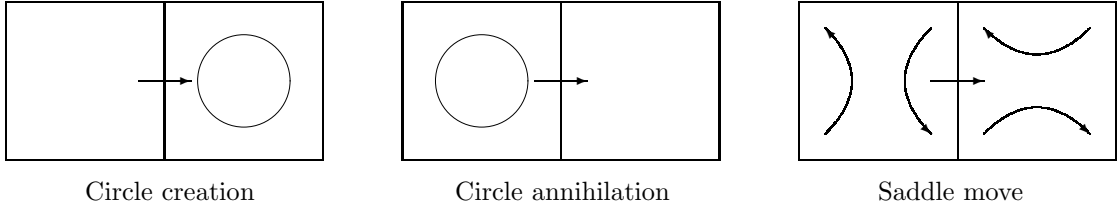
\begin{figure}[ht]

\setlength{\unitlength}{1pt}

\begin{picture}(420,75)(-210,-15)


\put(-210,0){\line(1,0){120}}

\put(-210,60){\line(1,0){120}}

\put(-210,0){\line(0,1){60}}

\put(-150,0){\line(0,1){60}}

\put(-90,0){\line(0,1){60}}

\put(-160,30){\vector(1,0){20}}

\put(-120,30){\circle{35}}

\put(-180,-15){\text{Circle creation}}


\put(-60,0){\line(1,0){120}}

\put(-60,60){\line(1,0){120}}

\put(-60,0){\line(0,1){60}}

\put(0,0){\line(0,1){60}}

\put(60,0){\line(0,1){60}}

\put(-10,30){\vector(1,0){20}}

\put(-30,30){\circle{35}}

\put(-40,-15){\text{Circle annihilation}}


\put(90,0){\line(1,0){120}}

\put(90,60){\line(1,0){120}}

\put(90,0){\line(0,1){60}}

\put(150,0){\line(0,1){60}}

\put(210,0){\line(0,1){60}}

\put(140,30){\vector(1,0){20}}

\qbezier(100,10)(120,30)(100,50)

\put(100,50){\vector(-1,1){0}}

\qbezier(140,10)(120,30)(140,50)

\put(140,10){\vector(1,-1){0}}

\qbezier(160,10)(180,30)(200,10)

\put(200,10){\vector(1,-1){0}}

\qbezier(160,50)(180,30)(200,50)

\put(160,50){\vector(-1,1){0}}

\put(120,-15){\text{Saddle move}}

\end{picture}

\caption{Graph Cobodisms}\label{graph-cobo}

\end{figure}

We call the three operations depicted in Figure \ref{graph-cobo} graph cobodisms, and, following \cite{KR1}, define the corresponding homomorphisms as below. 

Consider a circle $\bigcirc$. From Proposition \ref{fil-grade}, we have that $H_p^0(\bigcirc)=0$, and the filtered dimension of $H_p^1(\bigcirc)$ is equal to the graded dimension of $H_n^1(\bigcirc)$, which equals $[n]=q^{1-n}+q^{3-n}+\cdots+q^{n-3}+q^{n-1}$ from \cite{KR1}. So the lowest non-zero filtration of $H_p^1(\bigcirc)$ is $\fil^{1-n}H_p^1(\bigcirc)\cong \C$. This induces an injection $\iota:\C\left\langle 1 \right\rangle\rightarrow H_p(\bigcirc)$ by the composition $\C \xrightarrow{\cong} \fil^{1-n}H_p(\bigcirc) \rightarrow H_p(\bigcirc)$, where the second map is the standard inclusion, and $\C\left\langle 1 \right\rangle$ is a $\zed_2$-graded $\C$-linear space with $(\C\left\langle 1 \right\rangle)^0=0$ and $(\C\left\langle 1 \right\rangle)^1=\C$. $\iota$ induces the homomorphism of $H_p$ associated to circle creation, which is again denoted by $\iota$. More specifically, for any closed graph $\Gamma$, $\iota:H_p(\Gamma)\left\langle 1 \right\rangle\rightarrow H_p(\Gamma \sqcup \bigcirc)$ is given by the following composition:
\[
H_p(\Gamma)\left\langle 1 \right\rangle \xrightarrow{\cong} H_p(\Gamma)\otimes_{\C}\C\left\langle 1 \right\rangle \xrightarrow{\id\otimes\iota} H_p(\Gamma)\otimes_{\C}H_p(\bigcirc) \xrightarrow{\cong} H_p(\Gamma \sqcup \bigcirc).
\]
From the filtered dimension of $H_p^1(\bigcirc)$, we also have that $\fil^{n-1}H_p^1(\bigcirc)=H_p^1(\bigcirc)$, and $\fil^{n-1}H_p^1(\bigcirc)/\fil^{n-2}H_p^1(\bigcirc) \cong \C$. This induces a projection $\varepsilon:H_p(\bigcirc)\rightarrow\C\left\langle 1 \right\rangle$ given by the composition $H_p(\bigcirc)=\fil^{n-1}H_p(\bigcirc)\rightarrow\fil^{n-1}H_p(\bigcirc)/\fil^{n-2}H_p(\bigcirc)\xrightarrow{\cong}\C\left\langle 1 \right\rangle$, where the first map is the standard projection. $\varepsilon$ induces the homomorphism of $H_p$ associated to circle annihilation, which is again denoted by $\varepsilon$. More specifically, for any closed graph $\Gamma$, $\varepsilon:H_p(\Gamma \sqcup \bigcirc)\rightarrow H_p(\Gamma)\left\langle 1 \right\rangle$ is given by the following composition:
\[
H_p(\Gamma \sqcup \bigcirc) \xrightarrow{\cong} H_p(\Gamma)\otimes_{\C}H_p(\bigcirc) \xrightarrow{\id\otimes\varepsilon} H_p(\Gamma)\otimes_{\C}\C\left\langle 1 \right\rangle \xrightarrow{\cong} H_p(\Gamma)\left\langle 1 \right\rangle.
\]
Note that $\iota$ and $\varepsilon$ are homomorphisms of degree $1-n$.

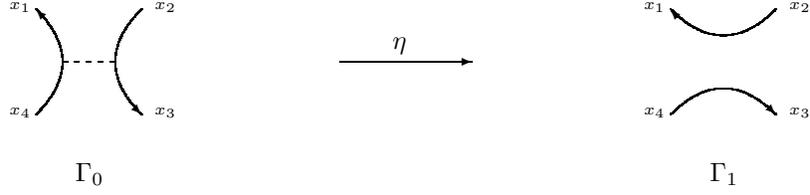
\begin{figure}[ht]

\setlength{\unitlength}{1pt}

\begin{picture}(420,75)(-210,-15)


\put(-5,35){$\eta$}

\put(-25,30){\vector(1,0){50}}


\qbezier(-140,10)(-120,30)(-140,50)

\put(-140,50){\vector(-1,1){0}}

\qbezier(-100,10)(-120,30)(-100,50)

\put(-100,10){\vector(1,-1){0}}

\multiput(-130,30)(4.5,0){5}{\line(1,0){2}}

\put(-150,50){\tiny{$x_1$}}

\put(-150,10){\tiny{$x_4$}}

\put(-95,50){\tiny{$x_2$}}

\put(-95,10){\tiny{$x_3$}}

\put(-125,-15){$\Gamma_0$}


\qbezier(100,10)(120,30)(140,10)

\put(140,10){\vector(1,-1){0}}

\qbezier(100,50)(120,30)(140,50)

\put(100,50){\vector(-1,1){0}}

\put(90,50){\tiny{$x_1$}}

\put(90,10){\tiny{$x_4$}}

\put(145,50){\tiny{$x_2$}}

\put(145,10){\tiny{$x_3$}}

\put(115,-15){$\Gamma_1$}

\end{picture}

\caption{Saddle Move}\label{saddle-move}

\end{figure}

Consider the graphs depicted in Figure \ref{saddle-move}. In the rest of this paper, we will use dashed lines to indicate where the saddle move is done. At the level of matrix factorization, the homomorphism associated to this saddle move is $\tilde{\eta}:C_p(\Gamma_0)\rightarrow C_p(\Gamma_1)\left\langle 1 \right\rangle$ defined by the vertical maps in the following diagram.
\[
\begin{CD}
\left.%
\begin{array}{c}
  R \\
  \oplus \\
  R\{2-2n\}
\end{array}%
\right. 
@>{\tiny{\left(%
\begin{array}{cc}
  \pi_{14} & x_3-x_2 \\
  \pi_{23} & x_4-x_1
\end{array}%
\right)}}>> 
\left.%
\begin{array}{c}
  R\{1-n\} \\
  \oplus \\
  R\{1-n\}
\end{array}%
\right. 
@>{\tiny{\left(%
\begin{array}{cc}
  x_1-x_4 & x_3-x_2 \\
  \pi_{23} & -\pi_{14} 
\end{array}%
\right)}}>> 
\left.%
\begin{array}{c}
  R \\
  \oplus \\
  R\{2-2n\}
\end{array}%
\right.  \\
@VV{\tiny{\left(%
\begin{array}{cc}
  e_{123}+e_{124}& 1 \\
  -e_{134}-e_{234} & 1 
\end{array}%
\right)}}V 
@VV{\tiny{\left(%
\begin{array}{cc}
  -1 & 1 \\
  -e_{123}-e_{234} & -e_{134}-e_{124} 
\end{array}%
\right)}}V 
@VV{\tiny{\left(%
\begin{array}{cc}
  e_{123}+e_{124}& 1 \\
  -e_{134}-e_{234} & 1 
\end{array}%
\right)},}V \\
\left.%
\begin{array}{c}
  R\{1-n\} \\
  \oplus \\
  R\{1-n\}
\end{array}%
\right.  
@>{\tiny{\left(%
\begin{array}{cc}
  x_2-x_1 & x_4-x_3 \\
  -\pi_{34} & \pi_{12} 
\end{array}%
\right)}}>> 
\left.%
\begin{array}{c}
  R \\
  \oplus \\
  R\{2-2n\}
\end{array}%
\right.  
@>{\tiny{\left(%
\begin{array}{cc}
  -\pi_{12} & x_4-x_3 \\
  -\pi_{34} & x_1-x_2 
\end{array}%
\right)}}>> 
\left.%
\begin{array}{c}
  R\{1-n\} \\
  \oplus \\
  R\{1-n\}
\end{array}%
\right. 
\end{CD}
\]
where the first row is $C_p(\Gamma_0)$, the second row is $C_p(\Gamma_1)\left\langle 1 \right\rangle$, and
\[
e_{ijk}=\frac{(x_k-x_j)p(x_i)+(x_i-x_k)p(x_j)+(x_j-x_i)p(x_k)}{2(x_i-x_j)(x_j-x_k)(x_k-x_i)}.
\]
It's straightforward to check that $\tilde{\eta}$ is a homomorphism of matrix factorizations of degree $n-1$. When the arcs involved are parts of two closed graphs, $\tilde{\eta}$ induces a homomorphism of cohomologies of degree $n-1$, which is denoted by $\eta$.

\begin{remark}\label{sign-problem}
As in \cite{KR1}, the homomorphism $\eta$ is defined only up to $\pm$ sign change. But this does not affect the arguments in the rest of this paper.
\end{remark}

\begin{figure}[ht]

\setlength{\unitlength}{1pt}

\begin{picture}(420,75)(-210,-15)

\put(-150,25){\Large{$\Phi$:}}


\put(-110,55){\vector(0,1){5}}

\qbezier(-110,55)(-110,45)(-100,45)

\qbezier(-100,45)(-90,45)(-90,30)

\qbezier(-90,30)(-90,15)(-100,15)

\qbezier(-100,15)(-110,15)(-110,5)

\put(-110,0){\line(0,1){5}}

\multiput(-100,15)(0,5){6}{\line(0,1){3}}

\put(-105,-15){$\Gamma$}

\put(-55,35){$\eta$}

\put(-70,30){\vector(1,0){40}}


\put(-10,0){\vector(0,1){60}}

\put(5,30){\circle{20}}

\put(-15,-15){$\Gamma\sqcup\bigcirc$}

\put(45,35){$\varepsilon$}

\put(30,30){\vector(1,0){40}}


\put(100,0){\vector(0,1){60}}

\put(95,-15){$\Gamma$}

\end{picture}

\caption{Homomorphism $\Phi$}\label{Phi}

\end{figure}
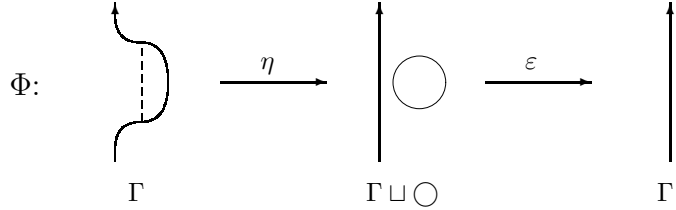

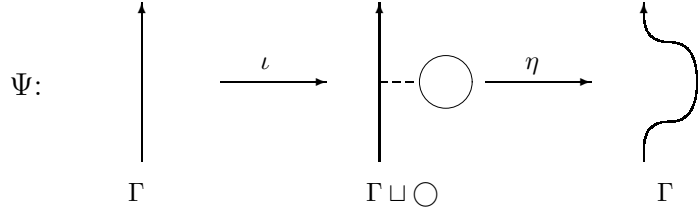
\begin{figure}[ht]

\setlength{\unitlength}{1pt}

\begin{picture}(420,75)(-210,-15)

\put(-150,25){\Large{$\Psi$:}}


\put(-100,0){\vector(0,1){60}}

\put(-55,35){$\iota$}

\put(-70,30){\vector(1,0){40}}

\put(-105,-15){$\Gamma$}


\put(-10,0){\vector(0,1){60}}

\put(15,30){\circle{20}}

\multiput(-10,30)(5,0){3}{\line(1,0){3}}

\put(-15,-15){$\Gamma\sqcup\bigcirc$}

\put(45,35){$\eta$}

\put(30,30){\vector(1,0){40}}


\put(90,55){\vector(0,1){5}}

\qbezier(90,55)(90,45)(100,45)

\qbezier(100,45)(110,45)(110,30)

\qbezier(110,30)(110,15)(100,15)

\qbezier(100,15)(90,15)(90,5)

\put(90,0){\line(0,1){5}}

\put(95,-15){$\Gamma$}

\end{picture}

\caption{Homomorphism $\Psi$}\label{Psi}

\end{figure}

\begin{lemma}\label{eta-iota-epsilon}
Let $\Gamma$ be a closed graph, and $\Phi$, $\Psi$ the homomorphisms from $H_p(\Gamma)$ to itself defined in Figures \ref{Phi}, \ref{Psi}. Then $\Phi$ and $\Psi$ are non-zero multiples of the identity map.
\end{lemma}
\begin{proof}
Denote by $\uparrow$ the part of $\Gamma$ depicted in Figures \ref{Phi} and \ref{Psi}. These figures define homomorphisms $\varphi:C_p(\uparrow)\rightarrow C_p(\uparrow)$ and $\psi:C_p(\uparrow)\rightarrow C_p(\uparrow)$. Clearly, $\Phi$ and $\Psi$ are homomorphisms on the cohomology induced by the $\zed_2$-chain maps $\varphi\otimes \id$ and $\psi\otimes\id$, where $\id$ is the identity map of the matrix factorization associated to the part of $\Gamma$ that is not depicted in Figures \ref{Phi} and \ref{Psi}.

From the definitions of $\eta$, $\iota$ and $\epsilon$, one can see that $\varphi$ and $\psi$ are both filtered homomorphisms of matrix factorizations. From Proposition \ref{hmf-tensor}, we know that 
\begin{eqnarray*}
\Hom_{HMF}(C_p(\uparrow),C_p(\uparrow)) & = & H^0(C_p(\uparrow)\otimes C_p(\uparrow)_{\bullet}) \\
& = & H^0(
\left(%
\begin{array}{cc}
  \pi_{ij} & x_i-x_j \\
  \pi_{ij} & x_j-x_i \\
\end{array}%
\right)\left\langle 1\right\rangle\{n-1\}) \\
& = & H^1(\bigcirc)\{n-1\}
\end{eqnarray*}
In particular, we have 
\[
\Hom_{hmf}(C_p(\uparrow),C_p(\uparrow)) = \fil^0 \Hom_{HMF}(C_p(\uparrow),C_p(\uparrow)) \cong \C. 
\]So any filtered homomorphism of matrix factorization from $C_p(\uparrow)$ to itself is homotopic to a scalar multiple of the identity map. To prove the lemma, we only need to check that $\varphi$ and $\psi$ are not homotopic to the zero map. We can do this by explicitly computing the homomorphisms $\hat{\varphi}$ and $\hat{\psi}$ defined in the Figure \ref{hat-phi}.

\begin{figure}[ht]

\setlength{\unitlength}{1pt}

\begin{picture}(420,100)(-210,-50)


\put(-150,25){\Large{$\hat{\varphi}$}:}

\put(-100,49){\line(0,1){2}}

\put(-102,52){$x_1$}

\put(-100,30){\circle{40}}

\multiput(-100,10)(0,5){8}{\line(0,1){3}} 

\put(-70,30){\vector(1,0){40}}

\put(-55,35){$\eta$}


\put(15,39){\line(0,1){2}}

\put(-15,39){\line(0,1){2}}

\put(-17,42){$x_1$}

\put(13,42){$x_2$}

\put(15,30){\circle{20}}

\put(-15,30){\circle{20}}

\put(30,30){\vector(1,0){40}}

\put(45,35){$\varepsilon$}


\put(100,49){\line(0,1){2}}

\put(98,52){$x_1$}

\put(100,30){\circle{40}}


\put(-100,-11){\line(0,1){2}}

\put(-102,-8){$x_1$}

\put(-150,-35){\Large{$\hat{\psi}$}:}

\put(-100,-30){\circle{40}}

\put(-70,-30){\vector(1,0){40}}

\put(-55,-25){$\iota$}


\put(15,-21){\line(0,1){2}}

\put(-15,-21){\line(0,1){2}}

\put(-17,-18){$x_1$}

\put(13,-18){$x_2$}

\put(15,-30){\circle{20}}

\put(-15,-30){\circle{20}}

\put(30,-30){\vector(1,0){40}}

\multiput(-5,-30)(5,0){2}{\line(1,0){3}}

\put(45,-25){$\eta$}


\put(100,-11){\line(0,1){2}}

\put(98,-8){$x_1$}

\put(100,-30){\circle{40}}

\end{picture}

\caption{Homomorphisms $\hat{\varphi}$ and $\hat{\psi}$}\label{hat-phi}

\end{figure}
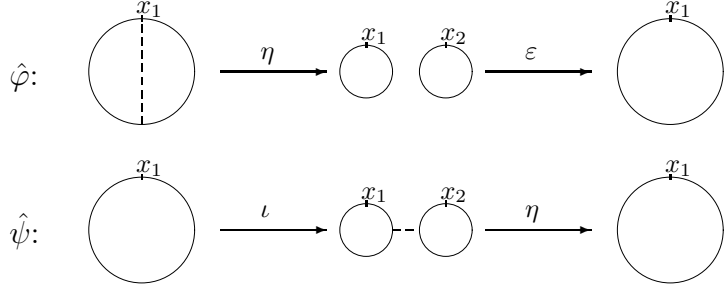

We know that 
\[
H_p(\bigcirc)\cong(\C[x_1]/(p'(x_1)))\{1-n\},  
\]
which has the basis $\{1,x_1,x_1^2,\cdots,x_1^{n-1}\}$, and 
\[
H_p(\bigcirc\sqcup\bigcirc) = H_p(\bigcirc)\otimes H_p(\bigcirc)\cong (\C[x_1,x_2]/(p'(x_1),p'(x_2)))\{2-2n\},
\]
which has the basis $\{x_1^kx_2^l~|~0\leq k,l \leq n-1\}$.

It is straightforward to check that, for $k=0,1,\cdots,n-1$,
\begin{eqnarray*}
\hat{\varphi}(x_1^k) & = & \epsilon(\eta(x_1^k)) = \epsilon(-2x_1^k(e_{112}+e_{122}))\\
                     & = & -\epsilon(x_1^k(\frac{\partial \pi_{12}}{\partial x_1}+\frac{\partial \pi_{12}}{\partial x_2})) \\
                     & = & -\epsilon(x_1^k(x_2^{n-1}+nx_2^{n-1}+ \text{ terms with power of } x_2 \text{ less than $n-1$})) \\
                     & = & -(n+1)\lambda x_1^k,
\end{eqnarray*}
where $\lambda$ is a non-zero scalar coming from the definition of $\epsilon$, and, similarly,
\[
\hat{\psi}(x_1^k) = \eta(\iota(x_1^k)) = \eta(\mu x_1^k) = \mu x_1^k,
\]
where $\mu$ is a non-zero scalar coming from the definition of $\iota$.

This shows that $\hat{\varphi}$ and $\hat{\psi}$ are non-zero homomorphisms, and hence completes the proof.
\end{proof}

\subsection{Khovanov-Rozansky homology}

\begin{figure}[ht]

\setlength{\unitlength}{1pt}

\begin{picture}(420,60)(-210,-15)


\put(-40,0){\vector(0,1){40}}

\put(-60,0){\vector(0,1){40}}

\put(-52,-15){\small{$\Gamma_0$}}

\put(-70,40){\tiny{$x_i$}}

\put(-35,40){\tiny{$x_j$}}

\put(-70,0){\tiny{$x_k$}}

\put(-35,0){\tiny{$x_l$}}


\multiput(-20,30)(7,0){6}{\line(1,0){5}}

\multiput(-20,10)(7,0){6}{\line(1,0){5}}

\put(20,30){\vector(1,0){0}}

\put(-20,10){\vector(-1,0){0}}


\put(47.5,10){\vector(1,1){0}}

\qbezier(40,0)(40,5)(47.5,10)

\qbezier(40,40)(40,35)(47.5,30)

\put(40,44){\vector(0,1){0}}

\put(52.5,10){\vector(-1,1){0}}

\qbezier(60,0)(60,5)(52.5,10)

\qbezier(60,40)(60,35)(52.5,30)

\put(60,44){\vector(0,1){0}}

\put(30,40){\tiny{$x_i$}}

\put(65,40){\tiny{$x_j$}}

\put(30,0){\tiny{$x_k$}}

\put(65,0){\tiny{$x_l$}}

\linethickness{5pt}

\put(50,10){\vector(0,1){20}}

\put(48,-15){\small{$\Gamma_1$}}

\end{picture}

\caption{$\Gamma_0$ and $\Gamma_1$}\label{maps}

\end{figure}
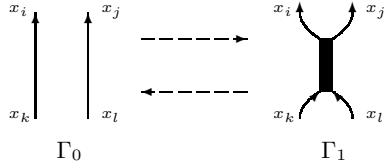

Let $\Gamma_0$ and $\Gamma_1$ be the two planar diagrams depicted
in Figure \ref{maps}. Then the matrix factorization
$C_p(\Gamma_0)$ is
\[
\left(%
\begin{array}{cc}
  \pi_{ik} & x_i-x_k \\
  \pi_{jl} & x_j-x_l \\
\end{array}%
\right)_{R}.
\]
Explicitly, this is
\[
\left[%
\begin{array}{c}
  R \\
  R\{2-2n\} \\
\end{array}%
\right] \xrightarrow{P_0}
\left[%
\begin{array}{c}
  R\{1-n\} \\
  R\{1-n\} \\
\end{array}%
\right] \xrightarrow{P_1}
\left[%
\begin{array}{c}
  R \\
  R\{2-2n\} \\
\end{array}%
\right],
\]
where
\[
P_0=\left(%
\begin{array}{cc}
  \pi_{ik} & x_j-x_l \\
  \pi_{jl} & -x_i+x_k \\
\end{array}%
\right), \hspace{1cm}
P_1=\left(%
\begin{array}{cc}
  x_i-x_k & x_j-x_l \\
  \pi_{jl} & -\pi_{ik} \\
\end{array}%
\right).
\]
The matrix factorization $C_p(\Gamma_1)$ is
\[
\left(%
\begin{array}{cc}
  u_{ijkl} & x_i+x_j-x_k-x_l \\
  v_{ijkl} & x_ix_j-x_kx_l \\
\end{array}%
\right)_{R}\{-1\}.
\]
Explicitly, this is
\[
\left[%
\begin{array}{c}
  R\{-1\} \\
  R\{3-2n\} \\
\end{array}%
\right] \xrightarrow{Q_0}
\left[%
\begin{array}{c}
  R\{-n\} \\
  R\{2-n\} \\
\end{array}%
\right] \xrightarrow{Q_1}
\left[%
\begin{array}{c}
  R\{-1\} \\
  R\{3-2n\} \\
\end{array}%
\right],
\]
where
\[
Q_0=\left(%
\begin{array}{cc}
  u_{ijkl} & x_ix_j-x_kx_l \\
  v_{ijkl} & -x_i-x_j+x_k+x_l \\
\end{array}%
\right), \hspace{.4cm}
Q_1=\left(%
\begin{array}{cc}
  x_i+x_j-x_k-x_l & x_ix_j-x_kx_l \\
  v_{ijkl} & -u_{ijkl} \\
\end{array}%
\right).
\]

Define $\chi_0:C_p(\Gamma_0)\rightarrow C_p(\Gamma_1)$ by
the matrices
\[
U_0=\left(%
\begin{array}{cc}
  x_k-x_j & 0 \\
  a_1 & 1 \\
\end{array}%
\right), \hspace{1cm} U_1=\left(%
\begin{array}{cc}
  x_k & -x_j \\
  -1 & 1 \\
\end{array}%
\right),
\]
and $\chi_1:C_p(\Gamma_1)\rightarrow C_p(\Gamma_0)$ by the
matrices
\[
V_0=\left(%
\begin{array}{cc}
  1 & 0 \\
  -a_1 & x_k-x_j \\
\end{array}%
\right), \hspace{1cm} V_1=\left(%
\begin{array}{cc}
  1 & x_j \\
  1 & x_k \\
\end{array}%
\right),
\]
where $a_1=-v_{ijkl}+(u_{ijkl}+x_iv_{ijkl}-\pi_{jl})/(x_i-x_k)$.  These are
homomorphisms of matrix factorizations of degree $1$, i.e., these
commute with the differential maps and raise the filtration by $1$. One can check that
\[
\chi_1\chi_0=(x_k-x_j)\id_{C_p(\Gamma_0)},
\hspace{1cm}
\chi_0\chi_1=(x_k-x_j)\id_{C_p(\Gamma_1)}.
\]

\begin{figure}[ht]

\setlength{\unitlength}{1pt}

\begin{picture}(420,120)(-210,-30)


\put(-80,10){\vector(0,1){40}}

\put(-100,10){\vector(0,1){40}}

\put(-95,0){$\Gamma_0$}


\put(-20,-30){\line(1,1){15}}

\put(5,-5){\vector(1,1){15}}

\put(20,-30){\vector(-1,1){40}}

\put(-5,-30){$P_-$}

\put(20,50){\line(-1,1){15}}

\put(-5,75){\vector(-1,1){15}}

\put(-20,50){\vector(1,1){40}}

\put(-5,50){$P_+$}

\multiput(30,70)(7,-3.5){6}{\line(2,-1){2.5}}

\put(70,50){\vector(2,-1){0}}

\put(50,65){\small{$1$}}

\multiput(-30,70)(-7,-3.5){6}{\line(-2,-1){2.5}}

\put(-70,50){\vector(-2,-1){0}}

\put(-55,65){\small{$0$}}

\multiput(30,-10)(7,3.5){6}{\line(2,1){2.5}}

\put(70,10){\vector(2,1){0}}

\put(50,-10){\small{$1$}}

\multiput(-30,-10)(-7,3.5){6}{\line(-2,1){2.5}}

\put(-70,10){\vector(-2,1){0}}

\put(-55,-10){\small{$0$}}


\put(87.5,20){\vector(1,1){0}}

\qbezier(80,10)(80,15)(87.5,20)

\qbezier(80,50)(80,45)(87.5,40)

\put(80,54){\vector(0,1){0}}

\put(92.5,20){\vector(-1,1){0}}

\qbezier(100,10)(100,15)(92.5,20)

\qbezier(100,50)(100,45)(92.5,40)

\put(100,54){\vector(0,1){0}}

\linethickness{5pt}

\put(90,20){\vector(0,1){20}}

\put(85,0){$\Gamma_1$}

\end{picture}

\caption{Resolutions of a crossing}\label{resolutions}

\end{figure}
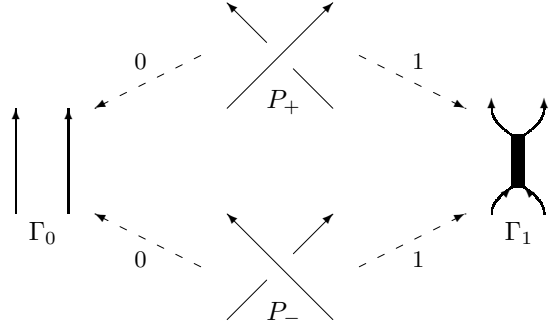

Let $D$ be an oriented link diagram. We put markings $\{x_1,\cdots,~x_m\}$
on $D$ so that none of the crossings is marked, and each arc
between two crossings has a marked point. For each crossing in
$D$, there are two ways to resolve it as shown in Figure
\ref{resolutions}. A resolution of $D$ is a planar diagram obtain
from $D$ by resolving all the crossings of $D$. Together with the
markings inherited from $D$, each resolution of $D$ is a planar
diagram satisfying the conditions in Subsection \ref{pd}.

For an arc $L^i_j$ from $x_j$ to $x_i$ that contains no crossings
and no other marked points, define $C_p(L^i_j)$ as above, and
consider it as the chain complex
\[
0 \rightarrow \underbrace{C_p(L^i_j)}_{0} \rightarrow 0,
\]
where $C_p(L^i_j)$ has cohomological degree $0$.

For a positive crossing $P_+$ in $D$, let $\Gamma_0$ and $\Gamma_1$ be the two
diagrams depicted in Figure \ref{resolutions}. Define $C_p(P_+)$
to be the chain complex
\[
0 \rightarrow \underbrace{C_p(\Gamma_1)\{n\}}_{-1} \xrightarrow{\chi_1}
\underbrace{C_p(\Gamma_0)\{n-1\}}_{0} \rightarrow 0,
\]
where $C_p(P_0)\{n-1\}$ has cohomological degree $0$, and
$C_p(P_1)\{n\}$ has cohomological degree $-1$.

For a negative crossing $P_-$ in $D$, define $C_p(P_-)$ to be the
chain complex
\[
0 \rightarrow \underbrace{C_p(\Gamma_0)\{1-n\}}_{0} \xrightarrow{\chi_0}
\underbrace{C_p(\Gamma_1)\{-n\}}_{1} \rightarrow 0,
\]
where $C_p(P_0)\{1-n\}$ has cohomological degree $0$, and
$C_p(P_1)\{-n\}$ has cohomological degree $1$.

Now define $C_p(D)$ to be the chain complex
\[
C_p(D) = (\bigotimes_{L^i_j}C_p(L^i_j)) ~\bigotimes~
(\bigotimes_{P}C_p(P)),
\]
where $L^i_j$ runs through all arcs in $D$ starting and ending in
marked points that contain no crossings and no other marked
points, and $P$ runs through all the crossings of $D$.

There are two differential maps on $C_p(D)$. One comes from the $0$-potential matrix factorization structure, which we denote by $d_{mf}$. The other comes from maps $\chi_0$, $\chi_1$, which we denote by $d_\chi$. These two differentials commute. 

Disregarding the differential $d_\chi$ on $C_p(D)$, we have a decomposition of filtered $0$-potential matrix factorization
\[
C_p(D)=\bigoplus_{\Gamma}C_p(\Gamma)\{p_\Gamma\},
\]
where $\Gamma$ runs through all resolutions of $D$, $p_\Gamma=(n-1)w+c_{\Gamma,+}-c_{\Gamma,-}$, in which $w$ is the writhe of $D$, and $c_{\Gamma,+}$ (resp. $c_{\Gamma,-}$) is the number of wide edges in $\Gamma$ from positive (resp. negative) crossings in $D$. For a given $\Gamma$, all elements of $C_p(\Gamma)$ have the same cohomological degree. The above decomposition gives us the following decomposition of filtered $\C$-linear spaces,
\[
H(C_p(D),d_{mf}) \cong \bigoplus_{\Gamma}H_p(\Gamma)\{p_\Gamma\}.
\]

$d_\chi$ induces a differential map on the cohomology $H(C_p(D),d_{mf})$, which we still denote by $d_\chi$. Define $H_p(D)$ to be the cohomology
\[
H_p(D)=H(H(C_p(D),d_{mf}),d_\chi).
\]
Note that $d_\chi$ is filtered under the quantum filtration. So $H_p(D)$ inherits the quantum filtration from $H(C_p(D),d_{mf})$, and $H_p(D)$ admits a $\zed\oplus\zed_2$-grading, where the $\zed$-grading is the cohomological grading, and the $\zed_2$-grading comes from the cyclic chain complex $(C_p(\Gamma), d_{mf})$. The $\zed_2$-grading is trivial in the sense that the cohomology is always concentrated on the $\zed_2$-degree equal to the $\zed_2$-number of components of $D$. Also, similar to \cite{KR1,KR2}, we have that $H_p(D)$ is independent of the markings on $D$ as long as none of the crossings is marked, and each arc between two crossings has a marked point. (See \cite{Lobb} for more details.)

As above, when $p(x)=x^{n+1}$, we use the notation $H_n(D)=H_{x^{n+1}}(D)$.

\subsection{The spectral sequence}\label{spec-def}
Now we generalize Gornik's construction in \cite{Gornik} to construct a spectral sequence converging to $H_p$ with $E_1$-term equal to $H_n$. 

Let $D$ be a link diagram. Consider the quantum filtration $\fil$ of the chain complex $(H(C_p(D),d_{mf}),d_{\chi})$. By construction, it is bounded from below, and, by Lemma \ref{m(p')}, it is also bounded from above. Let $\{E_k\}$ be the spectral sequence of this filtration, which, of course, converges to $H_p(D)$. In particular, we have  
\[
E_0^{i,j} = \fil^i H^{-i-j}(C_p(D),d_{mf})/ \fil^{i-1} H^{-i-j}(C_p(D),d_{mf}),
\] 
where $H^{l}(C,d)$ means the subspace of $H(C,d)$ of elements of cohomological degree $l$, and the differential on $E_0$ is the one induced by $d_{\chi}$, which we denote by $\hat{d}_{\chi}$. Denote by $d_{mf}^{(0)}$ and $d_{\chi}^{(0)}$ the top homogeneous part of $d_{mf}$ and $d_{\chi}$. By Proposition \ref{fil-grade}, there is an isomorphism
\[
\phi: H(C_p(D),d_{mf}^{(0)}) \xrightarrow{\cong} E_0.
\]
More specifically, for each pair $(i,j)$, $\phi$ is an isomorphism
\[
\phi: H^{i,j}(C_p(D),d_{mf}^{(0)}) \xrightarrow{\cong} E_0^{j,-i-j},
\]
where $H^{i,j}(C_p(D),d_{mf}^{(0)})$ is the subspace of $H(C_p(D),d_{mf}^{(0)})$ of elements of homological degree $i$ and quantum degree $j$. From the construction of $\phi$, it is straightforward to check that the following square commutes.
\[
\begin{CD}
H^{i,j}(C_p(D),d_{mf}^{(0)}) @>{d_{\chi}^{(0)}}>> H^{i+1,j}(C_p(D),d_{mf}^{(0)}) \\
@VV{\phi}V @VV{\phi}V \\
E_0^{j,-i-j} @>{\hat{d}_{\chi}}>> E_0^{j,-i-j-1} 
\end{CD}
\]
By definition, $E_1^{i,j}=H^{i+j}(E_0^{i,\ast}, \hat{d}_{\chi})$. So $\phi$ induces an isomorphism
\[
\hat{\phi}: H_n^{i,j}(D)=H^{i,j}(H^{\ast,j}(C_p(D),d_{mf}^{(0)}),d_{\chi}^{(0)}) \xrightarrow{\cong} E_1^{j,-i-2j}.
\]
In short, we have $E_1 \cong H_n(D)$.

\section{MOY Decompositions}\label{moy}

Khovanov and Rozansky's proof of the invariance of the $\mathfrak{sl}(n)$-cohomology under Reidemeister moves is based on several direct sum decompositions of the graph cohomology that correspond to Murakami, Ohtsuki and Yamada's recurrence relations of their graphic re-formulation of the $SU(n)$-HOMFLY polynomial in \cite{MOY}. Following Rasmussen \cite{Ras2}, we call these direct sum decompositions the MOY decompositions. In this section, we fix an integer $n\geq 2$ and a monic polynomial $p(x)=x^{n+1}+\sum_{i=0}^{n}c_i x^i$, and establish the MOY decompositions for the $H_p$ cohomology. 

\subsection{MOY decomposition I}

\begin{figure}[ht]

\setlength{\unitlength}{1pt}

\begin{picture}(420,65)(-210,-25)


\put(-90,-10){\vector(0,1){50}}

\put(-93,-25){$\Gamma_1$}

\put(-100,-10){\tiny{$x_2$}}

\put(-100,40){\tiny{$x_1$}}


\put(57.5,10){\vector(1,1){0}}

\qbezier(50,-10)(50,5)(57.5,10)

\qbezier(50,40)(50,25)(57,20)

\put(50,44){\vector(0,1){0}}

\qbezier(70,0)(70,5)(62.5,10)

\qbezier(70,30)(70,25)(62.5,20)

\qbezier(70,30)(70,40)(80,40)

\qbezier(80,40)(90, 40)(90,30)

\put(90,30){\vector(0,-1){30}}

\qbezier(70,0)(70,-10)(80,-10)

\qbezier(80,-10)(90, -10)(90,0)

\put(65,-25){$\Gamma$}

\put(40,-10){\tiny{$x_2$}}

\put(40,40){\tiny{$x_1$}}

\put(95,15){\tiny{$x_3$}}

\put(88.5,15){-}

\linethickness{5pt}

\put(60,10){\line(0,1){10}}

\end{picture}

\caption{$\Gamma_1$ and
$\Gamma$}\label{reduction1}

\end{figure}
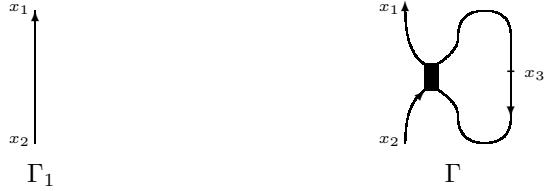

\begin{proposition}\label{moy1}
Let $\Gamma_1$ and $\Gamma$ be two closed graphs that are identical except in the part depicted in Figure \ref{reduction1}. Then there is a direct sum decomposition of filtered $\C$-linear space
\[
H_p(\Gamma)\left\langle 1\right\rangle \cong \bigoplus_{i=0}^{n-2}H_p(\Gamma_1)\{2-n+2i\}.
\]
\end{proposition}
\begin{proof}
The proof is very similar to that of Proposition 29 of \cite{KR1}. We only need to slightly generalize the transformations used there, and make sure that the generalized transformations preserve the filtration.

Define filtered homomorphisms
\[
\alpha:H_p(\Gamma_1)\{2-n\} \rightarrow H_p(\Gamma)\left\langle 1\right\rangle
\] 
and 
\[
\beta:H_p(\Gamma)\left\langle 1\right\rangle\rightarrow H_p(\Gamma_1)\{n-2\}
\]
by Figure \ref{moy1maps}.

\begin{figure}[ht]

\setlength{\unitlength}{1pt}

\begin{picture}(420,140)(-210,-100)


\put(-190,12){\Large{$\alpha$:}}

\put(-150,-10){\vector(0,1){50}}

\put(-153,-20){$\Gamma_1$}

\put(-160,-10){\tiny{$x_2$}}

\put(-160,40){\tiny{$x_1$}}

\put(-100,15){\vector(1,0){60}}

\put(-70,20){$\iota$}


\put(-20,-10){\vector(0,1){50}}

\put(0,0){\line(0,1){30}}

\qbezier(0,30)(0,40)(10,40)

\qbezier(10,40)(20, 40)(20,30)

\put(20,30){\vector(0,-1){30}}

\qbezier(0,0)(0,-10)(10,-10)

\qbezier(10,-10)(20, -10)(20,0)

\put(-10,-20){$\Gamma_2$}

\put(-30,-10){\tiny{$x_2$}}

\put(-30,40){\tiny{$x_1$}}

\put(25,15){\tiny{$x_3$}}

\put(18.5,15){-}

\put(40,15){\vector(1,0){60}}

\put(70,20){$\chi_0$}


\put(117.5,10){\vector(1,1){0}}

\qbezier(110,-10)(110,5)(117.5,10)

\qbezier(110,40)(110,25)(117,20)

\put(110,44){\vector(0,1){0}}

\qbezier(130,0)(130,5)(122.5,10)

\qbezier(130,30)(130,25)(122.5,20)

\qbezier(130,30)(130,40)(140,40)

\qbezier(140,40)(150, 40)(150,30)

\put(150,30){\vector(0,-1){30}}

\qbezier(130,0)(130,-10)(140,-10)

\qbezier(140,-10)(150, -10)(150,0)

\put(117,-20){$\Gamma$}

\put(100,-10){\tiny{$x_2$}}

\put(100,40){\tiny{$x_1$}}

\put(155,15){\tiny{$x_3$}}

\put(148.5,15){-}


\put(-190,-63){\Large{$\beta$:}}

\put(-152.5,-65){\vector(1,1){0}}

\qbezier(-160,-85)(-160,-70)(-152.5,-65)

\qbezier(-160,-35)(-160,-50)(-153,-55)

\put(-160,-31){\vector(0,1){0}}

\qbezier(-140,-75)(-140,-70)(-147.5,-65)

\qbezier(-140,-45)(-140,-50)(-147.5,-55)

\qbezier(-140,-45)(-140,-35)(-130,-35)

\qbezier(-130,-35)(-120,-35)(-120,-45)

\put(-120,-45){\vector(0,-1){30}}

\qbezier(-140,-75)(-140,-85)(-130,-85)

\qbezier(-130,-85)(-120, -85)(-120,-75)

\put(-153,-95){$\Gamma$}

\put(-170,-85){\tiny{$x_2$}}

\put(-170,-35){\tiny{$x_1$}}

\put(-115,-60){\tiny{$x_3$}}

\put(-121.5,-60){-}

\put(-100,-60){\vector(1,0){60}}

\put(-70,-55){$\chi_1$}


\put(-20,-85){\vector(0,1){50}}

\put(0,-75){\line(0,1){30}}

\qbezier(0,-45)(0,-35)(10,-35)

\qbezier(10,-35)(20, -35)(20,-45)

\put(20,-45){\vector(0,-1){30}}

\qbezier(0,-75)(0,-85)(10,-85)

\qbezier(10,-85)(20, -85)(20,-75)

\put(-10,-95){$\Gamma_2$}

\put(-30,-85){\tiny{$x_2$}}

\put(-30,-35){\tiny{$x_1$}}

\put(25,-60){\tiny{$x_3$}}

\put(18.5,-60){-}

\put(40,-60){\vector(1,0){60}}

\put(70,-55){$\varepsilon$}


\put(120,-85){\vector(0,1){50}}

\put(117,-95){$\Gamma_1$}

\put(110,-85){\tiny{$x_2$}}

\put(110,-35){\tiny{$x_1$}}

\linethickness{5pt}

\put(120,10){\line(0,1){10}}

\put(-150,-65){\line(0,1){10}}

\end{picture}

\caption{$\alpha$ and
$\beta$}\label{moy1maps}

\end{figure}
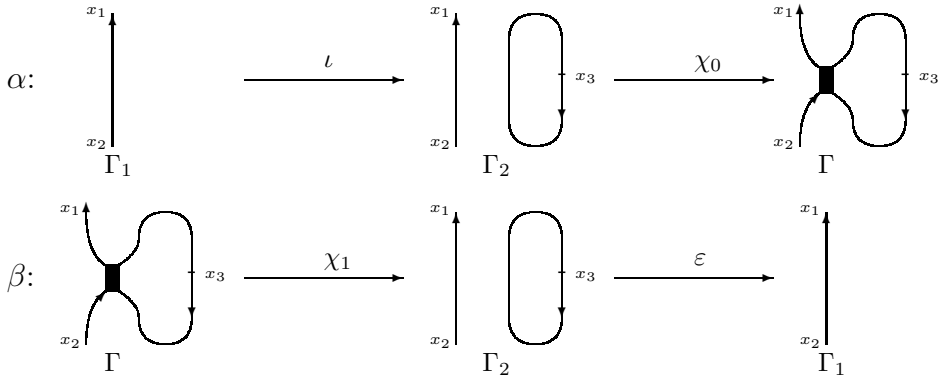

Then, for $0\leq i \leq n-2$, define a filtered homomorphism
\[
\alpha_i:H_p(\Gamma_1)\{2+2i-n\} \rightarrow H_p(\Gamma) \left\langle 1\right\rangle~\text{ by }~ \alpha_i=\sum_{k=0}^{i}\sum_{l=0}^{i-k}(n-k+1)c_{n-k+1}m(x_1^{i-l-k}x_3^l)\circ\alpha,
\]
where $c_k$ is the $k$-th coefficient of $p(x)$, in particular $c_{n+1}=1$, and $m(x_1^{i-l-k}x_3^l)$ is the homomorphism induced by multiplying $x_1^{i-j-k}x_3^j$. Similarly, define a filtered homomorphism
\[
\beta_i:H_p(\Gamma)\left\langle 1\right\rangle\rightarrow H_p(\Gamma_1)\{2+2i-n\} ~\text{ by }~ \beta_i=-\beta\circ m(x_3^{n-i-2}),
\]
where $m(x_3^{n-i-2})$ is the homomorphism induced by multiplying $x_3^{n-i-2}$.

We have, for $0\leq i,j \leq n-2$,
\begin{eqnarray*}
\beta_j\circ\alpha_i & = & -\sum_{k=0}^{i}\sum_{l=0}^{i-k}(n-k+1)c_{n-k+1}\beta\circ m(x_3^{n-j-2}) \circ m(x_1^{i-l-k}x_3^l)\circ\alpha \\
                     & = & -\sum_{k=0}^{i}\sum_{l=0}^{i-k}(n-k+1)c_{n-k+1}\varepsilon \circ \chi_1 \circ \chi_0 \circ m(x_3^{n-j-2}) \circ m(x_1^{i-l-k}x_3^l) \circ \iota \\
                     & = & -\sum_{k=0}^{i}\sum_{l=0}^{i-k} (n-k+1)c_{n-k+1}\varepsilon \circ m(x_1-x_3) \circ m(x_3^{n-j-2}) \circ m(x_1^{i-l-k}x_3^l) \circ \iota \\
                     & = & -\sum_{k=0}^{i} (n-k+1)c_{n-k+1}\varepsilon \circ m(x_3^{n-j-2}) \circ m(x_1^{i-k+1}-x_3^{i-k+1}) \circ \iota \\
                     & = & \sum_{k=0}^{i} (n-k+1)c_{n-k+1}\varepsilon \circ m(x_3^{n-j+i-k-1}-x_1^{i-k+1}x_3^{n-j-2})\circ \iota \\
                     & = & \sum_{k=0}^{i} (n-k+1)c_{n-k+1}\varepsilon \circ m(x_3^{n-j+i-k-1})\circ \iota \\
                     & = & \varepsilon \circ m(\sum_{k=0}^{i} (n-k+1)c_{n-k+1}x_3^{n-j+i-k-1})\circ \iota
\end{eqnarray*}

If $i<j$, then very term in the sum 
\[
\sum_{k=0}^{i} (n-k+1)c_{n-k+1}x_3^{n-j+i-k-1}
\] 
has power less than $n-1$. By the definition of $\varepsilon$ and $\iota$, we have $\beta_j\circ\alpha_i=0$. 

If $i>j$, then by Lemme \ref{m(p')}, we have 
\begin{eqnarray*}
m(\sum_{k=0}^{i} (n-k+1)c_{n-k+1}x_3^{n-j+i-k-1}) & = & m(x_3^{i-j-1}\sum_{k=0}^{i} (n-k+1)c_{n-k+1}x_3^{n-k}) \\
                                                  & = & m(x_3^{i-j-1}(p'(x_3)-\sum_{k=i+1}^{n} (n-k+1)c_{n-k+1}x_3^{n-k})) \\
                                                  & = & -m(x_3^{i-j-1}\sum_{k=i+1}^{n} (n-k+1)c_{n-k+1}x_3^{n-k}) \\
                                                  & = & -m(\sum_{k=i+1}^{n} (n-k+1)c_{n-k+1}x_3^{n-k+i-j-1}),
\end{eqnarray*}
in which every term has power less than $n-1$. This again implies $\beta_j\circ\alpha_i=0$. 

If $i=j$, then 
\[
\sum_{k=0}^{i} (n-k+1)c_{n-k+1}x_3^{n-j+i-k-1} = \sum_{k=0}^{i} (n-k+1)c_{n-k+1}x_3^{n-k-1},
\] 
whose leading term is $(n+1)x_3^{n-1}$, and all other terms have power less than $n-1$. So $\beta_i\circ\alpha_i=\lambda(n+1)\id$, where $\lambda$ is a constant determined by the definitions of $\varepsilon$ and $\iota$, and is independent of $i$.

Thus, we have proved
\[
\beta_j\circ\alpha_i=\lambda(n+1)\delta_{ij}\id_{H_p(\Gamma_1)\{2+2i-n\}}.
\]
Define 
\begin{eqnarray*}
\widetilde{\alpha} = \sum_{i=0}^{n-2}\alpha_i & : & \bigoplus_{i=0}^{n-2}H_p(\Gamma_1)\{2-n+2i\} \rightarrow H_p(\Gamma)\left\langle 1\right\rangle \\
\widetilde{\beta} = \frac{1}{\lambda(n+1)}\sum_{i=0}^{n-2}\beta_i & : & H_p(\Gamma)\left\langle 1\right\rangle \rightarrow \bigoplus_{i=0}^{n-2}H_p(\Gamma_1)\{2-n+2i\}.
\end{eqnarray*}
Then 
\[
\widetilde{\beta}\circ\widetilde{\alpha}=\id_{\bigoplus_{i=0}^{n-2}H_p(\Gamma_1)\{2-n+2i\}}.
\]
But, from Proposition \ref{fil-grade} and Proposition 29 of \cite{KR1}, we know that $\bigoplus_{i=0}^{n-2}H_p(\Gamma_1)\{2-n+2i\}$ and $H_p(\Gamma)\left\langle 1\right\rangle$ have that same filtered dimension. So $\widetilde{\beta}$ and $\widetilde{\alpha}$ are filtered isomorphisms that are inverses of each other.
\end{proof}

\subsection{MOY decomposition II}

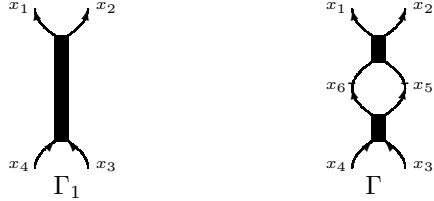
\begin{figure}[ht]

\setlength{\unitlength}{1pt}

\begin{picture}(420,70)(-210,-10)


\put(62.5,10){\vector(-1,1){0}}

\qbezier(70,0)(70,5)(62.5,10)

\qbezier(70,30)(70,25)(62.5,20)

\put(70,30){\vector(1,4){0}}

\put(50,30){\vector(-1,4){0}}

\put(57.5,10){\vector(1,1){0}}

\qbezier(50,0)(50,5)(57.5,10)

\qbezier(50,30)(50,25)(57.5,20)

\qbezier(70,30)(70,35)(62.5,40)

\qbezier(70,60)(70,55)(62.5,50)

\put(70,60){\vector(1,4){0}}

\put(50,60){\vector(-1,4){0}}

\qbezier(50,30)(50,35)(57.5,40)

\qbezier(50,60)(50,55)(57.5,50)

\put(55,-10){$\Gamma$}

\put(40,60){\tiny{$x_1$}}

\put(40,0){\tiny{$x_4$}}

\put(73,60){\tiny{$x_2$}}

\put(73,0){\tiny{$x_3$}}

\put(40,30){\tiny{$x_6$}}

\put(73,30){\tiny{$x_5$}}

\put(48.5,30){-}

\put(68.5,30){-}


\put(-62.5,10){\vector(1,1){0}}

\qbezier(-70,0)(-70,5)(-62.5,10)

\put(-57.5,10){\vector(-1,1){0}}

\qbezier(-50,0)(-50,5)(-57.5,10)

\qbezier(-70,60)(-70,55)(-62.5,50)

\put(-70,60){\vector(-1,4){0}}

\put(-50,60){\vector(1,4){0}}

\qbezier(-50,60)(-50,55)(-57.5,50)

\put(-63,-10){$\Gamma_1$}

\put(-80,60){\tiny{$x_1$}}

\put(-80,0){\tiny{$x_4$}}

\put(-47,60){\tiny{$x_2$}}

\put(-47,0){\tiny{$x_3$}}

\linethickness{5pt}

\put(-60,10){\line(0,1){40}}

\put(60,10){\line(0,1){10}}

\put(60,40){\line(0,1){10}}

\end{picture}

\caption{$\Gamma_1$ and
$\Gamma$}\label{reduction2}

\end{figure}

\begin{proposition}\label{moy2}
Let $\Gamma_1$ and $\Gamma$ be two closed graphs that are identical except in the part depicted in Figure \ref{reduction2}. Then there is a direct sum decomposition of filtered $\C$-linear space
\[
H_p(\Gamma) \cong H_p(\Gamma_1)\{1\} \oplus H_p(\Gamma_1)\{-1\}.
\]
\end{proposition}
\begin{proof}
The proof of MOY decomposition II in \cite{KR1} does not depend on the particular choice of the potential polynomial there. So this proposition is already proved in \cite{KR1}. We only need to check that the transformations in that proof preserve the quantum filtration. 

Let $R$ be the polynomial ring over $\C$ generated by marking on $\Gamma_1$. Then $\widetilde{R}=R[x_5,x_6]$ is the polynomial ring over $\C$ generated by marking on $\Gamma$. Let $s_1=x_5+x_6$, $s_2=x_5x_6$ and $\hat{R}=R[s_1,s_2]$. Then 
\[
\widetilde{R}=\hat{R}\oplus x_6 \cdot \hat{R}\cong \hat{R} \oplus \hat{R}\{2\}.
\]
The factor of $C_p(\Gamma)$ contributed by the part depicted in Figure \ref{reduction2} is
\[
\widetilde{M} =
\left(%
\begin{array}{cc}
  u_{1256} & x_1+x_2-s_1 \\
  v_{1256} & x_1x_2-s_2  \\
  u_{5634} & s_1-x_3-x_4 \\
  v_{5634} & s_2-x_3x_4  \\
\end{array}%
\right)_{\widetilde{R}}\{-2\}.
\]
Let 
\[
\hat{M} =
\left(%
\begin{array}{cc}
  u_{1256} & x_1+x_2-s_1 \\
  v_{1256} & x_1x_2-s_2  \\
  u_{5634} & s_1-x_3-x_4 \\
  v_{5634} & s_2-x_3x_4  \\
\end{array}%
\right)_{\hat{R}}.
\]
Then there is a decomposition of filtered matrix factorization
\[
\widetilde{M} = \hat{M}\{-2\} \oplus x_6\cdot\hat{M}\{-2\} \cong \hat{M}\{-2\} \oplus \hat{M}.
\]
Let $N$ be the matrix factorization associated to the part of $\Gamma$ outside Figure \ref{reduction2}. Using Proposition \ref{variable-exclusion}, we can exclude $s_1$ and $s_2$, which gives
\[
H(N\otimes_R\hat{M}) \cong 
H(N\otimes_R\left(%
\begin{array}{cc}
  u_{1234} & x_1+x_2-x_3-x_4 \\
  v_{1234} & x_1x_2-x_3x_4  \\
\end{array}%
\right)_{R})
=H_p(\Gamma_1)\{1\}.
\]
Thus
\[
H_p(\Gamma)= H(N\otimes_R\widetilde{M}) \cong H(N\otimes_R\hat{M})\{-2\}\oplus H(N\otimes_R\hat{M}) \cong H_p(\Gamma_1)\{-1\}\oplus H_p(\Gamma_1)\{1\}.
\]
\end{proof}

\subsection{MOY decomposition III}\label{subsection-moy3}

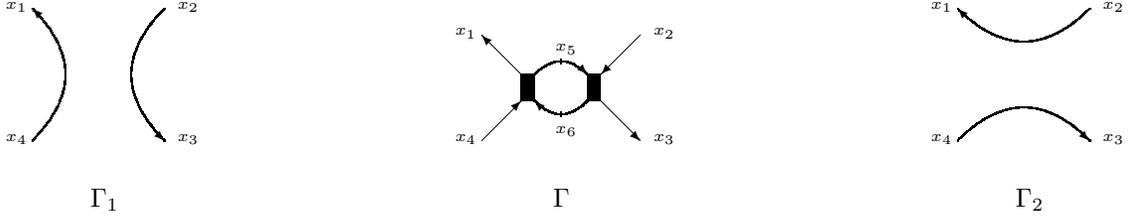
\begin{figure}[ht]

\setlength{\unitlength}{1pt}

\begin{picture}(420,80)(-210,-25)


\qbezier(-200,0)(-175,25)(-200,50)

\put(-200,50){\vector(-1,1){0}}

\qbezier(-150,0)(-175,25)(-150,50)

\put(-150,0){\vector(1,-1){0}}

\put(-178,-25){$\Gamma_1$}

\put(-210,50){\tiny{$x_1$}}

\put(-145,50){\tiny{$x_2$}}

\put(-210,0){\tiny{$x_4$}}

\put(-145,0){\tiny{$x_3$}}


\qbezier(-10,25)(0,35)(10,25)

\put(10,25){\vector(1,-1){0}}

\qbezier(-10,15)(0,5)(10,15)

\put(-10,15){\vector(-1,1){0}}

\put(-15,25){\vector(-1,1){15}}

\put(15,15){\vector(1,-1){15}}

\put(-30,0){\vector(1,1){15}}

\put(30,40){\vector(-1,-1){15}}

\put(-40,40){\tiny{$x_1$}}

\put(35,40){\tiny{$x_2$}}

\put(-40,0){\tiny{$x_4$}}

\put(35,0){\tiny{$x_3$}}

\put(-2,34){\tiny{$x_5$}}

\put(-2,3){\tiny{$x_6$}}

\put(0,29){\line(0,1){2}}

\put(0,9){\line(0,1){2}}

\put(-3,-25){$\Gamma$}


\qbezier(200,0)(175,25)(150,0)

\put(200,0){\vector(1,-1){0}}

\qbezier(150,50)(175,25)(200,50)

\put(150,50){\vector(-1,1){0}}

\put(172,-25){$\Gamma_2$}

\put(140,50){\tiny{$x_1$}}

\put(205,50){\tiny{$x_2$}}

\put(140,0){\tiny{$x_4$}}

\put(205,0){\tiny{$x_3$}}

\linethickness{5pt}

\put(-12.5,15){\line(0,1){10}}

\put(12.5,15){\line(0,1){10}}

\end{picture}

\caption{$\Gamma_1$, $\Gamma$ and
$\Gamma_2$}\label{reduction3}

\end{figure}

\begin{proposition}\label{moy3}
Let $\Gamma_1$ $\Gamma_2$ and $\Gamma$ be three closed graphs that are identical except in the part depicted in Figure \ref{reduction3}. Then there is a direct sum decomposition of filtered $\C$-linear space
\[
H_p(\Gamma) \cong H_p(\Gamma_2) \oplus (\bigoplus_{i=0}^{n-3}H_p(\Gamma_1)\left\langle 1 \right\rangle\{3-n+2i\}).
\]
\end{proposition}

Our proof of this decomposition is different from that in \cite{KR1}, which is partially implicit and depends on the homotopy splitting idempotents property. Instead, we give a completely explicit construction without using the homotopy splitting idempotents property. To start, we prove the following lemma.

\begin{figure}[ht]

\setlength{\unitlength}{1pt}

\begin{picture}(420,80)(-210,-25)

\qbezier(-10,25)(0,35)(10,25)

\put(10,25){\vector(1,-1){0}}

\qbezier(-10,15)(0,5)(10,15)

\put(-10,15){\vector(-1,1){0}}

\qbezier(-15,25)(0,60)(15,25)

\put(15,25){\vector(1,-2){0}}

\qbezier(-15,15)(0,-20)(15,15)

\put(-15,15){\vector(-1,2){0}}

\put(-2,45){\tiny{$x_1$}}

\put(-2,32){\tiny{$x_2$}}

\put(-2,-8){\tiny{$x_3$}}

\put(-2,4){\tiny{$x_4$}}

\put(0,41.5){\line(0,1){2}}

\put(0,29){\line(0,1){2}}

\put(0,9){\line(0,1){2}}

\put(0,-3.5){\line(0,1){2}}

\put(-3,-25){$\hat{\Gamma}$}

\linethickness{5pt}

\put(-12.5,15){\line(0,1){10}}

\put(12.5,15){\line(0,1){10}}

\end{picture}

\caption{$\hat{\Gamma}$}\label{hat-gamma}

\end{figure}
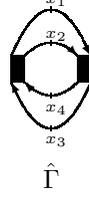

\begin{lemma}\label{H-hat-gamma}
Let $\hat{\Gamma}$ be the closed graph depicted in Figure \ref{hat-gamma}. Then $H_p(\hat{\Gamma})$ has a basis
\[
\{x_1^ix_2^jx_3^{\epsilon} ~|~ 0\leq i \leq n-1,~0\leq j\leq n-2,~\epsilon=0,1\},
\]
where the element $x_1^ix_2^jx_3^{\epsilon}$ has quantum degree $2i+2j+2\epsilon-2n+2$.
\end{lemma}

\begin{proof}
Let $R=\C[x_1,x_2]$, and $I$ be the ideal of $R$ generated by $\{u_{1212},v_{1212}\}$. Similar to the proof of Proposition \ref{moy3}, we have 
\[
C_p(\hat{\Gamma}) =
\left(%
\begin{array}{cc}
  u_{1234} & x_1+x_2-x_3-x_4 \\
  v_{1234} & x_1x_2-x_3x_4  \\
  u_{3412} & x_3+x_4-x_1-x_2 \\
  v_{3412} & x_3x_4-x_1x_2  \\
\end{array}%
\right)_{R[x_3,x_4]}\{-2\},
\]
which is quasi-isomorphic to
\[
\left(%
\begin{array}{cc}
  u_{1212} & 0 \\
  v_{1212} & 0  \\
\end{array}%
\right)_{R}\{-2\}
\oplus
x_3\cdot\left(%
\begin{array}{cc}
  u_{1212} & 0 \\
  v_{1212} & 0  \\
\end{array}%
\right)_{R} \{-2\}.
\]
It is easy to check that the cohomology of 
\[
\left(%
\begin{array}{cc}
  u_{1212} & 0 \\
  v_{1212} & 0  \\
\end{array}%
\right)_{R}
\]
is isomorphic to $(R/I)\{4-2n\}$. So $H_p(\hat{\Gamma})\cong (R/I \oplus x_3\cdot R/I)\{2-2n\}$. By Lemma \ref{m(p')}, we know that $p'(x_1)$ and $p'(x_2)$ act on $H_p(\hat{\Gamma})$ as zero homomorphisms. So $H_p(\hat{\Gamma})$ is generated as a $\C$-linear space by
\[
\{x_1^ix_2^jx_3^{\epsilon} ~|~ 0\leq i,j \leq n-1,~\epsilon=0,1\}.
\]
But $v_{1212}$ acts on $H_p(\hat{\Gamma})$ as the zero homomorphism. Note that, as elements of $H_p(\hat{\Gamma})$,
\[
0=v_{1212}= - (n+1)\sum_{l=0}^{n-1}x_1^{n-1-l}x_2^l + \text{ lower degree terms},
\]
which implies that, as elements of $H_p(\hat{\Gamma})$,
\[
x_2^{n-1}= - \sum_{l=0}^{n-2}x_1^{n-1-l}x_2^l + \text{ lower degree terms}.
\]
From this and $p'(x_1)=0$, one can see that $H_p(\hat{\Gamma})$ is generated as a $\C$-linear space by
\[
\{x_1^ix_2^jx_3^{\epsilon} ~|~ 0\leq i \leq n-1,~0\leq j\leq n-2,~\epsilon=0,1\}.
\]
Using MOY I and II decompositions and the filtered dimension of $H_p(\bigcirc)$, one can easily check that the filtered dimension of $H_p(\hat{\Gamma})$ is $[2][n][n-1]$, where $[k]=\frac{q^k-q^{-k}}{q-q^{-1}}$. This implies that the above set of elements is indeed a basis for $H_p(\hat{\Gamma})$.
\end{proof}

From the above proof, we know that $\fil^{1-2n}H_p(\hat{\Gamma})=0$, $\fil^{2-2n}H_p(\hat{\Gamma})\cong \C$, $\fil^{2n-2}H_p(\hat{\Gamma}) = H_p(\hat{\Gamma})$ and $\fil^{2n-2}H_p(\hat{\Gamma})/\fil^{2n-1}H_p(\hat{\Gamma}) \cong \C$. We define homomorphism
$\widetilde{\iota}: \C \rightarrow H_p(\hat{\Gamma})$ by the composition
\[
\C \xrightarrow{\cong} \fil^{2-2n}H_p(\hat{\Gamma}) \rightarrow H_p(\hat{\Gamma}),
\]
where the second map is the standard inclusion, and define homomorphism $\widetilde{\varepsilon}: H_p(\hat{\Gamma}) \rightarrow \C$ by the composition
\[
H_p(\hat{\Gamma})= \fil^{2n-2}H_p(\hat{\Gamma}) \rightarrow \fil^{2n-2}H_p(\hat{\Gamma})/\fil^{2n-1}H_p(\hat{\Gamma}) \xrightarrow{\cong} \C,
\]
where the first map is the standard projection. Note that $\widetilde{\iota}$ and $\widetilde{\varepsilon}$ are homomorphisms of degree $2-2n$, and induce homomorphisms between the $H_p$ cohomology of a graph and the $H_p$ cohomology of the disjoint union of that graph and $\hat{\Gamma}$.

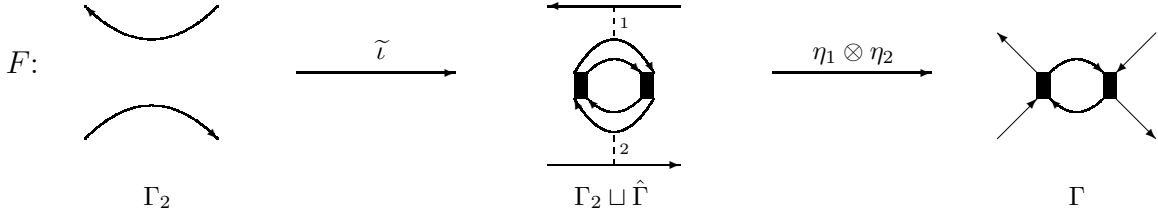
\begin{figure}[ht]

\setlength{\unitlength}{1pt}

\begin{picture}(420,80)(-210,-25)


\put(-230,25){\Large{$F$}:}

\qbezier(-150,0)(-175,25)(-200,0)

\put(-150,0){\vector(1,-1){0}}

\qbezier(-200,50)(-175,25)(-150,50)

\put(-200,50){\vector(-1,1){0}}

\put(-178,-25){$\Gamma_2$}

\put(-120,25){\vector(1,0){60}}

\put(-90,30){$\widetilde{\iota}$}


\qbezier(-10,25)(0,35)(10,25)

\put(10,25){\vector(1,-1){0}}

\qbezier(-10,15)(0,5)(10,15)

\put(-10,15){\vector(-1,1){0}}

\qbezier(-15,25)(0,50)(15,25)

\put(15,25){\vector(1,-2){0}}

\qbezier(-15,15)(0,-10)(15,15)

\put(-15,15){\vector(-1,2){0}}

\put(-15,-25){$\Gamma_2 \sqcup\hat{\Gamma}$}

\put(-25,-10){\vector(1,0){50}}

\put(25,50){\vector(-1,0){50}}

\multiput(0,-10)(0,4.5){3}{\line(0,1){2}}

\multiput(0,50)(0,-4.5){3}{\line(0,-1){2}}

\put(2,42){\tiny{$1$}}

\put(2,-6){\tiny{$2$}}

\put(60,25){\vector(1,0){60}}

\put(75,30){$\eta_1\otimes\eta_2$}


\qbezier(165,25)(175,35)(185,25)

\put(185,25){\vector(1,-1){0}}

\qbezier(165,15)(175,5)(185,15)

\put(165,15){\vector(-1,1){0}}

\put(160,25){\vector(-1,1){15}}

\put(190,15){\vector(1,-1){15}}

\put(145,0){\vector(1,1){15}}

\put(205,40){\vector(-1,-1){15}}

\put(172,-25){$\Gamma$}

\linethickness{5pt}

\put(162.5,15){\line(0,1){10}}

\put(187.5,15){\line(0,1){10}}

\put(-12.5,15){\line(0,1){10}}

\put(12.5,15){\line(0,1){10}}

\end{picture}

\caption{Definition of $F$, where $\eta_1$ and $\eta_2$ are homomorphisms associated to saddle moves $1$ and $2$.}\label{mapF}

\end{figure}

\begin{figure}[ht]

\setlength{\unitlength}{1pt}

\begin{picture}(420,80)(-210,-25)


\qbezier(200,0)(175,25)(150,0)

\put(200,0){\vector(1,-1){0}}

\qbezier(150,50)(175,25)(200,50)

\put(150,50){\vector(-1,1){0}}

\put(172,-25){$\Gamma_2$}


\qbezier(-10,25)(0,35)(10,25)

\put(10,25){\vector(1,-1){0}}

\qbezier(-10,15)(0,5)(10,15)

\put(-10,15){\vector(-1,1){0}}

\qbezier(-15,25)(0,50)(15,25)

\put(15,25){\vector(1,-2){0}}

\qbezier(-15,15)(0,-10)(15,15)

\put(-15,15){\vector(-1,2){0}}

\put(-15,-25){$\Gamma_2 \sqcup\hat{\Gamma}$}

\put(-25,-10){\vector(1,0){50}}

\put(25,50){\vector(-1,0){50}}

\put(60,25){\vector(1,0){60}}

\put(90,30){$\widetilde{\varepsilon}$}


\put(-230,25){\Large{$G$}:}

\qbezier(-185,25)(-175,35)(-165,25)

\put(-165,25){\vector(1,-1){0}}

\qbezier(-185,15)(-175,5)(-165,15)

\put(-185,15){\vector(-1,1){0}}

\put(-190,25){\vector(-1,1){25}}

\put(-160,15){\vector(1,-1){25}}

\put(-215,-10){\vector(1,1){25}}

\put(-135,50){\vector(-1,-1){25}}

\put(-178,-25){$\Gamma$}

\put(-105,30){$\eta_3\otimes\eta_4$}

\put(-120,25){\vector(1,0){60}}

\multiput(-205,40)(4,0){15}{\line(1,0){2}}

\multiput(-205,0)(4,0){15}{\line(1,0){2}}

\put(-175,42){\tiny{$3$}}

\put(-175,-7){\tiny{$4$}}

\linethickness{5pt}

\put(-187.5,15){\line(0,1){10}}

\put(-162.5,15){\line(0,1){10}}

\put(-12.5,15){\line(0,1){10}}

\put(12.5,15){\line(0,1){10}}

\end{picture}

\caption{Definition of $G$, where $\eta_3$ and $\eta_4$ are homomorphisms associated to saddle moves $3$ and $4$.}\label{mapG}

\end{figure}
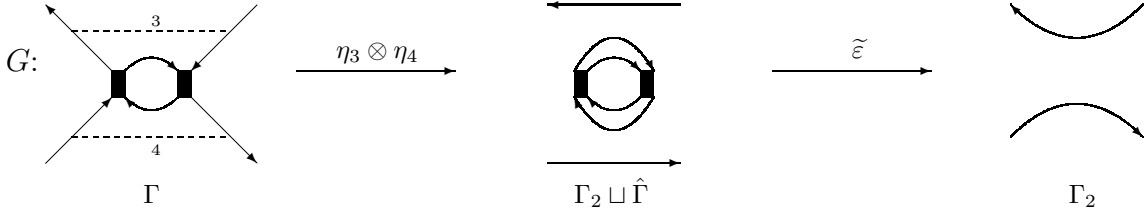

Now we define filtered homomorphisms $F:H_p(\Gamma_2)\rightarrow H_p(\Gamma)$ and $G:H_p(\Gamma) \rightarrow H_p(\Gamma_2)$ by Figures \ref{mapF} and \ref{mapG}.

\begin{lemma}\label{G-circ-F}
$G\circ F: H_p(\Gamma_2)\rightarrow H_p(\Gamma_2)$ is a non-zero multiple of the identity map.
\end{lemma}

\begin{proof}
Let $\Gamma'$ and $\Gamma_2'$ be the parts of $\Gamma$ and $\Gamma_2'$ depicted in Figure \ref{reduction3}. Then Figures \ref{mapF} and \ref{mapG} also define filtered homomorphisms $f:C_p(\Gamma_2')\rightarrow C_p(\Gamma')$ and $g: C_p(\Gamma') \rightarrow C_p(\Gamma_2')$. Moreover, $F$ and $G$ are homomorphisms of cohomologies induced by the chain maps $f\otimes\id$ and $g\otimes\id$, where $\id$ is the identity map of the matrix factorization associated to the part of $\Gamma$ not depicted in Figure \ref{reduction3}. By Proposition \ref{hmf-tensor}, it is straightforward to check that 
\[
\Hom_{HMF}(C_p(\Gamma_2'),C_p(\Gamma_2')) \cong H_p(\bigcirc~\bigcirc)\{2n-2\}\cong (H_p(\bigcirc) \otimes H_p(\bigcirc))\{2n-2\}.
\] 
In particular, $\Hom_{hmf}(\Gamma_2',\Gamma_2') \cong \C$, which means any filtered homomorphism from $C_p(\Gamma_2')$ to itself is homotopic to a multiple of the identity map. So $g\circ f$ is homotopic to a multiple of the identity map. Thus, we only need to check that $g\circ f$ is not homotopic to the zero map. We do this by closing $\Gamma_2'$ up to a circle and computing the homomorphism $\Phi$ on the cohomology induced by $g\circ f$ in the special case. Clearly, the homomorphism $\Phi$ is given by Figure \ref{map-phi-def}.

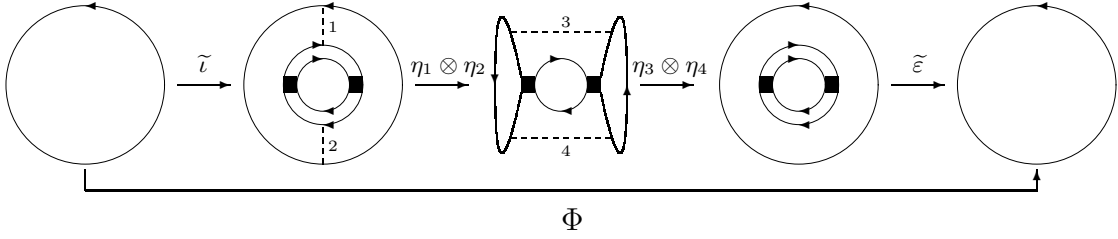
\begin{figure}[ht]

\setlength{\unitlength}{1pt}

\begin{picture}(420,85)(-210,-55)


\put(-180,0){\circle{60}}

\put(-180,30){\vector(-1,0){0}}

\put(-145,0){\vector(1,0){20}}

\put(-137,5){\small{$\widetilde{\iota}$}}


\put(-90,0){\circle{60}}

\put(-90,30){\vector(-1,0){0}}

\put(-90,0){\circle{20}}

\put(-90,10){\vector(1,0){0}}

\put(-90,-10){\vector(-1,0){0}}

\put(-90,0){\circle{30}}

\put(-90,15){\vector(1,0){0}}

\put(-90,-15){\vector(-1,0){0}}

\multiput(-90,15)(0,4){4}{\line(0,1){2}}

\multiput(-90,-30)(0,4){4}{\line(0,1){2}}

\put(-88,20){\tiny{$1$}}

\put(-88,-25){\tiny{$2$}}

\put(-55,0){\vector(1,0){20}}

\put(-57,5){\small{$\eta_1\otimes \eta_2$}}


\put(0,0){\circle{20}}

\put(0,10){\vector(1,0){0}}

\put(0,-10){\vector(-1,0){0}}

\qbezier(-15,3)(-25,50)(-25,0)

\qbezier(-15,-3)(-25,-50)(-25,0)

\put(-25,0){\vector(0,-1){0}}

\qbezier(15,3)(25,50)(25,0)

\qbezier(15,-3)(25,-50)(25,0)

\put(25,0){\vector(0,1){0}}

\multiput(-19,20)(4,0){10}{\line(1,0){2}}

\multiput(-19,-20)(4,0){10}{\line(1,0){2}}

\put(0,22){\tiny{$3$}}

\put(0,-27){\tiny{$4$}}

\put(30,0){\vector(1,0){20}}

\put(27,5){\small{$\eta_3\otimes \eta_4$}}


\put(90,0){\circle{60}}

\put(90,30){\vector(-1,0){0}}

\put(90,0){\circle{20}}

\put(90,10){\vector(1,0){0}}

\put(90,-10){\vector(-1,0){0}}

\put(90,0){\circle{30}}

\put(90,15){\vector(1,0){0}}

\put(90,-15){\vector(-1,0){0}}

\put(125,0){\vector(1,0){20}}

\put(133,5){\small{$\widetilde{\varepsilon}$}}


\put(180,0){\circle{60}}

\put(180,30){\vector(-1,0){0}}


\put(-180,-32){\line(0,-1){8}}

\put(-180,-40){\line(1,0){360}}

\put(180,-40){\vector(0,1){8}}

\put(0,-55){\Large{$\Phi$}}

\linethickness{5pt}

\put(-102.5,-3){\line(0,1){6}}

\put(-77.5,-3){\line(0,1){6}}

\put(-12.5,-3){\line(0,1){6}}

\put(12.5,-3){\line(0,1){6}}

\put(102.5,-3){\line(0,1){6}}

\put(77.5,-3){\line(0,1){6}}

\end{picture}

\caption{Definition of $\Phi$}\label{map-phi-def}

\end{figure}

Note that the composition $(\eta_3\otimes \eta_4)\circ(\eta_1\otimes \eta_2)$ is indeed a tensor product of homomorphisms. So we can change the order in the composition, and get the equivalent definition of $\Phi$ given in Figure \ref{map-phi-def2}.

\begin{figure}[ht]

\setlength{\unitlength}{1pt}

\begin{picture}(420,85)(-210,-55)


\put(-180,0){\circle{60}}

\put(-180,30){\vector(-1,0){0}}

\put(-145,0){\vector(1,0){20}}

\put(-137,5){\small{$\widetilde{\iota}$}}


\put(-90,0){\circle{60}}

\put(-90,30){\vector(-1,0){0}}

\put(-90,0){\circle{20}}

\put(-90,10){\vector(1,0){0}}

\put(-90,-10){\vector(-1,0){0}}

\qbezier(-105,3)(-90,50)(-75,3)

\put(-75,3){\vector(1,-4){0}}

\qbezier(-105,-3)(-90,-50)(-75,-3)

\put(-105,-3){\vector(-1,4){0}}

\multiput(-100,15)(4.5,0){5}{\line(1,0){2}}

\multiput(-100,-15)(4.5,0){5}{\line(1,0){2}}

\put(-92,15){\tiny{$3$}}

\put(-92,-20){\tiny{$4$}}

\put(-55,0){\vector(1,0){20}}

\put(-58,5){\small{$\eta_3\otimes \eta_4$}}


\put(0,0){\circle{60}}

\put(0,30){\vector(-1,0){0}}

\put(0,0){\circle{20}}

\put(0,10){\vector(1,0){0}}

\put(0,-10){\vector(-1,0){0}}

\put(0,0){\circle{30}}

\put(0,15){\vector(1,0){0}}

\put(0,-15){\vector(-1,0){0}}

\multiput(0,22.5)(0,4){2}{\line(0,1){2}}

\multiput(0,-22.5)(0,-4){2}{\line(0,-1){2}}

\put(2,23.5){\tiny{$1$}}

\put(2,-28){\tiny{$2$}}

\put(35,0){\vector(1,0){20}}

\put(31,5){\small{$\eta_1\otimes \eta_2$}}

\put(0,20){\oval(20,5)}

\put(0,-20){\oval(20,5)}


\put(90,0){\circle{60}}

\put(90,30){\vector(-1,0){0}}

\put(90,0){\circle{20}}

\put(90,10){\vector(1,0){0}}

\put(90,-10){\vector(-1,0){0}}

\put(90,0){\circle{30}}

\put(90,15){\vector(1,0){0}}

\put(90,-15){\vector(-1,0){0}}

\put(125,0){\vector(1,0){20}}

\put(133,5){\small{$\widetilde{\varepsilon}$}}


\put(180,0){\circle{60}}

\put(180,30){\vector(-1,0){0}}


\put(-180,-32){\line(0,-1){8}}

\put(-180,-40){\line(1,0){360}}

\put(180,-40){\vector(0,1){8}}

\put(0,-55){\Large{$\Phi$}}

\linethickness{5pt}

\put(-102.5,-3){\line(0,1){6}}

\put(-77.5,-3){\line(0,1){6}}

\put(-12.5,-3){\line(0,1){6}}

\put(12.5,-3){\line(0,1){6}}

\put(102.5,-3){\line(0,1){6}}

\put(77.5,-3){\line(0,1){6}}

\end{picture}

\caption{The second definition of $\Phi$}\label{map-phi-def2}

\end{figure}
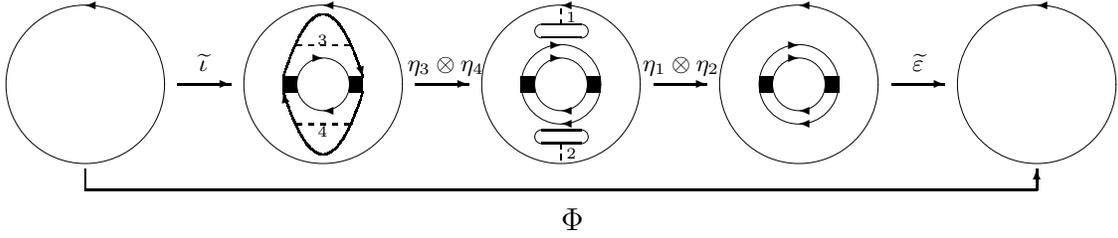

One can see that, in the second definition of $\Phi$, the homomorphisms $\eta_1\otimes \eta_2$ and $\widetilde{\varepsilon}$ commute. So we get the third equivalent definition of $\Phi$ as depicted in Figure \ref{map-phi-def3}.

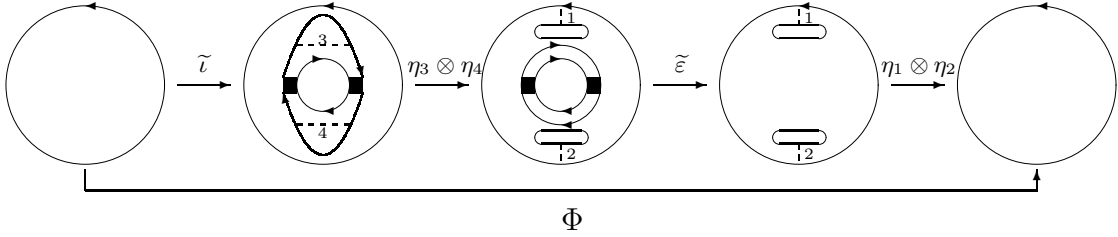
\begin{figure}[ht]

\setlength{\unitlength}{1pt}

\begin{picture}(420,85)(-210,-55)


\put(-180,0){\circle{60}}

\put(-180,30){\vector(-1,0){0}}

\put(-145,0){\vector(1,0){20}}

\put(-137,5){\small{$\widetilde{\iota}$}}


\put(-90,0){\circle{60}}

\put(-90,30){\vector(-1,0){0}}

\put(-90,0){\circle{20}}

\put(-90,10){\vector(1,0){0}}

\put(-90,-10){\vector(-1,0){0}}

\qbezier(-105,3)(-90,50)(-75,3)

\put(-75,3){\vector(1,-4){0}}

\qbezier(-105,-3)(-90,-50)(-75,-3)

\put(-105,-3){\vector(-1,4){0}}

\multiput(-100,15)(4.5,0){5}{\line(1,0){2}}

\multiput(-100,-15)(4.5,0){5}{\line(1,0){2}}

\put(-92,15){\tiny{$3$}}

\put(-92,-20){\tiny{$4$}}

\put(-55,0){\vector(1,0){20}}

\put(-58,5){\small{$\eta_3\otimes \eta_4$}}


\put(0,0){\circle{60}}

\put(0,30){\vector(-1,0){0}}

\put(0,0){\circle{20}}

\put(0,10){\vector(1,0){0}}

\put(0,-10){\vector(-1,0){0}}

\put(0,0){\circle{30}}

\put(0,15){\vector(1,0){0}}

\put(0,-15){\vector(-1,0){0}}

\multiput(0,22.5)(0,4){2}{\line(0,1){2}}

\multiput(0,-22.5)(0,-4){2}{\line(0,-1){2}}

\put(2,23.5){\tiny{$1$}}

\put(2,-28){\tiny{$2$}}

\put(0,20){\oval(20,5)}

\put(0,-20){\oval(20,5)}

\put(35,0){\vector(1,0){20}}

\put(43,5){\small{$\widetilde{\varepsilon}$}}


\put(90,0){\circle{60}}

\put(90,30){\vector(-1,0){0}}

\multiput(90,22.5)(0,4){2}{\line(0,1){2}}

\multiput(90,-22.5)(0,-4){2}{\line(0,-1){2}}

\put(92,23.5){\tiny{$1$}}

\put(92,-28){\tiny{$2$}}

\put(90,20){\oval(20,5)}

\put(90,-20){\oval(20,5)}

\put(125,0){\vector(1,0){20}}

\put(121,5){\small{$\eta_1\otimes \eta_2$}}


\put(180,0){\circle{60}}

\put(180,30){\vector(-1,0){0}}


\put(-180,-32){\line(0,-1){8}}

\put(-180,-40){\line(1,0){360}}

\put(180,-40){\vector(0,1){8}}

\put(0,-55){\Large{$\Phi$}}

\linethickness{5pt}

\put(-102.5,-3){\line(0,1){6}}

\put(-77.5,-3){\line(0,1){6}}

\put(-12.5,-3){\line(0,1){6}}

\put(12.5,-3){\line(0,1){6}}

\end{picture}

\caption{The third definition of $\Phi$}\label{map-phi-def3}

\end{figure}

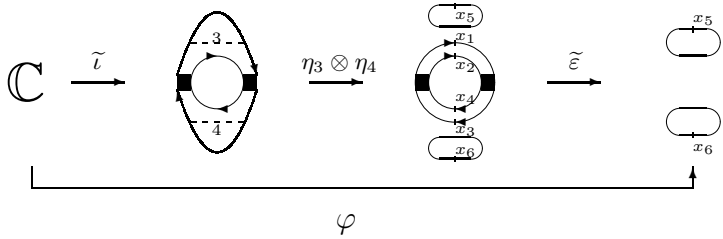
\begin{figure}[ht]

\setlength{\unitlength}{1pt}

\begin{picture}(420,85)(-210,-55)


\put(-170,-7){\Huge{$\C$}}

\put(-145,0){\vector(1,0){20}}

\put(-137,5){\small{$\widetilde{\iota}$}}


\put(-90,0){\circle{20}}

\put(-90,10){\vector(1,0){0}}

\put(-90,-10){\vector(-1,0){0}}

\qbezier(-105,3)(-90,50)(-75,3)

\put(-75,3){\vector(1,-4){0}}

\qbezier(-105,-3)(-90,-50)(-75,-3)

\put(-105,-3){\vector(-1,4){0}}

\multiput(-100,15)(4.5,0){5}{\line(1,0){2}}

\multiput(-100,-15)(4.5,0){5}{\line(1,0){2}}

\put(-92,15){\tiny{$3$}}

\put(-92,-20){\tiny{$4$}}

\put(-55,0){\vector(1,0){20}}

\put(-58,5){\small{$\eta_3\otimes \eta_4$}}


\put(0,0){\circle{20}}

\put(0,10){\vector(1,0){0}}

\put(0,-10){\vector(-1,0){0}}

\put(0,0){\circle{30}}

\put(0,15){\vector(1,0){0}}

\put(0,-15){\vector(-1,0){0}}

\put(0,25){\oval(20,8)}

\put(0,-25){\oval(20,8)}

\put(0,16){\tiny{$x_1$}}

\put(0,14){\line(0,1){2}}

\put(0,5){\tiny{$x_2$}}

\put(0,9){\line(0,1){2}}

\put(0,-19){\tiny{$x_3$}}

\put(0,-14){\line(0,-1){2}}

\put(0,-7){\tiny{$x_4$}}

\put(0,-9){\line(0,-1){2}}

\put(0,24){\tiny{$x_5$}}

\put(0,28){\line(0,1){2}}

\put(0,-27){\tiny{$x_6$}}

\put(0,-28){\line(0,-1){2}}

\put(35,0){\vector(1,0){20}}

\put(43,5){\small{$\widetilde{\varepsilon}$}}


\put(90,15){\oval(20,10)}

\put(90,-15){\oval(20,10)}

\put(90,23){\tiny{$x_5$}}

\put(90,19){\line(0,1){2}}

\put(90,-26){\tiny{$x_6$}}

\put(90,-19){\line(0,-1){2}}


\put(-160,-32){\line(0,-1){8}}

\put(-160,-40){\line(1,0){250}}

\put(90,-40){\vector(0,1){8}}

\put(-45,-55){\Large{$\varphi$}}

\linethickness{5pt}

\put(-102.5,-3){\line(0,1){6}}

\put(-77.5,-3){\line(0,1){6}}

\put(-12.5,-3){\line(0,1){6}}

\put(12.5,-3){\line(0,1){6}}

\end{picture}

\caption{Definition of $\varphi$}\label{map-small-phi}

\end{figure}

For now, we disregard the outer circle in Figure \ref{map-phi-def3}, and consider the homomorphism 
\[
\varphi=\widetilde{\varepsilon}\circ (\eta_3\otimes \eta_4) \circ \widetilde{\iota} : \C \rightarrow H_p(\bigcirc~\bigcirc)
\]
given by Figure \ref{map-small-phi}. This is a homomorphism induced by a chain map $\widetilde{\varphi}$ of degree $2-2n$. By Proposition \ref{hmf-tensor}, we have that 
\[
\Hom_{HMF}(\C,C_p(\bigcirc~\bigcirc))\cong H_p(\bigcirc~\bigcirc).
\] 
In particular, we have that 
\[
\fil^{2-2n}\Hom_{HMF}(\C,C_p(\bigcirc~\bigcirc)) \cong \C.
\] 
Denote by $\hat{\iota}:\C\rightarrow C_p(\bigcirc,\bigcirc)$ the composition of the two $\iota$-maps corresponding to the two circle creations. This is a non-zero element of $\fil^{2-2n}\Hom_{HMF}(\C,C_p(\bigcirc~\bigcirc))$. So $\widetilde{\varphi}$ is homotopic to a multiple of $\hat{\iota}$. 

Next, we prove that $\varphi\neq0$.
\begin{eqnarray*}
\varphi(1) & = & \lambda \cdot \widetilde{\varepsilon}\circ(\eta_3\otimes \eta_4)(1) \\
           & = & \lambda \cdot \widetilde{\varepsilon} ((-2)\cdot(e_{115}+e_{155})\cdot(-2)\cdot(e_{336}+e_{366})) \\
           & = & \lambda \cdot \widetilde{\varepsilon} ((\frac{\partial\pi_{15}}{\partial x_1}+\frac{\partial\pi_{15}}{\partial x_5}) \cdot (\frac{\partial\pi_{36}}{\partial x_3}+\frac{\partial\pi_{36}}{\partial x_6})) \\
           & = & \lambda \cdot \widetilde{\varepsilon}((n+1)^2\cdot x_1^{n-1} \cdot x_3^{n-1}),
\end{eqnarray*}
where $\lambda$ is a non-zero scalar given by the definition of $\widetilde{\iota}$, and, in the last step, we used the definition of $\widetilde{\varepsilon}$. As elements of $H_p(\hat{\Gamma})$, $x_1,~x_2,~x_3,~x_4$ satisfy $x_1+x_2=x_3+x_4$ and $x_1x_2=x_3x_4$. So 
\[
x_3^2-(x_1+x_2)x_3+x_1x_2=x_3^2-(x_3+x_4)x_3+x_3x_4 =(x_3-x_3)(x_3-x_4)=0.
\]
By induction, it's easy to show that 
\[
x_3^k = \frac{x_1^k-x_2^k}{x_1-x_2}x_3-x_1x_2\frac{x_1^{k-1}-x_2^{k-1}}{x_1-x_2}.
\]
Using this and the fact that $p'(x_1)=0$, one can easily check that
\[
x_1^{n-1} \cdot x_3^{n-1} = x_1^{n-1}x_2^{n-2}x_3+ \text{ terms with lower degrees}.
\] 
Thus, by Lemma \ref{H-hat-gamma} and the definition of $\widetilde{\varepsilon}$, 
\[
\widetilde{\varepsilon}(\cdot x_1^{n-1} \cdot x_3^{n-1}) = \widetilde{\varepsilon}(x_1^{n-1}x_2^{n-2}x_3)=\mu\neq0,
\]
where $\mu$ depends on the definition of $\widetilde{\varepsilon}$. Thus, $\varphi(1) = \lambda \mu(n+1)^2\neq0$. So $\varphi$ is a non-zero multiple of $\hat{\iota}$. Then, from Proposition \ref{eta-iota-epsilon}, we know that $\Phi$ is a non-zero multiple of the identity map, which implies that $G\circ F$ is a non-zero multiple of the identity map.
\end{proof}

Next, following the construction in \cite{KR1}, we define a pair of homomorphisms 
\[
\widetilde{\alpha}: \bigoplus_{i=0}^{n-3}H_p(\Gamma_1)\left\langle 1 \right\rangle\{3-n+2i\} \rightarrow H_p(\Gamma)
\]
and 
\[
\widetilde{\beta}: H_p(\Gamma) \rightarrow \bigoplus_{i=0}^{n-3}H_p(\Gamma_1)\left\langle 1 \right\rangle\{3-n+2i\}.
\]

\begin{figure}[ht]

\setlength{\unitlength}{1pt}

\begin{picture}(420,105)(-210,-50)


\qbezier(-200,0)(-175,25)(-200,50)

\put(-200,50){\vector(-1,1){0}}

\qbezier(-150,0)(-175,25)(-150,50)

\put(-150,0){\vector(1,-1){0}}

\put(-178,-25){$\Gamma_1$}

\put(-210,50){\tiny{$x_1$}}

\put(-145,50){\tiny{$x_2$}}

\put(-210,0){\tiny{$x_4$}}

\put(-145,0){\tiny{$x_3$}}

\put(-120,25){\vector(1,0){60}}

\put(-90,30){$\iota$}


\qbezier(-25,0)(0,25)(-25,50)

\put(-25,50){\vector(-1,1){0}}

\qbezier(25,0)(0,25)(25,50)

\put(25,0){\vector(1,-1){0}}

\put(-15,-25){$\Gamma_1\sqcup\bigcirc$}

\put(-35,50){\tiny{$x_1$}}

\put(30,50){\tiny{$x_2$}}

\put(-35,0){\tiny{$x_4$}}

\put(30,0){\tiny{$x_3$}}

\put(0,25){\oval(15,30)}

\put(-2,45){\tiny{$x_5$}}

\put(-2,3){\tiny{$x_6$}}

\put(0,39){\line(0,1){2}}

\put(0,10){\line(0,1){2}}

\put(60,25){\vector(1,0){60}}

\put(75,30){$\chi_0\otimes\chi_0$}


\qbezier(165,25)(175,35)(185,25)

\put(185,25){\vector(1,-1){0}}

\qbezier(165,15)(175,5)(185,15)

\put(165,15){\vector(-1,1){0}}

\put(160,25){\vector(-1,1){15}}

\put(190,15){\vector(1,-1){15}}

\put(145,0){\vector(1,1){15}}

\put(205,40){\vector(-1,-1){15}}

\put(135,40){\tiny{$x_1$}}

\put(210,40){\tiny{$x_2$}}

\put(135,0){\tiny{$x_4$}}

\put(210,0){\tiny{$x_3$}}

\put(173,34){\tiny{$x_5$}}

\put(173,3){\tiny{$x_6$}}

\put(175,29){\line(0,1){2}}

\put(175,9){\line(0,1){2}}

\put(172,-25){$\Gamma$}


\put(-175,-30){\line(0,-1){8}}

\put(-175,-38){\line(1,0){350}}

\put(175,-38){\vector(0,1){8}}

\put(0,-50){\Large{$\alpha$}}

\linethickness{5pt}

\put(162.5,15){\line(0,1){10}}

\put(187.5,15){\line(0,1){10}}

\end{picture}

\caption{Homomorphism $\alpha$, where the two $\chi_0$ maps are the one associated to the two wide edges in $\Gamma$.}\label{maps-alpha}
\end{figure}
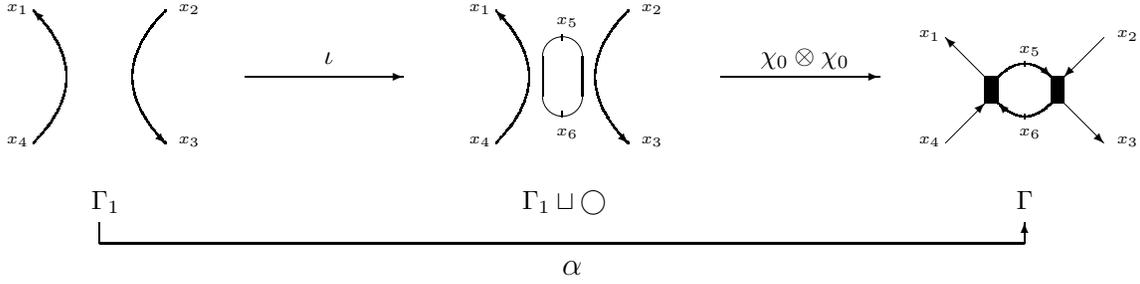

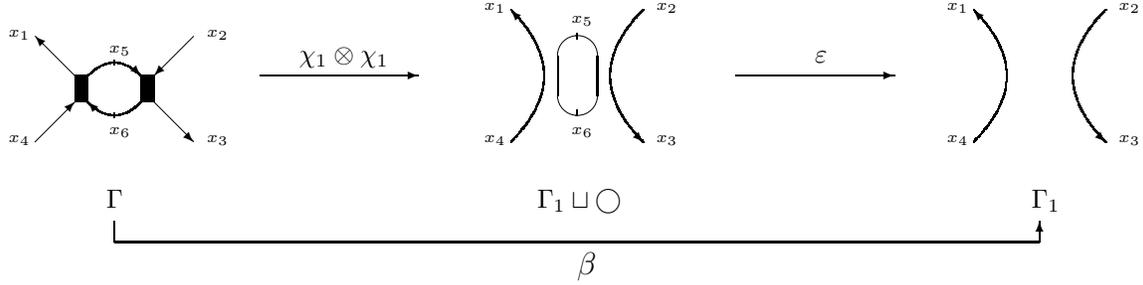
\begin{figure}[ht]

\setlength{\unitlength}{1pt}

\begin{picture}(420,105)(-210,-50)


\qbezier(150,0)(175,25)(150,50)

\put(150,50){\vector(-1,1){0}}

\qbezier(200,0)(175,25)(200,50)

\put(200,0){\vector(1,-1){0}}

\put(172,-25){$\Gamma_1$}

\put(140,50){\tiny{$x_1$}}

\put(205,50){\tiny{$x_2$}}

\put(140,0){\tiny{$x_4$}}

\put(205,0){\tiny{$x_3$}}


\qbezier(-25,0)(0,25)(-25,50)

\put(-25,50){\vector(-1,1){0}}

\qbezier(25,0)(0,25)(25,50)

\put(25,0){\vector(1,-1){0}}

\put(-15,-25){$\Gamma_1\sqcup\bigcirc$}

\put(-35,50){\tiny{$x_1$}}

\put(30,50){\tiny{$x_2$}}

\put(-35,0){\tiny{$x_4$}}

\put(30,0){\tiny{$x_3$}}

\put(0,25){\oval(15,30)}

\put(-2,45){\tiny{$x_5$}}

\put(-2,3){\tiny{$x_6$}}

\put(0,39){\line(0,1){2}}

\put(0,10){\line(0,1){2}}

\put(60,25){\vector(1,0){60}}

\put(90,30){$\varepsilon$}


\qbezier(-185,25)(-175,35)(-165,25)

\put(-165,25){\vector(1,-1){0}}

\qbezier(-185,15)(-175,5)(-165,15)

\put(-185,15){\vector(-1,1){0}}

\put(-190,25){\vector(-1,1){15}}

\put(-160,15){\vector(1,-1){15}}

\put(-205,0){\vector(1,1){15}}

\put(-145,40){\vector(-1,-1){15}}

\put(-215,40){\tiny{$x_1$}}

\put(-140,40){\tiny{$x_2$}}

\put(-215,0){\tiny{$x_4$}}

\put(-140,0){\tiny{$x_3$}}

\put(-177,34){\tiny{$x_5$}}

\put(-177,3){\tiny{$x_6$}}

\put(-175,29){\line(0,1){2}}

\put(-175,9){\line(0,1){2}}

\put(-178,-25){$\Gamma$}

\put(-120,25){\vector(1,0){60}}

\put(-105,30){$\chi_1\otimes\chi_1$}


\put(-175,-30){\line(0,-1){8}}

\put(-175,-38){\line(1,0){350}}

\put(175,-38){\vector(0,1){8}}

\put(0,-50){\Large{$\beta$}}

\linethickness{5pt}

\put(-187.5,15){\line(0,1){10}}

\put(-162.5,15){\line(0,1){10}}

\end{picture}

\caption{Homomorphism $\alpha$, where the two $\chi_1$ maps are the one associated to the two wide edges in $\Gamma$.}\label{maps-beta}
\end{figure}

First, we define by Figures \ref{maps-alpha} and \ref{maps-beta} the homomorphisms
\[
\alpha: H_p(\Gamma_1)\left\langle 1 \right\rangle\{3-n\} \rightarrow H_p(\Gamma)
\]
and 
\[
\beta: H_p(\Gamma) \rightarrow H_p(\Gamma_1)\left\langle 1 \right\rangle\{3-n\}.
\]

Then, for $0\leq i\leq n-3$, define 
\[
\alpha_i: H_p(\Gamma_1)\left\langle 1 \right\rangle\{n-3-2i\} \rightarrow H_p(\Gamma) 
\]
by
\[ 
\alpha_i = m(x_5^{n-3-i})\circ\alpha,
\]
and 
\[
\beta_i: H_p(\Gamma) \rightarrow H_p(\Gamma_1)\left\langle 1 \right\rangle\{n-3-2i\}
\]
by
\[ 
\beta_i = \beta\circ m(\sum_{k=0}^{i}(n-k+1)c_{n-k+1}\frac{(x_3-x_5)x_1^{i-k+2}+(x_5-x_1)x_3^{i-k+2}+(x_1-x_3)x_5^{i-k+2}}{(x_5-x_3)(x_1-x_5)(x_3-x_1)}),
\]
where the big fraction is in fact a polynomial, and $c_{n-k+1}$ is the $(n-k+1)$-th coefficient of $p(x)$. Similar to the computation in the proof of Proposition \ref{moy1}, we have that, for $0\leq i,j \leq n-3$,
\[
\beta_i\circ\alpha_j = \varepsilon\circ m(\sum_{k=0}^{i}(n-k+1)c_{n-k+1}x_5^{n-k+i-j-1}) \circ \iota = \nu(n+1)\delta_{ij} \id,
\]
where $\nu$ is a non-zero scalar coming from the definitions of $\iota$ and $\varepsilon$.
Now define $\widetilde{\alpha} = \sum_{i=0}^{n-3}\alpha_i$ and $\widetilde{\beta} = \frac{1}{\nu(n+1)} \sum_{i=0}^{n-3}\beta_i$. Then $\widetilde{\beta} \circ \widetilde{\alpha} = \id_{\bigoplus_{i=0}^{n-3}H_p(\Gamma_1)\left\langle 1 \right\rangle\{3-n+2i\}}$.

\begin{lemma}\label{alpha-beta-F-G}
$\widetilde{\beta} \circ F=0$ and $G \circ \widetilde{\alpha} =0$.
\end{lemma}
\begin{proof}
Let $\Gamma_1'$ and $\Gamma_2'$ be the parts of $\Gamma_1$ and $\Gamma_2$ depicted in Figure \ref{reduction3}. Then both $\widetilde{\beta} \circ F$ and $G \circ \widetilde{\alpha}$ are induced by filtered matrix factorization homomorphisms between $\bigoplus_{i=0}^{n-3}C_p(\Gamma_1')\left\langle 1 \right\rangle\{3-n+2i\}$ and $C_p(\Gamma_2')$. Using Proposition \ref{hmf-tensor}, it is easy to check that 
\[
\Hom_{HMF}(C_p(\Gamma_1'),C_p(\Gamma_2'))\cong \Hom_{HMF}(C_p(\Gamma_2'),C_p(\Gamma_1')) \cong H_p(\bigcirc)\{2n-2\}.
\]
In particular, this implies that 
\[
\Hom_{hmf}(\bigoplus_{i=0}^{n-3}C_p(\Gamma_1')\left\langle 1 \right\rangle\{3-n+2i\},C_p(\Gamma_2')) \cong \Hom_{hmf}(C_p(\Gamma_2'),\bigoplus_{i=0}^{n-3}C_p(\Gamma_1')\left\langle 1 \right\rangle\{3-n+2i\}) =0.
\]
So any filtered matrix factorization homomorphisms between $\bigoplus_{i=0}^{n-3}C_p(\Gamma_1')\left\langle 1 \right\rangle\{3-n+2i\}$ and $C_p(\Gamma_2')$ is homotopic to the zero map. This implies that $\widetilde{\beta} \circ F=0$ and $G \circ \widetilde{\alpha} =0$.
\end{proof}

\begin{proof}[Proof of Proposition \ref{moy3}]
Define
\begin{eqnarray*}
\widetilde{F}=(F,\widetilde{\alpha}) & : & H_p(\Gamma_2) \oplus (\bigoplus_{i=0}^{n-3}H_p(\Gamma_1)\left\langle 1 \right\rangle\{3-n+2i\}) \rightarrow H_p(\Gamma), \\
\widetilde{G}=(\kappa \cdot G,\widetilde{\beta})^{T} & : & H_p(\Gamma) \rightarrow H_p(\Gamma_2) \oplus (\bigoplus_{i=0}^{n-3}H_p(\Gamma_1)\left\langle 1 \right\rangle\{3-n+2i\}),
\end{eqnarray*}
where $\kappa\in \C$ satisfies $\kappa \cdot G\circ F=\id_{H_p(\Gamma_2)}$. Then, by the above arguments, we know that 
\[
\widetilde{G}\circ\widetilde{F}=\id_{H_p(\Gamma_2) \oplus (\bigoplus_{i=0}^{n-3}H_p(\Gamma_1)\left\langle 1 \right\rangle\{3-n+2i\})}.
\]
From Proposition 31 of \cite{KR1} and Proposition \ref{fil-grade}, we know that $H_p(\Gamma)$ and $H_p(\Gamma_2) \oplus (\bigoplus_{i=0}^{n-3}H_p(\Gamma_1)\left\langle 1 \right\rangle\{3-n+2i\})$ have the same filtered dimension. So $\widetilde{F}$ and $\widetilde{G}$ are filtered isomorphisms that are inverses of each other. This completes the proof.
\end{proof}

\subsection{MOY decomposition IV}

\begin{figure}[ht]

\setlength{\unitlength}{1pt}

\begin{picture}(420,120)(-210,-20)


\linethickness{.25pt}

\put(-210,30){\vector(0,1){30}}

\put(-170,0){\vector(0,1){30}}

\put(-202.5,10){\vector(1,1){0}}

\qbezier(-210,0)(-210,5)(-202.5,10)

\qbezier(-210,30)(-210,25)(-202.5,20)

\put(-197.5,10){\vector(-1,1){0}}

\qbezier(-190,0)(-190,5)(-197.5,10)

\qbezier(-190,30)(-190,25)(-197.5,20)

\qbezier(-190,30)(-190,35)(-182.5,40)

\qbezier(-190,60)(-190,55)(-182.5,50)

\put(-190,60){\vector(-1,4){0}}

\qbezier(-170,30)(-170,35)(-177.5,40)

\qbezier(-170,60)(-170,55)(-177.5,50)

\qbezier(-210,60)(-210,65)(-202.5,70)

\qbezier(-210,90)(-210,85)(-202.5,80)

\put(-210,90){\vector(-1,4){0}}

\qbezier(-190,60)(-190,65)(-197.5,70)

\qbezier(-190,90)(-190,85)(-197.5,80)

\put(-190,90){\vector(1,4){0}}

\put(-190,30){\vector(1,4){0}}

\put(-170,60){\vector(0,1){30}}

\put(-190,-20){$\Gamma_1$}

\put(-210,92){$x_6$}

\put(-190,92){$x_5$}

\put(-170,92){$x_4$}

\put(-210,-8){$x_3$}

\put(-190,-8){$x_2$}

\put(-170,-8){$x_1$}

\put(-208,42){$x_8$}

\put(-211,45){\line(1,0){2}}

\put(-188,57){$x_9$}

\put(-191,60){\line(1,0){2}}

\put(-188,27){$x_7$}

\put(-191,30){\line(1,0){2}}

\linethickness{5pt}

\put(-200,10){\line(0,1){10}}

\put(-180,40){\line(0,1){10}}

\put(-200,70){\line(0,1){10}}


\linethickness{.25pt}

\put(-90,0){\vector(0,1){90}}

\put(-50,0){\vector(0,1){30}}

\put(-70,0){\vector(0,1){30}}

\qbezier(-50,30)(-50,35)(-57.5,40)

\qbezier(-50,60)(-50,55)(-57.5,50)

\qbezier(-70,30)(-70,35)(-62.5,40)

\qbezier(-70,60)(-70,55)(-62.5,50)

\put(-50,60){\vector(0,1){30}}

\put(-70,60){\vector(0,1){30}}

\put(-70,-20){$\Gamma_2$}

\put(-90,92){$x_6$}

\put(-70,92){$x_5$}

\put(-50,92){$x_4$}

\put(-90,-8){$x_3$}

\put(-70,-8){$x_2$}

\put(-50,-8){$x_1$}

\linethickness{5pt}

\put(-60,40){\line(0,1){10}}


\linethickness{.25pt}

\put(90,30){\vector(0,1){30}}

\put(50,0){\vector(0,1){30}}

\put(82.5,10){\vector(-1,1){0}}

\qbezier(90,0)(90,5)(82.5,10)

\qbezier(90,30)(90,25)(82.5,20)

\put(77.5,10){\vector(1,1){0}}

\qbezier(70,0)(70,5)(77.5,10)

\qbezier(70,30)(70,25)(77.5,20)

\qbezier(70,30)(70,35)(62.5,40)

\qbezier(70,60)(70,55)(62.5,50)

\put(70,60){\vector(1,4){0}}

\qbezier(50,30)(50,35)(57.5,40)

\qbezier(50,60)(50,55)(57.5,50)

\qbezier(90,60)(90,65)(82.5,70)

\qbezier(90,90)(90,85)(82.5,80)

\put(90,90){\vector(1,4){0}}

\qbezier(70,60)(70,65)(77.5,70)

\qbezier(70,90)(70,85)(77.5,80)

\put(70,90){\vector(-1,4){0}}

\put(70,30){\vector(-1,4){0}}

\put(50,60){\vector(0,1){30}}

\put(70,-20){$\Gamma_3$}

\put(50,92){$x_6$}

\put(70,92){$x_5$}

\put(90,92){$x_4$}

\put(50,-8){$x_3$}

\put(70,-8){$x_2$}

\put(90,-8){$x_1$}

\put(92,42){$x_8$}

\put(89,45){\line(1,0){2}}

\put(72,57){$x_9$}

\put(69,60){\line(1,0){2}}

\put(72,27){$x_7$}

\put(69,30){\line(1,0){2}}

\linethickness{5pt}

\put(80,10){\line(0,1){10}}

\put(60,40){\line(0,1){10}}

\put(80,70){\line(0,1){10}}


\linethickness{.25pt}

\put(210,0){\vector(0,1){90}}

\put(170,0){\vector(0,1){30}}

\put(190,0){\vector(0,1){30}}

\qbezier(170,30)(170,35)(177.5,40)

\qbezier(170,60)(170,55)(177.5,50)

\qbezier(190,30)(190,35)(182.5,40)

\qbezier(190,60)(190,55)(182.5,50)

\put(170,60){\vector(0,1){30}}

\put(190,60){\vector(0,1){30}}

\put(190,-20){$\Gamma_4$}

\put(170,92){$x_6$}

\put(190,92){$x_5$}

\put(210,92){$x_4$}

\put(170,-8){$x_3$}

\put(190,-8){$x_2$}

\put(210,-8){$x_1$}

\linethickness{5pt}

\put(180,40){\line(0,1){10}}

\end{picture}

\caption{$\Gamma_1$, $\Gamma_2$,
$\Gamma_3$ and $\Gamma_4$}\label{reduction4}

\end{figure}
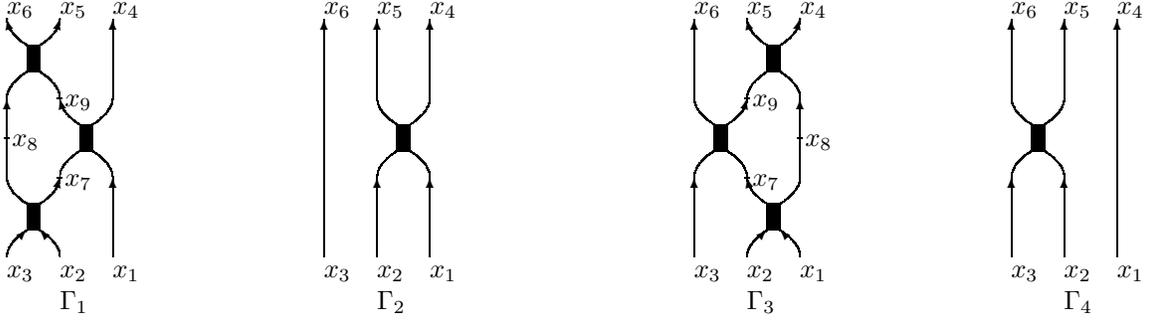

\begin{proposition}\label{moy4}
Let $\Gamma_1,\Gamma_2,\Gamma_3,\Gamma_4$ be four closed graphs that are identical except in the part depicted in Figure \ref{reduction4}. Then there is an isomorphism of filtered $\C$-linear spaces
\[
H_p(\Gamma_1) \oplus H_p(\Gamma_2) \cong H_p(\Gamma_3) \oplus H_p(\Gamma_4).
\]
\end{proposition}

Similar to \cite{KR1}, we consider the matrix factorization associated to the planar diagram $\Gamma$ that is identical to $\Gamma_1,\Gamma_2,\Gamma_3,\Gamma_4$ except in the part depicted in Figure \ref{3way}. 

\begin{figure}[ht]

\setlength{\unitlength}{1pt}

\begin{picture}(420,125)(-210,-10)


\put(-2.5,60){\vector(-1,1){40}}

\put(0,60){\vector(0,1){40}}

\put(2.5,60){\vector(1,1){40}}

\put(-45,102){$x_6$}

\put(-2,102){$x_5$}

\put(42.5,102){$x_4$}


\put(-42.5,0){\vector(1,1){40}}

\put(0,0){\vector(0,1){40}}

\put(42.5,0){\vector(-1,1){40}}

\put(-45,-8){$x_3$}

\put(-2,-8){$x_2$}

\put(42.5,-8){$x_1$}

\linethickness{5pt}

\put(0,40){\line(0,1){20}}

\end{picture}

\caption{$\Gamma$}\label{3way}

\end{figure}

Since $p(x)+p(y)+p(z)$ is symmetric in $x,y,z$, there is a unique three variable polynomial $h$ such that $h(x+y+z,xy+yz+zx,xyz)=p(x)+p(y)+p(z)$. Let 
\begin{eqnarray*}
s_1 & = & x_4+x_5+x_6, \\
s_2 & = & x_4x_5+x_5x_6+x_6x_4, \\
s_3 & = & x_4x_5x_6,\\
s_4 & = & x_1+x_2+x_3, \\
s_5 & = & x_1x_2+x_2x_3+x_3x_1, \\
s_6 & = & x_1x_2x_3.
\end{eqnarray*}
Then define
\begin{eqnarray*}
a_1 & = & \frac{h(s_1,s_2,s_3)-h(s_4,s_2,s_3)}{s_1-s_4}, \\
a_2 & = & \frac{h(s_4,s_2,s_3)-h(s_4,s_5,s_3)}{s_2-s_5}, \\
a_3 & = & \frac{h(s_4,s_5,s_3)-h(s_4,s_5,s_6)}{s_3-s_6}, \\
b_1 & = & s_1-s_4 = x_4+x_5+x_6-x_1-x_2-x_3, \\
b_2 & = & s_2-s_5 = x_4x_5+x_5x_6+x_6x_4-x_1x_2-x_2x_3-x_3x_1, \\
b_3 & = & s_3-s_6 = x_4x_5x_6-x_1x_2x_3.
\end{eqnarray*}
Let $R$ be the polynomial ring generated by the markings on $\Gamma$. The matrix factorization associated to the part of $\Gamma$ depicted in Figure \ref{3way} is 
\[
\Upsilon = 
\left(%
\begin{array}{cc}
  a_1 & b_1\\
  a_2 & b_2 \\
  a_3 & b_3 \\
\end{array}%
\right)_{R}\{-3\},
\]
which has potential $p(x_4)+p(x_5)+p(x_6)-p(x_1)-p(x_2)-p(x_3)$. Thus, $C_p(\Gamma)$ has potential $0$, and is therefore a  $\zed_2$-complex. Let $H_p(\Gamma)$ be its cohomology. Proposition \ref{moy4} follows from the lemma below.

\begin{lemma}\label{Upsilon-decomp}
\begin{eqnarray*}
H_p(\Gamma_1) & = & H_p(\Gamma) \oplus H_p(\Gamma_4), \\
H_p(\Gamma_3) & = & H_p(\Gamma) \oplus H_p(\Gamma_2).
\end{eqnarray*}
\end{lemma}
\begin{proof}
We only prove that
\[
H_p(\Gamma_3) = H_p(\Gamma) \oplus H_p(\Gamma_2).
\]
The other decomposition follows similarly. The proof is due to Rasmussen \cite{Ras2}, which is more explicit than Khovanov and Rozansky's original approach in \cite{KR1}. We only need to check the transformations in Rasmussen's proof preserve the quantum filtration.

Let $M$ be the matrix factorization associated to the part of $\Gamma_3$ shown in Figure \ref{reduction4}, and $N$ the matrix factorization associated to the part of $\Gamma_3$ outside Figure \ref{reduction4}. Let $R=\C[x_1,\cdots,x_9]$. Then 
\[
M=
\left(%
\begin{array}{cc}
  u_{7821} & x_7+x_8-x_2-x_1\\
  v_{7821} & x_7x_8-x_2x_1 \\
  u_{6937} & x_6+x_9-x_3-x_7\\
  v_{6937} & x_6x_9-x_3x_7 \\
  u_{5498} & x_5+x_4-x_9-x_8\\
  v_{5498} & x_5x_4-x_9x_8 \\
\end{array}%
\right)_R\{-3\},
\]
where
\begin{eqnarray*}
u_{ijkl} & = & u(x_i,x_j,x_k,x_l), \\
v_{ijkl} & = & v(x_i,x_j,x_k,x_l).
\end{eqnarray*}
By Corollary \ref{row-op}, we have that

\[
M \cong
\left(%
\begin{array}{cc}
  u_{7821}+x_2v_{7821} & x_7+x_8-x_2-x_1\\
  v_{7821} & (x_2-x_8)(x_2-x_7) \\
  u_{6937}+x_3v_{6937} & x_6+x_9-x_3-x_7\\
  v_{6937} & (x_3-x_9)(x_3-x_6) \\
  u_{5498}+x_9v_{5498} & x_5+x_4-x_9-x_8\\
  v_{5498} & (x_9-x_4)(x_9-x_5) \\
\end{array}%
\right)_R\{-3\}.
\]
Then, using Proposition \ref{variable-exclusion}, we exclude variables $x_1$, $x_7$ and $x_8$, and get
\begin{eqnarray*}
H_p(\Gamma_3)
& = & H(N\otimes M) \\
& \cong & H(N \otimes \left(%
\begin{array}{cc}
  v_{7821} & (x_2+x_9-x_4-x_5)(x_2+x_3-x_6-x_9) \\
  v_{6937} & (x_3-x_9)(x_3-x_6) \\
  v_{5498} & (x_9-x_4)(x_9-x_5) \\
\end{array}%
\right)_{R'})\{-3\},
\end{eqnarray*}
where $R'=R/I$ and $I=(x_7+x_8-x_2-x_1,x_6+x_9-x_3-x_7,x_5+x_4-x_9-x_8)$. Note that $R'$ is isomorphic to the polynomial ring over $\C$ generated by all the marking on $\Gamma_3$ except $x_1$, $x_7$ and $x_8$. Now apply Proposition \ref{variable-exclusion} to the last row to eliminate $x_9^2$, which gives
\[
H_p(\Gamma_3) \cong H(N \otimes \left(%
\begin{array}{cc}
  v_{7821} & (x_3-x_6)x_9+x_4x_5+(x_2-x_4-x_5)(x_2+x_3-x_6) \\
  v_{6937} & (x_3-x_9)(x_3-x_6) \\
\end{array}%
\right)_{R''})\{-3\},
\]
where $R''=R'/(x_9-x_4)(x_9-x_5)R'$. Applying Corollary \ref{row-op} again, we get
\begin{eqnarray*}
& &\left(%
 \begin{array}{cc}
  v_{7821} & (x_3-x_6)x_9+x_4x_5+(x_2-x_4-x_5)(x_2+x_3-x_6) \\
  v_{6937} & (x_3-x_9)(x_3-x_6) \\
\end{array}%
\right)_{R''} \\
& \cong &
\left(%
\begin{array}{cc}
  v_{7821} & (x_3-x_6)x_3+x_4x_5+(x_2-x_4-x_5)(x_2+x_3-x_6) \\
  v_{6937}-v_{7821} & (x_3-x_9)(x_3-x_6) \\
\end{array}%
\right)_{R''} = M''.
\end{eqnarray*}
Denote by $M''$ the last matrix factorization above. Let 
\begin{eqnarray*}
b_1 & = & (x_3-x_6)x_3+x_4x_5+(x_2-x_4-x_5)(x_2+x_3-x_6), \\
b_2 & = & (x_3-x_9)(x_3-x_6).
\end{eqnarray*}
Then $M''$ is
\[
\left.%
\begin{array}{c}
  R'' \\
  \oplus \\
  R''\{6-2n\}
\end{array}%
\right.
\xrightarrow{\left(%
\begin{array}{cc}
  v_{7821} & b_2 \\
  v_{6931}-v_{7821} & -b_1 
\end{array}%
\right)}
\left.%
\begin{array}{c}
  R''\{3-n\} \\
  \oplus \\
  R''\{3-n\}
\end{array}%
\right.
\xrightarrow{\left(%
\begin{array}{cc}
  b_1 & b_2 \\
  v_{6931}-v_{7821}& -v_{7821} 
\end{array}%
\right)}
\left.%
\begin{array}{c}
  R'' \\
  \oplus \\
  R''\{6-2n\}
\end{array}%
\right..
\]
Let $R_0=\C[x_2\cdots,x_6]$. Then 
\[
R'' = R_0 \oplus (x_9-x_3)R_0 =R_0 \oplus (x_9+x_3-x_4-x_5)R_0.
\]
Using these two decompositions of $R''$, we rewrite $M''$ as
{\tiny
\[
\left.%
\begin{array}{c}
  R_0 \\
  \oplus \\
  (x_9-x_3)R_0 \\
  \oplus \\
  R_0\{6-2n\} \\
  \oplus \\
  (x_9+x_3-x_4-x_5)R_0\{6-2n\}
\end{array}%
\right.
\xrightarrow{\left(%
\begin{array}{cc}
  A_1 & B_2 \\
  A_2 & -B_1 
\end{array}%
\right)}
\left.%
\begin{array}{c}
  R_0\{3-n\} \\
  \oplus \\
  (x_9-x_3)R_0\{3-n\} \\
  \oplus \\
  R_0\{3-n\} \\
  \oplus \\
  (x_9+x_3-x_4-x_5)R_0\{3-n\}
\end{array}%
\right.
\xrightarrow{\left(%
\begin{array}{cc}
  B_1' & B_2' \\
  A_2'& -A_1' 
\end{array}%
\right)}
\left.%
\begin{array}{c}
  R_0 \\
  \oplus \\
  (x_9-x_3)R_0 \\
  \oplus \\
  R_0\{6-2n\} \\
  \oplus \\
  (x_9+x_3-x_4-x_5)R_0\{6-2n\}
\end{array}%
\right.,
\]}where $A_1$, $A_2$, $B_1$, $B_2$, $A_1'$, $A_2'$, $B_1'$, $B_2'$ are $2\times2$ matrices with entries in $R_0$ representing corresponding linear maps (multiplications) in the previous diagram. Clearly, among these matrices, $X$ and $X'$ represent the same linear map under different basis. Since $x_9$ does not appear in $b_1$, one can see
\[
B_1=B_1'=
\left(%
\begin{array}{cc}
  b_1 & 0 \\
  0& b_1 
\end{array}%
\right).
\]
It is also easy to check that
\[
B_2=B_2'=
\left(%
\begin{array}{cc}
  0 & y \\
  z& 0 
\end{array}%
\right),
\]
where $y=(x_3-x_4)(x_3-x_5)(x_3-x_6)$ and $z=x_6-x_3$.

Using the equation $x_9^2=(x_4+x_5)x_9-x_4x_5$ in $R''$ repeatedly, we find polynomials $\alpha,\beta\in R_0$ with $\deg \alpha \leq 2n-2$, $\deg \beta \leq 2n-4$ and $v_{7821}=\alpha+(x_9-x_3)\beta$ in $R''$. Then
\[
A_1=
\left(%
\begin{array}{cc}
  \alpha & -(x_3-x_4)(x_3-x_5)\beta \\
  \beta & \alpha + (x_4+x_5-2x_3)\beta 
\end{array}%
\right)
=
\left(%
\begin{array}{cc}
  \alpha & \frac{y\beta}{z} \\
  \beta & \alpha + q\beta \end{array}%
\right),
\]
where $q=x_4+x_5-2x_3$, $y$ and $z$ are defined above. After changing the basis, we have 
\[
A_1'=\left(%
\begin{array}{cc}
  1 & q \\
  0 & 1
\end{array}%
\right)
\left(%
\begin{array}{cc}
  \alpha & \frac{y\beta}{z} \\
  \beta & \alpha + q\beta 
\end{array}%
\right)
\left(%
\begin{array}{cc}
  1 & -q \\
  0 & 1
\end{array}%
\right)
=
\left(%
\begin{array}{cc}
  \alpha+q\beta & \frac{y\beta}{z} \\
  \beta & \alpha
\end{array}%
\right).
\]
Similarly, there are polynomials $\gamma,\delta\in R_0$ with $\deg \gamma \leq 2n-2$, $\deg \delta \leq 2n-4$ and $v_{6931}-v_{7821}=\gamma+(x_9+x_3-x_4-x_5)\delta$ in $R''$. Then
\[
A_2=A_2'=
\left(%
\begin{array}{cc}
  \gamma & q\gamma+\frac{y\delta}{z} \\
  \delta & \gamma
\end{array}%
\right).
\]
Now use the fact that
\[
\left(%
\begin{array}{cc}
  A_1 & B_2 \\
  A_2 & -B_1 
\end{array}%
\right)
\left(%
\begin{array}{cc}
  B_1' & B_2' \\
  A_2'& -A_1' 
\end{array}%
\right)
=
\left(%
\begin{array}{cccc}
  w & 0&0&0 \\
  0 & w&0&0 \\
  0 & 0&w&0 \\
  0 & 0&0&w \\
\end{array}%
\right),
\]
where $w=p(x_4)+p(x_5)+p(x_6)-p(x_2)-p(x_3)-p(x_4+x_5+x_6-x_2-x_3)\in R_0$. Comparing the $(2,1)$-entries on both sides, we get $b_1\beta+z\gamma=0$ in $R_0$. But $b_1$ and $z$ are relatively prime. So there exists a polynomial $\rho\in R_0$ with $\deg\rho=\deg\beta-2$ such that $\beta=z\rho$ and, therefore, $\gamma=-b_1\rho$.

Define
\begin{eqnarray*}
h_0 & : & \left.%
\begin{array}{c}
  R_0 \\
  \oplus \\
  (x_9-x_3)R_0 \\
  \oplus \\
  R_0\{6-2n\} \\
  \oplus \\
  (x_9+x_3-x_4-x_5)R_0\{6-2n\}
\end{array}%
\right.
\rightarrow
\left.%
\begin{array}{c}
  R_0 \\
  \oplus \\
  (x_9-x_3)R_0 \\
  \oplus \\
  R_0\{6-2n\} \\
  \oplus \\
  (x_9+x_3-x_4-x_5)R_0\{6-2n\}
\end{array}%
\right., \\
h_1 & : &
\left.%
\begin{array}{c}
  R_0\{3-n\} \\
  \oplus \\
  (x_9-x_3)R_0\{3-n\} \\
  \oplus \\
  R_0\{3-n\} \\
  \oplus \\
  (x_9+x_3-x_4-x_5)R_0\{3-n\}
\end{array}%
\right.
\rightarrow
\left.%
\begin{array}{c}
  R_0\{3-n\} \\
  \oplus \\
  (x_9-x_3)R_0\{3-n\} \\
  \oplus \\
  R_0\{3-n\} \\
  \oplus \\
  (x_9+x_3-x_4-x_5)R_0\{3-n\}
\end{array}%
\right.
\end{eqnarray*}
by 
\begin{eqnarray*}
h_0 & = & \left(%
\begin{array}{cccc}
  0 & 0&0&0 \\
  0 & 0&0&0 \\
  -\rho & 0&0&0 \\
  0 & -\rho&0&0 \\
\end{array}%
\right), \\
h_1 & = & 0.
\end{eqnarray*}
Apply Lemma \ref{general-twist} on $M''$ and $h_0$, $h_1$. We get $M''\cong\hat{M}$, where $\hat{M}$ is the filtered matrix factorization
\[
\left.%
\begin{array}{c}
  R_0 \\
  \oplus \\
  (x_9-x_3)R_0 \\
  \oplus \\
  R_0\{6-2n\} \\
  \oplus \\
  (x_9+x_3-x_4-x_5)R_0\{6-2n\}
\end{array}%
\right.
\xrightarrow{\hat{d}_0}
\left.%
\begin{array}{c}
  R_0\{3-n\} \\
  \oplus \\
  (x_9-x_3)R_0\{3-n\} \\
  \oplus \\
  R_0\{3-n\} \\
  \oplus \\
  (x_9+x_3-x_4-x_5)R_0\{3-n\}
\end{array}%
\right.
\xrightarrow{\hat{d}_1}
\left.%
\begin{array}{c}
  R_0 \\
  \oplus \\
  (x_9-x_3)R_0 \\
  \oplus \\
  R_0\{6-2n\} \\
  \oplus \\
  (x_9+x_3-x_4-x_5)R_0\{6-2n\}
\end{array}%
\right.
\]
with
\begin{eqnarray*}
\hat{d}_0 & = & \left(%
\begin{array}{cccc}
  \alpha & 0&0&y \\
  0 & \alpha+q\beta&z&0 \\
  0 & q\gamma+\frac{y\delta}{z}&-b_1&0 \\
  \delta & 0&0&-b_1 \\
\end{array}%
\right) \\
\hat{d}_1 & = & \left(%
\begin{array}{cccc}
  b_1 & 0&0&y \\
  0 & b_1&z&0 \\
  0 & q\gamma+\frac{y\delta}{z}&-(\alpha+q\beta)&0 \\
  \delta & 0&0&-\alpha \\
\end{array}%
\right).
\end{eqnarray*}
So $\hat{M}\cong M_1\oplus M_2$, where $M_1$ is
{\tiny \[
\left.%
\begin{array}{c}
  R_0 \\
  \oplus \\
  (x_9+x_3-x_4-x_5)R_0\{6-2n\}
\end{array}%
\right.
\xrightarrow{\left(%
\begin{array}{cc}
  \alpha &y \\
  \delta &-b_1 \\
\end{array}%
\right)}
\left.%
\begin{array}{c}
  R_0\{3-n\} \\
  \oplus \\
  (x_9+x_3-x_4-x_5)R_0\{3-n\}
\end{array}%
\right.
\xrightarrow{\left(%
\begin{array}{cc}
  b_1 & y \\
  \delta & -\alpha \\
\end{array}%
\right)}
\left.%
\begin{array}{c}
  R_0 \\
  \oplus \\
  (x_9+x_3-x_4-x_5)R_0\{6-2n\}
\end{array}%
\right.,
\]} and $M_2$ is 
{\tiny \[
\left.%
\begin{array}{c}
  (x_9-x_3)R_0 \\
  \oplus \\
  R_0\{6-2n\} \\
\end{array}%
\right.
\xrightarrow{\left(%
\begin{array}{cc}
 \alpha+q\beta&z \\
 q\gamma+\frac{y\delta}{z}&-b_1 \\
\end{array}%
\right)}
\left.%
\begin{array}{c}
  (x_9-x_3)R_0\{3-n\} \\
  \oplus \\
  R_0\{3-n\} \\
\end{array}%
\right.
\xrightarrow{\left(%
\begin{array}{cc}
  b_1&z \\
  q\gamma+\frac{y\delta}{z}&-(\alpha+q\beta) \\
\end{array}%
\right)}
\left.%
\begin{array}{c}
  (x_9-x_3)R_0 \\
  \oplus \\
  R_0\{6-2n\} \\
\end{array}%
\right..
\]}
It is easy to see that 
\begin{eqnarray*}
M_1 & \cong & \left(%
\begin{array}{cc}
  \alpha& b_1 \\
  \delta & y \\
\end{array}%
\right)_{R_0}, \\
M_2 & \cong & \left(%
\begin{array}{cc}
  \alpha+q\beta & b_1 \\
  q\gamma+\frac{y\delta}{z}& z \\
\end{array}%
\right)_{R_0}\{2\},
\end{eqnarray*}
both of which have potential $w=p(x_4)+p(x_5)+p(x_6)-p(x_2)-p(x_3)-p(x_4+x_5+x_6-x_2-x_3)\in R_0$. And we have that
\[
H_p(\Gamma_3)\cong H(N\otimes M_1)\{-3\} \oplus H(N\otimes M_2)\{-3\}.
\]

Applying Corollary \ref{row-op} and Lemma \ref{negative-differential} on $M_2$, we get
\[
M_2 \cong \left(%
\begin{array}{cc}
  f & (x_2-x_4)(x_2-x_5) \\
  g& x_3-x_6 \\
\end{array}%
\right)_{R_0}\{2\}
\cong
\left(%
\begin{array}{cc}
  f & (x_2-x_4)(x_2-x_5) \\
  -g& x_6-x_3 \\
\end{array}%
\right)_{R_0}\{2\},
\]
where $f=\alpha+q\beta$ and $g=q\gamma+\frac{y\delta}{z}+(x_2+x_3-x_4-x_5)(\alpha+q\beta)$. Then $\deg f \leq 2n-2$ and $\deg g \leq 2n$. By Corollary \ref{row-op}, the matrix factorization associated to the part of $\Gamma_2$ depicted in Figure \ref{reduction4} is
\[
\left(%
\begin{array}{cc}
  u_{5421} & x_5+x_4-x_2-x_1 \\
  v_{5421} & x_5x_4-x_2x_1 \\
  \pi_{63} & x_6-x_3
\end{array}%
\right)_{R_0[x_1]}\{-1\} \cong
\left(%
\begin{array}{cc}
  u_{5421}+x_2v_{5421} & x_5+x_4-x_2-x_1 \\
  v_{5421} & (x_2-x_4)(x_2-x_5) \\
  \pi_{63} & x_6-x_3
\end{array}%
\right)_{R_0[x_1]}\{-1\}.
\]
We apply Proposition \ref{variable-exclusion} to exclude $x_1$, and get $H_p(\Gamma_2) \cong H(N\otimes M_2')$, where 
\[
M_2'=\left(%
\begin{array}{cc}
  v_{5421} & (x_2-x_4)(x_2-x_5) \\
  \pi_{63} & x_6-x_3
\end{array}%
\right)_{R_0}\{-1\},
\]
which also have the potential $w=p(x_4)+p(x_5)+p(x_6)-p(x_2)-p(x_3)-p(x_4+x_5+x_6-x_2-x_3)\in R_0$. By Lemma 3.10 of \cite{Ras2}, since $(x_2-x_4)(x_2-x_5)$ and $x_6-x_3$ are relatively prime, there exists a $k\in R_0$ such that 
\begin{eqnarray*}
f& = &v_{5421}+k(x_6-x_3) \\
-g& = &\pi_{63}-k(x_2-x_4)(x_2-x_5).
\end{eqnarray*}
It is easy to see that $\deg k \leq 2n-4$. By Corollary \ref{twist}, we have $M_2\cong M_2'\{3\}$, and therefore
\[
H_p(\Gamma_2) \cong H(N\otimes M_2') \cong H(N\otimes M_2)\{-3\}.
\]

Define $x_1=x_4+x_5+x_6-x_2-x_3$, and $s_1,\cdots,s_6$ as above. One can check that $b_1=s_2-s_5$. So, by Corollary \ref{row-op} and Lemma \ref{negative-differential}, 
\[
M_1 \cong \left(%
\begin{array}{cc}
  \alpha& s_2-s_5 \\
  \delta & y \\
\end{array}%
\right)_{R_0}
\cong
\left(%
\begin{array}{cc}
  \alpha -x_3\delta & s_2-s_5 \\
  \delta & s_6-s_3 \\
\end{array}%
\right)_{R_0}
\cong
\left(%
\begin{array}{cc}
  \alpha -x_3\delta & s_2-s_5 \\
  -\delta & s_3-s_6 \\
\end{array}%
\right)_{R_0},
\]
where $\deg (\alpha-x_3\delta) \leq 2n-2$, $\deg\delta\leq2n-4$. We exclude $x_1$ from $C_p(\Gamma)$ by Proposition \ref{variable-exclusion}, and get
\[
H_p(\Gamma) = H(N\otimes \Upsilon) \cong H(N \otimes M_1')\{-3\},
\]
where
\[
M_1' = \left(%
\begin{array}{cc}
  a_2 & s_2-s_5 \\
  a_3 & s_3-s_6 \\
\end{array}%
\right)_{R_0},
\]
in which $a_2$ and $a_3$ are polynomials defined above with $\deg a_2=2n-2$ and $\deg a_3 = 2n-4$. Again, by Lemma 3.10 of \cite{Ras2}, since $s_2-s_5$ and $s_3-s_6$ are relatively prime, there exits $l\in R_0$ such that
\begin{eqnarray*}
\alpha -x_3\delta & = & a_2 + l(s_3-s_6) \\
-\delta & = & a_3-l(s_2-s_5),
\end{eqnarray*}
which clearly implies that $\deg l \leq 2n-8$. So, by Corollary \ref{twist}, $M_1\cong M_1'$, and therefore
\[
H_p(\Gamma) \cong H(N \otimes M_1)\{-3\}.
\]
Thus,
\[
H_p(\Gamma_3)\cong H(N\otimes M_1)\{-3\} \oplus H(N\otimes M_2)\{-3\} \cong H_p(\Gamma)\oplus H_p(\Gamma_2).
\]
\end{proof}

\section{Invariance}\label{reidemeister}

\subsection{Filtered chain complexes}
Recall that $\fil$ is a filtration of a $\C$-linear chain complex $(C,d)$ of $\fil$ is a filtration of the underlying $\C$-linear space $C$, and $d$ is a filtered linear map under $\fil$. The cohomology of a filtered chain complex inherits the filtration. A chain map $f$ is a filtered isomorphism (resp. quasi-isomorphism) of filtered chain complexes if $f$ is an isomorphism (resp. quasi-isomorphism) of chain complexes, and $f$ and $f^{-1}$ (resp. $f_\ast$ and $f_\ast^{-1}$) are both filtered homomorphisms. In this case, we said that the two filtered chain complexes are isomorphic (resp. quasi-isomorphic).

Again, for the rest of this section, we fix an $n\geq2$ and a polynomial
\[
p(x)=x^{n+1}+\sum_{i=0}^{n}c_i x^i \in \C[x].
\]
Clearly, for any link diagram $D$, the quantum filtration is a filtration of the chain complex $(H(C_p(D),d_{mf}),d_\chi)$. With the above MOY decompositions in hand, we prove Theorem \ref{fil-inv} by showing $H_p$ is invariant under Reidemeister moves. The proof is close to that in \cite{KR1}, in which we will make repeated use of the following Gaussian Elimination Lemma, which is the filtered version of Lemma 4.2 of \cite{Bar-fast}.

\begin{lemma}[Gaussian Elimination, \cite{Bar-fast}]\label{gaussian-elimination}
Consider a filtered double chain complex $\text{I}$, that is $\text{I}$ is filtered, and both differentials of $\text{I}$ are filtered $\C$-linear maps. Assume that
\[
\text{I}=\cdots\rightarrow C\xrightarrow{\left(%
\begin{array}{c}
  \alpha\\
  \beta \\
\end{array}%
\right)}
\left.%
\begin{array}{c}
  A\\
  \oplus \\
  D
\end{array}%
\right.
\xrightarrow{
\left(%
\begin{array}{cc}
  \phi & \delta\\
  \gamma & \varepsilon \\
\end{array}%
\right)}
\left.%
\begin{array}{c}
  B\\
  \oplus \\
  E
\end{array}%
\right.
\xrightarrow{
\left(%
\begin{array}{cc}
  \mu & \nu\\
\end{array}%
\right)} F \rightarrow \cdots,
\]
where the arrows give one of the differentials of $\text{I}$, and each Latin letter represents a filtered chain complex, which give the other differential of $\text{I}$. If $\phi:A\rightarrow B$ is also an isomorphism of filtered $\C$-linear spaces, then $I$ is isomorphic to the filtered double chain complex
\[
\text{II}=
\cdots\rightarrow C\xrightarrow{\left(%
\begin{array}{c}
  0\\
  \beta \\
\end{array}%
\right)}
\left.%
\begin{array}{c}
  A\\
  \oplus \\
  D
\end{array}%
\right.
\xrightarrow{
\left(%
\begin{array}{cc}
  \phi & 0\\
  0 & \varepsilon-\gamma\phi^{-1}\delta \\
\end{array}%
\right)}
\left.%
\begin{array}{c}
  B\\
  \oplus \\
  E
\end{array}%
\right.
\xrightarrow{
\left(%
\begin{array}{cc}
  0 & \nu\\
\end{array}%
\right)} F \rightarrow \cdots.
\]
In particular, the filtered total complex of $\text{I}$ is quasi-isomorphic to the filtered total complex of the filtered double chain complex
\[
\text{III}=
\cdots\rightarrow C \xrightarrow{\beta} D
\xrightarrow{\varepsilon-\gamma\phi^{-1}\delta} E\xrightarrow{\nu} F \rightarrow \cdots.
\]
\end{lemma}
\begin{proof}
Define maps $f:\text{I}\rightarrow\text{II}$ and $g:\text{II}\rightarrow\text{I}$ by
\[
\begin{CD}
\cdots @>>> C @>{\left(%
\begin{array}{c}
  \alpha\\
  \beta \\
\end{array}%
\right)}>> \left.%
\begin{array}{c}
  A\\
  \oplus \\
  D
\end{array}%
\right.  
@>{\left(%
\begin{array}{cc}
  \phi & \delta\\
  \gamma & \varepsilon \\
\end{array}%
\right)}>> 
\left.%
\begin{array}{c}
  B\\
  \oplus \\
  E
\end{array}%
\right.
@>{
\left(%
\begin{array}{cc}
  \mu & \nu\\
\end{array}%
\right)}>> F @>>> \cdots \\
@VV{\id}V @VV{\id}V @VV{\left(%
\begin{array}{cc}
  \id & \phi^{-1}\delta\\
  0 & \id \\
\end{array}%
\right)}V @VV{\left(%
\begin{array}{cc}
  \id & 0\\
  -\gamma\phi^{-1} & \id \\
\end{array}%
\right)}V @VV{\id}V @VV{\id}V\\
\cdots @>>> C @>{\left(%
\begin{array}{c}
  0\\
  \beta \\
\end{array}%
\right)}>> \left.%
\begin{array}{c}
  A\\
  \oplus \\
  D
\end{array}%
\right. @>{\left(%
\begin{array}{cc}
  \phi & 0\\
  0 & \varepsilon-\gamma\phi^{-1}\delta \\
\end{array}%
\right)}>> \left.%
\begin{array}{c}
  B\\
  \oplus \\
  E
\end{array}%
\right.
@>{
\left(%
\begin{array}{cc}
  0 & \nu\\
\end{array}%
\right)}>> F @>>> \cdots \\
@VV{\id}V @VV{\id}V @VV{\left(%
\begin{array}{cc}
  \id & -\phi^{-1}\delta\\
  0 & \id \\
\end{array}%
\right)}V @VV{\left(%
\begin{array}{cc}
  \id & 0\\
  \gamma\phi^{-1} & \id \\
\end{array}%
\right)}V @VV{\id}V @VV{\id}V \\
\cdots @>>> C @>{\left(%
\begin{array}{c}
  \alpha\\
  \beta \\
\end{array}%
\right)}>> \left.%
\begin{array}{c}
  A\\
  \oplus \\
  D
\end{array}%
\right.  
@>{\left(%
\begin{array}{cc}
  \phi & \delta\\
  \gamma & \varepsilon \\
\end{array}%
\right)}>> 
\left.%
\begin{array}{c}
  B\\
  \oplus \\
  E
\end{array}%
\right.
@>{
\left(%
\begin{array}{cc}
  \mu & \nu\\
\end{array}%
\right)}>> F @>>> \cdots \\
\end{CD}
\]
Using that the two differentials of $\text{I}$ anti-commute, it is straightforward to check that $f$ and $g$ are isomorphisms of filtered double chain complexes. Since the cohomology of the total chain complex of $0\rightarrow A \xrightarrow{\phi} B\rightarrow0$ is $0$, the standard projection and inclusion induce the quasi-isomorphisms between the (filtered) total complexes of $\text{II}$ and $\text{III}$.
\end{proof}

\subsection{Reidemeister Move I}\label{R-Move-I}

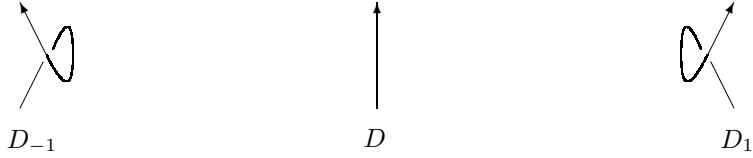
\begin{figure}[ht]

\setlength{\unitlength}{1pt}

\begin{picture}(420,60)(-210,-20)


\put(0,0){\vector(0,1){40}}

\put(-5,-15){$D$}


\put(-125,20){\vector(-1,2){10}}

\qbezier(-125,20)(-115,0)(-115,20)

\qbezier(-115,20)(-115,40)(-122.5,22.5)

\put(-135,0){\line(1,2){8.75}}

\put(-140,-15){$D_{-1}$}


\put(125,20){\vector(1,2){10}}

\qbezier(125,20)(115,0)(115,20)

\qbezier(115,20)(115,40)(122.5,22.5)

\put(135,0){\line(-1,2){8.75}}

\put(130,-15){$D_{1}$}

\end{picture}

\caption{Reidemeister move I}\label{Rmove1}

\end{figure}

\begin{proposition}\label{reidemeister1}
Assume that $D$, $D_{-1}$ and $D_1$ are link diagrams that are identical except in the part depicted in Figure \ref{Rmove1}. Then there are isomorphisms
\[
H_p(D) \cong H_p(D_{-1}) \cong H_p(D_1)
\]
that preserve both the cohomological grading and the quantum filtration.
\end{proposition}

\begin{figure}[ht]

\setlength{\unitlength}{1pt}

\begin{picture}(420,60)(-210,-20)


\put(-150,-10){\vector(0,1){50}}

\put(-153,-20){$D$}

\put(-160,-10){\tiny{$x_2$}}

\put(-160,40){\tiny{$x_1$}}


\put(-20,-10){\vector(0,1){50}}

\put(0,0){\line(0,1){30}}

\qbezier(0,30)(0,40)(10,40)

\qbezier(10,40)(20, 40)(20,30)

\put(20,30){\vector(0,-1){30}}

\qbezier(0,0)(0,-10)(10,-10)

\qbezier(10,-10)(20, -10)(20,0)

\put(-10,-20){$\Gamma_0$}

\put(-30,-10){\tiny{$x_2$}}

\put(-30,40){\tiny{$x_1$}}

\put(25,15){\tiny{$x_3$}}

\put(18.5,15){-}


\put(117.5,10){\vector(1,1){0}}

\qbezier(110,-10)(110,5)(117.5,10)

\qbezier(110,40)(110,25)(117,20)

\put(110,44){\vector(0,1){0}}

\qbezier(130,0)(130,5)(122.5,10)

\qbezier(130,30)(130,25)(122.5,20)

\qbezier(130,30)(130,40)(140,40)

\qbezier(140,40)(150, 40)(150,30)

\put(150,30){\vector(0,-1){30}}

\qbezier(130,0)(130,-10)(140,-10)

\qbezier(140,-10)(150, -10)(150,0)

\put(117,-20){$\Gamma_1$}

\put(100,-10){\tiny{$x_2$}}

\put(100,40){\tiny{$x_1$}}

\put(155,15){\tiny{$x_3$}}

\put(148.5,15){-}

\linethickness{5pt}

\put(120,10){\line(0,1){10}}

\end{picture}

\caption{$D$, $\Gamma_0$ and
$\Gamma_1$}\label{Rmove1-decomp}

\end{figure}
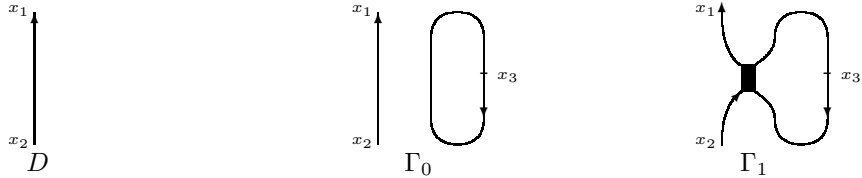

\begin{proof}
Let $\Gamma_0$ and $\Gamma_1$ be planar diagrams that are identical to $D$ except in the part depicted in Figure \ref{Rmove1-decomp}. Then $(H(C_p(D_{-1}),d_{mf}),d_\chi)$ is the total complex of
\[
0\rightarrow H(C_p(\Gamma_0),d_{mf})\{1-n\} \xrightarrow{\chi_0} H(C_p(\Gamma_1),d_{mf})\{-n\} \rightarrow 0,
\]
where $\chi_0$ is the map associated to the wide edge in $\Gamma_1$, and the other differential is from resolving all other crossings of $D_{-1}$. Since 
\[
H_p(\bigcirc)=\C[x_3]/(p(x_3))\{1-n\}=\left\langle 1,x_3,x_3^2,\cdots,x_3^{n-1}\right\rangle\{1-n\},
\] 
there is a decomposition of filtered chain complex
\begin{eqnarray*}
& & H(C_p(\Gamma_0),d_{mf}) \{1-n\} \\ & = & H(C_p(D),d_{mf})\otimes_\C \left\langle 1,x_3,\cdots,x_3^{n-2}\right\rangle\{2-2n\} \oplus H(C_p(D),d_{mf})\otimes_\C \left\langle x_3^{n-1}\right\rangle\{2-2n\} \\ 
& \cong & H(C_p(D),d_{mf})\otimes_\C \left\langle 1,x_3,\cdots,x_3^{n-2}\right\rangle\{2-2n\} \oplus H(C_p(D),d_{mf}).
\end{eqnarray*}
From the proof of Proposition \ref{moy1},
\[
\chi_0|_{H(C_p(D),d_{mf})\otimes_\C \left\langle 1,x_3,\cdots,x_3^{n-2}\right\rangle\{2-2n\}}:H(C_p(D),d_{mf})\otimes_\C \left\langle 1,x_3,\cdots,x_3^{n-2}\right\rangle\{2-2n\}\rightarrow H(C_p(\Gamma_1),d_{mf})\{-n\}
\]
is an isomorphism of filtered chain complexes. Applying the Gaussian Elimination Lemma (Lemma \ref{gaussian-elimination}), we get 
\[
H_p(D_{-1}) \cong H_p(D).
\]

Similarly, $(H(C_p(D_1),d_{mf}),d_\chi)$ is the total complex of
\[
0\rightarrow H(C_p(\Gamma_1),d_{mf})\{n\} \xrightarrow{\chi_1} H(C_p(\Gamma_0),d_{mf})\{n-1\} \rightarrow 0.
\]
There is a decomposition of filtered chain complex
\begin{eqnarray*}
& & H(C_p(\Gamma_0),d_{mf}) \{n-1\} \\ & = & H(C_p(D),d_{mf})\otimes_\C \left\langle 1 \right\rangle \oplus H(C_p(D),d_{mf})\otimes_\C \left\langle x_3,\cdots,x_3^{n-1}\right\rangle\\ 
& \cong & H(C_p(D),d_{mf}) \oplus H(C_p(D),d_{mf})\otimes_\C \left\langle x_3,\cdots,x_3^{n-1}\right\rangle.
\end{eqnarray*}
Again, from the proof of Proposition \ref{moy1},
\[
\chi_1: H(C_p(\Gamma_1),d_{mf})\{n\} \rightarrow H(C_p(D),d_{mf})\otimes_\C \left\langle x_3,\cdots,x_3^{n-1}\right\rangle
\]
is an isomorphism of filtered chain complexes. Applying the Gaussian Elimination Lemma (Lemma \ref{gaussian-elimination}), we get 
\[
H_p(D_1) \cong H_p(D).
\]
\end{proof}

\subsection{Reidemeister Move II$_a$}\label{R-Move-II-a}

\begin{figure}[ht]

\setlength{\unitlength}{1pt}

\begin{picture}(420,60)(-210,-20)


\put(-10,0){\vector(0,1){40}}

\put(10,0){\vector(0,1){40}}

\put(-5,-15){$D$}


\qbezier(-145,0)(-105,20)(-145,40)

\put(-145,40){\vector(-2,1){0}}

\put(-125,0){\line(-2,1){8}}

\qbezier(-137,6)(-155,20)(-137,34)

\put(-133,36){\vector(2,1){8}}

\put(-140,-15){$D_1$}


\qbezier(145,0)(105,20)(145,40)

\put(145,40){\vector(2,1){0}}

\put(125,0){\line(2,1){8}}

\qbezier(137,6)(155,20)(137,34)

\put(133,36){\vector(-2,1){8}}

\put(130,-15){$D_2$}

\end{picture}

\caption{Reidemeister move II$_a$}\label{Rmove2a}

\end{figure}
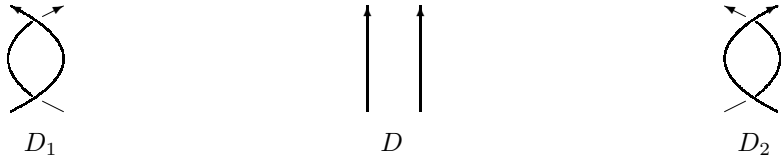

\begin{proposition}\label{reidemeister2a}
Assume that $D$, $D_1$ and $D_2$ are link diagrams that are identical except in the part depicted in Figure \ref{Rmove2a}. Then there are isomorphisms
\[
H_p(D) \cong H_p(D_1) \cong H_p(D_2)
\]
that preserve both the cohomological grading and the quantum filtration.
\end{proposition}

\begin{figure}[ht]

\setlength{\unitlength}{1pt}

\begin{picture}(420,160)(-210,-10)

\linethickness{.25pt}


\put(-70,80){\vector(0,1){60}}

\put(-50,80){\vector(0,1){60}}

\put(-100,105){$\Gamma_{00}$}

\put(-80,140){\tiny{$x_1$}}

\put(-80,80){\tiny{$x_4$}}

\put(-47,140){\tiny{$x_2$}}

\put(-47,80){\tiny{$x_3$}}


\linethickness{.25pt}

\put(62.5,120){\vector(-3,2){0}}

\qbezier(70,80)(70,120)(62.5,120)

\put(57.5,120){\vector(3,2){0}}

\qbezier(50,80)(50,120)(57.5,120)

\qbezier(50,140)(50,135)(57.5,130)

\qbezier(70,140)(70,135)(62.5,130)

\put(70,140){\vector(1,4){0}}

\put(50,140){\vector(-1,4){0}}

\put(20,105){$\Gamma_{10}$}

\put(40,140){\tiny{$x_1$}}

\put(40,80){\tiny{$x_4$}}

\put(73,140){\tiny{$x_2$}}

\put(73,80){\tiny{$x_3$}}

\linethickness{5pt}

\put(60,120){\line(0,1){10}}


\linethickness{.25pt}

\put(62.5,10){\vector(-1,1){0}}

\qbezier(70,0)(70,5)(62.5,10)

\qbezier(70,30)(70,25)(62.5,20)

\put(70,30){\vector(1,4){0}}

\put(50,30){\vector(-1,4){0}}

\put(57.5,10){\vector(1,1){0}}

\qbezier(50,0)(50,5)(57.5,10)

\qbezier(50,30)(50,25)(57.5,20)

\qbezier(70,30)(70,35)(62.5,40)

\qbezier(70,60)(70,55)(62.5,50)

\put(70,60){\vector(1,4){0}}

\put(50,60){\vector(-1,4){0}}

\qbezier(50,30)(50,35)(57.5,40)

\qbezier(50,60)(50,55)(57.5,50)

\put(20,25){$\Gamma_{11}$}

\put(40,60){\tiny{$x_1$}}

\put(40,0){\tiny{$x_4$}}

\put(73,60){\tiny{$x_2$}}

\put(73,0){\tiny{$x_3$}}

\put(40,30){\tiny{$x_6$}}

\put(73,30){\tiny{$x_5$}}

\put(48.5,30){-}

\put(68.5,30){-}

\linethickness{5pt}

\put(60,10){\line(0,1){10}}

\put(60,40){\line(0,1){10}}


\linethickness{.25pt}

\put(-62.5,10){\vector(1,1){0}}

\qbezier(-70,0)(-70,5)(-62.5,10)

\put(-57.5,10){\vector(-1,1){0}}

\qbezier(-50,0)(-50,5)(-57.5,10)

\qbezier(-70,60)(-70,20)(-62.5,20)

\put(-70,60){\vector(0,0){0}}

\put(-50,60){\vector(0,1){0}}

\qbezier(-50,60)(-50,20)(-57.5,20)

\put(-100,25){$\Gamma_{01}$}

\put(-80,60){\tiny{$x_1$}}

\put(-80,0){\tiny{$x_4$}}

\put(-47,60){\tiny{$x_2$}}

\put(-47,0){\tiny{$x_3$}}

\linethickness{5pt}

\put(-60,10){\line(0,1){10}}

\end{picture}

\caption{Local diagrams related to Reidemeister move
II$_a$}\label{Rmove2a-decomp}
\end{figure}
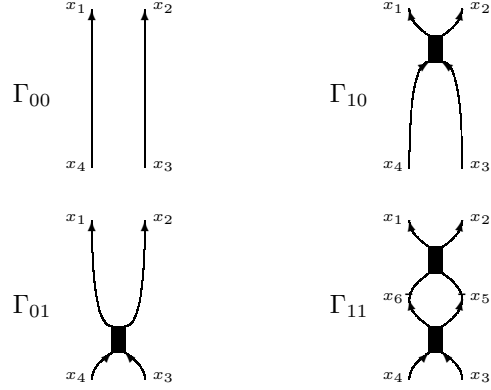

Consider the planar diagrams $\Gamma_{00},\Gamma_{01},\Gamma_{10},\Gamma_{11}$, which are identical to $D$ except in the part depicted in Figure \ref{Rmove2a-decomp}. Clearly, $\Gamma_{01}=\Gamma_{10}$. From Proposition \ref{moy2}, we have that 
\[
H(C_p(\Gamma_{11}),d_{mf}) \cong H(C_p(\Gamma_{10}),d_{mf})\{-1\} \oplus H(C_p(\Gamma_{01}),d_{mf})\{1\}.
\]
Let $P:H(C_p(\Gamma_{11}),d_{mf})\rightarrow H(C_p(\Gamma_{01}),d_{mf})\{1\}$ be the corresponding projection, and $J:H(C_p(\Gamma_{10}),d_{mf})\{-1\} \rightarrow H(C_p(\Gamma_{11}),d_{mf})$ the corresponding inclusion.

\begin{lemma}\label{Rmove2a-lemma}
Let $\chi_0:H(C_p(\Gamma_{01}),d_{mf})\rightarrow H(C_p(\Gamma_{11}),d_{mf})$ be the $\chi_0$-map associated to the upper wide edge in $\Gamma_{11}$ and $\chi_1:H(C_p(\Gamma_{11}),d_{mf})\rightarrow H(C_p(\Gamma_{10}),d_{mf})$ the $\chi_1$-map associated to the lower wide edge in $\Gamma_{11}$. (Both maps are of degree $1$.) Then the homomorphisms
\begin{eqnarray*}
P\circ\chi_0 & : & H(C_p(\Gamma_{01}),d_{mf})\{1\}\rightarrow H(C_p(\Gamma_{01}),d_{mf})\{1\}, \\
\chi_1 \circ J & : & H(C_p(\Gamma_{10}),d_{mf})\{-1\} \rightarrow H(C_p(\Gamma_{10}),d_{mf})\{-1\}
\end{eqnarray*}
are isomorphisms of filtered chain complexes, where the differentials of these chain complexes are from resolving all the crossing of these diagrams (outside Figure \ref{Rmove2a-decomp}).
\end{lemma}

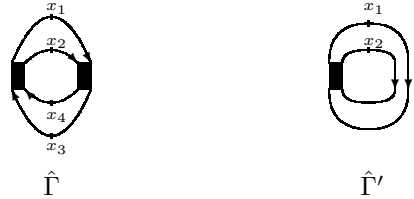
\begin{figure}[ht]

\setlength{\unitlength}{1pt}

\begin{picture}(420,80)(-210,-25)


\qbezier(-70,25)(-60,35)(-50,25)

\put(-50,25){\vector(1,-1){0}}

\qbezier(-70,15)(-60,5)(-50,15)

\put(-70,15){\vector(-1,1){0}}

\qbezier(-75,25)(-60,60)(-45,25)

\put(-45,25){\vector(1,-2){0}}

\qbezier(-75,15)(-60,-20)(-45,15)

\put(-75,15){\vector(-1,2){0}}

\put(-62,45){\tiny{$x_1$}}

\put(-62,32){\tiny{$x_2$}}

\put(-62,-8){\tiny{$x_3$}}

\put(-62,4){\tiny{$x_4$}}

\put(-60,41.5){\line(0,1){2}}

\put(-60,29){\line(0,1){2}}

\put(-60,9){\line(0,1){2}}

\put(-60,-3.5){\line(0,1){2}}

\put(-63,-25){$\hat{\Gamma}$}


\qbezier(45,25)(45,40)(60,40)

\qbezier(60,40)(75,40)(75,25)

\qbezier(50,25)(50,30)(60,30)

\qbezier(60,30)(70,30)(70,25)

\put(70,25){\vector(0,-1){10}}

\put(75,25){\vector(0,-1){10}}

\qbezier(50,15)(50,10)(60,10)

\qbezier(60,10)(70,10)(70,15)

\qbezier(45,15)(45,0)(60,0)

\qbezier(60,0)(75,0)(75,15)

\put(58,45){\tiny{$x_1$}}

\put(58,32){\tiny{$x_2$}}

\put(60,41){\line(0,-1){2}}

\put(60,29){\line(0,1){2}}

\put(57,-25){$\hat{\Gamma}'$}

\linethickness{5pt}

\put(-72.5,15){\line(0,1){10}}

\put(-47.5,15){\line(0,1){10}}

\put(47.5,15){\line(0,1){10}}

\end{picture}

\caption{$\hat{\Gamma}$ and $\hat{\Gamma}'$}\label{hat-gamma2}

\end{figure}

\begin{proof}
Let $\hat{\Gamma}_{00},\hat{\Gamma}_{01},\hat{\Gamma}_{10},\hat{\Gamma}_{11}$ be the local open graphs depicted in Figure \ref{Rmove2a-decomp}. We claim that, when interpreted as homomorphisms of the $C_p$ complexes of these local open graphs,
\begin{eqnarray*}
P\circ\chi_0& : & C_p(\hat{\Gamma}_{01})\{1\}\rightarrow C_p(\hat{\Gamma}_{01})\{1\} \\
\chi_1\circ J & : & C_p(\Gamma_{10})\{-1\} \rightarrow C_p(\Gamma_{10})\{-1\}
\end{eqnarray*} 
are homotopic to non-zero multiples of the identity maps, which implies the lemma. First, observe that both homomorphisms are filtered. By Proposition \ref{hmf-tensor}, 
\[
\Hom_{HMF}(C_p(\hat{\Gamma}_{01}),C_p(\hat{\Gamma}_{01})) \cong \Hom_{HMF}(C_p(\hat{\Gamma}_{10}),C_p(\hat{\Gamma}_{10})) \cong H_p(\hat{\Gamma})\{2n-2\},
\]
where $\hat{\Gamma}$ is depicted in Figure \ref{hat-gamma2}. But, by Lemma \ref{H-hat-gamma}, the filtered dimension of $H_p(\hat{\Gamma})$ is $[2][n][n-1]$, where $[k]=\frac{q^k-q^{-k}}{q-q^{-1}}$. This implies that 
\[
\Hom_{hmf}(C_p(\hat{\Gamma}_{01}),C_p(\hat{\Gamma}_{01})) \cong \Hom_{hmf}(C_p(\hat{\Gamma}_{10}),C_p(\hat{\Gamma}_{10})) \cong \C,
\]
which means any filtered homomorphisms of matrix factorizations from $C_p(\hat{\Gamma}_{01})~(=C_p(\hat{\Gamma}_{10}))$ to itself is homotopic to a multiple of the identity map. So, to prove the lemma, we only need to show that the above homomorphisms are not homotopic to the $0$ map. We do this by showing that, after closing $\hat{\Gamma}_{11}$ to $\hat{\Gamma}$ and $\hat{\Gamma}_{01}$, $\hat{\Gamma}_{10}$ to $\hat{\Gamma}'$ in Figure \ref{hat-gamma2}, the above homomorphisms induce non-zero maps of the cohomology. According to Lemma \ref{H-hat-gamma}, $H_p(\hat{\Gamma})$ has the basis 
\[
\{x_1^ix_2^jx_3^{\epsilon} ~|~ 0\leq i \leq n-1,~0\leq j\leq n-2,~\epsilon=0,1\},
\]
where the element $x_1^ix_2^jx_3^{\epsilon}$ has quantum degree $2i+2j+2\epsilon-2n+2$. By a similar computation, $H_p(\hat{\Gamma}')$ has the basis
\[
\{x_1^ix_2^j ~|~ 0\leq i \leq n-1,~0\leq j\leq n-2\},
\]
where the element $x_1^ix_2^j$ has quantum degree $2i+2j-2n+3$.

By the construction of the decomposition $H_p(\hat{\Gamma})\cong H_p(\hat{\Gamma}')\{-1\}\oplus H_p(\hat{\Gamma}')\{1\}$, we can see that, when using the above bases, 
\[
J: H_p(\hat{\Gamma}')\{-1\} \rightarrow H_p(\hat{\Gamma})
\] 
is given by $J(x_1^ix_2^j)=x_1^ix_2^j$, and 
\[
P:H_p(\hat{\Gamma})\rightarrow H_p(\hat{\Gamma}')\{1\}
\] 
is given by $P(x_1^ix_2^j)=0$, $P(x_1^ix_2^jx_3)=x_1^ix_2^j$.

By Remark \ref{direct-sum-variable-exclusion} and the definition of $\chi_1$, it is easy to check that $\chi_1(x_1^ix_2^j)=x_1^ix_2^j$. So $\chi_1\circ J(x_1^ix_2^j)=x_1^ix_2^j\neq0$. This means $\chi_1\circ J$ induces a non-zero homomorphism on the cohomology. 

Similarly, by Remark \ref{direct-sum-variable-exclusion} and the definition of $\chi_0$, one can check that 
\[
\chi_0(x_1^ix_2^j) = (x_1-x_4)x_1^ix_2^j = (x_3-x_2)x_1^ix_2^j = x_1^ix_2^jx_3 - x_1^ix_2^{j+1}.
\]
Then it is easy to see that $P\circ\chi_0(x_1^ix_2^j)=P(x_1^ix_2^jx_3 - x_1^ix_2^{j+1})=P(x_1^ix_2^jx_3)=x_1^ix_2^j\neq0$. (The case when $j=n-2$ needs slightly more work since $x_1^ix_2^{n-1}$ is not a base element. We can use $v_{1212}=0=p'(x_1)$ to convert it into a base element. See the proof of Lemma \ref{H-hat-gamma}.) So $P\circ\chi_0$ also induces a non-zero homomorphism on the cohomology. 
\end{proof}

\begin{proof}[Proof of Proposition \ref{reidemeister2a}]
$(H(C_p(D_1),d_{mf}),d_\chi)$ is the total complex of
\[
0\rightarrow H(C_p(\Gamma_{01}),d_{mf})\{1\} \xrightarrow{\left(%
\begin{array}{c}
  f_1\\
  f_0 \\
\end{array}%
\right)}
\left.%
\begin{array}{c}
  H(C_p(\Gamma_{11}),d_{mf})\\
  \oplus \\
  H(C_p(\Gamma_{00}),d_{mf})
\end{array}%
\right.
\xrightarrow{\left(%
\begin{array}{cc}
  g_1 & g_0
\end{array}%
\right)}
H(C_p(\Gamma_{10}),d_{mf})\{-1\} \rightarrow 0.
\]
By the above decomposition of $H(C_p(\Gamma_{11}),d_{mf})$, this is isomorphic to

\[
0\rightarrow H(C_p(\Gamma_{01}),d_{mf})\{1\} \xrightarrow{\left(%
\begin{array}{c}
  P\circ \chi_0\\
  f_{11} \\
  f_0 \\
\end{array}%
\right)}
\left.%
\begin{array}{c}
  H(C_p(\Gamma_{01}),d_{mf})\{1\} \\
  \oplus \\
  H(C_p(\Gamma_{10}),d_{mf})\{-1\}\\
  \oplus \\
  H(C_p(\Gamma_{00}),d_{mf})
\end{array}%
\right.
\xrightarrow{\left(%
\begin{array}{ccc}
  g_{10} & \chi_1\circ J & g_0
\end{array}%
\right)}
H(C_p(\Gamma_{10}),d_{mf})\{-1\} \rightarrow 0.
\]
Using Lemma \ref{Rmove2a-lemma}, we apply Lemma \ref{gaussian-elimination} to the above chain complex twice, which gives a quasi-isomorphism of filtered chain complexes
\[
(H(C_p(D_1),d_{mf}),d_\chi) \rightarrow (H(C_p(\Gamma_{00}),d_{mf}),d_\chi) = (H(C_p(D),d_{mf}),d_\chi).
\]
Hence, $H_p(D)\cong H_p(D_1)$. $H_p(D)\cong H_p(D_2)$ follows similarly.
\end{proof}

\subsection{Reidemeister Move II$_b$}\label{R-Move-II-b}

\begin{figure}[ht]

\setlength{\unitlength}{1pt}

\begin{picture}(420,60)(-210,-20)


\put(120,10){\vector(1,0){40}}

\put(160,30){\vector(-1,0){40}}

\put(135,-15){$D$}


\qbezier(-160,30)(-140,-10)(-120,30)

\put(-160,30){\vector(-1,2){0}}

\put(-160,10){\line(1,1){6}}

\qbezier(-150,20)(-140,30)(-130,20)

\put(-126,16){\vector(1,-1){6}}

\put(-145,-15){$D_1$}

\end{picture}

\caption{Reidemeister move II$_b$}\label{Rmove2b}

\end{figure}
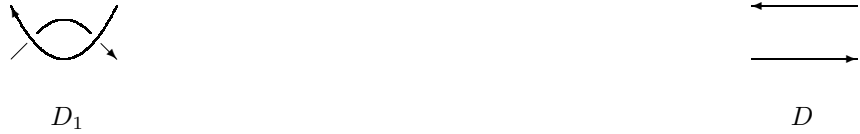

\begin{proposition}\label{reidemeister2b}
Assume that $D$ and $D_1$ are link diagrams that are identical except in the part depicted in Figure \ref{Rmove2b}. Then there is an isomorphism
\[
H_p(D) \cong H_p(D_1)
\]
that preserves both the cohomological grading and the quantum filtration.
\end{proposition}

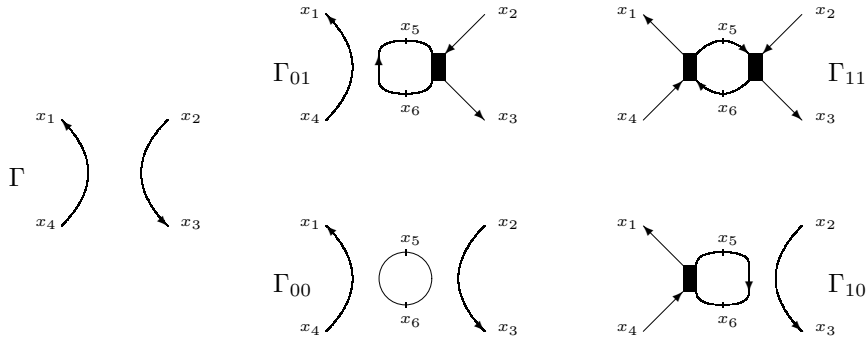
\begin{figure}[ht]

\setlength{\unitlength}{1pt}

\begin{picture}(420,160)(-210,-10)

\linethickness{.25pt}


\linethickness{.25pt}

\put(-45,95){\vector(1,-1){15}}

\put(-30,120){\vector(-1,-1){15}}

\qbezier(-50,105)(-50,110)(-60,110)

\qbezier(-60,110)(-70,110)(-70,105)

\qbezier(-50,95)(-50,90)(-60,90)

\qbezier(-60,90)(-70,90)(-70,95)

\put(-70,95){\vector(0,1){10}}

\put(-90,120){\vector(-1,1){0}}

\qbezier(-90,80)(-70,100)(-90,120)

\put(-100,120){\tiny{$x_1$}}

\put(-25,120){\tiny{$x_2$}}

\put(-100,80){\tiny{$x_4$}}

\put(-25,80){\tiny{$x_3$}}

\put(-62,114){\tiny{$x_5$}}

\put(-62,83){\tiny{$x_6$}}

\put(-60,109){\line(0,1){2}}

\put(-60,89){\line(0,1){2}}

\put(-110,95){$\Gamma_{01}$}

\linethickness{5pt}

\put(-47.5,95){\line(0,1){10}}


\linethickness{.25pt}

\qbezier(50,105)(60,115)(70,105)

\put(70,105){\vector(1,-1){0}}

\qbezier(50,95)(60,85)(70,95)

\put(50,95){\vector(-1,1){0}}

\put(45,105){\vector(-1,1){15}}

\put(75,95){\vector(1,-1){15}}

\put(30,80){\vector(1,1){15}}

\put(90,120){\vector(-1,-1){15}}

\put(20,120){\tiny{$x_1$}}

\put(95,120){\tiny{$x_2$}}

\put(20,80){\tiny{$x_4$}}

\put(95,80){\tiny{$x_3$}}

\put(58,114){\tiny{$x_5$}}

\put(58,83){\tiny{$x_6$}}

\put(60,109){\line(0,1){2}}

\put(60,89){\line(0,1){2}}

\put(100,95){$\Gamma_{11}$}

\linethickness{5pt}

\put(72.5,95){\line(0,1){10}}

\put(47.5,95){\line(0,1){10}}


\linethickness{.25pt}

\put(45,25){\vector(-1,1){15}}

\put(30,0){\vector(1,1){15}}

\qbezier(50,25)(50,30)(60,30)

\qbezier(60,30)(70,30)(70,25)

\qbezier(50,15)(50,10)(60,10)

\qbezier(60,10)(70,10)(70,15)

\put(70,25){\vector(0,-1){10}}

\put(90,0){\vector(1,-1){0}}

\qbezier(90,00)(70,20)(90,40)

\put(20,40){\tiny{$x_1$}}

\put(95,40){\tiny{$x_2$}}

\put(20,0){\tiny{$x_4$}}

\put(95,0){\tiny{$x_3$}}

\put(58,34){\tiny{$x_5$}}

\put(58,3){\tiny{$x_6$}}

\put(60,29){\line(0,1){2}}

\put(60,9){\line(0,1){2}}

\put(100,15){$\Gamma_{10}$}

\linethickness{5pt}

\put(47.5,15){\line(0,1){10}}


\linethickness{.25pt}

\put(-30,0){\vector(1,-1){0}}

\qbezier(-30,40)(-50,20)(-30,0)

\put(-90,40){\vector(-1,1){0}}

\qbezier(-90,0)(-70,20)(-90,40)

\put(-60,20){\circle{20}}

\put(-100,40){\tiny{$x_1$}}

\put(-25,40){\tiny{$x_2$}}

\put(-100,0){\tiny{$x_4$}}

\put(-25,0){\tiny{$x_3$}}

\put(-62,34){\tiny{$x_5$}}

\put(-62,3){\tiny{$x_6$}}

\put(-60,29){\line(0,1){2}}

\put(-60,9){\line(0,1){2}}

\put(-110,15){$\Gamma_{00}$}


\linethickness{.25pt}

\put(-150,40){\vector(1,-1){0}}

\qbezier(-150,80)(-170,60)(-150,40)

\put(-190,80){\vector(-1,1){0}}

\qbezier(-190,40)(-170,60)(-190,80)

\put(-200,80){\tiny{$x_1$}}

\put(-145,80){\tiny{$x_2$}}

\put(-200,40){\tiny{$x_4$}}

\put(-145,40){\tiny{$x_3$}}

\put(-210,55){$\Gamma$}

\end{picture}

\caption{Local diagrams related to Reidemeister move
II$_b$}\label{Rmove2b-decomp}
\end{figure}

\begin{proof}
Let $\Gamma,\Gamma_{00},\Gamma_{01},\Gamma_{10},\Gamma_{11}$ be planar diagrams that are identical to $D$ except in the part depicted in Figure \ref{Rmove2b-decomp}. Then $(H(C_p(D_1),d_{mf}),d_\chi)$ is the total complex of
\[
0\rightarrow H(C_p(\Gamma_{01}),d_{mf})\{1\} \xrightarrow{\left(%
\begin{array}{c}
  f_1\\
  f_0 \\
\end{array}%
\right)}
\left.%
\begin{array}{c}
  H(C_p(\Gamma_{11}),d_{mf})\\
  \oplus \\
  H(C_p(\Gamma_{00}),d_{mf})
\end{array}%
\right.
\xrightarrow{\left(%
\begin{array}{cc}
  g_1 & g_0
\end{array}%
\right)}
H(C_p(\Gamma_{10}),d_{mf})\{-1\} \rightarrow 0.
\]

Note that $g_0:H(C_p(\Gamma_{00}),d_{mf}) \rightarrow H(C_p(\Gamma_{10}),d_{mf})\{-1\}$ is induced by the $\chi_0$-map associated to the left wide edge. By the proof of Proposition \ref{moy1}, 
\[
H(C_p(\Gamma_{00}),d_{mf}) \cong H(C_p(\Gamma_{10}),d_{mf})\{-1\} \oplus H(C_p(\Gamma),d_{mf})\{n-1\},
\] 
and the composition
\[
H(C_p(\Gamma_{10}),d_{mf})\{-1\} \xrightarrow{J} H(C_p(\Gamma_{00}),d_{mf}) \xrightarrow{\chi_0} H(C_p(\Gamma_{10}),d_{mf})\{-1\},
\]
is an isomorphism of filtered chain complexes, where $J$ is the inclusion from the above decomposition, and $\chi_0$ is the $\chi_0$-map associated to the left wide edge. By Lemma \ref{gaussian-elimination}, $(H(C_p(D_1),d_{mf}),d_\chi)$ is quasi-isomorphic to the total complex of
\[
0\rightarrow
H(C_p(\Gamma_{01}),d_{mf})\{1\}
\xrightarrow{\left(%
\begin{array}{c}
  f_1 \\
  \varepsilon\circ f_0
\end{array}%
\right)}
\left.%
\begin{array}{c}
  H(C_p(\Gamma_{11}),d_{mf})\\
  \oplus \\
  H(C_p(\Gamma),d_{mf})\{n-1\}
\end{array}%
\right.
\rightarrow 0,
\]
where $\varepsilon:H(C_p(\Gamma_{00}),d_{mf}) \rightarrow H(C_p(\Gamma),d_{mf})\{n-1\}$ is defined in Subsection \ref{graph-cobodisms}, and is clearly the projection map from the above decomposition of $H(C_p(\Gamma_{00}),d_{mf})$. 

Again, by Proposition \ref{moy1}, 
\[
H(C_p(\Gamma_{01}),d_{mf})\{1\} \cong \bigoplus_{i=0}^{n-2} H(C_p(\Gamma),d_{mf})\{3-n+2i\}.
\]
Note that $f_0$ is induced by the $\chi_1$-map associated to the right wide edge. From the proof of Proposition \ref{moy1}, the composition
\[
H(C_p(\Gamma),d_{mf})\{n-1\} \xrightarrow{J'} H(C_p(\Gamma_{01}),d_{mf})\{1\} \xrightarrow{\varepsilon\circ f_0} H(C_p(\Gamma),d_{mf})\{n-1\},
\]
is an isomorphism of filtered chain complexes, where $J': H(C_p(\Gamma),d_{mf})\{n-1\} \rightarrow H(C_p(\Gamma_{01}),d_{mf})\{1\}$ is the inclusion from the above decomposition of $H(C_p(\Gamma_{01}),d_{mf})\{1\}$. Applying Lemma \ref{gaussian-elimination} again, we get that 
$(H(C_p(D_1),d_{mf}),d_\chi)$ is quasi-isomorphic to the total complex of
\[
0\rightarrow \bigoplus_{i=0}^{n-3} H(C_p(\Gamma),d_{mf})\{3-n+2i\} \xrightarrow{f_1\circ J''} H(C_p(\Gamma_{11}),d_{mf})\rightarrow 0,
\]
where 
\[
J'': \bigoplus_{i=0}^{n-3} H(C_p(\Gamma),d_{mf})\{3-n+2i\} \rightarrow H(C_p(\Gamma_{01}),d_{mf})\{1\}
\] 
is the projection from the above decomposition of $H(C_p(\Gamma_{01}),d_{mf})\{1\}$.

By Proposition \ref{moy3},
\[
H(C_p(\Gamma_{11}),d_{mf}) \cong H(C_p(D),d_{mf}) \oplus (\bigoplus_{i=0}^{n-3} H(C_p(\Gamma),d_{mf})\{3-n+2i\}).
\]
Note that $f_1: H(C_p(\Gamma_{01}),d_{mf})\{1\} \rightarrow H(C_p(\Gamma_{11}),d_{mf})$ is induced by the $\chi_0$-map of the left wide edge. From the proof Proposition \ref{moy3}, the composition
\[
\bigoplus_{i=0}^{n-3} H(C_p(\Gamma),d_{mf})\{3-n+2i\} \xrightarrow{f_1\circ J''} H(C_p(\Gamma_{11}),d_{mf}) \xrightarrow{P'} \bigoplus_{i=0}^{n-3} H(C_p(\Gamma),d_{mf})\{3-n+2i\},
\]
where $P': H(C_p(\Gamma_{11}),d_{mf}) \rightarrow\bigoplus_{i=0}^{n-3} H(C_p(\Gamma),d_{mf})\{3-n+2i\}$ is the projection map from the above decomposition of $H(C_p(\Gamma_{11}),d_{mf})$, is an isomorphism of filtered chain complexes. Then, by Lemma \ref{gaussian-elimination}, $(H(C_p(D_1),d_{mf}),d_\chi)$ is quasi-isomorphic to the total complex of
\[
0 \rightarrow H(C_p(D),d_{mf}) \rightarrow 0,
\]
which is $(H(C_p(D),d_{mf}),d_\chi)$. Thus, $H_p(D_1)\cong H_p(D)$.
\end{proof}

\subsection{Reidemeister Move III}\label{R-Move-III}

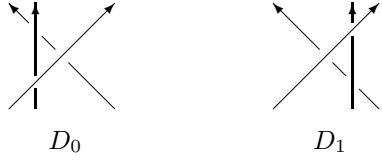
\begin{figure}[ht]

\setlength{\unitlength}{1pt}

\begin{picture}(420,60)(-210,-15)

\put(-70,0){\vector(1,1){40}}

\put(-60,0){\line(0,1){7.5}}

\put(-60,12.5){\vector(0,1){27.5}}

\put(-30,0){\line(-1,1){17.5}}

\put(-52.5,22.5){\line(-1,1){5}}

\put(-62.5,32.5){\vector(-1,1){7.5}}

\put(-55,-15){$D_0$}

\put(30,0){\vector(1,1){40}}

\put(60,0){\line(0,1){27.5}}

\put(60,32.5){\vector(0,1){7.5}}

\put(70,0){\line(-1,1){7.5}}

\put(57.5,12.5){\line(-1,1){5}}

\put(47.5,22.5){\vector(-1,1){17.5}}

\put(45,-15){$D_1$}

\end{picture}

\caption{Reidemeister move III}\label{Rmove3}

\end{figure}

\begin{proposition}\label{reidemeister3}
Assume that $D_0$ and $D_1$ are link diagrams that are identical except in the part depicted in Figure \ref{Rmove3}. Then there is an isomorphism
\[
H_p(D_0) \cong H_p(D_1)
\]
that preserves both the cohomological grading and the quantum filtration.
\end{proposition}

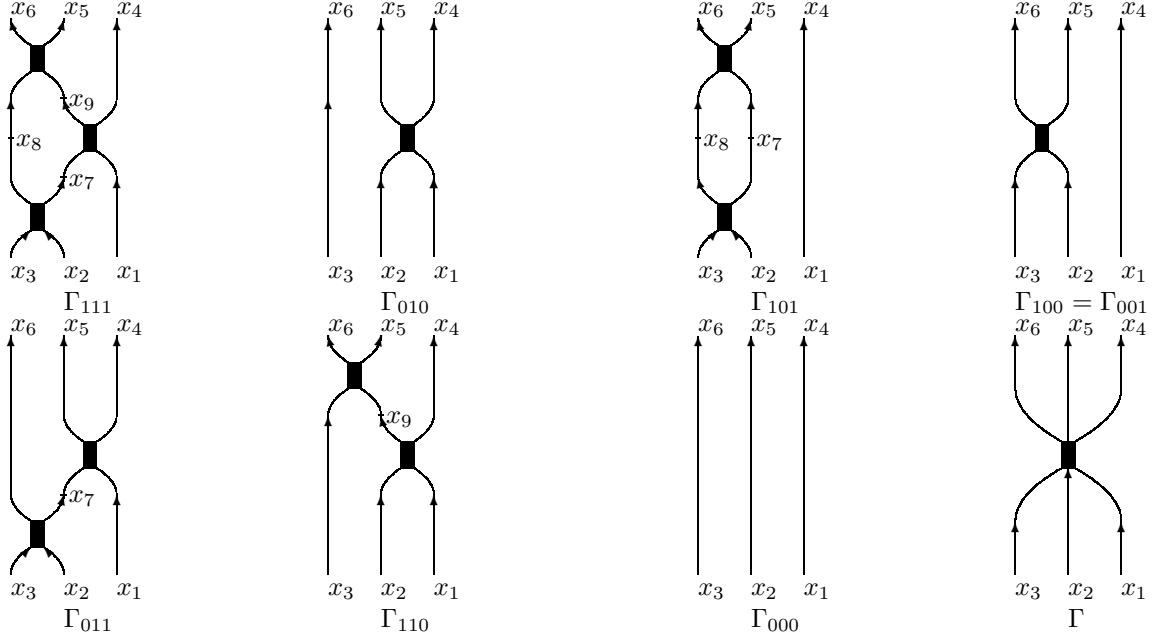
\begin{figure}[ht]

\setlength{\unitlength}{1pt}

\begin{picture}(420,240)(-210,-140)


\linethickness{.25pt}

\put(-210,30){\vector(0,1){30}}

\put(-170,0){\vector(0,1){30}}

\put(-202.5,10){\vector(1,1){0}}

\qbezier(-210,0)(-210,5)(-202.5,10)

\qbezier(-210,30)(-210,25)(-202.5,20)

\put(-197.5,10){\vector(-1,1){0}}

\qbezier(-190,0)(-190,5)(-197.5,10)

\qbezier(-190,30)(-190,25)(-197.5,20)

\qbezier(-190,30)(-190,35)(-182.5,40)

\qbezier(-190,60)(-190,55)(-182.5,50)

\put(-190,60){\vector(-1,4){0}}

\qbezier(-170,30)(-170,35)(-177.5,40)

\qbezier(-170,60)(-170,55)(-177.5,50)

\qbezier(-210,60)(-210,65)(-202.5,70)

\qbezier(-210,90)(-210,85)(-202.5,80)

\put(-210,90){\vector(-1,4){0}}

\qbezier(-190,60)(-190,65)(-197.5,70)

\qbezier(-190,90)(-190,85)(-197.5,80)

\put(-190,90){\vector(1,4){0}}

\put(-190,30){\vector(1,4){0}}

\put(-170,60){\vector(0,1){30}}

\put(-190,-20){$\Gamma_{111}$}

\put(-210,92){$x_6$}

\put(-190,92){$x_5$}

\put(-170,92){$x_4$}

\put(-210,-8){$x_3$}

\put(-190,-8){$x_2$}

\put(-170,-8){$x_1$}

\put(-208,42){$x_8$}

\put(-211,45){\line(1,0){2}}

\put(-188,57){$x_9$}

\put(-191,60){\line(1,0){2}}

\put(-188,27){$x_7$}

\put(-191,30){\line(1,0){2}}

\linethickness{5pt}

\put(-200,10){\line(0,1){10}}

\put(-180,40){\line(0,1){10}}

\put(-200,70){\line(0,1){10}}


\linethickness{.25pt}

\put(-90,0){\vector(0,1){90}}

\put(-50,0){\vector(0,1){30}}

\put(-70,0){\vector(0,1){30}}

\qbezier(-50,30)(-50,35)(-57.5,40)

\qbezier(-50,60)(-50,55)(-57.5,50)

\qbezier(-70,30)(-70,35)(-62.5,40)

\qbezier(-70,60)(-70,55)(-62.5,50)

\put(-50,60){\vector(0,1){30}}

\put(-70,60){\vector(0,1){30}}

\put(-70,-20){$\Gamma_{010}$}

\put(-90,92){$x_6$}

\put(-70,92){$x_5$}

\put(-50,92){$x_4$}

\put(-90,-8){$x_3$}

\put(-70,-8){$x_2$}

\put(-50,-8){$x_1$}

\linethickness{5pt}

\put(-60,40){\line(0,1){10}}


\linethickness{.25pt}

\put(90,0){\vector(0,1){90}}

\put(70,30){\vector(0,1){30}}

\put(62.5,10){\vector(-1,1){0}}

\qbezier(70,0)(70,5)(62.5,10)

\qbezier(70,30)(70,25)(62.5,20)

\put(57.5,10){\vector(1,1){0}}

\qbezier(50,0)(50,5)(57.5,10)

\qbezier(50,30)(50,25)(57.5,20)

\put(50,30){\vector(0,1){30}}

\qbezier(70,60)(70,65)(62.5,70)

\qbezier(70,90)(70,85)(62.5,80)

\put(70,90){\vector(1,4){0}}

\qbezier(50,60)(50,65)(57.5,70)

\qbezier(50,90)(50,85)(57.5,80)

\put(50,90){\vector(-1,4){0}}

\put(50,30){\vector(-1,4){0}}

\put(70,-20){$\Gamma_{101}$}

\put(50,92){$x_6$}

\put(70,92){$x_5$}

\put(90,92){$x_4$}

\put(50,-8){$x_3$}

\put(70,-8){$x_2$}

\put(90,-8){$x_1$}

\put(52,42){$x_8$}

\put(49,45){\line(1,0){2}}

\put(72,42){$x_7$}

\put(69,45){\line(1,0){2}}

\linethickness{5pt}

\put(60,10){\line(0,1){10}}

\put(60,70){\line(0,1){10}}


\linethickness{.25pt}

\put(210,0){\vector(0,1){90}}

\put(170,0){\vector(0,1){30}}

\put(190,0){\vector(0,1){30}}

\qbezier(170,30)(170,35)(177.5,40)

\qbezier(170,60)(170,55)(177.5,50)

\qbezier(190,30)(190,35)(182.5,40)

\qbezier(190,60)(190,55)(182.5,50)

\put(170,60){\vector(0,1){30}}

\put(190,60){\vector(0,1){30}}

\put(170,-20){$\Gamma_{100}=\Gamma_{001}$}

\put(170,92){$x_6$}

\put(190,92){$x_5$}

\put(210,92){$x_4$}

\put(170,-8){$x_3$}

\put(190,-8){$x_2$}

\put(210,-8){$x_1$}

\linethickness{5pt}

\put(180,40){\line(0,1){10}}


\linethickness{.25pt}

\put(-210,-90){\vector(0,1){60}}

\put(-190,-60){\vector(0,1){30}}

\put(-170,-120){\vector(0,1){30}}

\put(-202.5,-110){\vector(1,1){0}}

\qbezier(-210,-120)(-210,-115)(-202.5,-110)

\qbezier(-210,-90)(-210,-95)(-202.5,-100)

\put(-197.5,-110){\vector(-1,1){0}}

\qbezier(-190,-120)(-190,-115)(-197.5,-110)

\qbezier(-190,-90)(-190,-95)(-197.5,-100)

\qbezier(-190,-90)(-190,-85)(-182.5,-80)

\qbezier(-190,-60)(-190,-65)(-182.5,-70)

\qbezier(-170,-90)(-170,-85)(-177.5,-80)

\qbezier(-170,-60)(-170,-65)(-177.5,-70)

\put(-190,-90){\vector(1,4){0}}

\put(-170,-60){\vector(0,1){30}}

\put(-190,-140){$\Gamma_{011}$}

\put(-210,-28){$x_6$}

\put(-190,-28){$x_5$}

\put(-170,-28){$x_4$}

\put(-210,-128){$x_3$}

\put(-190,-128){$x_2$}

\put(-170,-128){$x_1$}

\put(-188,-93){$x_7$}

\put(-191,-90){\line(1,0){2}}

\linethickness{5pt}

\put(-200,-110){\line(0,1){10}}

\put(-180,-80){\line(0,1){10}}


\linethickness{.25pt}

\put(-90,30){\vector(0,1){30}}

\put(-50,0){\vector(0,1){30}}

\qbezier(-70,-90)(-70,-85)(-62.5,-80)

\qbezier(-70,-60)(-70,-65)(-62.5,-70)

\put(-70,-60){\vector(-1,4){0}}

\qbezier(-50,-90)(-50,-85)(-57.5,-80)

\qbezier(-50,-60)(-50,-65)(-57.5,-70)

\qbezier(-90,-60)(-90,-55)(-82.5,-50)

\qbezier(-90,-30)(-90,-35)(-82.5,-40)

\put(-90,-30){\vector(-1,4){0}}

\qbezier(-70,-60)(-70,-55)(-77.5,-50)

\qbezier(-70,-30)(-70,-35)(-77.5,-40)

\put(-70,-30){\vector(1,4){0}}

\put(-50,-60){\vector(0,1){30}}

\put(-90,-120){\vector(0,1){60}}

\put(-70,-120){\vector(0,1){30}}

\put(-50,-120){\vector(0,1){30}}

\put(-70,-140){$\Gamma_{110}$}

\put(-90,-28){$x_6$}

\put(-70,-28){$x_5$}

\put(-50,-28){$x_4$}

\put(-90,-128){$x_3$}

\put(-70,-128){$x_2$}

\put(-50,-128){$x_1$}

\put(-68,-63){$x_9$}

\put(-71,-60){\line(1,0){2}}

\linethickness{5pt}

\put(-60,-80){\line(0,1){10}}

\put(-80,-50){\line(0,1){10}}


\linethickness{.25pt}

\put(50,-120){\vector(0,1){90}}

\put(70,-120){\vector(0,1){90}}

\put(90,-120){\vector(0,1){90}}

\put(70,-140){$\Gamma_{000}$}

\put(50,-28){$x_6$}

\put(70,-28){$x_5$}

\put(90,-28){$x_4$}

\put(50,-128){$x_3$}

\put(70,-128){$x_2$}

\put(90,-128){$x_1$}


\linethickness{.25pt}

\qbezier(187.5,-70)(170,-60)(170,-50)

\put(170,-50){\vector(0,1){20}}

\put(190,-70){\vector(0,1){40}}

\qbezier(192.5,-70)(210,-60)(210,-50)

\put(210,-50){\vector(0,1){20}}

\qbezier(187.5,-80)(170,-90)(170,-100)

\put(170,-120){\vector(0,1){20}}

\put(190,-120){\vector(0,1){40}}

\qbezier(192.5,-80)(210,-90)(210,-100)

\put(210,-120){\vector(0,1){20}}

\put(190,-140){$\Gamma$}

\put(170,-28){$x_6$}

\put(190,-28){$x_5$}

\put(210,-28){$x_4$}

\put(170,-128){$x_3$}

\put(190,-128){$x_2$}

\put(210,-128){$x_1$}

\linethickness{5pt}

\put(190,-80){\line(0,1){10}}

\end{picture}

\caption{Local diagrams related to Reidemeister move
III}\label{Rmove3-decomp}

\end{figure}

\begin{proof}
Let $\Gamma$ and $\Gamma_{ijk}$, $i,j,k\in \{0,1\}$, be the planar diagrams identical to $D_0$ except in the part depicted in Figure \ref{Rmove3-decomp}. Then $(H(C_p(D_0),d_{mf}),d_\chi)$ is the total complex of
{\tiny
\[
0\rightarrow H(C_p(\Gamma_{111}),d_{mf})\{3n\} 
\xrightarrow{d_0}  
\left.%
\begin{array}{c}
  H(C_p(\Gamma_{011}),d_{mf})\{3n-1\}\\
  \oplus \\
  H(C_p(\Gamma_{101}),d_{mf})\{3n-1\}\\
  \oplus \\
  H(C_p(\Gamma_{110}),d_{mf})\{3n-1\}\\
\end{array}%
\right.
\xrightarrow{d_1} 
\left.%
\begin{array}{c}
  H(C_p(\Gamma_{001}),d_{mf})\{3n-2\}\\
  \oplus \\
  H(C_p(\Gamma_{010}),d_{mf})\{3n-2\}\\
  \oplus \\
  H(C_p(\Gamma_{100}),d_{mf})\{3n-2\}\\
\end{array}%
\right. 
\xrightarrow{d_2}
H(C_p(\Gamma_{000}),d_{mf})\{3n-3\}
\rightarrow 0,
\]} where 
\begin{eqnarray*}
d_0 & = &
\left(%
\begin{array}{c}
  \eta_1\\
  \eta_2\\
  \eta_3\\
\end{array}%
\right), \\
d_1 & = & 
\left(%
\begin{array}{ccc}
  \xi_{11} & \xi_{12} & 0 \\
  \xi_{21} &  0       & \xi_{23} \\
  0        & \xi_{32} & \xi_{33}  \\
\end{array}%
\right), \\
d_2 & = & 
\left(%
\begin{array}{ccc}
  \zeta_1 & \zeta_2 & \zeta_3 \\
\end{array}%
\right)
\end{eqnarray*}
are filtered homomorphisms of chain complexes. 

By Propositions \ref{moy2} and \ref{moy4},
\begin{eqnarray*}
H(C_p(\Gamma_{111}),d_{mf}) & \cong & H(C_p(\Gamma),d_{mf}) \oplus H(C_p(\Gamma_{100}),d_{mf}), \\
H(C_p(\Gamma_{101}),d_{mf}) & \cong & H(C_p(\Gamma_{100}),d_{mf})\{-1\} \oplus H(C_p(\Gamma_{100}),d_{mf})\{1\}.
\end{eqnarray*}
Let $J_{\pm1}:H(C_p(\Gamma_{100}),d_{mf})\{\pm1\} \rightarrow H(C_p(\Gamma_{101}),d_{mf})$ and $P_{\pm1}:H(C_p(\Gamma_{101}),d_{mf}) \rightarrow H(C_p(\Gamma_{100}),d_{mf})\{\pm1\}$ be the inclusion and projection from the above decomposition of $H(C_p(\Gamma_{101}),d_{mf})$. Then there is a filtered homomorphism $J': H(C_p(\Gamma_{100}),d_{mf}) \rightarrow H(C_p(\Gamma_{111}),d_{mf})$ given by the composition
\[
H(C_p(\Gamma_{100}),d_{mf}) \xrightarrow{J_{-1}} H(C_p(\Gamma_{101}),d_{mf}) \xrightarrow{\chi_0} H(C_p(\Gamma_{111}),d_{mf}),
\]
where $\chi_0$ is induced by the $\chi_0$-map associated to the wide edge on the right side. By Proposition \ref{hmf-tensor},
\[
\Hom_{hmf}(C_p(\Gamma_{100}),C_p(\Gamma_{100})) \cong \Hom_{hmf}(C_p(\Gamma_{100}),C_p(\Gamma_{111})) \cong \C.
\]
So $J'$ is a multiple of the inclusion map $J: H(C_p(\Gamma_{100}),d_{mf}) \rightarrow H(C_p(\Gamma_{111}),d_{mf})$ from the above decomposition of $H(C_p(\Gamma_{111}),d_{mf})$. Let $\chi_1:H(C_p(\Gamma_{111}),d_{mf}) \rightarrow H(C_p(\Gamma_{101}),d_{mf})$ be induced by the $\chi_1$-map associated to the wide edge on the right side. Using $\chi_1\circ\chi_0=(x_7-x_4)\id_{H(C_p(\Gamma_{101}),d_{mf})}$, it is easy to check that $P_{1}\circ\chi_1 \circ J': H(C_p(\Gamma_{100}),d_{mf}) \rightarrow H(C_p(\Gamma_{100}),d_{mf})$ is a non-zero multiple of the identity map. This implies that $J'\neq 0$ and $P_1\circ\eta_2\circ J:H(C_p(\Gamma_{100}),d_{mf}) \rightarrow H(C_p(\Gamma_{100}),d_{mf})$ is an isomorphism of filtered chain complexes. By Lemma \ref{gaussian-elimination}, $(H(C_p(D_0),d_{mf}),d_\chi)$ is quasi-isomorphic to the total complex of
{\tiny
\[
0\rightarrow H(C_p(\Gamma),d_{mf})\{3n\} 
\rightarrow  
\left.%
\begin{array}{c}
  H(C_p(\Gamma_{011}),d_{mf})\{3n-1\}\\
  \oplus \\
  H(C_p(\Gamma_{100}),d_{mf})\{3n-2\}\\
  \oplus \\
  H(C_p(\Gamma_{110}),d_{mf})\{3n-1\}\\
\end{array}%
\right.
\rightarrow 
\left.%
\begin{array}{c}
  H(C_p(\Gamma_{001}),d_{mf})\{3n-2\}\\
  \oplus \\
  H(C_p(\Gamma_{010}),d_{mf})\{3n-2\}\\
  \oplus \\
  H(C_p(\Gamma_{100}),d_{mf})\{3n-2\}\\
\end{array}%
\right. 
\rightarrow
H(C_p(\Gamma_{000}),d_{mf})\{3n-3\}
\rightarrow 0.
\]}

By Lemma \ref{Rmove2a-lemma}, the homomorphism $\xi_{32}\circ J_{-1}: H(C_p(\Gamma_{100}),d_{mf}) \rightarrow H(C_p(\Gamma_{100}),d_{mf})$ is an isomorphism of filtered chain complexes. So, by Lemma \ref{gaussian-elimination}, $(H(C_p(D_0),d_{mf}),d_\chi)$ is quasi-isomorphic to the total complex of
{\tiny
\[
0\rightarrow H(C_p(\Gamma),d_{mf})\{3n\} 
\xrightarrow{d_0'}  
\left.%
\begin{array}{c}
  H(C_p(\Gamma_{011}),d_{mf})\{3n-1\}\\
  \oplus \\
  H(C_p(\Gamma_{110}),d_{mf})\{3n-1\}\\
\end{array}%
\right.
\xrightarrow{d_1'} 
\left.%
\begin{array}{c}
  H(C_p(\Gamma_{001}),d_{mf})\{3n-2\}\\
  \oplus \\
  H(C_p(\Gamma_{010}),d_{mf})\{3n-2\}\\
\end{array}%
\right. 
\xrightarrow{d_2'}
H(C_p(\Gamma_{000}),d_{mf})\{3n-3\}
\rightarrow 0,
\]} where
\begin{eqnarray*}
d_0' & = &
\left(%
\begin{array}{c}
  f_1\\
  f_2\\
\end{array}%
\right), \\
d_1' & = & 
\left(%
\begin{array}{cc}
  f_3 & f_5 \\
  f_4 &  f_6 \\
\end{array}%
\right), \\
d_2' & = & 
\left(%
\begin{array}{cc}
  f_7 & f_8 \\
\end{array}%
\right).
\end{eqnarray*}
Similarly, $(H(C_p(D_1),d_{mf}),d_\chi)$ is quasi-isomorphic to the total complex of
{\tiny
\[
0\rightarrow H(C_p(\Gamma),d_{mf})\{3n\} 
\xrightarrow{d_0''}  
\left.%
\begin{array}{c}
  H(C_p(\Gamma_{011}),d_{mf})\{3n-1\}\\
  \oplus \\
  H(C_p(\Gamma_{110}),d_{mf})\{3n-1\}\\
\end{array}%
\right.
\xrightarrow{d_1''} 
\left.%
\begin{array}{c}
  H(C_p(\Gamma_{001}),d_{mf})\{3n-2\}\\
  \oplus \\
  H(C_p(\Gamma_{010}),d_{mf})\{3n-2\}\\
\end{array}%
\right. 
\xrightarrow{d_2''}
H(C_p(\Gamma_{000}),d_{mf})\{3n-3\}
\rightarrow 0,
\]} where
\begin{eqnarray*}
d_0'' & = &
\left(%
\begin{array}{c}
  f_1'\\
  f_2'\\
\end{array}%
\right), \\
d_1'' & = & 
\left(%
\begin{array}{cc}
  f_3' & f_5' \\
  f_4' &  f_6' \\
\end{array}%
\right), \\
d_2'' & = & 
\left(%
\begin{array}{cc}
  f_7' & f_8' \\
\end{array}%
\right).
\end{eqnarray*}
Here, $f_i$ and $f_i'$ are filtered homomorphisms of filtered chain complexes for $i=1,2,\cdots,8$.

In the rest of this proof, we assume $n\geq 3$. (The proof is easier when $n=2$. In fact, when $n=2$, the $H_p$ is either the Khovanov homology defined in \cite{K1} or Lee's variant defined in \cite{Lee2}, the invariance of which have been established.)

Consider the filtered homomorphisms \[
f_i,f_i': H(C_p(Z_1),d_{mf}) \rightarrow H(C_p(Z_2),d_{mf})\{-1\}
\] 
from the above chain complexes. Let $\hat{Z}_1$ and $\hat{Z}_2$ be the parts of $Z_1$ and $Z_2$ visible in Figure \ref{Rmove3-decomp}. (These are open graphs.) Then $f_i$ and $f_i'$ are induced by filtered matrix factorization homomorphisms 
\[
\widetilde{f}_i, \widetilde{f}_i':C_p(\hat{Z}_1)\rightarrow C_p(\hat{Z}_2)\{-1\},
\] 
in the sense that $f_i$ and $f_i'$ are the homomorphisms of cohomologies induced by chain maps 
\[
\id\otimes\widetilde{f}_i,~\id\otimes\widetilde{f}_i': (C_p(Z_1),d_{mf}) \rightarrow (C_p(Z_2),d_{mf})\{-1\},
\] 
where $\id$ is the identity map of the matrix factorization associated to the part of $Z_1$ outside Figure \ref{Rmove3-decomp}. Using Proposition \ref{hmf-tensor}, one can check that
\[
\Hom_{hmf}(C_p(\hat{Z}_1),C_p(\hat{Z}_2)) \cong \C.
\]
(Compare to Lemma 26 of \cite{KR1}.) So $\widetilde{f}_i$ is homotopic to a multiple of $\widetilde{f}_i'$, and, hence, $f_i$ is a multiple of $f_i'$. And, using the argument of the proof of Lemma 27 of \cite{KR1}, we get $f_i\circ f_j\neq0$, $f_i'\circ f_j'\neq0$, where the domain of $f_i$ is the codomain of $f_j$. Thus, for $i=1,2,\cdots,8$, there exists an $\lambda_i\in \C\setminus \{0\}$ such that $f_i'=\lambda_if_i$. And, using $d'\circ d'=d''\circ d''=0$, we have
\begin{eqnarray*} 
\lambda_1\lambda_3 = \lambda_2\lambda_5, & & \lambda_1\lambda_4 = \lambda_2\lambda_6, \\
\lambda_3\lambda_7 = \lambda_4\lambda_8, & & \lambda_5\lambda_7 = \lambda_6\lambda_8.
\end{eqnarray*}
Now define a map $F$ between the two reduced chain complexes by
\begin{eqnarray*}
H(C_p(\Gamma),d_{mf})\{3n\} & \xrightarrow{\id} & H(C_p(\Gamma),d_{mf})\{3n\}, \\
\left.%
\begin{array}{c}
  H(C_p(\Gamma_{011}),d_{mf})\{3n-1\}\\
  \oplus \\
  H(C_p(\Gamma_{110}),d_{mf})\{3n-1\}\\
\end{array}%
\right.
& \xrightarrow{\left(%
\begin{array}{cc}
  \lambda_1 & 0 \\
  0 &  \lambda_2 \\
\end{array}%
\right)} &
\left.%
\begin{array}{c}
  H(C_p(\Gamma_{011}),d_{mf})\{3n-1\}\\
  \oplus \\
  H(C_p(\Gamma_{110}),d_{mf})\{3n-1\}\\
\end{array}%
\right., \\
\left.%
\begin{array}{c}
  H(C_p(\Gamma_{001}),d_{mf})\{3n-2\}\\
  \oplus \\
  H(C_p(\Gamma_{010}),d_{mf})\{3n-2\}\\
\end{array}%
\right. & \xrightarrow{\left(%
\begin{array}{cc}
  \lambda_1\lambda_3 & 0 \\
  0 &  \lambda_1\lambda_4 \\
\end{array}%
\right)} &
\left.%
\begin{array}{c}
  H(C_p(\Gamma_{001}),d_{mf})\{3n-2\}\\
  \oplus \\
  H(C_p(\Gamma_{010}),d_{mf})\{3n-2\}\\
\end{array}%
\right., \\
H(C_p(\Gamma_{000}),d_{mf})\{3n-3\} & \xrightarrow{\lambda_1\lambda_3\lambda_7\id} & H(C_p(\Gamma_{000}),d_{mf})\{3n-3\}.
\end{eqnarray*}
It is straightforward to check that $F$ is an isomorphism of filtered chain complexes. Thus, 
\[
H_p(D_0) \cong H_p(D_1).
\]
\end{proof}

Clearly, Theorem \ref{fil-inv} follows from Propositions \ref{reidemeister1}, \ref{reidemeister2a}, \ref{reidemeister2b} and \ref{reidemeister3}.

\subsection{Invariance of the spectral sequence $\{E_k\}$}

Next, we prove Theorem \ref{spec} by showing that the spectral sequence $\{E_k\}$ defined in Subsection \ref{spec-def} is invariant under Reidemeister moves. 

\begin{proof}[Proof of Theorem \ref{spec}]
Let $D_1$ and $D_2$ be two link diagrams that differ by a Reidemeister move. In the proofs of Propositions \ref{reidemeister1}, \ref{reidemeister2a}, \ref{reidemeister2b} and \ref{reidemeister3}, we constructed an quasi-isomorphism of filtered chain complexes
\[
f: (H(C_p(D_1),d_{mf}),d_\chi) \rightarrow (H(C_p(D_2),d_{mf}),d_\chi).
\]
By the definition of $\{E_k\}$, $f$ induces $\zed\oplus\zed$-graded chain maps
\[
f_k: E_k(D_1) \rightarrow E_k(D_2).
\]
Let $\hat{d}_{mf}$ and $\hat{d}_\chi$ be the top homogeneous parts of $d_{mf}$ and $d_\chi$. Then the Khovanov-Rozansky $\mathfrak{sl}(n)$-cohomology $H_n$ with the potential polynomial $\hat{p}(x)=x^{n+1}$ is defined by
\[
H_n(D_i) = (H(C_p(D_i),\hat{d}_{mf}),\hat{d}_\chi).
\]
Specially, the above construction gives an quasi-isomorphism of $\zed\oplus\zed$-graded chain complexes
\[
\hat{f}: (H(C_p(D_1),\hat{d}_{mf}),\hat{d}_\chi) \rightarrow (H(C_p(D_2),\hat{d}_{mf}),\hat{d}_\chi).
\]
From Subsection \ref{spec-def}, for $i=1,2$, there is an isomorphism of $\zed\oplus\zed$-graded chain complexes
\[
\phi_i: (H(C_p(D_i),\hat{d}_{mf}),\hat{d}_\chi) \rightarrow E_0(D_i).
\]
It is straightforward to check that the following square commutes.
\[
\begin{CD}
(H(C_p(D_1),\hat{d}_{mf}),\hat{d}_\chi) @>{\hat{f}}>> (H(C_p(D_2),\hat{d}_{mf}),\hat{d}_\chi) \\
@VV{\phi_1}V @VV{\phi_2}V \\
E_0(D_1) @>{f_0}>> E_0(D_2) 
\end{CD}
\]
Thus, $f_0$ is also an quasi-isomorphism. So $f_1$ is an isomorphism of $\zed\oplus\zed$-graded chain complexes, which implies that, for $k\geq1$, $f_k$ is an isomorphism of $\zed\oplus\zed$-graded chain complexes. (See e.g. Mapping Lemma 5.2.4 of \cite{Chuck-Weibel}.)
\end{proof}

\subsection{Link cobodisms}\label{linkcobodisms}

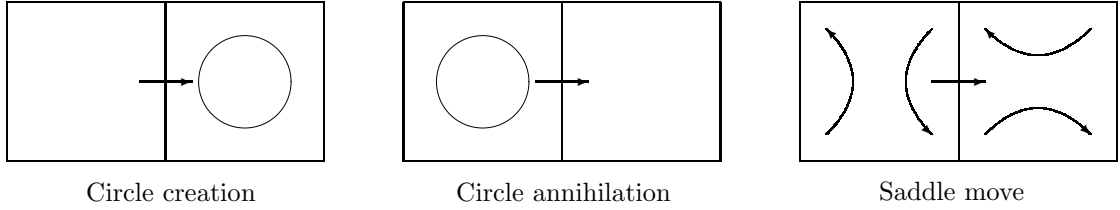
\begin{figure}[ht]

\setlength{\unitlength}{1pt}

\begin{picture}(420,75)(-210,-15)


\put(-210,0){\line(1,0){120}}

\put(-210,60){\line(1,0){120}}

\put(-210,0){\line(0,1){60}}

\put(-150,0){\line(0,1){60}}

\put(-90,0){\line(0,1){60}}

\put(-160,30){\vector(1,0){20}}

\put(-120,30){\circle{35}}

\put(-180,-15){\text{Circle creation}}


\put(-60,0){\line(1,0){120}}

\put(-60,60){\line(1,0){120}}

\put(-60,0){\line(0,1){60}}

\put(0,0){\line(0,1){60}}

\put(60,0){\line(0,1){60}}

\put(-10,30){\vector(1,0){20}}

\put(-30,30){\circle{35}}

\put(-40,-15){\text{Circle annihilation}}


\put(90,0){\line(1,0){120}}

\put(90,60){\line(1,0){120}}

\put(90,0){\line(0,1){60}}

\put(150,0){\line(0,1){60}}

\put(210,0){\line(0,1){60}}

\put(140,30){\vector(1,0){20}}

\qbezier(100,10)(120,30)(100,50)

\put(100,50){\vector(-1,1){0}}

\qbezier(140,10)(120,30)(140,50)

\put(140,10){\vector(1,-1){0}}

\qbezier(160,10)(180,30)(200,10)

\put(200,10){\vector(1,-1){0}}

\qbezier(160,50)(180,30)(200,50)

\put(160,50){\vector(-1,1){0}}

\put(120,-15){\text{Saddle move}}

\end{picture}

\caption{Morse moves}\label{morse-moves}

\end{figure}

A cobodism from link $L_0$ to link $L_1$ is a properly smoothly embedded oriented surface $S$ in $S^3\times[0,1]$ such that $\partial S = (-L_0\times\{0\}) \sqcup (L_1\times\{1\})$, where $-L_0$ is $L_0$ with the opposite orientation. A link cobodism is said to be elementary if it corresponds to either a Reidemeister move or a Morse move (depicted in Figure \ref{morse-moves}). Any link cobodism admits a movie presentation, that is, a decomposition into elementary cobodisms. Two movie presentations present the same cobodism if one can be changed into the other by the movie moves depicted in Figures 59, 60 in \cite{KR1}. (See \cite{CS1,CS2} for more details.)

If a link cobodism $S$ from $L_0$ to $L_1$ corresponds to a Reidemeister move, then define 
\[
\widetilde{\Psi}_S: H(C_p(L_0),d_{mf}) \rightarrow H(C_p(L_1),d_{mf})
\]
to be the quasi-isomorphism constructed in Subsections \ref{R-Move-I}-\ref{R-Move-III}. If $S$ corresponds to a Morse move, then define 
\[
\widetilde{\Psi}_S: H(C_p(L_0),d_{mf}) \rightarrow H(C_p(L_1),d_{mf}) \left\langle v_1(S) \right\rangle \{(n-1)\chi(S)\}
\] 
to be the homomorphism induced by the corresponding graph cobodism defined in Subsection \ref{graph-cobodisms}, where $v_1(S)$ is the sum of numbers of components of $L_0$, $L_1$. If $(S_1,S_2,\cdots,S_m)$ is a movie presentation of a cobodism $S$ from $L_0$ to $L_1$, where each $S_i$ is an elementary link cobodism, then define 
\[
\widetilde{\Psi}_{(S_1,S_2,\cdots,S_m)} = \widetilde{\Psi}_{S_m}\circ\cdots\circ \widetilde{\Psi}_{S_2} \circ  \widetilde{\Psi}_{S_1}: H(C_p(L_0),d_{mf}) \rightarrow H(C_p(L_1),d_{mf}) \left\langle v_1(S) \right\rangle \{(n-1)\chi(S)\}.
\] 
$\widetilde{\Psi}_{(S_1,S_2,\cdots,S_m)}$ is a chain map that preserves the cohomological grading and the quantum filtration. So, it induces a homomorphism 
\[
\Psi_{(S_1,S_2,\cdots,S_m),k}: E_k(L_0) \rightarrow E_k(L_1),  
\]
which shifts the first $\zed$-grading by $-(n-1)\chi(S)$.

Using results from \cite{KR1}, we prove Theorem \ref{cobo-homo} by showing that, when $k\geq1$, up to multiplication by a non-zero scalar, $\Psi_{(S_1,S_2,\cdots,S_m),k}$ depends only on the cobodism $S$ not on the choice of the movies presentation $(S_1,S_2,\cdots,S_m)$.

\begin{proof}[Proof of Theorem \ref{cobo-homo}]
Let $(S_1,S_2,\cdots,S_m)$ and $(S_1',S_2',\cdots,S_l')$ be two movie presentations of a link cobodism $S$ from $L_0$ to $L_1$. Consider the potential polynomial $x^{n+1}$. The above construction gives two filtered chain maps
\begin{eqnarray*}
\widetilde{\Psi}_{(S_1,S_2,\cdots,S_m)}^{(0)} & : & (H(C_p(L_0),d_{mf}^{(0)}),d_\chi^{(0)}) \rightarrow (H(C_p(L_1),d_{mf}^{(0)}),d_\chi^{(0)}) \left\langle v_1(S) \right\rangle \{(n-1)\chi(S)\}, \\
\widetilde{\Psi}_{(S_1',S_2',\cdots,S_l')}^{(0)} & : & (H(C_p(L_0),d_{mf}^{(0)}),d_\chi^{(0)}) \rightarrow (H(C_p(L_1),d_{mf}^{(0)}),d_\chi^{(0)}) \left\langle v_1(S) \right\rangle \{(n-1)\chi(S)\},
\end{eqnarray*}
where $d_{mf}^{(0)}$ and $d_\chi^{(0)}$ are the top homogeneous parts of differentials $d_{mf}$ and $d_\chi$ of $C_p(L_i)$. Note that the cohomologies of these chain complexes are the $\mathfrak{sl}(n)$-cohomologies of $L_0$ and $L_1$, and the homomorphisms induced by these chain maps are exactly the homomorphisms of the $\mathfrak{sl}(n)$-cohomologies induced by these movies presentations. The following squares commute.
\[
\begin{CD}
(H(C_p(L_0),d_{mf}^{(0)}),d_\chi^{(0)}) @>{\widetilde{\Psi}_{(S_1,S_2,\cdots,S_m)}^{(0)}}>> (H(C_p(L_1),d_{mf}^{(0)}),d_\chi^{(0)}) \\
@VV{\phi_0}V @VV{\phi_1}V \\
E_0(L_0) @>{\Psi_{(S_1,S_2,\cdots,S_m),0}}>> E_0(L_1) 
\end{CD},
\]
\[
\begin{CD}
(H(C_p(L_0),d_{mf}^{(0)}),d_\chi^{(0)}) @>{\widetilde{\Psi}_{(S_1',S_2',\cdots,S_l')}^{(0)}}>> (H(C_p(L_1),d_{mf}^{(0)}),d_\chi^{(0)}) \\
@VV{\phi_0}V @VV{\phi_1}V \\
E_0(L_0) @>{\Psi_{(S_1',S_2',\cdots,S_l'),0}}>> E_0(L_1) 
\end{CD},
\]
where $\phi_i$ is the isomorphism from Subsection \ref{spec-def}. Thus, $\Psi_{(S_1,S_2,\cdots,S_m),1}$ and $\Psi_{(S_1',S_2',\cdots,S_l'),1}$ are the homomorphisms of the cohomologies induced by $\widetilde{\Psi}_{(S_1,S_2,\cdots,S_m)}^{(0)}$ and $\widetilde{\Psi}_{(S_1',S_2',\cdots,S_l')}^{(0)}$. By Proposition 37 of \cite{KR1}, there is a $\lambda\in\C\setminus\{0\}$ such that $\Psi_{(S_1,S_2,\cdots,S_m),1}=\lambda\Psi_{(S_1',S_2',\cdots,S_l'),1}$. Hence,  for $k\geq1$,
\[
\Psi_{(S_1,S_2,\cdots,S_m),k}=\lambda\Psi_{(S_1',S_2',\cdots,S_l'),k}.
\]
\end{proof}

\begin{remark}
Clearly, the chain map $\widetilde{\Psi}_{(S_1,S_2,\cdots,S_m)}$ also induces a homomorphism from $H_p(L_0)$ to $H_p(L_1)$. I expect that, up to multiplication by a non-zero scalar, this homomorphism is also independent of the choice of the movie presentation $(S_1,S_2,\cdots,S_m)$. But the proof is probably harder since we can not use the above shortcut in this situation.
\end{remark}

\section{Slice Euler Characteristic and Generalized Rasmussen Invariants}\label{rasmussen}

The goal of this section is to prove Theorem \ref{ras-genus}. Our proof is based on Gornik's computations in \cite{Gornik} and Rasmussen's arguments in \cite{Ras1}. Since Gornik's computations can be generalized to any monic polynomial $p(x)$ of degree $n+1$ whose derivative $p'(x)$ has $n$ distinct roots, our proof can be carried out for any of such potentials. But, for simplicity, we work only with the potential polynomial $x^{n+1}-(n+1)x$ and leave the general case to the readers. In the rest of this section, we fix $n\geq 2$ and
\[
p(x)=x^{n+1}-(n+1)x.
\]

\begin{lemma}[\cite{Gornik}]\label{gornik-circle}
If $\Gamma$ is a round circle, then $H_p(\Gamma)\cong\C[x]/(x^n-1)$, as a $\C$-linear space, has a basis $\{f_0(x),f_1(x),\cdots,f_{n-1}(x)\}$, where
\[
f_k(x) = \sum_{l=0}^{n-1} x^l e^{-\frac{2kl\pi}{n}i} \in \C[x]/(x^n-1),
\]
and $i$ is the imaginary unit.
\end{lemma}

The next 2 lemmas describe the effects of graph cobodisms on $H_p$.

\begin{lemma}\label{gornik-epsilon-iota}
If $\Gamma$ is a round circle, then
\begin{eqnarray*}
\iota & : & \C \rightarrow H_p(\Gamma), \\
\varepsilon & : & H_p(\Gamma) \rightarrow \C
\end{eqnarray*}
are given by 
\begin{eqnarray*}
\iota(1) & = & \lambda \sum_{k=0}^{n-1} f_k(x), \\
\varepsilon (f_k(x)) & = & \mu e^{\frac{2k\pi}{n}i},
\end{eqnarray*}
where $\lambda,\mu\in\C\setminus\{0\}$ come from the definitions of $\iota, \varepsilon$, and $f_k(x)$ is the basis element defined in Lemma \ref{gornik-circle}.
\end{lemma}
\begin{proof}
Straightforward.
\end{proof}

\begin{figure}[ht]

\setlength{\unitlength}{1pt}

\begin{picture}(420,45)(-210,-15)


\put(-80,15){\circle{30}}

\put(-40,15){\circle{30}}

\put(-96,15){\line(1,0){2}}

\put(-105,15){$x_1$}

\put(-26,15){\line(1,0){2}}

\put(-23,15){$x_2$}

\put(-62,16){\tiny{$1$}}

\multiput(-65,15)(4,0){3}{\line(1,0){2}}

\put(-60,-15){$\Gamma$}


\put(10,23){$\eta_1$}

\put(0,20){\vector(1,0){30}}

\put(10,3){$\eta_2$}

\put(30,10){\vector(-1,0){30}}


\put(60,15){\circle{30}}

\put(61,14){\tiny{$2$}}

\multiput(60,0)(0,4){8}{\line(0,1){2}}

\put(44,15){\line(1,0){2}}

\put(39,15){$x$}

\put(60,-15){$\Gamma_0$}

\end{picture}

\caption{Saddle moves}\label{gornik-saddle-move}

\end{figure}
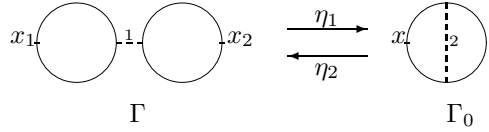

\begin{lemma}\label{gornik-circle-saddle}
Let $\Gamma$ and $\Gamma_0$ be the two graphs depicted in Figure \ref{gornik-saddle-move}. Then $\Gamma$ has the basis $\{f_k(x_1)f_l(x_2)~|~ k,l=0,1,\cdots,n-1\}$, and $\Gamma$ has the basis $\{f_k(x)~|~ k=0,1,\cdots,n-1\}$. If $\eta_1$ and $\eta_2$ are the two saddle moves depicted in Figure \ref{gornik-saddle-move}, then
\begin{eqnarray*}
\eta_1(f_k(x_1) f_l(x_2)) & = & n \delta_{kl} f_k(x), \\
\eta_2(f_k(x)) & = & n e^{-\frac{2k\pi}{n}i}  f_k(x_1) f_k(x_2).
\end{eqnarray*}
\end{lemma}
\begin{proof}
By the definition of $\eta_1$,  
\begin{eqnarray*}
\eta_1(f_k(x_1) f_l(x_2)) & = & f_k(x)f_l(x)~(\in \C[x]/(x^n-1)\cong H_p(\Gamma_0)), \\
                                 & = & (\sum_{r=0}^{n-1} x^r e^{-\frac{2kr\pi}{n}i})(\sum_{q=0}^{n-1} x^q e^{-\frac{2lq\pi}{n}i}), \\
                                 & = & \sum_{m=0}^{n-1}x^m e^{-\frac{2km\pi}{n}i} (\sum_{r=0}^m e^{-\frac{2r(k-l)\pi}{n}i}) + \sum_{m=0}^{n-2}x^{m+n} e^{-\frac{2k(m+n)\pi}{n}i} (\sum_{r=m+1}^{n-1}e^{-\frac{2r(k-l)\pi}{n}i}), \\
                                 & = & (\sum_{m=0}^{n-1}x^m e^{-\frac{2km\pi}{n}i}) (\sum_{r=0}^{n-1} e^{-\frac{2r(k-l)\pi}{n}i}), \\
                                 & = & n\delta_{kl} f_k(x).
\end{eqnarray*}

Note that, for $p(x)=x^{n+1}-(n+1)x$, 
\[
e_{ijk}=\frac{(x_k-x_j)p(x_i)+(x_i-x_k)p(x_j)+(x_j-x_i)p(x_k)}{2(x_i-x_j)(x_j-x_k)(x_k-x_i)} = \frac{1}{2}\sum_{a+b+c=n-1}x_i^a x_j^b x_k^c.
\]
By definition of $\eta_2$,
\begin{eqnarray*}
\eta_2(f_k(x)) & = & -2(e_{122}+e_{112})f_k(x_1), \\
               & = & n(\sum_{c=0}^{n-1}x_1^{n-1-c}x_2^c)(\sum_{r=0}^{n-1} x_1^r e^{-\frac{2kr\pi}{n}i}), \\
               & = & n\sum_{m=0}^{n-1} x_1^m e^{-\frac{2km\pi}{n}i} (\sum_{r=0}^m x_2^{n-1-m+r} e^{-\frac{2k(r-m)\pi}{n}i}) + n \sum_{m=0}^{n-2}x_1^{m+n}e^{-\frac{2km\pi}{n}i} (\sum_{r=m+1}^{n-1}x_2^{r-m-1} e^{-\frac{2k(r-m)\pi}{n}i}), \\
               & = & n e^{-\frac{2k\pi}{n}i}  (\sum_{m=0}^{n-1} x_1^m e^{-\frac{2km\pi}{n}i}) (\sum_{q=0}^{n-1} x_2^q e^{-\frac{2kq\pi}{n}i}), \\
               & = & n e^{-\frac{2k\pi}{n}i} f_k(x_1)f_k(x_2).
\end{eqnarray*}
\end{proof}

Next, we recall Gornik's definition of admissible states of a graph. Let $\Gamma$ be a closed graph whose edges are $e_1,\cdots,e_m$. We mark each edge $e_j$ of $\Gamma$ by a single variable $x_j$. Then $H_p(\Gamma)$ is a $\C[x_1,\cdots, x_m]$-module. By Lemma \ref{m(p')}, the action of $\C[x_1,\cdots, x_m]$ on $H_p(\Gamma)$ factors through $\C[x_1,\cdots, x_m]/(x_1^n-1,\cdots,x_m^n-1)$. A state of $\Gamma$ is a function $\varphi:\{e_1,\cdots,e_m\}\rightarrow \{0,1,\cdots,n-1\}$. A state $\varphi$ of $\Gamma$ is called admissible if, at every (double) wide edge in $\Gamma$, $\varphi(e_1)\neq\varphi(e_2)$ and $\{\varphi(e_1),\varphi(e_2)\}=\{\varphi(e_3),\varphi(e_4)\}$ where $e_1$, $e_2$ are the two entering regular edges at this wide edge, and $e_3$, $e_4$ are the two exiting regular edges See Figure 6 of \cite{Gornik}.

For a state $\varphi$ of $\Gamma$, define 
\[
Q_\varphi = \prod_{j=1}^m f_{\varphi(e_j)}(x_j) \in \C[x_1,\cdots, x_m]/(x_1^n-1,\cdots,x_m^n-1).
\]

\begin{theorem}[Theorem 4, \cite{Gornik}]\label{admissible-states}
For any state $\varphi$ of $\Gamma$ and $\alpha\in H_p(\Gamma)$, we have
\[
\alpha\in Q_\varphi H_p(\Gamma) ~ \Longleftrightarrow ~ x_j\alpha = e^{\frac{2\varphi(e_j)\pi}{n}i}\alpha, ~\forall~j=1,\cdots,m.
\]
If a state $\varphi$ is not admissible, then $Q_\varphi H_p(\Gamma)=0$. For each admissible state $\varphi$ of $\Gamma$, the $\C$-linear space $Q_\varphi H_p(\Gamma)$ is $1$-dimensional, and
\[
H_p(\Gamma) = \bigoplus_\varphi Q_\varphi H_p(\Gamma),
\]
where $\varphi$ runs through all admissible states of $\Gamma$. 
\end{theorem}

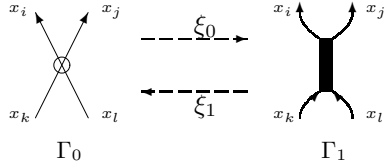
\begin{figure}[ht]

\setlength{\unitlength}{1pt}

\begin{picture}(420,60)(-210,-15)


\put(-40,0){\vector(-1,2){20}}

\put(-60,0){\vector(1,2){20}}

\put(-50,20){\circle{6}}

\put(-52,-15){\small{$\Gamma_0$}}

\put(-70,40){\tiny{$x_i$}}

\put(-35,40){\tiny{$x_j$}}

\put(-70,0){\tiny{$x_k$}}

\put(-35,0){\tiny{$x_l$}}


\put(0,32){$\xi_0$}

\multiput(-20,30)(7,0){6}{\line(1,0){5}}

\multiput(-20,10)(7,0){6}{\line(1,0){5}}

\put(0,2){$\xi_1$}

\put(20,30){\vector(1,0){0}}

\put(-20,10){\vector(-1,0){0}}


\put(47.5,10){\vector(1,1){0}}

\qbezier(40,0)(40,5)(47.5,10)

\qbezier(40,40)(40,35)(47.5,30)

\put(40,44){\vector(0,1){0}}

\put(52.5,10){\vector(-1,1){0}}

\qbezier(60,0)(60,5)(52.5,10)

\qbezier(60,40)(60,35)(52.5,30)

\put(60,44){\vector(0,1){0}}

\put(30,40){\tiny{$x_i$}}

\put(65,40){\tiny{$x_j$}}

\put(30,0){\tiny{$x_k$}}

\put(65,0){\tiny{$x_l$}}

\linethickness{5pt}

\put(50,10){\vector(0,1){20}}

\put(48,-15){\small{$\Gamma_1$}}

\end{picture}

\caption{A virtual crossing and a wide edge}\label{virtual-crossing}

\end{figure}

To compute the $H_p$-cohomology of a link, Gornik \cite{Gornik} introduced a new type of singularities: virtual crossings, which is depicted on the left side of Figure \ref{virtual-crossing}. Consider the two graphs in Figure \ref{maps}. The matrix factorization
$C_p(\Gamma_0)$ is
\[
\left(%
\begin{array}{cc}
  \pi_{il} & x_i-x_l \\
  \pi_{jk} & x_j-x_k \\
\end{array}%
\right)_{R}.
\]
Explicitly, this is
\[
\left[%
\begin{array}{c}
  R \\
  R\{2-2n\} \\
\end{array}%
\right] \xrightarrow{P_0}
\left[%
\begin{array}{c}
  R\{1-n\} \\
  R\{1-n\} \\
\end{array}%
\right] \xrightarrow{P_1}
\left[%
\begin{array}{c}
  R \\
  R\{2-2n\} \\
\end{array}%
\right],
\]
where
\[
P_0=\left(%
\begin{array}{cc}
  \pi_{il} & x_j-x_k \\
  \pi_{jk} & -x_i+x_l \\
\end{array}%
\right), \hspace{1cm}
P_1=\left(%
\begin{array}{cc}
  x_i-x_l & x_j-x_k \\
  \pi_{jk} & -\pi_{il} \\
\end{array}%
\right).
\]
The matrix factorization $C_p(\Gamma_1)$ is
\[
\left(%
\begin{array}{cc}
  u_{ijkl} & x_i+x_j-x_k-x_l \\
  v_{ijkl} & x_ix_j-x_kx_l \\
\end{array}%
\right)_{R}\{-1\}.
\]
Explicitly, this is
\[
\left[%
\begin{array}{c}
  R\{-1\} \\
  R\{3-2n\} \\
\end{array}%
\right] \xrightarrow{Q_0}
\left[%
\begin{array}{c}
  R\{-n\} \\
  R\{2-n\} \\
\end{array}%
\right] \xrightarrow{Q_1}
\left[%
\begin{array}{c}
  R\{-1\} \\
  R\{3-2n\} \\
\end{array}%
\right],
\]
where
\[
Q_0=\left(%
\begin{array}{cc}
  u_{ijkl} & x_ix_j-x_kx_l \\
  v_{ijkl} & -x_i-x_j+x_k+x_l \\
\end{array}%
\right), \hspace{.4cm}
Q_1=\left(%
\begin{array}{cc}
  x_i+x_j-x_k-x_l & x_ix_j-x_kx_l \\
  v_{ijkl} & -u_{ijkl} \\
\end{array}%
\right).
\]

Gornik also defined $\xi_0:C_p(\Gamma_0)\rightarrow C_p(\Gamma_1)$ by
the matrices
\[
U_0=\left(%
\begin{array}{cc}
  x_l-x_j & 0 \\
  a_1 & 1 \\
\end{array}%
\right), \hspace{1cm} U_1=\left(%
\begin{array}{cc}
  x_l & -x_j \\
  -1 & 1 \\
\end{array}%
\right),
\]
and $\xi_1:C_p(\Gamma_1)\rightarrow C_p(\Gamma_0)$ by the
matrices
\[
V_0=\left(%
\begin{array}{cc}
  1 & 0 \\
  -a_1 & x_l-x_j \\
\end{array}%
\right), \hspace{1cm} V_1=\left(%
\begin{array}{cc}
  1 & x_j \\
  1 & x_l \\
\end{array}%
\right),
\]
where $a_1=-v_{ijkl}+(u_{ijkl}+x_iv_{ijkl}-\pi_{jk})/(x_i-x_l)$.  These are
homomorphisms of matrix factorizations of degree $1$, i.e., these
commute with the differential maps and raise the filtration by $1$. One can check that
\[
\xi_1\xi_0=(x_l-x_j)\id_{C_p(\Gamma_0)},
\hspace{1cm}
\xi_0\xi_1=(x_l-x_j)\id_{C_p(\Gamma_1)}.
\]

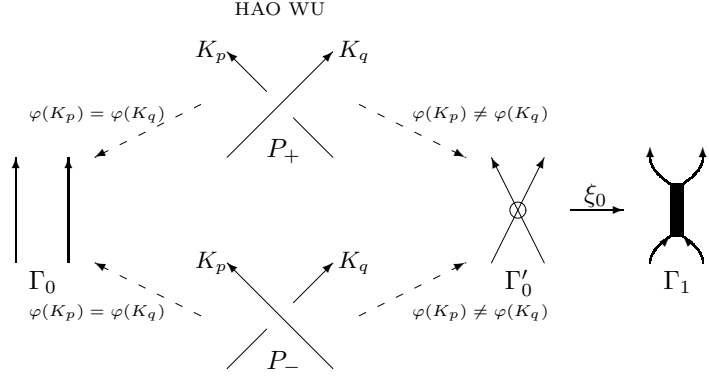
\begin{figure}[ht]

\setlength{\unitlength}{1pt}

\begin{picture}(420,120)(-210,-30)


\put(-80,10){\vector(0,1){40}}

\put(-100,10){\vector(0,1){40}}

\put(-95,0){$\Gamma_0$}


\put(22,8){\small{$K_q$}}

\put(-32,8){\small{$K_p$}}

\put(22,88){\small{$K_q$}}

\put(-32,88){\small{$K_p$}}

\put(-20,-30){\line(1,1){15}}

\put(5,-5){\vector(1,1){15}}

\put(20,-30){\vector(-1,1){40}}

\put(-5,-30){$P_-$}

\put(20,50){\line(-1,1){15}}

\put(-5,75){\vector(-1,1){15}}

\put(-20,50){\vector(1,1){40}}

\put(-5,50){$P_+$}

\multiput(30,70)(7,-3.5){6}{\line(2,-1){2.5}}

\put(70,50){\vector(2,-1){0}}

\put(50,65){\tiny{$\varphi(K_p)\neq\varphi(K_q)$}}

\multiput(-30,70)(-7,-3.5){6}{\line(-2,-1){2.5}}

\put(-70,50){\vector(-2,-1){0}}

\put(-95,65){\tiny{$\varphi(K_p)=\varphi(K_q)$}}

\multiput(30,-10)(7,3.5){6}{\line(2,1){2.5}}

\put(70,10){\vector(2,1){0}}

\put(50,-10){\tiny{$\varphi(K_p)\neq\varphi(K_q)$}}

\multiput(-30,-10)(-7,3.5){6}{\line(-2,1){2.5}}

\put(-70,10){\vector(-2,1){0}}

\put(-95,-10){\tiny{$\varphi(K_p)=\varphi(K_q)$}}


\put(80,10){\vector(1,2){20}}

\put(100,10){\vector(-1,2){20}}

\put(90,30){\circle{6}}

\put(85,0){$\Gamma_0'$}

\put(115,33){$\xi_0$}

\put(110,30){\vector(1,0){20}}

\put(147.5,20){\vector(1,1){0}}

\qbezier(140,10)(140,15)(147.5,20)

\qbezier(140,50)(140,45)(147.5,40)

\put(140,54){\vector(0,1){0}}

\put(152.5,20){\vector(-1,1){0}}

\qbezier(160,10)(160,15)(152.5,20)

\qbezier(160,50)(160,45)(152.5,40)

\put(160,54){\vector(0,1){0}}

\linethickness{5pt}

\put(150,20){\vector(0,1){20}}

\put(145,0){$\Gamma_1$}

\end{picture}

\caption{Resolutions of a crossing}\label{gornik-resolutions}

\end{figure}

Consider a link diagram $D$. Let $\{K_1,\cdots,K_m\}$ be the set of components of the link represented by $D$. A state of $D$ is a function 
\[
\varphi: \{K_1,\cdots,K_m\} \rightarrow \{0,1,\cdots, n-1\}.
\]
Gornik explicitly constructed a basis for $H_p(D)$ in \cite{Gornik}. Given a state $\varphi$ of $D$, one resolves every crossing in $D$ into either pair of parallel arcs or a virtual crossing by the rules given in Figure \ref{gornik-resolutions}. This gives a planar diagram $D_\varphi^0$ with only virtual crossings. Every component $K_j^0$ of $D_\varphi^0$ is a circle with only virtual crossings, whose $H_p$ cohomology is clearly isomorphic to that of a round circle, i.e. $H_p(K_i^0)\cong \C[x_i]/(x_i^n-1)$. So 
\[
H_p(D_\varphi^0)\cong H_p(K_1^0)\otimes_\C\cdots\otimes_\C H_p(K_l^0) \cong \C[x_1,\cdots,x_l]/(x_1^n-1,\cdots,x_l^n-1),
\] 
where $l$ is the number of components of $D_\varphi^0$. The state $\varphi$ of $D$ induces a state $\varphi^0$ of $D_\varphi^0$. Define $f_\varphi^0=f_{\varphi^0(K_1^0)}(x_1)\cdots f_{\varphi^0(K_l^0)}(x_l) \in H_p(D_\varphi^0)$. Now change each virtual crossing in $D_\varphi^0$ into a wide edge. This gives a graph $\Gamma_\varphi$, which is a complete resolution of $D$. Also, $\varphi$ induces an admissible state $\hat{\varphi}$ of $\Gamma_\varphi$. Clearly, the $\xi_0$ maps associated to these virtual crossings induce a homomorphism 
\[
F_\varphi: H_p(D_\varphi^0) \rightarrow H_p(\Gamma_\varphi) ~(\subset H(C_p(D),d_{mf})).
\] 
Define $f_\varphi = F_\varphi (f_\varphi^0)$. It is straightforward to check that $f_\varphi \in Q_{\hat{\varphi}} H_p(\Gamma_\varphi)$.

\begin{theorem}[\cite{Gornik}]\label{gornik-basis}
For every state $\varphi$ of $D$, $f_\varphi$ is a cocycle in $(H(C_p(D),d_{mf}),d_\chi)$. Moreover, 
\[
\{[f_\varphi]~|~\varphi \text{ is a state of } D\}
\] 
is a $\C$-linear basis for $H_p(D)$. In particular, $\dim_\C H_p(D) = n^m$, where $m$ is the number of components of the link represented by $D$.
\end{theorem}

Next, we want to understand how elementary cobodisms, i.e. Reidemeister moves and Morse moves, affect this basis. The effects of Morse moves are easy to understand. The following proposition follows easily from Lemmas \ref{gornik-epsilon-iota} and \ref{gornik-circle-saddle}.

\begin{proposition}\label{gornik-morse}
Let $S$ be an elementary link cobodism from $D$ to $D'$, and $\Psi:H_p(D)\rightarrow H_p(D')$ the homomorphism induces by $S$.
\begin{enumerate}
	\item If $S$ corresponds to a circle creation, then any state $\varphi$ of $D$ induces $n$ states $\varphi_j$, $j=0,1,\cdots,n-1$, of $D'$, where $\varphi_j$ agrees with $\varphi$ on components of $D$, and has value $j$ on the new component, and $\Psi([f_{\varphi}])=\lambda \sum_{j=0}^{n-1} [f_{\varphi_j}]$ for some $\lambda\in\C\setminus\{0\}$.
	      	
	\item If $S$ corresponds to a circle annihilation, then any state $\varphi$ of $D$ induces a state $\varphi'$ of $D'$ which agrees with $\varphi$ on components of $D'$, and $\Psi([f_{\varphi}])=\mu_\varphi [f_{\varphi'}]$ for some $\mu_\varphi\in\C\setminus\{0\}$.

  \item Assume that $S$ corresponds to a saddle move. If the values of a state $\varphi$ of $D$ on the two strands involved in the move are different, then $\Psi([f_{\varphi}])=0$. If the values of a state $\varphi$ of $D$ on the two strands involved in the move are equal, including the case when the two strands belong to the same component, then $\varphi$ induces a state $\varphi'$ of $D'$ which agrees with $\varphi$ on the unchanged components, and takes the common value on the changed component(s), and $\Psi([f_{\varphi}])=\nu_\varphi [f_{\varphi'}]$ for some $\nu_\varphi\in\C\setminus\{0\}$.
\end{enumerate}
\end{proposition}

\begin{proof}
Note that the homomorphism induced by a Morse move commutes with the $\xi_0$-homomorphisms used in the definition of $f_\varphi$. So we can do the Morse move before we change the virtual crossings into wide edges. Clearly, Lemmas \ref{gornik-epsilon-iota} and \ref{gornik-circle-saddle} are applicable to planar diagrams with only virtual crossings. And the proposition follows.
\end{proof}

It is harder to determine the effects of Reidemeister moves on the basis of $H_p$ given in Theorem \ref{gornik-basis}. We start by studying the admissible states of some special graphs.

\begin{figure}[ht]

\setlength{\unitlength}{1pt}

\begin{picture}(420,80)(-210,-25)

\qbezier(-10,25)(0,35)(10,25)

\put(10,25){\vector(1,-1){0}}

\qbezier(-10,15)(0,5)(10,15)

\put(-10,15){\vector(-1,1){0}}

\qbezier(-15,25)(0,60)(15,25)

\put(15,25){\vector(1,-2){0}}

\qbezier(-15,15)(0,-20)(15,15)

\put(-15,15){\vector(-1,2){0}}

\put(-2,45){\tiny{$x_1$}}

\put(-2,32){\tiny{$x_2$}}

\put(-2,-8){\tiny{$x_3$}}

\put(-2,4){\tiny{$x_4$}}

\put(0,41.5){\line(0,1){2}}

\put(0,29){\line(0,1){2}}

\put(0,9){\line(0,1){2}}

\put(0,-3.5){\line(0,1){2}}

\put(-3,-25){$\hat{\Gamma}$}

\linethickness{5pt}

\put(-12.5,15){\line(0,1){10}}

\put(12.5,15){\line(0,1){10}}

\end{picture}

\caption{$\hat{\Gamma}$}\label{gornik-hat-gamma}

\end{figure}
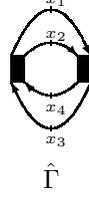

First we consider the admissible states of $\hat{\Gamma}$ in Figure \ref{gornik-hat-gamma}. For $k,l \in \{0,1,\cdots,n-1\}$ with $k\neq l$, define $\varphi_{kl}$ to be the state of $\hat{\Gamma}$ with $\varphi_{kl}(e_1)=\varphi_{kl}(e_3)= k$, $\varphi_{kl}(e_2)=\varphi_{kl}(e_4)= l$, and $\varphi_{kl}'$ to be the state with $\varphi_{kl}'(e_1)=\varphi_{kl}'(e_4)= k$, $\varphi_{kl}'(e_2)=\varphi_{kl}'(e_3)= l$. These states are all the admissible states of $\hat{\Gamma}$. From the proofs of Proposition \ref{moy2} and Lemma \ref{H-hat-gamma}, we have that 
\[
H_p(\hat{\Gamma})\cong \C[x_1,x_2,x_3,x_4]/(u_{1212},v_{1212},x_1+x_2-x_3-x_4,x_1x_2-x_3x_4).
\] 
Define 
\begin{eqnarray*}
f_{\varphi_{kl}} & = & f_k(x_1)f_l(x_2)f_k(x_3)f_l(x_4) \in \C[x_1,x_2,x_3,x_4]/(u_{1212},v_{1212},x_1+x_2-x_3-x_4,x_1x_2-x_3x_4) (\cong H_p(\hat{\Gamma})), \\
f_{\varphi_{kl}'} & = & f_k(x_1)f_l(x_2)f_l(x_3)f_k(x_4) \in \C[x_1,x_2,x_3,x_4]/(u_{1212},v_{1212},x_1+x_2-x_3-x_4,x_1x_2-x_3x_4) (\cong H_p(\hat{\Gamma})).
\end{eqnarray*}

\begin{lemma}\label{gornik-hat-gamma-comp}
$f_{\varphi_{kl}}$ is a none-zero element of $Q_{\varphi_{kl}} H_p(\hat{\Gamma})$, and $f_{\varphi_{kl}'}$ is a none-zero element of $Q_{\varphi_{kl}'} H_p(\hat{\Gamma})$. In particular, $\exists \lambda\in\C\setminus\{0\}$ such that $\widetilde{\varepsilon}(f_{\varphi_{kl}}) = -\lambda e^{\frac{2(k+l)\pi}{n}i}$ and $\widetilde{\varepsilon}(f_{\varphi_{kl}'}) = \lambda e^{\frac{2(k+l)\pi}{n}i}$ $\forall ~k,l \in \{0,1,\cdots,n-1\}$ with $k\neq l$, where the homomorphism $\widetilde{\varepsilon}: H_p(\hat{\Gamma}) \rightarrow \C$ is defined in Subsection \ref{subsection-moy3}.
\end{lemma}
\begin{proof}
By Lemma \ref{m(p')}, $x_j^n-1=0$ in $H_p(\hat{\Gamma})$ for $j=1,2,3,4$. From this, it is easy to check that $x_jf_{\varphi_{kl}} = e^{\frac{2\varphi_{kl}(e_j)\pi}{n}i} f_{\varphi_{kl}}$ and $x_jf_{\varphi_{kl}'} = e^{\frac{2\varphi_{kl}'(e_j)\pi}{n}i} f_{\varphi_{kl}'}$. So $f_{\varphi_{kl}}$ and $f_{\varphi_{kl}'}$ are elements of $Q_{\varphi_{kl}} H_p(\hat{\Gamma})$ and $Q_{\varphi_{kl}'} H_p(\hat{\Gamma})$. Since $x_1,~x_2,~x_3,~x_4$ satisfy $x_1+x_2=x_3+x_4$ and $x_1x_2=x_3x_4$, we have
\[
x_3^2-(x_1+x_2)x_3+x_1x_2=x_3^2-(x_3+x_4)x_3+x_3x_4 =(x_3-x_3)(x_3-x_4)=0.
\]
By induction, it's easy to show that 
\[
x_3^q = \frac{x_1^q-x_2^q}{x_1-x_2}x_3-x_1x_2\frac{x_1^{q-1}-x_2^{q-1}}{x_1-x_2}.
\]
Thus, for $k,l \in \{0,1,\cdots,n-1\}$ with $k\neq l$,
\begin{eqnarray*}
x_3^q f_k(x_1)f_l(x_2) & = & (\frac{x_1^q-x_2^q}{x_1-x_2}x_3-x_1x_2\frac{x_1^{q-1}-x_2^{q-1}}{x_1-x_2}) f_k(x_1)f_l(x_2) \\
& = & (\frac{e^{\frac{2qk\pi}{n}i}-e^{\frac{2ql\pi}{n}i}}{e^{\frac{2k\pi}{n}i}-e^{\frac{2l\pi}{n}i}}x_3-e^{\frac{2(k+l)\pi}{n}i}\frac{e^{\frac{2(q-1)k\pi}{n}i}-e^{\frac{2(q-1)l\pi}{n}i}}{e^{\frac{2k\pi}{n}i}-e^{\frac{2l\pi}{n}i}}) f_k(x_1)f_l(x_2).
\end{eqnarray*}
So, for $k,l \in \{0,1,\cdots,n-1\}$ with $k\neq l$,
\begin{eqnarray*}
f_{\varphi_{kl}} & = & f_k(x_1)f_l(x_2)f_k(x_3)f_l(x_4) \\
                 & = & (\sum_{q=0}^{n-1}e^{-\frac{2qk\pi}{n}i}x_3^q) f_k(x_1)f_l(x_2)f_l(x_4) \\
                 & = & (\sum_{q=0}^{n-1}e^{-\frac{2qk\pi}{n}i} (\frac{e^{\frac{2qk\pi}{n}i}-e^{\frac{2ql\pi}{n}i}}{e^{\frac{2k\pi}{n}i}-e^{\frac{2l\pi}{n}i}}x_3-e^{\frac{2(k+l)\pi}{n}i}\frac{e^{\frac{2(q-1)k\pi}{n}i}-e^{\frac{2(q-1)l\pi}{n}i}}{e^{\frac{2k\pi}{n}i}-e^{\frac{2l\pi}{n}i}})) f_k(x_1)f_l(x_2)f_l(x_4) \\
                 & = & (\sum_{q=0}^{n-1} (\frac{1-e^{\frac{2q(l-k)\pi}{n}i}}{e^{\frac{2k\pi}{n}i}-e^{\frac{2l\pi}{n}i}}x_3-e^{\frac{2l\pi}{n}i}\frac{1-e^{\frac{2(q-1)(l-k)\pi}{n}i}}{e^{\frac{2k\pi}{n}i}-e^{\frac{2l\pi}{n}i}})) f_k(x_1)f_l(x_2)f_l(x_4) \\
                 & = & \frac{n}{e^{\frac{2k\pi}{n}i}-e^{\frac{2l\pi}{n}i}}(x_3-e^{\frac{2l\pi}{n}i}) f_k(x_1)f_l(x_2)f_l(x_4) \\
                 & \cdots & (\text{repeat the above steps for } x_4) \\
                 & = & \frac{n}{e^{\frac{2k\pi}{n}i}-e^{\frac{2l\pi}{n}i}}(x_3-e^{\frac{2l\pi}{n}i})  \frac{n}{e^{\frac{2l\pi}{n}i}-e^{\frac{2k\pi}{n}i}}(x_4-e^{\frac{2k\pi}{n}i})  f_k(x_1)f_l(x_2) \\
                 & = & -\frac{n^2}{(e^{\frac{2k\pi}{n}i}-e^{\frac{2l\pi}{n}i})^2} (x_3-e^{\frac{2l\pi}{n}i})(x_1+x_2-x_3-e^{\frac{2k\pi}{n}i})  f_k(x_1)f_l(x_2) \\
                 & = & -\frac{n^2}{(e^{\frac{2k\pi}{n}i}-e^{\frac{2l\pi}{n}i})^2} (x_3-e^{\frac{2l\pi}{n}i}) (e^{\frac{2k\pi}{n}i}+e^{\frac{2l\pi}{n}i}-x_3-e^{\frac{2k\pi}{n}i})  f_k(x_1)f_l(x_2) \\
                 & = & \frac{n^2}{(e^{\frac{2k\pi}{n}i}-e^{\frac{2l\pi}{n}i})^2} (x_3-e^{\frac{2l\pi}{n}i})^2  f_k(x_1)f_l(x_2) \\
                 & = & \frac{n^2}{(e^{\frac{2k\pi}{n}i}-e^{\frac{2l\pi}{n}i})^2} (x_3^2-2e^{\frac{2l\pi}{n}i}x_3+e^{\frac{4l\pi}{n}i})  f_k(x_1)f_l(x_2) \\
                 & = & \frac{n^2}{(e^{\frac{2k\pi}{n}i}-e^{\frac{2l\pi}{n}i})^2} ((e^{\frac{2k\pi}{n}i}+e^{\frac{2l\pi}{n}i})x_3- e^{\frac{2(l+k)\pi}{n}i} -2e^{\frac{2l\pi}{n}i}x_3+e^{\frac{4l\pi}{n}i})  f_k(x_1)f_l(x_2) \\
                 & = & \frac{n^2}{e^{\frac{2k\pi}{n}i}-e^{\frac{2l\pi}{n}i}} (x_3- e^{\frac{2l\pi}{n}i}) f_k(x_1)f_l(x_2).
\end{eqnarray*}
Similarly, for $k,l \in \{0,1,\cdots,n-1\}$ with $k\neq l$,
\[
f_{\varphi_{kl}'} = -\frac{n^2}{e^{\frac{2k\pi}{n}i}-e^{\frac{2l\pi}{n}i}} (x_3- e^{\frac{2k\pi}{n}i}) f_k(x_1)f_l(x_2).
\]

According to Lemma \ref{H-hat-gamma}, $H_p(\hat{\Gamma})$, as a $\C$-linear space, has the basis 
\[
\{x_1^ix_2^jx_3^{\epsilon} ~|~ 0\leq i \leq n-1,~0\leq j\leq n-2,~\epsilon=0,1\},
\]
where the element $x_1^ix_2^jx_3^{\epsilon}$ has quantum degree $2i+2j+2\epsilon-2n+2$. By the definition of $\widetilde{\varepsilon}$, we have that $\widetilde{\varepsilon}(x_1^{n-1}x_2^{n-2}x_3)=\lambda'\in\C\setminus\{0\}$, and the value of $\widetilde{\varepsilon}$ on any other element in this basis is $0$. Since $p(x)=x^{n+1}-(n+1)x$, we have
\[
v_{1212}= - (n+1)\sum_{l=0}^{n-1}x_1^{n-1-l}x_2^l,
\]
which implies that, as elements of $H_p(\hat{\Gamma})$,
\[
x_2^{n-1}=-\sum_{l=0}^{n-2}x_1^{n-1-l}x_2^l.
\]
Using this and that $x_1^n=1$ in $H_p(\hat{\Gamma})$, one can check that, for any $k,l \in \{0,1,\cdots,n-1\}$ with $k\neq l$,
\[
f_{\varphi_{kl}} = \frac{n^2}{e^{\frac{2k\pi}{n}i}-e^{\frac{2l\pi}{n}i}} (x_3- e^{\frac{2l\pi}{n}i}) f_k(x_1)f_l(x_2) = -n^2 e^{\frac{2(k+l)\pi}{n}i} x_1^{n-1}x_2^{n-2}x_3 + \text{ lower degree terms}.
\]
Let $\lambda=\lambda'n^2 ~(\in\C\setminus\{0\})$. Then $\widetilde{\varepsilon}(f_{\varphi_{kl}}) = -\lambda e^{\frac{2(k+l)\pi}{n}i}$ and, similarly, $\widetilde{\varepsilon}(f_{\varphi_{kl}'}) = \lambda e^{\frac{2(k+l)\pi}{n}i}$ $\forall ~k,l \in \{0,1,\cdots,n-1\}$ with $k\neq l$. 
\end{proof}

To understand the effect of Reidemeister move III, we need to consider graphs with both double wide edges, i.e. the usual wide edges, and triple wide edges, i.e. wide edges with three regular edges entering in one end and three regular edges exiting at the other end. (See Figure \ref{3way}.) A state for such a graph is still a function from the set of all regular edges to $\{0,1,\cdots,n-1\}$. We say that a state $\varphi$ of such a graph $\Gamma$ is admissible if both of the following statements are true.

\begin{enumerate}
	\item At every double wide edge in $\Gamma$, $\varphi(e_1)\neq\varphi(e_2)$ and $\{\varphi(e_1),\varphi(e_2)\}=\{\varphi(e_3),\varphi(e_4)\}$ 
	      where $e_1$, $e_2$ are the two entering regular edges at this wide edge, and $e_3$, $e_4$ are the two exiting regular edges; 
	      	
	\item At every double wide edge in $\Gamma$, $\varphi(e_1),\varphi(e_2),\varphi(e_3)$ are pairwise distinct and
	      $\{\varphi(e_1),\varphi(e_2),\varphi(e_3)\} = \{\varphi(e_4),\varphi(e_5),\varphi(e_6)\}$ 
	      where $e_1$, $e_2$ and $e_3$ are the three entering regular edges at this wide edge, and $e_4$, $e_5$ and $e_6$ are the three exiting regular 
	      edges. 
\end{enumerate}

For a graph $\Gamma$ with double and triple wide edges, we mark each regular edge $e_j$ by a single variable $x_j$. For a state $\varphi$ of $\Gamma$, define 
\[
Q_\varphi = \prod_{j=1}^m f_{\varphi(e_j)}(x_j) \in \C[x_1,\cdots, x_m]/(x_1^n-1,\cdots,x_m^n-1).
\]
The following is a generalization of Theorem \ref{admissible-states} (Theorem 4 of \cite{Gornik}).

\begin{proposition}\label{general-admissible-states}
Let $\Gamma$ be a closed graph with double and triples wide edges. 

\begin{enumerate}
	\item For any state $\varphi$ of $\Gamma$ and $\alpha\in H_p(\Gamma)$,
\[
\alpha\in Q_\varphi H_p(\Gamma) ~ \Longleftrightarrow ~ x_j\alpha = e^{\frac{2\varphi(e_j)\pi}{n}i}\alpha, ~\forall~j.
\]
	\item If a state $\varphi$ is not admissible, then $Q_\varphi H_p(\Gamma)=0$. For each admissible state $\varphi$ of $\Gamma$, the $\C$-linear space $Q_\varphi H_p(\Gamma)$ is $1$-dimensional.
	\item \[
H_p(\Gamma) = \bigoplus_\varphi Q_\varphi H_p(\Gamma),
\]
where $\varphi$ runs through all admissible states of $\Gamma$. 
\end{enumerate}
\end{proposition}

\begin{figure}[ht]

\setlength{\unitlength}{1pt}

\begin{picture}(420,120)(-210,-20)


\linethickness{.25pt}

\put(-110,30){\vector(0,1){30}}

\put(-70,0){\vector(0,1){30}}

\put(-102.5,10){\vector(1,1){0}}

\qbezier(-110,0)(-110,5)(-102.5,10)

\qbezier(-110,30)(-110,25)(-102.5,20)

\put(-97.5,10){\vector(-1,1){0}}

\qbezier(-90,0)(-90,5)(-97.5,10)

\qbezier(-90,30)(-90,25)(-97.5,20)

\qbezier(-90,30)(-90,35)(-82.5,40)

\qbezier(-90,60)(-90,55)(-82.5,50)

\put(-90,60){\vector(-1,4){0}}

\qbezier(-70,30)(-70,35)(-77.5,40)

\qbezier(-70,60)(-70,55)(-77.5,50)

\qbezier(-110,60)(-110,65)(-102.5,70)

\qbezier(-110,90)(-110,85)(-102.5,80)

\put(-110,90){\vector(-1,4){0}}

\qbezier(-90,60)(-90,65)(-97.5,70)

\qbezier(-90,90)(-90,85)(-97.5,80)

\put(-90,90){\vector(1,4){0}}

\put(-90,30){\vector(1,4){0}}

\put(-70,60){\vector(0,1){30}}

\put(-90,-20){$\Gamma'$}

\put(-110,92){$x_6$}

\put(-90,92){$x_5$}

\put(-70,92){$x_4$}

\put(-110,-8){$x_3$}

\put(-90,-8){$x_2$}

\put(-70,-8){$x_1$}

\put(-108,42){$x_8$}

\put(-111,45){\line(1,0){2}}

\put(-88,57){$x_9$}

\put(-91,60){\line(1,0){2}}

\put(-88,27){$x_7$}

\put(-91,30){\line(1,0){2}}

\linethickness{5pt}

\put(-100,10){\line(0,1){10}}

\put(-80,40){\line(0,1){10}}

\put(-100,70){\line(0,1){10}}


\linethickness{.25pt}

\qbezier(87.5,50)(70,60)(70,70)

\put(70,70){\vector(0,1){20}}

\put(90,50){\vector(0,1){40}}

\qbezier(92.5,50)(110,60)(110,70)

\put(110,70){\vector(0,1){20}}

\qbezier(87.5,40)(70,30)(70,20)

\put(70,0){\vector(0,1){20}}

\put(90,0){\vector(0,1){40}}

\qbezier(92.5,40)(110,30)(110,20)

\put(110,0){\vector(0,1){20}}

\put(90,-20){$\Gamma$}

\put(70,92){$x_6$}

\put(90,92){$x_5$}

\put(110,92){$x_4$}

\put(70,-8){$x_3$}

\put(90,-8){$x_2$}

\put(110,-8){$x_1$}

\linethickness{5pt}

\put(90,40){\line(0,1){10}}


\linethickness{.25pt}

\put(210,0){\vector(0,1){90}}

\put(170,0){\vector(0,1){30}}

\put(190,0){\vector(0,1){30}}

\qbezier(170,30)(170,35)(177.5,40)

\qbezier(170,60)(170,55)(177.5,50)

\qbezier(190,30)(190,35)(182.5,40)

\qbezier(190,60)(190,55)(182.5,50)

\put(170,60){\vector(0,1){30}}

\put(190,60){\vector(0,1){30}}

\put(170,-20){$\Gamma''$}

\put(170,92){$x_6$}

\put(190,92){$x_5$}

\put(210,92){$x_4$}

\put(170,-8){$x_3$}

\put(190,-8){$x_2$}

\put(210,-8){$x_1$}

\linethickness{5pt}

\put(180,40){\line(0,1){10}}

\end{picture}

\caption{The replacement of a triple wide edge}\label{replacing-triple}

\end{figure}
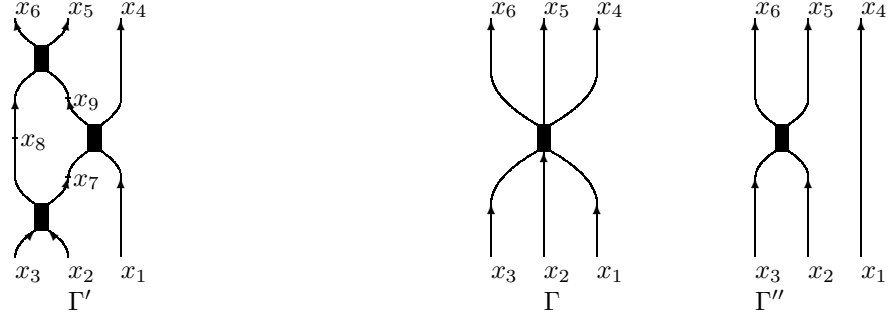

\begin{proof}
The proof of (1) is simple and identical to that of the corresponding statement in Theorem 4 of \cite{Gornik}. And, similar to Gornik's proof of that theorem, it is easy to see that (3) follows from (1) and (2). So we only need to prove (2). We do this by inducting on the number of triple wide edges in $\Gamma$. If $\Gamma$ has no triple wide edges, then the statement is true from Gornik's theorem. Assume the statement is true for graphs with $m-1$ triple edges, and $\Gamma$ has $m$ triple edges. Let $\varphi$ be a state of $\Gamma$. Let $E$ be a triple edge of $\Gamma$. Replacing $E$ by the three-double-wide-edge configuration in Figure \ref{replacing-triple}, we get a graph $\Gamma'$ with $m-1$ triple wide edges. By Proposition \ref{moy3}, there is a decomposition of $R$-modules
\[
H_p(\Gamma') = H_p(\Gamma) \oplus H_p(\Gamma''),
\]
where $R$ is the polynomial ring generated by all the markings on $\Gamma$. Let $J: H_p(\Gamma) \rightarrow H_p(\Gamma')$ be the inclusion from this decomposition, and $P: H_p(\Gamma') \rightarrow H_p(\Gamma'')$ the projection from this decomposition. Then $J$ is injective. Let $\varphi_{jkl}$ be the state of $\Gamma'$ that agrees with $\varphi$ on the edges inherited from $\Gamma$ and $\varphi_{jkl}(e_7)=j$, $\varphi_{jkl}(e_8)=k$, $\varphi_{jkl}(e_9)=l$. Using (1), it is easy to check that
\[
J(Q_\varphi H_p(\Gamma)) \subset \bigoplus_{j,k,l} Q_{\varphi_{jkl}} H_p(\Gamma').
\]

If $\varphi(e_1),\varphi(e_2),\varphi(e_3)$ are not pairwise distinct, then we can assume that $\varphi(e_2)=\varphi(e_3)$. This can be achieved by creating virtual crossings just outside the visible part of $\Gamma$ in Figure \ref{replacing-triple} before replacing the triple wide edge with the three double wide edges, which do not affect our argument. It is easy to see that none of $\varphi_{jkl}$ is admissible. So, $Q_\varphi H_p(\Gamma)=0$. Similarly, if $\varphi(e_4),\varphi(e_5),\varphi(e_6)$ are not pairwise distinct, then $Q_\varphi H_p(\Gamma)=0$.

Assume $\varphi_{jkl}$ is admissible for some $j,k,l$. Then 
\begin{eqnarray*}
\{j,k\} & =& \{\varphi(e_2),\varphi(e_3)\}, \\ 
l\in \{j,\varphi(e_1)\} & \subset & \{\varphi(e_1),\varphi(e_2),\varphi(e_3)\}.
\end{eqnarray*} 
Thus, 
\[
\{\varphi(e_4),\varphi(e_5),\varphi(e_6)\} \subset \{j,k,l,\varphi(e_1)\}\subset\{\varphi(e_1),\varphi(e_2),\varphi(e_3)\}.
\] 
Similarly, 
\[
\{\varphi(e_1),\varphi(e_2),\varphi(e_3)\} \subset \{\varphi(e_4),\varphi(e_5),\varphi(e_6)\}.
\]
This implies that $Q_\varphi H_p(\Gamma)=0$ if $\{\varphi(e_1),\varphi(e_2),\varphi(e_3)\} \neq \{\varphi(e_4),\varphi(e_5),\varphi(e_6)\}$. Now we can conclude that $Q_\varphi H_p(\Gamma)=0$ if $\varphi$ is not admissible.

Next, assume $\varphi$ is admissible. After possibly introducing some virtual crossings before replacing the triple wide edge with the three double wide edges, which does not affect our result, we can assume that $\varphi(e_1)=\varphi(e_6)$, $\varphi(e_2)=\varphi(e_5)$, $\varphi(e_3)=\varphi(e_4)$. Then $\varphi_{jkl}$ is admissible if and only if $l=\varphi(e_1)=\varphi(e_6)$, $k=\varphi(e_2)=\varphi(e_5)$, $j=\varphi(e_3)=\varphi(e_4)$. So, $J(Q_\varphi H_p(\Gamma)) \subset Q_{\varphi_{jkl}} H_p(\Gamma')$, where $j,k,l$ are as above. From the proof of Proposition \ref{reidemeister3}, one can see that $P$ factors through the $\chi_1$-map associated to the right wide edge in $\Gamma'$. By Gornik's computation, we have that $P(Q_{\varphi_{jkl}} H_p(\Gamma'))=0$. So $Q_{\varphi_{jkl}} H_p(\Gamma') \subset J(H_p(\Gamma))$. But $H_p(\Gamma) = \sum_{\psi}Q_{\psi} H_p(\Gamma)$, where $\psi$ runs through all states of $\Gamma$. From the above argument, we know that, if $\psi\neq\varphi$, then $Q_{\psi} H_p(\Gamma)$ is either $0$, or mapped by $J$ into some $Q_{\psi'} H_p(\Gamma')$, where $\psi'$ is a state of $\Gamma'$ not equal to $\varphi_{jkl}$. This implies that $Q_{\varphi_{jkl}} H_p(\Gamma') \subset J(Q_\varphi H_p(\Gamma))$. So $J(Q_\varphi H_p(\Gamma)) = Q_{\varphi_{jkl}} H_p(\Gamma')$. But $J$ is injective. This means $Q_\varphi H_p(\Gamma) \cong Q_{\varphi_{jkl}} H_p(\Gamma')$, which is $1$-dimensional by induction hypothesis.
\end{proof}

With Lemma \ref{gornik-hat-gamma-comp} and Proposition \ref{general-admissible-states} in hand, we are ready to determine the effects of Reidemeister moves.

\begin{proposition}\label{gornik-reidemeister}
If $D$ and $D'$ are two link diagrams that differ by a Reidemeister move, then there is a canonical 1-1 correspondence between the states of $D$ and states of $D'$ such that, if $\varphi$ and $\varphi'$ are states of $D$ and $D'$ corresponding to each other, and $\Psi:H_p(D)\rightarrow H_p(D')$ is the isomorphism induced by the Reidemeister move, then $\exists~ \lambda_\varphi \in \C\setminus \{0\}$ with $\Psi([f_\varphi])=\lambda_\varphi[f_{\varphi'}]$.
\end{proposition}
\begin{proof}
Since $D$ and $D'$ differ by a Reidemeister move, there is natural 1-1 correspondence between the components of $D$ and $D'$. This induces the 1-1 correspondence between the states of $D$ and $D'$, i.e. corresponding states have equal values on corresponding components. To show that this correspondence is compatible with the isomorphism induced by the Reidemeister move, we need to discuss the type of the Reidemeister move involved.

Reidemeister Move I. (i) Assume $D'$ is $D$ plus a positive kink. Then $\Psi$ is induced by the composition of the two homomorphisms in Figure \ref{gornik-Rmove1a}, where the homomorphism $J$ is the inclusion of $H(C_p(D_0),d_{mf})$ into $H(C_p(D'),d_{mf})$. Let $K$ be the component of $D$ containing the visible strand. Write $k=\varphi(K)$. For $j=0,1\cdots,n-1$, let $\varphi_j$ be the state of $D_0$ that agrees with $\varphi$ on components of $D$, and has value $j$ on the new circle. Then $f_{\varphi'}=J(f_{\varphi_k})$. By Lemma \ref{gornik-epsilon-iota}, there exits $\lambda\neq0$ such that
\[
\iota(f_\varphi)= \lambda f_\varphi \cdot (\sum_{j=0}^{n-1} f_j(x)) = \lambda\sum_{j=0}^{n-1} f_{\varphi_j}. 
\]
From \cite{Gornik}, $J_\ast([f_{\varphi_j}])=0$ if $j\neq k$. (See, e.g., Figure 6 of \cite{Gornik}.) Also $\Psi([f_\varphi]) = J_\ast \circ \iota_\ast ([f_\varphi]) = \lambda[f_{\varphi'}]$.

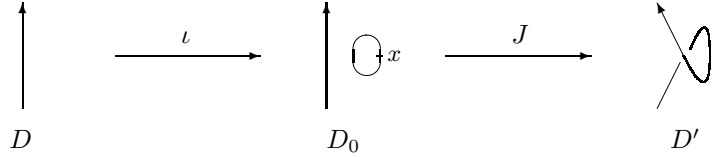
\begin{figure}[ht]

\setlength{\unitlength}{1pt}

\begin{picture}(420,60)(-210,-20)


\put(-125,0){\vector(0,1){40}}

\put(-130,-15){$D$}

\put(-90,20){\vector(1,0){55}}

\put(-65,25){$\iota$}


\put(-10,0){\vector(0,1){40}}

\put(5,20){\oval(10,15)}

\put(9,20){\line(1,0){2}}

\put(13,18){\small{$x$}}

\put(-10,-15){$D_0$}

\put(35,20){\vector(1,0){55}}

\put(60,25){$J$}


\put(125,20){\vector(-1,2){10}}

\qbezier(125,20)(135,0)(135,20)

\qbezier(135,20)(135,40)(127.5,22.5)

\put(115,0){\line(1,2){8.75}}

\put(120,-15){$D'$}

\end{picture}

\caption{Positive kink}\label{gornik-Rmove1a}

\end{figure}

(ii) Assume $D'$ is $D$ plus a negative kink. Then $\Psi^{-1}$ is induced by the composition of the homomorphisms in Figure \ref{gornik-Rmove1b}, where $P$ is the projection from $H(C_p(D'),d_{mf})$ onto $H(C_p(D_0),d_{mf})$. Define $k$ and $\varphi_j$ as above. Then $P(f_{\varphi'})= f_{\varphi_k}$. And, by Lemma \ref{gornik-epsilon-iota}, $\varepsilon(f_{\varphi_k})= \mu' f_\varphi$ for some $\mu'\in\C\setminus\{0\}$. Thus, $\Psi^{-1}([f_{\varphi'}]) = \mu' [f_\varphi]$. 

\begin{figure}[ht]

\setlength{\unitlength}{1pt}

\begin{picture}(420,60)(-210,-20)


\put(-125,20){\vector(1,2){10}}

\qbezier(-125,20)(-135,0)(-135,20)

\qbezier(-135,20)(-135,40)(-127.5,22.5)

\put(-115,0){\line(-1,2){8.75}}

\put(-130,-15){$D'$}

\put(-90,20){\vector(1,0){55}}

\put(-65,25){$P$}


\put(-10,0){\vector(0,1){40}}

\put(5,20){\oval(10,15)}

\put(9,20){\line(1,0){2}}

\put(13,18){\small{$x$}}

\put(-10,-15){$D_0$}

\put(35,20){\vector(1,0){55}}

\put(60,25){$\varepsilon$}


\put(125,0){\vector(0,1){40}}

\put(120,-15){$D$}

\end{picture}

\caption{Negative kink}\label{gornik-Rmove1b}

\end{figure}
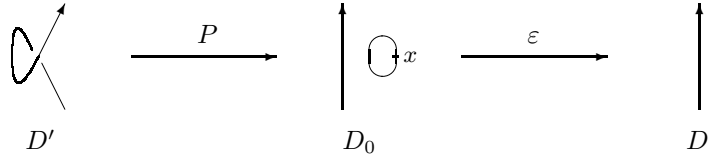

\begin{figure}[ht]

\setlength{\unitlength}{1pt}

\begin{picture}(420,160)(-210,-10)

\linethickness{.25pt}


\qbezier(-185,50)(-145,70)(-185,90)

\put(-185,90){\vector(-2,1){0}}

\put(-165,50){\line(-2,1){8}}

\qbezier(-177,56)(-195,70)(-177,84)

\put(-173,86){\vector(2,1){8}}

\put(-180,35){$D'$}


\put(-70,80){\vector(0,1){60}}

\put(-50,80){\vector(0,1){60}}

\put(-120,105){$D=\Gamma_{00}$}

\put(-80,140){\tiny{$x_1$}}

\put(-80,80){\tiny{$x_4$}}

\put(-47,140){\tiny{$x_2$}}

\put(-47,80){\tiny{$x_3$}}


\linethickness{.25pt}

\put(62.5,120){\vector(-3,2){0}}

\qbezier(70,80)(70,120)(62.5,120)

\put(57.5,120){\vector(3,2){0}}

\qbezier(50,80)(50,120)(57.5,120)

\qbezier(50,140)(50,135)(57.5,130)

\qbezier(70,140)(70,135)(62.5,130)

\put(70,140){\vector(1,4){0}}

\put(50,140){\vector(-1,4){0}}

\put(20,105){$\Gamma_{10}$}

\put(40,140){\tiny{$x_1$}}

\put(40,80){\tiny{$x_4$}}

\put(73,140){\tiny{$x_2$}}

\put(73,80){\tiny{$x_3$}}

\linethickness{5pt}

\put(60,120){\line(0,1){10}}


\linethickness{.25pt}

\put(62.5,10){\vector(-1,1){0}}

\qbezier(70,0)(70,5)(62.5,10)

\qbezier(70,30)(70,25)(62.5,20)

\put(70,30){\vector(1,4){0}}

\put(50,30){\vector(-1,4){0}}

\put(57.5,10){\vector(1,1){0}}

\qbezier(50,0)(50,5)(57.5,10)

\qbezier(50,30)(50,25)(57.5,20)

\qbezier(70,30)(70,35)(62.5,40)

\qbezier(70,60)(70,55)(62.5,50)

\put(70,60){\vector(1,4){0}}

\put(50,60){\vector(-1,4){0}}

\qbezier(50,30)(50,35)(57.5,40)

\qbezier(50,60)(50,55)(57.5,50)

\put(20,25){$\Gamma_{11}$}

\put(40,60){\tiny{$x_1$}}

\put(40,0){\tiny{$x_4$}}

\put(73,60){\tiny{$x_2$}}

\put(73,0){\tiny{$x_3$}}

\put(40,30){\tiny{$x_6$}}

\put(73,30){\tiny{$x_5$}}

\put(48.5,30){-}

\put(68.5,30){-}

\linethickness{5pt}

\put(60,10){\line(0,1){10}}

\put(60,40){\line(0,1){10}}


\linethickness{.25pt}

\put(-62.5,10){\vector(1,1){0}}

\qbezier(-70,0)(-70,5)(-62.5,10)

\put(-57.5,10){\vector(-1,1){0}}

\qbezier(-50,0)(-50,5)(-57.5,10)

\qbezier(-70,60)(-70,20)(-62.5,20)

\put(-70,60){\vector(0,0){0}}

\put(-50,60){\vector(0,1){0}}

\qbezier(-50,60)(-50,20)(-57.5,20)

\put(-100,25){$\Gamma_{01}$}

\put(-80,60){\tiny{$x_1$}}

\put(-80,0){\tiny{$x_4$}}

\put(-47,60){\tiny{$x_2$}}

\put(-47,0){\tiny{$x_3$}}

\linethickness{5pt}

\put(-60,10){\line(0,1){10}}

\end{picture}

\caption{Reidemeister move II$_a$}\label{gornik-Rmove2a}
\end{figure}

Reidemeister move II$_a$. Assume $D'$ is $D$ with a pair of canceling crossings created by a Reidemeister II$_a$ move. Recall that $\Psi^{-1}$ is induced by Gaussian elimination which eliminates the chain complexes associated to graphs $\Gamma_{01}$ $\Gamma_{10}$ and $\Gamma_{11}$ in Figure \ref{gornik-Rmove2a}. Let the components of $D$ containing the two visible strands be $K_1$ and $K_2$, and $K_1'$, $K_2'$ the corresponding components of $D'$. Consider a state $\varphi$ of $D$. Write $k_1=\varphi(K_1)$ and $k_2=\varphi(K_2)$. Then $k_1=\varphi'(K_1')$ and $k_2=\varphi'(K_2')$. If $k_1=k_2$ and $f_{\varphi'}$ is an element of $H(C_p(\Gamma_{00}),d_{mf})$. By the definition of $\Psi^{-1}$, it's clear that $\Psi^{-1}([f_{\varphi'}])=f_{\varphi}$. 

Next, we consider the case $k_1\neq k_2$. In this case, $f_{\varphi'}$ is an element of $H(C_p(\Gamma_{11}),d_{mf})$. From Gornik's computation, see e.g. Figure 6 of \cite{Gornik}, we have $\chi_1(f_{\varphi})=0$, where $\chi_1$ is the $\chi_1$-map associated to the lower wide edge. By Lemma \ref{Rmove2a-lemma}, this implies that $f_{\varphi} \in H(C_p(\Gamma_{01}),d_{mf})\{1\}$ in the decomposition
\[
H(C_p(\Gamma_{11}),d_{mf}) \cong H(C_p(\Gamma_{10}),d_{mf})\{-1\} \oplus H(C_p(\Gamma_{01}),d_{mf})\{1\}.
\]
Let $\chi_0$ be the $\chi_0$-map associated to the upper wide edge. Then $\chi_0$ admits the decomposition depicted in Figure \ref{gornik-Rmove2a-twist}, where the vertical homomorphism is the obvious isomorphism. Let $\varphi_{01}$ be the admissible state of $\Gamma_{01}$ that takes values $k_1$ and $k_2$ in the part of $\Gamma_{01}$ depicted in Figure \ref{gornik-Rmove2a-twist} as shown, and agrees with $\varphi'$ everywhere else. From the decomposition of $\chi_0$, we have $P\circ\chi_0(f_{\varphi_{01}})=f_{\varphi'}$, where $f_{\varphi_{01}}$ is the element of $H(C_p(\Gamma_{01}),d_{mf})$ associated to $\varphi_{01}$, and $P$ is the project of $H(C_p(\Gamma_{11}),d_{mf})$ onto $H(C_p(\Gamma_{01}),d_{mf})\{1\}$ from the above decomposition of $H(C_p(\Gamma_{11}),d_{mf})$. Recall that $P\circ\chi_0$ is the isomorphism we used to eliminated the two copies of $H(C_p(\Gamma_{01}),d_{mf})\{1\}$ in $H(C_p(D'),d_{mf})$. Thus, in the process of Gaussian Elimination, $f_{\varphi'}$ is mapped to $\chi_1'(f_{\varphi_{01}})\in H(C_p(\Gamma_{00}),d_{mf})$, where $\chi_1'$ is the $\chi_1$-map associated to the lower wide edge. It is clear that $\chi_1'(f_{\varphi_{01}})$ is a non-zero multiple of $f_{\varphi}$. (See Figure 6 of \cite{Gornik}.) This implies that $\Psi^{-1}([f_{\varphi'}])=\lambda[f_{\varphi}]$ for some $\lambda\in\C\setminus\{0\}$.

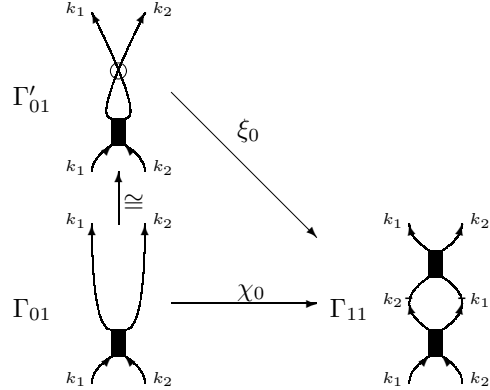
\begin{figure}[ht]

\setlength{\unitlength}{1pt}

\begin{picture}(420,160)(-210,-10)

\linethickness{.25pt}


\linethickness{.25pt}

\put(-40,110){\vector(1,-1){55}}

\put(-15,93){$\xi_0$}

\put(-62.5,90){\vector(1,1){0}}

\qbezier(-70,80)(-70,85)(-62.5,90)

\put(-57.5,90){\vector(-1,1){0}}

\qbezier(-50,80)(-50,85)(-57.5,90)

\qbezier(-50,140)(-70,100)(-62.5,100)

\put(-70,140){\vector(-1,2){0}}

\put(-50,140){\vector(1,2){0}}

\qbezier(-70,140)(-50,100)(-57.5,100)

\put(-60,118){\circle{6}}

\put(-100,105){$\Gamma_{01}'$}

\put(-80,140){\tiny{$k_1$}}

\put(-80,80){\tiny{$k_1$}}

\put(-47,140){\tiny{$k_2$}}

\put(-47,80){\tiny{$k_2$}}

\linethickness{5pt}

\put(-60,90){\line(0,1){10}}


\linethickness{.25pt}

\put(62.5,10){\vector(-1,1){0}}

\qbezier(70,0)(70,5)(62.5,10)

\qbezier(70,30)(70,25)(62.5,20)

\put(70,30){\vector(1,4){0}}

\put(50,30){\vector(-1,4){0}}

\put(57.5,10){\vector(1,1){0}}

\qbezier(50,0)(50,5)(57.5,10)

\qbezier(50,30)(50,25)(57.5,20)

\qbezier(70,30)(70,35)(62.5,40)

\qbezier(70,60)(70,55)(62.5,50)

\put(70,60){\vector(1,4){0}}

\put(50,60){\vector(-1,4){0}}

\qbezier(50,30)(50,35)(57.5,40)

\qbezier(50,60)(50,55)(57.5,50)

\put(20,25){$\Gamma_{11}$}

\put(40,60){\tiny{$k_1$}}

\put(40,0){\tiny{$k_1$}}

\put(73,60){\tiny{$k_2$}}

\put(73,0){\tiny{$k_2$}}

\put(40,30){\tiny{$k_2$}}

\put(73,30){\tiny{$k_1$}}

\put(48.5,30){-}

\put(68.5,30){-}

\linethickness{5pt}

\put(60,10){\line(0,1){10}}

\put(60,40){\line(0,1){10}}


\linethickness{.25pt}

\put(-60,60){\vector(0,1){20}}

\put(-58,65){$\cong$}

\put(-40,30){\vector(1,0){55}}

\put(-15,33){$\chi_0$}

\put(-62.5,10){\vector(1,1){0}}

\qbezier(-70,0)(-70,5)(-62.5,10)

\put(-57.5,10){\vector(-1,1){0}}

\qbezier(-50,0)(-50,5)(-57.5,10)

\qbezier(-70,60)(-70,20)(-62.5,20)

\put(-70,60){\vector(0,0){0}}

\put(-50,60){\vector(0,1){0}}

\qbezier(-50,60)(-50,20)(-57.5,20)

\put(-100,25){$\Gamma_{01}$}

\put(-80,60){\tiny{$k_1$}}

\put(-80,0){\tiny{$k_1$}}

\put(-47,60){\tiny{$k_2$}}

\put(-47,0){\tiny{$k_2$}}

\linethickness{5pt}

\put(-60,10){\line(0,1){10}}

\end{picture}

\caption{Decomposing the $\chi_0$-map}\label{gornik-Rmove2a-twist}
\end{figure}

\begin{figure}[ht]

\setlength{\unitlength}{1pt}

\begin{picture}(420,160)(-210,-10)

\linethickness{.25pt}


\linethickness{.25pt}

\put(-45,95){\vector(1,-1){15}}

\put(-30,120){\vector(-1,-1){15}}

\qbezier(-50,105)(-50,110)(-60,110)

\qbezier(-60,110)(-70,110)(-70,105)

\qbezier(-50,95)(-50,90)(-60,90)

\qbezier(-60,90)(-70,90)(-70,95)

\put(-70,95){\vector(0,1){10}}

\put(-90,120){\vector(-1,1){0}}

\qbezier(-90,80)(-70,100)(-90,120)

\put(-100,120){\tiny{$x_1$}}

\put(-25,120){\tiny{$x_2$}}

\put(-25,80){\tiny{$x_3$}}

\put(-62,114){\tiny{$x_5$}}

\put(-60,109){\line(0,1){2}}

\put(-60,89){\line(0,1){2}}

\put(-110,95){$\Gamma_{01}$}

\linethickness{5pt}

\put(-47.5,95){\line(0,1){10}}


\linethickness{.25pt}

\qbezier(50,105)(60,115)(70,105)

\put(70,105){\vector(1,-1){0}}

\qbezier(50,95)(60,85)(70,95)

\put(50,95){\vector(-1,1){0}}

\put(45,105){\vector(-1,1){15}}

\put(75,95){\vector(1,-1){15}}

\put(30,80){\vector(1,1){15}}

\put(90,120){\vector(-1,-1){15}}

\put(20,120){\tiny{$x_1$}}

\put(95,120){\tiny{$x_2$}}

\put(20,80){\tiny{$x_4$}}

\put(95,80){\tiny{$x_3$}}

\put(58,114){\tiny{$x_5$}}

\put(58,83){\tiny{$x_6$}}

\put(60,109){\line(0,1){2}}

\put(60,89){\line(0,1){2}}

\put(100,95){$\Gamma_{11}$}

\linethickness{5pt}

\put(72.5,95){\line(0,1){10}}

\put(47.5,95){\line(0,1){10}}


\linethickness{.25pt}

\put(45,25){\vector(-1,1){15}}

\put(30,0){\vector(1,1){15}}

\qbezier(50,25)(50,30)(60,30)

\qbezier(60,30)(70,30)(70,25)

\qbezier(50,15)(50,10)(60,10)

\qbezier(60,10)(70,10)(70,15)

\put(70,25){\vector(0,-1){10}}

\put(90,0){\vector(1,-1){0}}

\qbezier(90,00)(70,20)(90,40)

\put(20,40){\tiny{$x_1$}}

\put(95,40){\tiny{$x_2$}}

\put(20,0){\tiny{$x_4$}}

\put(58,34){\tiny{$x_5$}}

\put(60,29){\line(0,1){2}}

\put(60,9){\line(0,1){2}}

\put(100,15){$\Gamma_{10}$}

\linethickness{5pt}

\put(47.5,15){\line(0,1){10}}


\linethickness{.25pt}

\put(-30,0){\vector(1,-1){0}}

\qbezier(-30,40)(-50,20)(-30,0)

\put(-90,40){\vector(-1,1){0}}

\qbezier(-90,0)(-70,20)(-90,40)

\put(-60,20){\circle{20}}

\put(-100,40){\tiny{$x_1$}}

\put(-25,40){\tiny{$x_2$}}

\put(-62,34){\tiny{$x_5$}}

\put(-60,29){\line(0,1){2}}

\put(-60,9){\line(0,1){2}}

\put(-110,15){$\Gamma_{00}$}


\linethickness{.25pt}

\put(-150,0){\vector(1,-1){0}}

\qbezier(-150,40)(-170,20)(-150,0)

\put(-190,40){\vector(-1,1){0}}

\qbezier(-190,0)(-170,20)(-190,40)

\put(-200,40){\tiny{$x_1$}}

\put(-145,40){\tiny{$x_2$}}

\put(-210,15){$\Gamma$}


\linethickness{.25pt}

\put(-150,80){\vector(1,-1){0}}

\qbezier(-150,120)(-170,100)(-190,120)

\put(-190,120){\vector(-1,1){0}}

\qbezier(-190,80)(-170,100)(-150,80)

\put(-200,120){\tiny{$x_1$}}

\put(-145,80){\tiny{$x_3$}}

\put(-210,95){$D$}


\qbezier(180,110)(160,70)(140,110)

\put(140,110){\vector(-1,2){0}}

\put(174,96){\vector(1,-1){6}}

\qbezier(170,100)(160,110)(150,100)

\put(146,96){\line(-1,-1){6}}

\put(190,95){$D'$}

\put(135,115){\tiny{$x_1$}}

\put(180,85){\tiny{$x_3$}}

\end{picture}

\caption{Reidemeister move II$_b$}\label{gornik-Rmove2b-decomp}
\end{figure}
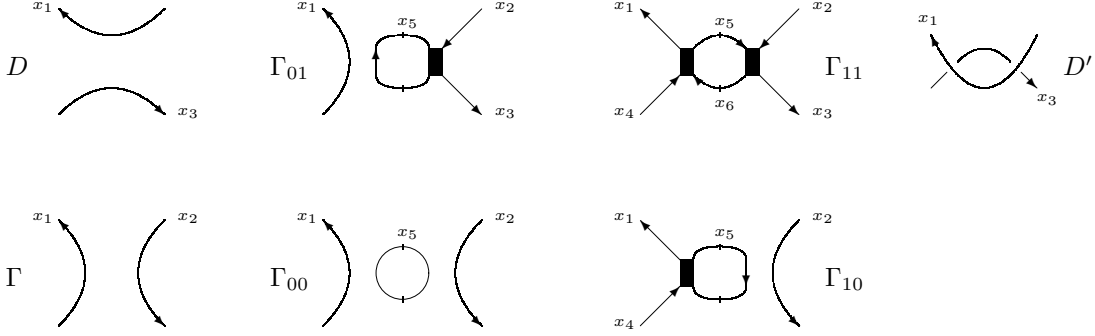

Reidemeister move II$_b$. Assume $D'$ is $D$ with a pair of canceling crossings created by a Reidemeister II$_b$ move. Then $\Psi^{-1}$ is induced by the Gaussian elimination that eliminates $\Gamma_{00}$, $\Gamma_{01}$, $\Gamma_{10}$ and the $\bigoplus_{i=0}^{n-3} H(C_p(\Gamma),d_{mf})\{3-n+2i\}$ component of $H(C_p(\Gamma_{11}),d_{mf})$ in the decomposition
\[
H(C_p(\Gamma_{11}),d_{mf}) \cong H(C_p(D),d_{mf}) \oplus (\bigoplus_{i=0}^{n-3} H(C_p(\Gamma),d_{mf})\{3-n+2i\}),
\]
where $\Gamma$, $\Gamma_{00}$, $\Gamma_{01}$, $\Gamma_{10}$ and $\Gamma_{11}$ are depicted in Figure \ref{gornik-Rmove2b-decomp}. Let the components of $D$ containing the two visible strands be $K_1$ and $K_3$, and $K_1'$, $K_3'$ the corresponding components of $D'$. Consider a state $\varphi$ of $D$. Write $k=\varphi(K_1)$ and $l=\varphi(K_3)$. Then $k=\varphi'(K_1')$ and $l=\varphi'(K_3')$. 

If $k\neq l$, then $\varphi'$ induces an admissible state $\hat{\varphi}'$ of $\Gamma_{11}$ such that $\hat{\varphi}'(e_1)=\hat{\varphi}'(e_2)=\hat{\varphi}'(e_6)=k$, $\hat{\varphi}'(e_3)=\hat{\varphi}'(e_4)=\hat{\varphi}'(e_5)=l$, where $e_j$ is the edge marked by $x_j$. The projection 
\[
G: H(C_p(\Gamma_{11}),d_{mf}) \rightarrow H(C_p(D),d_{mf})
\]
from the above decomposition of $H(C_p(\Gamma_{11}),d_{mf})$ is defined in Figure \ref{mapG}. It is the composition of two saddle moves followed by a $\widetilde{\varepsilon}$-map. By Lemmas \ref{gornik-circle-saddle} and \ref{gornik-hat-gamma-comp},
\[
G(f_{\varphi'})=\lambda' \widetilde{\varepsilon}(f_{\varphi'_{kl}}) f_{\varphi} = \lambda'' f_{\varphi},
\]
where $f_{\varphi'_{kl}}$ is the admissible state of $\hat{\Gamma}$ defined in Lemma \ref{gornik-hat-gamma-comp}, and $\lambda'$, $\lambda''$ are non-zero complex numbers. By its construction in Subsection \ref{subsection-moy3}, the projection 
\[
\widetilde{\beta} : H(C_p(\Gamma_{11}),d_{mf}) \rightarrow \bigoplus_{i=0}^{n-3} H(C_p(\Gamma),d_{mf})\{3-n+2i\}
\]
from the above decomposition of $H(C_p(\Gamma_{11}),d_{mf})$ factors through the $\chi_1$-maps of the two wide edges. From Gornik's computation (see Figure 6 of \cite{Gornik}), $\widetilde{\beta}(f_{\varphi'})=0$. This means that $f_{\varphi'}\in H(C_p(D),d_{mf})$ in the above decomposition of $H(C_p(\Gamma_{11}),d_{mf})$ and is a non-zero multiple of $f_{\varphi}$, which implies that $\Psi^{-1}([f_{\varphi'}])=\lambda''[f_{\varphi}]$.

If $k= l$, then $\varphi'$ induces an admissible state $\hat{\varphi}'$ of $\Gamma_{00}$ such that $\hat{\varphi}'(e_1) = \hat{\varphi}'(e_2) = \hat{\varphi}'(e_5) = k$, where $e_j$ is the edge marked by $x_j$. Recall that, in the first step of the Gaussian elimination, we eliminated $H(C_p(\Gamma_{10}),d_{mf})$ and the $\bigoplus_{i=0}^{n-2} H(C_p(\Gamma),d_{mf})\{1-n+2i\}$ component of $H(C_p(\Gamma_{00}),d_{mf})$ in the decomposition
\[
H(C_p(\Gamma_{00}),d_{mf}) = \bigoplus_{i=0}^{n-1} H(C_p(\Gamma),d_{mf})\{1-n+2i\}.
\]
So, after this reduction, $f_{\varphi'}$ is mapped to $\varepsilon(f_{\varphi'})$ in the remaining component $H(C_p(\Gamma),d_{mf})\{n-1\}= x_5^{n-1}H(C_p(\Gamma),d_{mf})\{1-n\}$ of $H(C_p(\Gamma_{00}),d_{mf})$. In the decomposition
\[
H(C_p(\Gamma_{01}),d_{mf})\{1\} = \bigoplus_{i=0}^{n-2} H(C_p(\Gamma),d_{mf})\{3-n+2i\} = \bigoplus_{i=0}^{n-2} x_5^{i}H(C_p(\Gamma),d_{mf})\{3-n\},
\]
the top component $x_5^{n-2}H(C_p(\Gamma),d_{mf})\{3-n\}$ is isomorphic to the remaining component $x_5^{n-1}H(C_p(\Gamma),d_{mf})\{1-n\}$ of $H(C_p(\Gamma_{00}),d_{mf})$ via the $\chi_1$-map of the right wide edge. This isomorphism is used to eliminate these two components. Let $\hat{\varphi}'_j$ be the state of $\Gamma_{01}$ that agrees with $\hat{\varphi}'$ everywhere except $\hat{\varphi}'_j(e_5)=j$. The $\hat{\varphi}'_j$ is admissible if and only if $j\neq k$. Use the fact that
\[
x_5^{n-2}=\frac{1}{n}\sum_{j=0}^{n-1} e^{-\frac{4j\pi}{n}i}f_j(x_5),
\]
it is easy to check that $\sum_{j\neq k} e^{-\frac{4j\pi}{n}i}f_{\hat{\varphi}'_j} \in x_5^{n-2}H(C_p(\Gamma),d_{mf})\{3-n\}$, and 
\[
\varepsilon(f_{\varphi'}) = \mu \chi_1(\sum_{j\neq k} e^{-\frac{4j\pi}{n}i}f_{\hat{\varphi}'_j}),
\]
where $\mu\in\C\setminus\{0\}$, and $j$ runs through $\{0,1,\cdots,n-1\}\setminus\{k\}$. So, after the second elimination, $f_{\varphi'}$ is mapped to $-\mu\chi_0(\sum_{j\neq k} e^{-\frac{4j\pi}{n}i}f_{\hat{\varphi}'_j}) \in H(C_p(\Gamma_{11}),d_{mf})$, where $\chi_0$ is the $\chi_0$-map associated to the left wide edge. Use the method in the proof of Lemma \ref{alpha-beta-F-G}, one can check that $\widetilde{\beta}\circ\chi_0|_{x_5^{n-2}H(C_p(\Gamma),d_{mf})\{3-n\}}=0$. Then, $\chi_0(\sum_{j\neq k} e^{-\frac{4j\pi}{n}i}f_{\hat{\varphi}'_j}) \in H(C_p(D),d_{mf})$ in the above decomposition of $H(C_p(\Gamma_{11}),d_{mf})$. More specifically, in this decomposition, $\chi_0(\sum_{j\neq k} e^{-\frac{4j\pi}{n}i}f_{\hat{\varphi}'_j})$ corresponds to the element
\[
G(\chi_0(\sum_{j\neq k} e^{-\frac{4j\pi}{n}i}f_{\hat{\varphi}'_j})) \in H(C_p(D),d_{mf}).
\]
Denote by $\eta_3,\eta_4$ the homomorphisms associated to the two saddle moves in the definition of $G$. Then
\begin{eqnarray*}
G(\chi_0(\sum_{j\neq k} e^{-\frac{4j\pi}{n}i}f_{\hat{\varphi}'_j}))  
& = & \widetilde{\varepsilon}\circ\eta_3\circ\eta_4\circ \chi_0(\sum_{j\neq k} e^{-\frac{4j\pi}{n}i}f_{\hat{\varphi}'_j}) \\
& = & \widetilde{\varepsilon}\circ\chi_0\circ\eta_3\circ\eta_4(\sum_{j\neq k} e^{-\frac{4j\pi}{n}i}f_{\hat{\varphi}'_j}) \\
& = & \mu'\widetilde{\varepsilon}\circ\chi_0(\sum_{j\neq k} e^{-\frac{4j\pi}{n}i} f_{kj}) f_{\varphi},
\end{eqnarray*}
where $\mu'\in\C\setminus\{0\}$, $f_{kj}=f_k(x_1)f_j(x_5) \in H_p(\hat{\Gamma}')$, and $\hat{\Gamma}'$ is in depicted in Figure \ref{gornik-hat-gamma2}. 

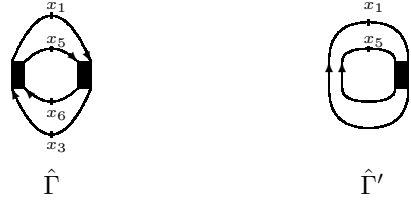
\begin{figure}[ht]

\setlength{\unitlength}{1pt}

\begin{picture}(420,80)(-210,-25)


\qbezier(-70,25)(-60,35)(-50,25)

\put(-50,25){\vector(1,-1){0}}

\qbezier(-70,15)(-60,5)(-50,15)

\put(-70,15){\vector(-1,1){0}}

\qbezier(-75,25)(-60,60)(-45,25)

\put(-45,25){\vector(1,-2){0}}

\qbezier(-75,15)(-60,-20)(-45,15)

\put(-75,15){\vector(-1,2){0}}

\put(-62,45){\tiny{$x_1$}}

\put(-62,32){\tiny{$x_5$}}

\put(-62,-8){\tiny{$x_3$}}

\put(-62,4){\tiny{$x_6$}}

\put(-60,41.5){\line(0,1){2}}

\put(-60,29){\line(0,1){2}}

\put(-60,9){\line(0,1){2}}

\put(-60,-3.5){\line(0,1){2}}

\put(-63,-25){$\hat{\Gamma}$}


\qbezier(45,25)(45,40)(60,40)

\qbezier(60,40)(75,40)(75,25)

\qbezier(50,25)(50,30)(60,30)

\qbezier(60,30)(70,30)(70,25)

\put(50,15){\vector(0,1){10}}

\put(45,15){\vector(0,1){10}}

\qbezier(50,15)(50,10)(60,10)

\qbezier(60,10)(70,10)(70,15)

\qbezier(45,15)(45,0)(60,0)

\qbezier(60,0)(75,0)(75,15)

\put(58,45){\tiny{$x_1$}}

\put(58,32){\tiny{$x_5$}}

\put(60,41){\line(0,-1){2}}

\put(60,29){\line(0,1){2}}

\put(57,-25){$\hat{\Gamma}'$}

\linethickness{5pt}

\put(-72.5,15){\line(0,1){10}}

\put(-47.5,15){\line(0,1){10}}

\put(72.5,15){\line(0,1){10}}

\end{picture}

\caption{$\hat{\Gamma}$ and $\hat{\Gamma}'$}\label{gornik-hat-gamma2}

\end{figure}

It is easy to see that
\begin{eqnarray*}
H_p(\hat{\Gamma}) & \cong & \C[x_1,x_3,x_5]/(u_{1515},v_{1515},(x_3-x_1)(x_3-x_5)), \\
H_p(\hat{\Gamma}') & \cong & \C[x_1,x_5]/(u_{1515},v_{1515}).
\end{eqnarray*}
By Remark \ref{direct-sum-variable-exclusion} and the definition of $\chi_0$, we know that $\chi_0=m(x_3-x_5):H_p(\hat{\Gamma}')\{1\}\rightarrow H_p(\hat{\Gamma})$. Thus, from Lemma \ref{gornik-hat-gamma-comp},
\begin{eqnarray*}
\chi_0(\sum_{j\neq k} e^{-\frac{4j\pi}{n}i} f_{kj}) & = & \sum_{j\neq k}e^{-\frac{4j\pi}{n}i}(x_3-x_5)f_k(x_1)f_j(x_5) \\
& = & \sum_{j\neq k}e^{-\frac{4j\pi}{n}i}(x_3-e^{\frac{2j\pi}{n}i})f_k(x_1)f_j(x_5) \\
& = & \frac{1}{n^2} \sum_{j\neq k} e^{-\frac{4j\pi}{n}i}(e^{\frac{2k\pi}{n}i}-e^{\frac{2j\pi}{n}i})f_{\varphi_{kj}},
\end{eqnarray*}
where $f_{\varphi_{kj}}$ is defined in Lemma \ref{gornik-hat-gamma-comp}. From that lemma, 
\begin{eqnarray*}
\widetilde{\varepsilon}\circ\chi_0(\sum_{j\neq k} e^{-\frac{4j\pi}{n}i} f_{kj}) & = & \widetilde{\varepsilon}(\frac{1}{n^2} \sum_{j\neq k} e^{-\frac{4j\pi}{n}i}(e^{\frac{2k\pi}{n}i}-e^{\frac{2j\pi}{n}i})f_{\varphi_{kj}}) \\
& = & -\frac{\lambda}{n^2} \sum_{j\neq k} e^{-\frac{4j\pi}{n}i}(e^{\frac{2k\pi}{n}i}-e^{\frac{2j\pi}{n}i}) e^{\frac{2(k+j)\pi}{n}i} \\
& = & -\frac{\lambda}{n^2} \sum_{j\neq k} (e^{\frac{2(2k-j)\pi}{n}i}-e^{\frac{2k\pi}{n}i}) \\
& = & \frac{\lambda}{n}\cdot e^{\frac{2k\pi}{n}i} \neq 0,
\end{eqnarray*}
where $\lambda$ is the non-zero scalar from Lemma \ref{gornik-hat-gamma-comp}.

Thus, $G(\chi_0(\sum_{j\neq k} e^{-\frac{4j\pi}{n}i}f_{\hat{\varphi}'_j}))$ is a non-zero multiple of $f_{\varphi}$, which implies that $\Psi^{-1}([f_{\varphi'}])=\lambda' [f_{\varphi}]$ for some $\lambda' \in \C\setminus\{0\}$.

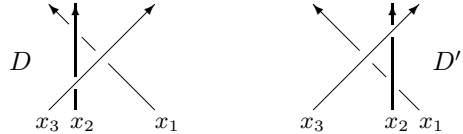
\begin{figure}[ht]

\setlength{\unitlength}{1pt}

\begin{picture}(420,52)(-210,-7)

\put(-70,0){\vector(1,1){40}}

\put(-60,0){\line(0,1){7.5}}

\put(-60,12.5){\vector(0,1){27.5}}

\put(-30,0){\line(-1,1){17.5}}

\put(-52.5,22.5){\line(-1,1){5}}

\put(-62.5,32.5){\vector(-1,1){7.5}}

\put(-62,-7){\small{$x_2$}}

\put(-75,-7){\small{$x_3$}}

\put(-30,-7){\small{$x_1$}}

\put(-85,15){$D$}

\put(30,0){\vector(1,1){40}}

\put(60,0){\line(0,1){27.5}}

\put(60,32.5){\vector(0,1){7.5}}

\put(70,0){\line(-1,1){7.5}}

\put(57.5,12.5){\line(-1,1){5}}

\put(47.5,22.5){\vector(-1,1){17.5}}

\put(75,15){$D'$}

\put(57,-7){\small{$x_2$}}

\put(70,-7){\small{$x_1$}}

\put(25,-7){\small{$x_3$}}

\end{picture}

\caption{Reidemeister move III}\label{gornik-Rmove3}

\end{figure}

Reidemeister move III. Assume $D$ and $D'$ deffer by a Reidemeister move III as depicted in Figure \ref{gornik-Rmove3}. Let $\varphi$ be a state of $D$, and $\varphi(K_1)=j$, $\varphi(K_2)=k$, $\varphi(K_3)=l$, where $K_j$ is the component of $D$ marked by $x_j$. Then $\varphi'$ is the state of $D'$ with $\varphi'(K_1')=j$, $\varphi'(K_2')=k$, $\varphi'(K_3')=l$, where $K_j'$ is the component of $D$ corresponding to $K_j$.

\begin{figure}[ht]

\setlength{\unitlength}{1pt}

\begin{picture}(420,240)(-210,-140)


\linethickness{.25pt}

\put(-210,30){\vector(0,1){30}}

\put(-170,0){\vector(0,1){30}}

\put(-202.5,10){\vector(1,1){0}}

\qbezier(-210,0)(-210,5)(-202.5,10)

\qbezier(-210,30)(-210,25)(-202.5,20)

\put(-197.5,10){\vector(-1,1){0}}

\qbezier(-190,0)(-190,5)(-197.5,10)

\qbezier(-190,30)(-190,25)(-197.5,20)

\qbezier(-190,30)(-190,35)(-182.5,40)

\qbezier(-190,60)(-190,55)(-182.5,50)

\put(-190,60){\vector(-1,4){0}}

\qbezier(-170,30)(-170,35)(-177.5,40)

\qbezier(-170,60)(-170,55)(-177.5,50)

\qbezier(-210,60)(-210,65)(-202.5,70)

\qbezier(-210,90)(-210,85)(-202.5,80)

\put(-210,90){\vector(-1,4){0}}

\qbezier(-190,60)(-190,65)(-197.5,70)

\qbezier(-190,90)(-190,85)(-197.5,80)

\put(-190,90){\vector(1,4){0}}

\put(-190,30){\vector(1,4){0}}

\put(-170,60){\vector(0,1){30}}

\put(-190,-20){$\Gamma_{111}$}

\put(-210,92){$x_6$}

\put(-190,92){$x_5$}

\put(-170,92){$x_4$}

\put(-210,-8){$x_3$}

\put(-190,-8){$x_2$}

\put(-170,-8){$x_1$}

\put(-208,42){$x_8$}

\put(-211,45){\line(1,0){2}}

\put(-188,57){$x_9$}

\put(-191,60){\line(1,0){2}}

\put(-188,27){$x_7$}

\put(-191,30){\line(1,0){2}}

\linethickness{5pt}

\put(-200,10){\line(0,1){10}}

\put(-180,40){\line(0,1){10}}

\put(-200,70){\line(0,1){10}}


\linethickness{.25pt}

\put(-90,0){\vector(0,1){90}}

\put(-50,0){\vector(0,1){30}}

\put(-70,0){\vector(0,1){30}}

\qbezier(-50,30)(-50,35)(-57.5,40)

\qbezier(-50,60)(-50,55)(-57.5,50)

\qbezier(-70,30)(-70,35)(-62.5,40)

\qbezier(-70,60)(-70,55)(-62.5,50)

\put(-50,60){\vector(0,1){30}}

\put(-70,60){\vector(0,1){30}}

\put(-70,-20){$\Gamma_{010}$}

\put(-90,92){$x_6$}

\put(-70,92){$x_5$}

\put(-50,92){$x_4$}

\put(-90,-8){$x_3$}

\put(-70,-8){$x_2$}

\put(-50,-8){$x_1$}

\linethickness{5pt}

\put(-60,40){\line(0,1){10}}


\linethickness{.25pt}

\put(90,0){\vector(0,1){90}}

\put(70,30){\vector(0,1){30}}

\put(62.5,10){\vector(-1,1){0}}

\qbezier(70,0)(70,5)(62.5,10)

\qbezier(70,30)(70,25)(62.5,20)

\put(57.5,10){\vector(1,1){0}}

\qbezier(50,0)(50,5)(57.5,10)

\qbezier(50,30)(50,25)(57.5,20)

\put(50,30){\vector(0,1){30}}

\qbezier(70,60)(70,65)(62.5,70)

\qbezier(70,90)(70,85)(62.5,80)

\put(70,90){\vector(1,4){0}}

\qbezier(50,60)(50,65)(57.5,70)

\qbezier(50,90)(50,85)(57.5,80)

\put(50,90){\vector(-1,4){0}}

\put(50,30){\vector(-1,4){0}}

\put(70,-20){$\Gamma_{101}$}

\put(50,92){$x_6$}

\put(70,92){$x_5$}

\put(90,92){$x_4$}

\put(50,-8){$x_3$}

\put(70,-8){$x_2$}

\put(90,-8){$x_1$}

\put(52,42){$x_8$}

\put(49,45){\line(1,0){2}}

\put(72,42){$x_7$}

\put(69,45){\line(1,0){2}}

\linethickness{5pt}

\put(60,10){\line(0,1){10}}

\put(60,70){\line(0,1){10}}


\linethickness{.25pt}

\put(210,0){\vector(0,1){90}}

\put(170,0){\vector(0,1){30}}

\put(190,0){\vector(0,1){30}}

\qbezier(170,30)(170,35)(177.5,40)

\qbezier(170,60)(170,55)(177.5,50)

\qbezier(190,30)(190,35)(182.5,40)

\qbezier(190,60)(190,55)(182.5,50)

\put(170,60){\vector(0,1){30}}

\put(190,60){\vector(0,1){30}}

\put(170,-20){$\Gamma_{100}=\Gamma_{001}$}

\put(170,92){$x_6$}

\put(190,92){$x_5$}

\put(210,92){$x_4$}

\put(170,-8){$x_3$}

\put(190,-8){$x_2$}

\put(210,-8){$x_1$}

\linethickness{5pt}

\put(180,40){\line(0,1){10}}


\linethickness{.25pt}

\put(-210,-90){\vector(0,1){60}}

\put(-190,-60){\vector(0,1){30}}

\put(-170,-120){\vector(0,1){30}}

\put(-202.5,-110){\vector(1,1){0}}

\qbezier(-210,-120)(-210,-115)(-202.5,-110)

\qbezier(-210,-90)(-210,-95)(-202.5,-100)

\put(-197.5,-110){\vector(-1,1){0}}

\qbezier(-190,-120)(-190,-115)(-197.5,-110)

\qbezier(-190,-90)(-190,-95)(-197.5,-100)

\qbezier(-190,-90)(-190,-85)(-182.5,-80)

\qbezier(-190,-60)(-190,-65)(-182.5,-70)

\qbezier(-170,-90)(-170,-85)(-177.5,-80)

\qbezier(-170,-60)(-170,-65)(-177.5,-70)

\put(-190,-90){\vector(1,4){0}}

\put(-170,-60){\vector(0,1){30}}

\put(-190,-140){$\Gamma_{011}$}

\put(-210,-28){$x_6$}

\put(-190,-28){$x_5$}

\put(-170,-28){$x_4$}

\put(-210,-128){$x_3$}

\put(-190,-128){$x_2$}

\put(-170,-128){$x_1$}

\put(-188,-93){$x_7$}

\put(-191,-90){\line(1,0){2}}

\linethickness{5pt}

\put(-200,-110){\line(0,1){10}}

\put(-180,-80){\line(0,1){10}}


\linethickness{.25pt}

\put(-90,30){\vector(0,1){30}}

\put(-50,0){\vector(0,1){30}}

\qbezier(-70,-90)(-70,-85)(-62.5,-80)

\qbezier(-70,-60)(-70,-65)(-62.5,-70)

\put(-70,-60){\vector(-1,4){0}}

\qbezier(-50,-90)(-50,-85)(-57.5,-80)

\qbezier(-50,-60)(-50,-65)(-57.5,-70)

\qbezier(-90,-60)(-90,-55)(-82.5,-50)

\qbezier(-90,-30)(-90,-35)(-82.5,-40)

\put(-90,-30){\vector(-1,4){0}}

\qbezier(-70,-60)(-70,-55)(-77.5,-50)

\qbezier(-70,-30)(-70,-35)(-77.5,-40)

\put(-70,-30){\vector(1,4){0}}

\put(-50,-60){\vector(0,1){30}}

\put(-90,-120){\vector(0,1){60}}

\put(-70,-120){\vector(0,1){30}}

\put(-50,-120){\vector(0,1){30}}

\put(-70,-140){$\Gamma_{110}$}

\put(-90,-28){$x_6$}

\put(-70,-28){$x_5$}

\put(-50,-28){$x_4$}

\put(-90,-128){$x_3$}

\put(-70,-128){$x_2$}

\put(-50,-128){$x_1$}

\put(-68,-63){$x_9$}

\put(-71,-60){\line(1,0){2}}

\linethickness{5pt}

\put(-60,-80){\line(0,1){10}}

\put(-80,-50){\line(0,1){10}}


\linethickness{.25pt}

\put(50,-120){\vector(0,1){90}}

\put(70,-120){\vector(0,1){90}}

\put(90,-120){\vector(0,1){90}}

\put(70,-140){$\Gamma_{000}$}

\put(50,-28){$x_6$}

\put(70,-28){$x_5$}

\put(90,-28){$x_4$}

\put(50,-128){$x_3$}

\put(70,-128){$x_2$}

\put(90,-128){$x_1$}


\linethickness{.25pt}

\qbezier(187.5,-70)(170,-60)(170,-50)

\put(170,-50){\vector(0,1){20}}

\put(190,-70){\vector(0,1){40}}

\qbezier(192.5,-70)(210,-60)(210,-50)

\put(210,-50){\vector(0,1){20}}

\qbezier(187.5,-80)(170,-90)(170,-100)

\put(170,-120){\vector(0,1){20}}

\put(190,-120){\vector(0,1){40}}

\qbezier(192.5,-80)(210,-90)(210,-100)

\put(210,-120){\vector(0,1){20}}

\put(190,-140){$\Gamma$}

\put(170,-28){$x_6$}

\put(190,-28){$x_5$}

\put(210,-28){$x_4$}

\put(170,-128){$x_3$}

\put(190,-128){$x_2$}

\put(210,-128){$x_1$}

\linethickness{5pt}

\put(190,-80){\line(0,1){10}}

\end{picture}

\caption{Local diagrams related to Reidemeister move
III}\label{gornik-Rmove3-decomp}

\end{figure}
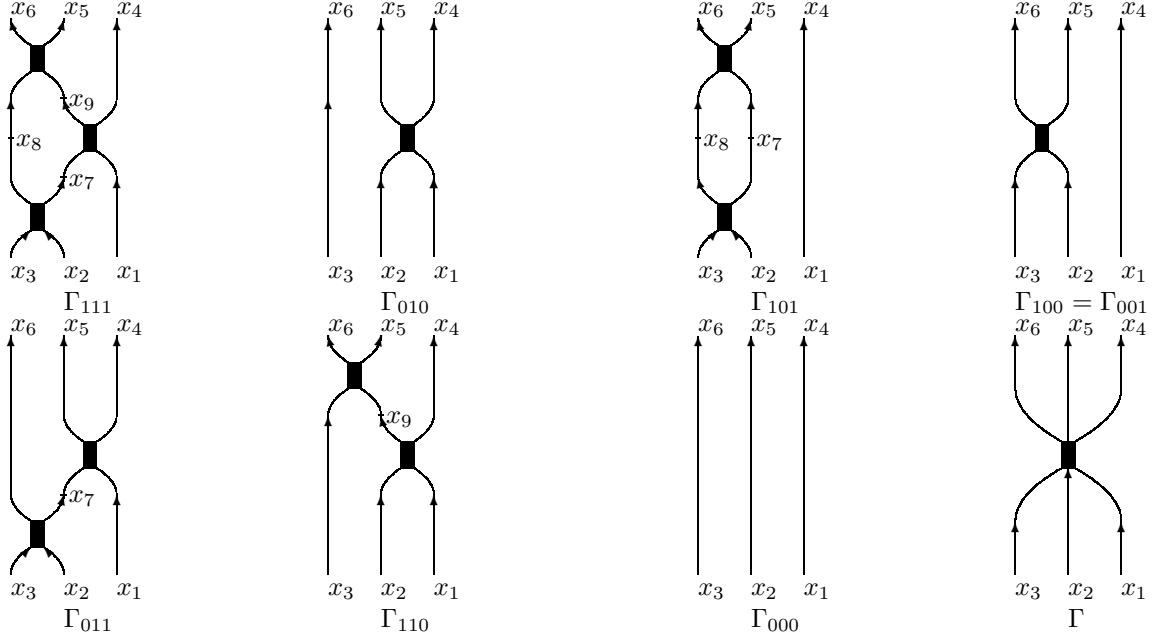

Recall that $\Psi:H_p(D)\rightarrow H_p(D')$ is induced by the Gaussian eliminations that reduce $(H(C_p(D),d_{mf}),d_\chi)$ and $(H(C_p(D'),d_{mf}),d_\chi)$ to chain complexes $(C,d')$ and $(C,d'')$ of the form
{\tiny
\[
0\rightarrow H(C_p(\Gamma),d_{mf})\{3n\} 
\rightarrow  
\left.%
\begin{array}{c}
  H(C_p(\Gamma_{011}),d_{mf})\{3n-1\}\\
  \oplus \\
  H(C_p(\Gamma_{110}),d_{mf})\{3n-1\}\\
\end{array}%
\right.
\rightarrow 
\left.%
\begin{array}{c}
  H(C_p(\Gamma_{001}),d_{mf})\{3n-2\}\\
  \oplus \\
  H(C_p(\Gamma_{010}),d_{mf})\{3n-2\}\\
\end{array}%
\right. 
\rightarrow
H(C_p(\Gamma_{000}),d_{mf})\{3n-3\}
\rightarrow 0,
\]}
and the isomorphism between $(C,d')$ and $(C,d'')$ is given by 
\begin{eqnarray*}
H(C_p(\Gamma),d_{mf})\{3n\} & \xrightarrow{\id} & H(C_p(\Gamma),d_{mf})\{3n\}, \\
\left.%
\begin{array}{c}
  H(C_p(\Gamma_{011}),d_{mf})\{3n-1\}\\
  \oplus \\
  H(C_p(\Gamma_{110}),d_{mf})\{3n-1\}\\
\end{array}%
\right.
& \xrightarrow{\left(%
\begin{array}{cc}
  \lambda_1 & 0 \\
  0 &  \lambda_2 \\
\end{array}%
\right)} &
\left.%
\begin{array}{c}
  H(C_p(\Gamma_{011}),d_{mf})\{3n-1\}\\
  \oplus \\
  H(C_p(\Gamma_{110}),d_{mf})\{3n-1\}\\
\end{array}%
\right., \\
\left.%
\begin{array}{c}
  H(C_p(\Gamma_{001}),d_{mf})\{3n-2\}\\
  \oplus \\
  H(C_p(\Gamma_{010}),d_{mf})\{3n-2\}\\
\end{array}%
\right. & \xrightarrow{\left(%
\begin{array}{cc}
  \lambda_1\lambda_3 & 0 \\
  0 &  \lambda_1\lambda_4 \\
\end{array}%
\right)} &
\left.%
\begin{array}{c}
  H(C_p(\Gamma_{001}),d_{mf})\{3n-2\}\\
  \oplus \\
  H(C_p(\Gamma_{010}),d_{mf})\{3n-2\}\\
\end{array}%
\right., \\
H(C_p(\Gamma_{000}),d_{mf})\{3n-3\} & \xrightarrow{\lambda_1\lambda_3\lambda_7\id} & H(C_p(\Gamma_{000}),d_{mf})\{3n-3\},
\end{eqnarray*}
where $\lambda_1, \cdots, \lambda_7$ are non-zero scalars.

If $j=k=l$, then $f_\varphi$ and $f_{\varphi'}$ are multiples of the same non-zero element of $H(C_p(\Gamma_{000}),d_{mf})\{3n-3\}$. So they are not affected by the Gaussian elimination, and the above isomorphism maps $f_\varphi$ to a non-zero multiple of $f_{\varphi'}$. Thus, $\Psi([f_\varphi])=\lambda[f_{\varphi'}]$, where $\lambda\neq0$.

If $j=k\neq l$, then $f_\varphi$ and $f_{\varphi'}$ are multiples of the same non-zero element of $H(C_p(\Gamma_{011}),d_{mf})\{3n-1\}$. So they are not affected by the Gaussian elimination, and the above isomorphism maps $f_\varphi$ to a non-zero multiple of $f_{\varphi'}$. Thus, $\Psi([f_\varphi])=\lambda[f_{\varphi'}]$, where $\lambda\neq0$. The case when $j\neq k= l$ is similar.

If $j=l\neq k$, then $f_\varphi$ is a non-zero element of $H(C_p(\Gamma_{101}),d_{mf})\{3n-1\}$. Let $\varphi_{100}$ be the admissible state of $\Gamma_{100}$ that agrees with $\varphi$ outside the visible part of $\Gamma_{100}$ in Figure \ref{gornik-Rmove3-decomp} with $\varphi_{100}(e_1)=\varphi_{100}(e_3)=\varphi_{100}(e_4)=\varphi_{100}(e_6)=j$ and $\varphi_{100}(e_2)=\varphi_{100}(e_5)=k$, where $e_c$ is the edge marked by $x_c$. Let $f_{100}$ be a non-zero element of $Q_{\varphi_{100}} H_p(\Gamma_{100})$. Then, similar to the proof of Reidemeister move II$_a$, we know that $f_\varphi$ is mapped to a non-zero multiple of 
\[
\chi_1\circ J(f_{100}) + \chi_1'\circ J(f_{100}) \in  H(C_p(\Gamma_{011}),d_{mf})\{3n-1\} \oplus H(C_p(\Gamma_{110}),d_{mf})\{3n-1\}
\] 
after the Gaussian elimination, where $\chi_1$ and $\chi_1'$ are the $\chi_1$-maps associated to the upper and lower wide edges of $\Gamma_{111}$ and $J:H(C_p(\Gamma_{100}),d_{mf})\{3n\} \rightarrow H(C_p(\Gamma_{111}),d_{mf})\{3n\}$ is the inclusion map from the decomposition
\[
H(C_p(\Gamma_{111}),d_{mf})\{3n\} \cong H(C_p(\Gamma),d_{mf})\{3n\} \oplus H(C_p(\Gamma_{100}),d_{mf})\{3n\}.
\]
Let $\varphi_{011,c}$ be the state of $\Gamma_{011}$ agrees with $\varphi$ outside the visible part of $\Gamma_{011}$ in Figure \ref{gornik-Rmove3-decomp} with $\varphi_{011,c}(e_1)=\varphi_{011,c}(e_3)=\varphi_{011,c}(e_4)=\varphi_{011,c}(e_6)=j$, $\varphi_{011,c}(e_2)=\varphi_{011,c}(e_5)=k$ and $\varphi_{011,c}(e_7)=c$, and $\varphi_{110,c}$ be the state of $\Gamma_{110}$ agrees with $\varphi$ outside the visible part of $\Gamma_{110}$ in Figure \ref{gornik-Rmove3-decomp} with $\varphi_{110,c}(e_1)=\varphi_{110,c}(e_3)=\varphi_{110,c}(e_4)=\varphi_{110,c}(e_6)=j$, $\varphi_{110,c}(e_2)=\varphi_{110,c}(e_5)=k$ and $\varphi_{110,c}(e_9)=c$. Clearly, $\varphi_{011,c}$ is admissible if and only if $c=k$; $\varphi_{110,c}$ is admissible if and only if $c=k$. From Theorem \ref{admissible-states},
\begin{eqnarray*}
\chi_1\circ J(f_{100}) & \in & Q_{\varphi_{011,k}} H(C_p(\Gamma_{011}),d_{mf})\{3n-1\},  \\
\chi_1'\circ J(f_{100})& \in & Q_{\varphi_{110,k}} H(C_p(\Gamma_{110}),d_{mf})\{3n-1\}.
\end{eqnarray*}
Let $f_{011}$ and $f_{110}$ be non-zero elements of the $1$-dimensional spaces $H(C_p(\Gamma_{011}),d_{mf})\{3n-1\}$ and $Q_{\varphi_{110,k}} H(C_p(\Gamma_{110}),d_{mf})\{3n-1\}$. Then $\chi_1\circ J(f_{100}) + \chi_1'\circ J(f_{100})=\mu_1 f_{011} + \mu_2 f_{110}$, where $\mu_1,\mu_2\in \C$, and $(\mu_1,\mu_2)\neq0$ since $\chi_1\circ J(f_{100}) + \chi_1'\circ J(f_{100})$ is a cocycle of $(C,d')$ and $[\chi_1\circ J(f_{100}) + \chi_1'\circ J(f_{100})]$ corresponds to $[f_\varphi]\neq0$ in the Gaussian elimination. From Gornik's computation, see Figure 6 of \cite{Gornik}, we have $d'(f_{011})\neq0$, $d'(f_{110})\neq0$. This implies that $\mu_1\neq 0$ and $\mu_2\neq0$. Altogether, we have that, after the Gaussian elimination, $f_\varphi$ is mapped to $\mu_1 f_{011} + \mu_2 f_{110}$, where $\mu_1\neq 0$ and $\mu_2\neq0$. Similarly, after the Gaussian elimination, $f_{\varphi'}$ is mapped to $\nu_1 f_{011} + \nu_2 f_{110}$, where $\nu_1\neq 0$ and $\nu_2\neq0$. Now use the fact that $d'(\mu_1 f_{011} + \mu_2 f_{110}) =d''(\nu_1 f_{011} + \nu_2 f_{110})=0$ and
\[
\left(%
\begin{array}{cc}
  \lambda_1\lambda_3 & 0 \\
  0 &  \lambda_1\lambda_4 \\
\end{array}%
\right)
\circ d' = 
d'' \circ \left(%
\begin{array}{cc}
  \lambda_1 & 0 \\
  0 &  \lambda_2 \\
\end{array}%
\right),
\]
it is easy to see that $\lambda_1\mu_1 f_{011} + \lambda_2\mu_2 f_{110}$ is a non-zero multiple of $\nu_1 f_{011} + \nu_2 f_{110}$. This implies that $\Psi([f_\varphi])=\lambda[f_{\varphi'}]$ for some $\lambda\in \C\setminus\{0\}$.

If $j, k, l$ are pairwise distinct, then $\varphi$ induces an admissible state of $\Gamma_{111}$, and $f_\varphi\in H(C_p(\Gamma_{111}),d_{mf})\{3n\}$. By Gornik's computation, see Figure 6 of \cite{Gornik}, $\chi_1''(f_\varphi)=0$, where $\chi_1''$ is the $\chi_1$-map associated to the right wide edge in $\Gamma_{111}$. But the projection $P'':H(C_p(\Gamma_{111}),d_{mf})\{3n\}\rightarrow H(C_p(\Gamma_{100}),d_{mf})\{3n\}$ from the decomposition 
\[
H(C_p(\Gamma_{111}),d_{mf})\{3n\} \cong H(C_p(\Gamma),d_{mf})\{3n\} \oplus H(C_p(\Gamma_{100}),d_{mf})\{3n\}
\]
factors through $\chi_1'':H(C_p(\Gamma_{111}),d_{mf})\{3n\}\rightarrow H(C_p(\Gamma_{101}),d_{mf})\{3n\}$. This means that $f_\varphi\in H(C_p(\Gamma),d_{mf})\{3n\}$ in the above decomposition. Let $\hat{\varphi}$ be the admissible state of $\Gamma$ induced by $\varphi$. Then $f_\varphi$ is a non-zero element of $Q_{\hat{\varphi}}H(C_p(\Gamma),d_{mf})\{3n\}$. Similarly, $f_{\varphi'}$ is a non-zero element of $Q_{\hat{\varphi}}H(C_p(\Gamma),d_{mf})\{3n\}$. According to Proposition \ref{general-admissible-states}, $Q_{\hat{\varphi}}H(C_p(\Gamma),d_{mf})\{3n\}$ is $1$-dimensional. So $f_{\varphi'}$ is a non-zero multiple of $f_\varphi$, which implies that $\Psi([f_\varphi])=\lambda[f_{\varphi'}]$ for some $\lambda\in \C\setminus\{0\}$.
\end{proof}

We are now ready to generalize Rasmussen's arguments in \cite{Ras1}. Let $S$ be a cobodism from link $D$ to link $D'$, and $(S_1,\cdots,S_m)$ a movie presentation of $S$, where each $S_j$ is an elementary link cobodism. Let $D_{j-1}$ and $D_j$ be the initial and terminal ends of $S_j$. In particular, $D_0=D$ and $D_m=D'$. Define $\overline{S}_j=S_1\cup\cdots\cup S_j$, and $\overline{\Psi}_j=\Psi_{S_j}\circ\cdots\circ\Psi_{S_1}$. $\overline{S}_j$ is a cobodism from $D=D_0$ to $D_j$. A component of $D_j$ is called evolved if it is a boundary component of a component of $\overline{S}_j$ that contains at least one component of $D$. A component of $D_j$ is called created if it is not evolved. For $k=0,1,\cdots,n-1$, let $\varphi^{(k)}$ be the state of $D$ whose value is $k$ on every component of $D$. A state $\psi$ of $D_j$ is said to be compatible with $\varphi^{(k)}$ if any two components of $D\sqcup D_j$ contained in the same component of $\overline{S}_j$ have the same value under the function $\varphi^{(k)} \sqcup \psi$. In particular, if $\psi$ is compatible with $\varphi^{(k)}$, then its value is $k$ on every evolved component of $D_j$. The following is a generalization of Proposition 4.1 of \cite{Ras1}.

\begin{proposition}[Compare with Proposition 4.1 of \cite{Ras1}]\label{gornik-rasmussen-wu}
If $S$ has no closed components, then
\[
\overline{\Psi}_j([f_{\varphi^{(k)}}]) = \sum_{\psi} \lambda_{\psi,j} [f_\psi],
\]
where $\psi$ runs through all states of $D_j$ compatible with $\varphi^{(k)}$, and each $\lambda_{\psi,j}$ is a non-zero scalar.
\end{proposition}

\begin{proof}
Define $\overline{\Psi}_0=\id:H_p(D)\rightarrow H_p(D)$. We prove the proposition by induction on $j$. when $j=0$, the proposition is trivially true. Assume the proposition is true for $j-1$. Consider $\overline{\Psi}_j = \Psi_{S_j}\circ \overline{\Psi}_{j-1}$. 

If $S_j$ corresponds to a Reidemeister move, then, by Proposition \ref{gornik-reidemeister}, the proposition is clearly true for $j$.

If $S_j$ corresponds to a circle creation, then each state $\psi$ of $D_{j-1}$ induces $n$ states $\psi_l$ of $D_j$, where $l=0,1,\cdots,n-1$, and $\psi_l$ agrees with $\psi$ on components inherited from $D_{j-1}$, and takes value $l$ on the new component. Clearly, if $\psi$ is compatible with $\varphi^{(k)}$, so is $\psi_l$, $l=0,1,\cdots,n-1$. Moreover, any state of $D_j$ compatible with $\varphi^{(k)}$ is induced by a state of $D_{j-1}$ compatible with $\varphi^{(k)}$. By Proposition \ref{gornik-morse}, there is a non-zero scalar $\lambda$ such that 
\[
\Psi_{S_j}([f_\psi]) = \lambda \sum_{l=0}^{n-1} [f_{\psi_l}].
\]
Thus,
\[
\overline{\Psi}_j([f_{\varphi^{(k)}}]) = \sum_{\psi} \sum_{l=0}^{n-1} \lambda \lambda_{\psi,j-1} [f_{\psi_l}],
\]
where $\psi$ runs through all states of $D_{j-1}$ compatible with $\varphi^{(k)}$, and $\lambda_\psi\neq0$ $\forall~ \psi$. Note that each state of $D_j$ compatible with $\varphi^{(k)}$ appears exactly once in the above sum. The proposition for $j$ follows in this case.

If $S_j$ corresponds to a circle annihilation, then every state $\psi$ of $D_{j-1}$ induces a unique state $\psi'$ of $D_j$, which agrees with $\psi$ on the unchanged components. Clearly, if a state $\psi$ of $D_{j-1}$ is compatible with $\varphi^{(k)}$, then the state $\psi'$ of $D_j$ induced by $\psi$ is also compatible with $\varphi^{(k)}$. Moreover, any state of $D_j$ compatible with $\varphi^{(k)}$ is induced by a state of $D_{j-1}$ compatible with $\varphi^{(k)}$. If two different states of $D_{j-1}$ compatible with $\varphi^{(k)}$ induce the same state of $D_j$, then the circle being annihilated must be the only boundary component of a component of $\overline{S}_{j-1}$. Which means $\overline{S}_j$ has a closed component. This contradicts the assumption that $S$ has no closed components. Thus, different states of $D_{j-1}$ induce different states $D_j$. By Proposition \ref{gornik-morse}, for any state $\psi$ of $D_{j-1}$, there exists a non-zero scalar $\mu_{\psi}$ such that $\Psi_{S_j}([f_\psi]) = \mu_{\psi} [f_{\psi'}]$. So
\[
\overline{\Psi}_j([f_{\varphi^{(k)}}]) = \Psi_{S_j}(\sum_{\psi} \lambda_{\psi,j-1} [f_{\psi_l}]) = \sum_{\psi'} \mu_{\psi}\lambda_{\psi,j-1} [f_{\psi'}],
\]
where $\psi$ runs through all states of $D_{j-1}$ compatible with $\varphi^{(k)}$, and $\psi'$ runs through all states of $D_j$ compatible with $\varphi^{(k)}$. The proposition for $j$ follows in this case.

Finally, consider the case when $S_j$ corresponds to a saddle move. Let $K_1$ and $K_2$ be the components involved in the move. (It's possible that $K_1=K_2$.) Let $\psi$ be a state of $D_{j-1}$. If $\psi(K_1)\neq\psi(K_2)$, then $\Psi_{S_j}([f_{\psi}])=0$. If $\psi$ is compatible with $\varphi^{(k)}$, and $\psi(K_1)=\psi(K_2)$, then $\psi$ induces a unique state $\psi'$ of $D_j$ compatible with $\varphi^{(k)}$, which agrees with $\psi$ on the unchanged components, and takes the value $\psi(K_1)=\psi(K_2)$ on the changed components. Moreover, any state of $D_j$ compatible with $\varphi^{(k)}$ is induced by a unique state $\psi$ of $D_{j-1}$ compatible with $\varphi^{(k)}$ with $\psi(K_1)=\psi(K_2)$. By Proposition \ref{gornik-morse}, for any state $\psi$ of $D_{j-1}$ with $\psi(K_1)=\psi(K_2)$, there exists a non-zero scalar $\nu_{\psi}$ such that $\Psi_{S_j}([f_\psi]) = \nu_{\psi} [f_{\psi'}]$, where $\psi'$ is the state of $D_j$ induced by $\psi$. So
\[
\overline{\Psi}_j([f_{\varphi^{(k)}}]) = \Psi_{S_j}(\sum_{\psi} \lambda_{\psi,j-1} [f_{\psi_l}]) = \sum_{\psi'} \nu_{\psi}\lambda_{\psi,j-1} [f_{\psi'}],
\]
where $\psi$ runs through all states of $D_{j-1}$ compatible with $\varphi^{(k)}$, and $\psi'$ runs through all states of $D_j$ compatible with $\varphi^{(k)}$. The proposition for $j$ follows in this case.
\end{proof}

Theorem \ref{ras-genus} follows from Proposition \ref{gornik-rasmussen-wu}.

\begin{proof}[Proof of Theorem \ref{ras-genus}]
Assume that $L$ is a link in $S^3$, and $S$ is an oriented compact surface with no closed components smoothly embedded in $D^4$ bounded by $-L$ and $\chi(S)=\chi_s(L)$. After possibly a small perturbation of $S$ that fixes $L$, we delete a small disk in the interior of $S$, and get a link cobodism $S^0$ from $L$ to the unknot $\bigcirc$. Note that $\chi(S^0)=\chi(S)-1$. Let $(S_1,\cdots,S_m)$ be a movie presentation of $S^0$, and $\Psi=\Psi_{S_m}\circ\cdots\circ\Psi_{S_1}$. For $k=0,1,\cdots,n-1$, let $\varphi^{(k)}$ be the state of $L$ that takes value $k$ on every component of $L$, and $\varphi_0^{(k)}$ the state of $\bigcirc$ that takes value $k$ on the only component of $\bigcirc$. Since $S$ has no closed components, $\bigcirc$ must be evolved from the components of $L$. This means that $\varphi_0^{(k)}$ is the only state of $\bigcirc$ compatible with $\varphi^{(k)}$. By Proposition \ref{gornik-rasmussen-wu}, $\Psi([f_{\varphi^{(k)}}])=\lambda_k [f_{\varphi_0^{(k)}}]$, where $\lambda_k \in \C\setminus\{0\}$. But $\{[f_{\varphi_0^{(0)}}], [f_{\varphi_0^{(1)}}], \cdots, [f_{\varphi_0^{(n-1)}}]\}$ is a basis for $H_p(\bigcirc)$. So $\Psi:H_p(L)\rightarrow H_p(\bigcirc)$ is surjective. Let $[f]$ be a non-zero element of $H_p(\bigcirc)$ of degree $n-1=g^{max}_p(\bigcirc)$, then there exists a non-zero element $[\widetilde{f}]$ of $H_p(L)$ such that $\Psi([\widetilde{f}])=[f]$. But $\Psi\in\fil^{-(n-1)\chi(S^0)}\Hom(H_p(L),H_p(\bigcirc))$. So 
\[
n-1 = g^{max}_p(\bigcirc) = \deg ([f]) \leq \deg ([\widetilde{f}]) - (n-1)\chi(S^0) \leq g^{max}_p(L) - (n-1)\chi(S^0),
\]
and hence 
\[
g^{max}_p(L) \geq (n-1)\chi(S^0)+ g^{max}_p(\bigcirc) = (n-1) \chi(S) = (n-1)\chi_s(L).
\]

Next, consider the special situation when $L=K$ is a knot. In this case, $\Psi$ is an isomorphism since $\Psi:H_p(K)\rightarrow H_p(\bigcirc)$ is surjective and $\dim H_p(K) = \dim H_p(\bigcirc) =n$. Let $[f']$ be a non-zero element of $H_p(K)$ with $\deg ([f']) =g^{min}_p(K)$. Then $\Psi([f'])\neq0$, and 
\[
1-n = g^{min}_p(\bigcirc) \leq \deg (\Psi([f'])) \leq -(n-1)\chi(S^0)+ \deg ([f']) = -(n-1)\chi(S^0) +g^{min}_p(K).
\]
So,
\[
g^{min}_p(K) \geq (n-1)\chi(S^0) + g^{min}_p(\bigcirc) = (n-1)(\chi(S^0)-1) = (n-1)(\chi(S)-2).
\]
Note that, in this case, $S$ has only one boundary component and $\chi(S)=1-2g_s(K)$, where $g_s(K)$ is the slice genus of $K$. Therefore,
\[
s_n(K) = \frac{1}{2}(g^{max}_p(K)+g^{min}_p(K)) \geq (n-1)(\chi(S)-1) = -2(n-1)g_s(K).
\]
Now reverse the orientation of $S^0$, we get a cobodism $-S^0$ of genus $g$ from the unknot $\bigcirc$ to $K$. A movie presentation of $-S^0$ induces a homomorphism $\Phi:H_p(\bigcirc)\rightarrow H_p(K)$. Then $\Phi$ is also an isomorphism, and $\Phi\in \fil^{-(n-1)\chi(S^0)}\Hom(H_p(\bigcirc),H_p(K))$. We can argue as above to show that
\begin{eqnarray*}
n-1 & = & g^{max}_p(\bigcirc) \geq (n-1)\chi(S^0) + g^{max}_p(K), \\
1-n & = & g^{min}_p(\bigcirc) \geq (n-1)\chi(S^0) + g^{min}_p(K). 
\end{eqnarray*}
So,
\[
s_n(K) = \frac{1}{2}(g^{max}_p(K)+g^{min}_p(K)) \leq -(n-1)\chi(S^0) = -(n-1)(\chi(S)-1) = 2(n-1)g_s(K).
\]
Thus, 
\[
|s_n(K)|\leq 2(n-1) g_s(K).
\]
\end{proof}

\bibliographystyle{amsplain}

\end{document}